\documentclass[11pt]{article}
\usepackage{latexsym,amsfonts,amsmath,epsfig,amssymb}
\usepackage[colorlinks,linkcolor=red,anchorcolor=blue,urlcolor=red,citecolor=blue]{hyperref}
\usepackage{graphicx,times}
\usepackage{mathrsfs}
\usepackage{caption}
\usepackage{enumerate}
\usepackage{subfigure}
\usepackage{graphicx}
\usepackage{epstopdf}
\usepackage{multirow}
\usepackage{booktabs}
\usepackage{cite}
\usepackage{color}
\usepackage{algorithm}
\usepackage{algorithmicx}
\usepackage{algpseudocode}
\usepackage[numbers,sort&compress]{natbib}
\usepackage{mathtools}
\usepackage{array}
\usepackage{booktabs} 
\usepackage{cases}
\usepackage{stmaryrd}
\usepackage{appendix}
\graphicspath{{figure/}}

\usepackage{array}
\newcommand{\PreserveBackslash}[1]{\let\temp=\\#1\let\\=\temp}
\newcolumntype{C}[1]{>{\PreserveBackslash\centering}p{#1}}
\newcolumntype{R}[1]{>{\PreserveBackslash\raggedleft}p{#1}}
\newcolumntype{L}[1]{>{\PreserveBackslash\raggedright}p{#1}}

\newcommand{\qed}{\hfill $\Box$}

\newtheorem{theorem}{Theorem}

\newtheorem{Prop}{Proposition}

\newtheorem{lemma}{Lemma}

\newtheorem{assumption}{Assumption}
\newtheorem{remark}{Remark}

\newtheorem{definition}{Definition}

\newcounter{nextauthor}
\def\mathrm{\mbox}

\begin{document}
	

\title{\Large {\bf Sparse Tucker Decomposition and Graph Regularization for  High-Dimensional Time Series Forecasting}}
\author{Sijia Xia\footnotemark[1],  \ Michael K. Ng\footnotemark[2],  \ and \ Xiongjun Zhang\footnotemark[3]}

\renewcommand{\thefootnote}{\fnsymbol{footnote}}

\footnotetext[1]{School of Mathematics and Statistics,
	Central China Normal University, Wuhan 430079, China (e-mail: sijiax@mails.ccnu.edu.cn).
}
\footnotetext[2]{Department of Mathematics, Hong Kong Baptist University,
	Kowloon Tong, Hong Kong (e-mail: michael-ng@hkbu.edu.hk).
	The research of this author was supported
	in part by GDSTC: Guangdong and Hong Kong Universities “1+1+1” Joint Research Collaboration Scheme UICR0800008-24,
National Key Research and Development Program of China under Grant 2024YFE0202900,
RGC GRF 12300125 and Joint NSFC and RGC N-HKU769/21.}
\footnotetext[3]{School of Mathematics and Statistics, and Key Laboratory of Nonlinear Analysis \& Applications (Ministry of Education),
	Central China Normal University, Wuhan 430079, China (e-mail: xjzhang@ccnu.edu.cn).
	The research of this author  was supported in part by the National Natural Science Foundation of China under Grant No. 12171189, Hubei Provincial Natural Science Foundation of
	China under Grant No. 2025AFB966, and  Fundamental
	Research Funds for the Central Universities under Grant No. CCNU24ai002.
}

\renewcommand{\thefootnote}{\arabic{footnote}}


\maketitle \vspace*{0mm}
\begin{center}
	\begin{minipage}{5.5in}
		$${\bf Abstract}$$

Existing methods of  
vector autoregressive model 
for multivariate time series analysis
make use of  		
low-rank matrix approximation or Tucker decomposition
to reduce the dimension of the over-parameterization issue.
		In this paper, we propose a sparse Tucker decomposition method with graph regularization for high-dimensional vector autoregressive time series.
		By stacking the time-series transition matrices into a third-order tensor,
the sparse Tucker decomposition is employed to characterize important 
interactions within the transition third-order tensor and 
reduce the number of parameters.
		Moreover, the graph regularization is employed to measure the local consistency of the response, predictor and temporal factor matrices in the vector autoregressive model.
The two proposed regularization techniques can be shown to more accurate 
parameters estimation.
		A non-asymptotic error bound of the estimator of the proposed method is established,
which is lower than those of the existing matrix or tensor based methods.
		A proximal alternating linearized minimization algorithm is designed to solve the resulting model and its global convergence is established under very mild conditions.
		Extensive numerical experiments on synthetic data and real-world datasets are carried out to verify the superior performance of the proposed method over existing state-of-the-art
		methods.
	\end{minipage}
\end{center}

\begin{center}
	\begin{minipage}{6.1in}
	\end{minipage}
\end{center}

\begin{center}
	\begin{minipage}{6.1in}
	\end{minipage}
\end{center}

\newpage
\section{Introduction}

With the  rapid development of information technology,
high-dimensional time series data have emerged across various fields in the current data-abundant environment,
such as finance \cite{tsay2005analysis}, economics \cite{sims1980macroeconomics}, ecology \cite{hampton2013quantifying},
 and meteorology \cite{duchon2012time}.
Classical vector autoregressive
(VAR) models are widely employed for modeling multivariate time
series data due to the flexible ability  for capturing cross-variable temporal dynamics, which have been applied in a wide range of areas, such as neuroimaging \cite{gorrostieta2012investigating}, signal processing \cite{basu2019low},  traffic state estimation \cite{chen2025forecasting}.
Consider an $m$-dimensional zero-mean VAR model of order $p$ (denoted by VAR($p$)) in the following form
\begin{equation}\label{test0}
	\begin{aligned}
		\mathbf{y}_{t}
		&=\mathbf{W}_1 \mathbf{y}_{t-1}+\mathbf{W}_2 \mathbf{y}_{t-2}+\cdots+\mathbf{W}_p\mathbf{y}_{t-p}+ \boldsymbol{\varepsilon}_t,\quad 1\leq t\leq T,
	\end{aligned}
\end{equation}
where $\{\mathbf{y}_{t}\}$ is the observed data with
$\mathbf{y}_{t}=(y_{1t}, y_{2t},\ldots,y_{mt})^H\in\mathbb{R}^{m}$,
$\mathbf{W}_i\in\mathbb{R}^{m \times m}, i=1,2,\ldots,p$, are the transition matrices,
$p$ is the lag order of the VAR model, $T$ denotes the sample size,
$\boldsymbol{\varepsilon}_t$ denotes independent and identically distributed (i.i.d.) error with $\boldsymbol{\varepsilon}_t = (\varepsilon_{1t},\varepsilon_{2t},\ldots, \varepsilon_{mt})^H\in\mathbb{R}^{m}, \mathbb{E}(\boldsymbol{\varepsilon}_t) = 0$, and var$(\boldsymbol{\varepsilon}_t) < \infty$.

In real-world applications, the dimension $m$ in (\ref{test0}) is often very large, 
which implies that the number of coefficient
parameters (i.e., $m^2p$) is large.
Hence an unrestricted VAR($p$)
model is likely to encounter the difficulty of over-parameterization,
where the number of parameters (e.g., coefficient matrices in time-varying VAR) inevitably exceeds the
number of observations.
As a result,  the corresponding VAR methods  cannot provide
reliable estimates nor accurate forecasts without further restrictions \cite{de2008forecasting}.
Estimation consistency of high-dimensional VAR models is achievable under  certain regularity conditions about the transition matrices \cite{wang2023rate}. For example,
if the coefficient matrices have an unobserved low-dimensional structure, such as sparsity
or low-rankness, the structure-inducing regularization methods, including least absolute shrinkage and selection operator (Lasso) \cite{Basu2015}
and nuclear norm penalty \cite{Negahban2011} give consistent estimates under the Gaussian assumption of the time series.

On the other hand, the
number of parameters in $\mathbf{W}_i$ increases quadratically with the dimension $m$, which makes it difficult
to apply VAR models to high-dimensional data. To overcome it, a commonly used approach is
to assume sparsity in parameter matrices, and many sparsity-imposing or inducing methods
can then be employed for estimation and variable selection, including $\ell_1$ regularization \cite{Basu2015, han2015direct, kock2025data, wong2020lasso},
weighted $\ell_1$ regularization \cite{wang2023regularized}, weakly sparse constraint on the transition matrix \cite{han2015direct}, and nonconcave penalization method based on $\ell_1$ norm \cite{zhu2020nonconcave}.
However, unlike the linear regression for Lasso, the time series  data exhibit temporal and cross-sectional dependence,
while the $\ell_1$ regularization based methods neglect this property,
which will seriously affect the accuracy of estimators of the sparse regularization based methods.

To address the over-parameterization issue of the VAR
model and the data dependence for time series data, some  
reduced-rank regression based methods were proposed and studied, see
 \cite{Negahban2011, basu2019low, raskutti2019convex, miao2023high} and references therein.
In these papers, the transition matrices are low-rank, which can reduce the number of 
parameters.
For example, by forming the transition matrices into a larger matrix, Basu et al. \cite{basu2019low} proposed a low-rank plus sparse estimation for high-dimensional VAR models, where  only order-one was discussed in this model.
 Samadi et al. \cite{samadi2024reduced} proposed a reduced-rank envelope
 VAR model by combining the envelope model into the reduced-rank
 VAR model to extract relevant information from complex data efficiently.
Additionally, Reinsel et al. \cite{reinsel1983some} proposed an autoregressive index model
based on the low-rankness assumption on the stacking matrices of $\mathbf{W}_i^H$.
 However, the reduced-rank model based matrix methods only consider the low-rankness of the specially stacking matrices, which can only reduce the dimensionality of these parameter matrices along one direction.

By reformulating  the transition matrices into a third-order tensor, Wang et al. \cite{Wang2022}
developed a tensor Tucker decomposition method for VAR time series modeling, where the transition tensor is decomposed
into a core tensor and three sparse factor matrices, and the low-rankness can be explored along different directions.
However, the above tensor based method did not consider the local patterns of the time series data.
Bahadori et al. \cite{bahadori2014fast} proposed
a unified low-rank tensor learning framework for multivariate spatio-temporal analysis by using
low Tucker rank constraint and spatial Laplacian regularization, which can be applied for VAR time series.
Moreover,  Harris et al. \cite{harris2021time} proposed a time-varying autoregressive model
by incorporating CANDECOMP/PARAFAC (CP)  \cite{hitchcock1927expression} decomposition and smoothness priors over time for multivariate time series.
However, the CP rank is NP-hard to determine in general \cite{hillar2013most}.
Besides, there is no theoretical guarantee in the above two work.

\subsection{Our Proposal}

In this paper, we propose a sparse Tucker decomposition method with graph regularization for high-dimensional  VAR time series.
Specifically, by stacking the transition matrices of (\ref{test0}) into a third-order tensor,
the sparse Tucker decomposition is employed to explore the low-rankness of the transition tensor along different dimensions,
which can reduce the parameters greatly for small Tucker rank.
Moreover, the sparsity is imposed on the core tensor, which can reduce the parameters further
and select significant variables along the response, predictor, and temporal factors for VAR time series.
Besides, the graph regularization on the factor matrices is utilized to characterize the local patterns
of the transition tensor, which is  capable of preserving the intrinsic manifold structure of data in the response, predictor and temporal factors.
Then a non-asymptotic error bound of the estimator of the proposed model is established under some conditions,
which is smaller than those of the existing sparse and low-rank based methods.
A proximal alternating linearized minimization (PALM) algorithm is developed to solve the resulting model
and its global convergence is established under very mild conditions.
Numerical experiments on simulated data and real-world datasets
substantiate the superiority of the proposed method compared with other competition methods.

The remaining parts of this paper are organized as follows. In the next section, some notations and preliminaries about tensors are given.
In Section \ref{SPTDMGR}, we propose a sparse Tucker decomposition with graph regularization method for high-dimensional VAR time series.
Moreover, 
a non-asymptotic error bound of the estimator of the proposed model is established.
In Section \ref{PALMG}, a PALM algorithm is designed to solve the resulting model and its global convergence is established.
In Section \ref{ExperM}, some numerical experiments are conducted to demonstrate the effectiveness of the proposed method.
Finally, the concluding remarks are given in Section \ref{CondRem}.
All technical proofs of the lemmas and main results
are deferred to the Appendix.

\subsection{Preliminaries}

The basic symbols and notations used throughout this paper  are summarized in Table \ref{tableNota}.

\begin{table}[htbp]
	\caption{Summary of the notations.}\label{tableNota}
	\footnotesize
	\centering
	\begin{tabular}{l|c}
		\toprule
		Notations & Description \\
		\midrule
		$a$/$\bf a$/$\bf A$/$\mathcal A$ & Scalars/Vectors/Matrices/Tensors \\
        ${\bf I}_{n}$ & The identity matrix with  size $n\times n$ \\
		$\mathcal A_{ijl}$ &  The $(i,j,l)$-th element of $\mathcal A$\\
		$\cdot^{H}$ & The conjugate transpose operator \\
        $\cdot^{-1}$ & The inverse operator \\
		$\textup{tr}(\cdot)$ & The trace of a matrix \\
		$\mathcal A^{\langle i\rangle}$ & The $i$-th frontal slice of $\mathcal A$ \\
		$\mathcal A_{(i)}$ &  The mode-$i$ unfolding of $\mathcal A$ \\
		$\|{\bf a}\|_2$ &  The $\ell_2$ norm of a vector ${\bf a}$ \\
		$\sigma_j({\bf A})$ & The $j$-th largest singular value of $\bf A$
		\\
        $\|\cdot\|_0$ & The $\ell_0$ norm defined as the number of nonzero entries		\\
		$\|{\bf A}\|_{*}$ & The nuclear norm of ${\bf A}$
		defined as $\|{\bf A}\|_*:=\sum_{j=1}^{\min\{n_1,n_2\}}\sigma_j({\bf A})$	\\
		$\|{\bf A}\|$ & The spectral norm of ${\bf A}$ defined as $\|{\bf A}\|:= \sigma_1({\bf A})$
		\\
		$\|{\bf A}\|_F$ &  The Frobenius norm of $\bf A$ defined as $\|{\bf A}\|_F:=\sqrt{\textup{Tr}({\bf A}^H\bf A)}$ \\
		$\langle \mathcal{A},\mathcal{B} \rangle$
		& The inner product of two tensors defined as $\langle \mathcal{A},\mathcal{B} \rangle:=\sum_{i=1}^{n_3}\textup{Tr}((\mathcal A^{\langle i\rangle})^H\mathcal B^{\langle i\rangle})$ \\
		$\operatorname{vec}(\mathcal A)$  & Vectorize a tensor $\mathcal A$ by the lexicographical order \\
		$\|\mathcal A\|_\infty$ & Tensor $\ell_\infty$ norm of $\mathcal A$ defined as $\|\mathcal A\|_{\infty}:=\max|\mathcal A_{ijl}|$ \\
     $\|\mathcal A\|_1$ & Tensor $\ell_1$ norm of $\mathcal A$ defined as $\|\mathcal A\|_{1}:=\sum_{i,j,l}|\mathcal A_{ijl}|$ \\
		$\|\mathcal{A}\|_F$ &  Tensor Frobenius norm of $\mathcal A$ defined as $\|\mathcal{A}\|_F:=\sqrt{\langle\mathcal{A},\mathcal{A} \rangle }$ \\
		$\delta_D(\cdot)$ & The indicator function of a set $D$ with $\delta_D(x)=0$ if $x\in D$, otherwise $+\infty$ \\
		\bottomrule
	\end{tabular}
\end{table}

%


For any matrix $\mathbf{A}\in \mathbb{C}^{n_1 \times n_2}, \mathbf{B}\in
\mathbb{C}^{n_3 \times n_4}$, their Kronecker product \cite[Section 1.3.6]{golub2013matrix} is defined as
$$
\mathbf{A} \otimes \mathbf{B}=\left(\begin{array}{cccc}
	a_{11} \mathbf{B} & a_{12} \mathbf{B} & \cdots & a_{1 n_2} \mathbf{B} \\
	a_{21} \mathbf{B} & a_{22} \mathbf{B} & \cdots & a_{2 n_2} \mathbf{B} \\
	\vdots & \vdots & & \vdots \\
	a_{n_1 1} \mathbf{B} & a_{n_1 2} \mathbf{B} & \cdots & a_{n_1 n_2} \mathbf{B}
\end{array}\right).
$$
For any two tensors $\mathcal{X} \in \mathbb{R}^{n_1 \times n_2 \times n_3}$ and $\mathcal{Y} \in \mathbb{R}^{n_4 \times n_5 \times n_6}$, their tensor outer product, denoted by $\mathcal{X} \circ \mathcal{Y} \in \mathbb{R}^{n_1 \times n_2 \times n_3\times n_4 \times n_5 \times n_6}$, is defined as
$
(\mathcal{X} \circ \mathcal{Y})_{i_1 i_2 i_3 i_4 i_5 i_6}=\mathcal{X}_{i_1 i_2 i_3} \mathcal{Y}_{i_4 i_5 i_6},
$
where $1 \leq i_j \leq n_j, j=1,\ldots, 6$.
The $k$-mode product of a tensor $\mathcal{X}\in \mathbb{R}^{n_1 \times \cdots \times n_d}$ with a matrix $\mathbf{A}\in \mathbb{R}^{J\times n_k}$ is denoted  by
$\mathcal{X}\times_k \mathbf{A}$ and is of size $n_1\times \cdots \times n_{k-1}\times J\times n_{k+1}\times \cdots\times n_d$ \cite{kolda2009tensor}.
The point-wise manner is defined as
$$
(\mathcal{X}\times_k \mathbf{A})_{i_1\cdots i_{k-1}ji_{k+1}\cdots i_d}=\sum_{i_k=1}^{n_k}\mathcal{X}_{i_1\cdots i_d}\mathbf{A}_{ji_k}.
$$
The mode-$k$ matricization of
a tensor $\mathcal{X}\in \mathbb{R}^{n_1 \times \cdots \times n_d}$ is denoted by $\mathcal{X}_{(k)}$
and arranges the mode-$k$ fibers to be
the columns of the resulting matrix.
Here the tensor element $(i_1, i_2,\ldots,i_d)$ maps to the matrix element $(i_k, j)$, where
$
j=1+\sum_{m=1,m \neq k}^d (i_m-1)L_m  \ \textup{with} \  L_m=\prod_{s=1, s\neq k}^{m-1}n_s.
$
Moreover, $\mathcal{Y}=\mathcal{X}\times_k \mathbf{A}$ is equivalent to $\mathcal{Y}_{(k)}=\mathbf{A}\mathcal{X}_{(k)}$.

Next we recall the definition of subgradient of a general function, which is given as follows.

\begin{definition}\cite[Definition 8.3]{Rockafellar1998}\label{thm12}
	Consider a function $f:\mathbb{R}^n\rightarrow \mathbb{R}\cup\{+\infty\}$ and any
	$\mathbf{x}\in\textup{dom}(f)$. For a vector
	$\mathbf{v}\in\mathbb{R}^n$, one says
	
	(a)
	$\mathbf{v}$ is a regular subgradient of $f$ at $\mathbf{x}$, written as $\hat{\partial} f(\mathbf{x})$, if
	$$
	\liminf _{\substack{\mathbf{y} \neq \mathbf{x} \\ \mathbf{y} \rightarrow \mathbf{x}}} \frac{1}{\|\mathbf{x}-\mathbf{y}\|_2}\left[f(\mathbf{y})-f(\mathbf{x})-\langle \mathbf{v}, \mathbf{y}-\mathbf{x}\rangle\right] \geq 0.
	$$
	
	(b)
	$\mathbf{v}$ is a subgradient of $f$ at $\mathbf{x}$, written as ${\partial} f(\mathbf{x})$, if there exist
	sequences $\mathbf{x}^{k} \rightarrow \mathbf{x}, f(\mathbf{x}^{k}) \rightarrow f(\mathbf{x}), \mathbf{v}^{k} \in \hat{\partial} f(\mathbf{x}^{k}) $ with $\mathbf{v}^{k}\rightarrow \mathbf{v}$.
\end{definition}

%
%

\section{Sparse Tucker Decomposition  with Graph Regularization}\label{SPTDMGR}

In this section, we present the proposed method. 
Let each frontal slice of $\mathcal{W}\in \mathbb{R}^{m\times m \times p}$ be  $\mathcal{W}^{\langle i\rangle}=\mathbf{W}_i, i=1,2,\ldots,p$.
Then problem (\ref{test0}) can be rewritten as
\begin{equation}\label{test2}
\begin{aligned}
\mathbf{y}_{t}
&=\left(\mathbf{W}_1, \mathbf{W}_2, \ldots, \mathbf{W}_{p}\right)\left(\begin{array}{c}
\mathbf{y}_{t-1} \\
\mathbf{y}_{t-2} \\
\vdots \\
\mathbf{y}_{t-p}
\end{array}\right)+\boldsymbol{\varepsilon}_t\\
&=\mathcal{W}_{(1)}\mathbf{x}_{t}+\boldsymbol{\varepsilon}_t,
\end{aligned}
\end{equation}
where $\mathbf{x}_{t}=(\mathbf{y}_{t-1}^{H},\mathbf{y}_{t-2}^{H},\ldots,
\mathbf{y}_{t-p}^{H})^{H}\in \mathbb{R}^{mp \times 1}$.
 Therefore, the Tucker decomposition
of $\mathcal{W}$ can be employed to explore the low-rankness of different modes in the VAR model,
i.e.,
\begin{equation}\label{tucker}
\mathcal{W}=\mathcal{G} \times_1 \mathbf{A}_1 \times_2 \mathbf{A}_2 \times_3 \mathbf{A}_3
=:\llbracket \mathcal{G}; \mathbf{A}_1, \mathbf{A}_2, \mathbf{A}_3\rrbracket,
\end{equation}
where $\mathcal{G} \in \mathbb{R}^{r_1 \times r_2 \times r_3 }$ is the core tensor, $\mathbf{A}_1\in\mathbb{R}^{ m \times r_1}, \mathbf{A}_2\in\mathbb{R}^{ m \times r_2}$
and
$\mathbf{A}_3\in\mathbb{R}^{ p \times r_3}$ are the response, predictor and temporal factor matrices, respectively.
Consequently, problem (\ref{test2}) can be expressed as
$$
\begin{aligned}
\mathbf{y}_{t}
=\mathcal{W}_{(1)}\mathbf{x}_{t}+\boldsymbol{\varepsilon}_t
=(\mathcal{G} \times_1 \mathbf{A}_1 \times_2 \mathbf{A}_2 \times_3 \mathbf{A}_3)
_{(1)}\mathbf{x}_{t}+\boldsymbol{\varepsilon}_t.
\end{aligned}
$$
Moreover, in order to obtain a unique Tucker decomposition of $\mathcal{W}$,
we consider the column orthogonal constraints on the factor matrices \cite{kolda2009tensor}, i.e.,  
$$
\mathbf{A}_i^{H}\mathbf{A}_i=\mathbf{I}_{r_i}, i=1,2,3.
$$ 
Note that 
$$
\mathbf{A}_1^H\mathbf{y}_{t}
=\mathcal{G}_{(1)}(\mathbf{A}_3\otimes\mathbf{A}_2)^H\mathbf{x}_{t}+\mathbf{A}_1^H\boldsymbol{\varepsilon}_t.
$$ 
Let
 $\mathbf{A}_1^H\mathbf{y}_{t}:=\mathbf{f}_t=(f_{1t},f_{2t},\ldots, f_{r_1t})^H\in\mathbb{R}^{r_1}$, which denotes $r_1$ response factors across $m$ variables of $\mathbf{y}_{t}$.
Additionally, $f_{it}=\sum_{j=1}^m(\mathbf{A}_1)_{ji}y_{jt}$ is the $i$-th response factor.
 If $(\mathbf{A}_1)_{ji}=0$, then $f_{it}$ is irrelevant to $y_{jt}$.
 In this case, $\mathbf{A}_1$ can be interpreted as the loadings of the response factors. 
 Similarly, $\mathbf{A}_2$ and $\mathbf{A}_3$ can be interpreted as the loadings of the predictor and temporal factors, respectively, see \cite{Wang2022} for more discussions. 

\begin{figure}[h]
	\centering
	\subfigure[CCD]{%
		\includegraphics[width=0.52\textwidth]{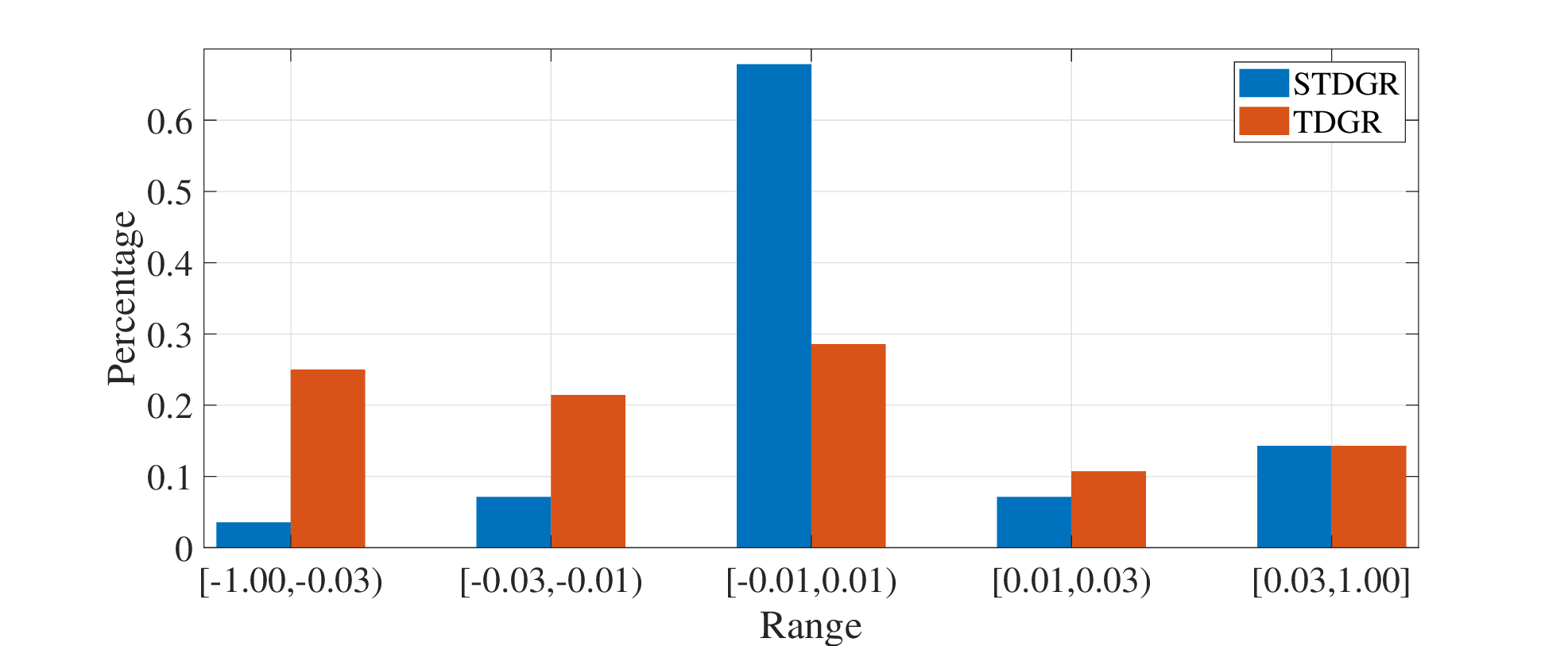}%
	}%
	\hspace{-7mm}%
	\subfigure[PeMS04]{%
		\includegraphics[width=0.52\textwidth]{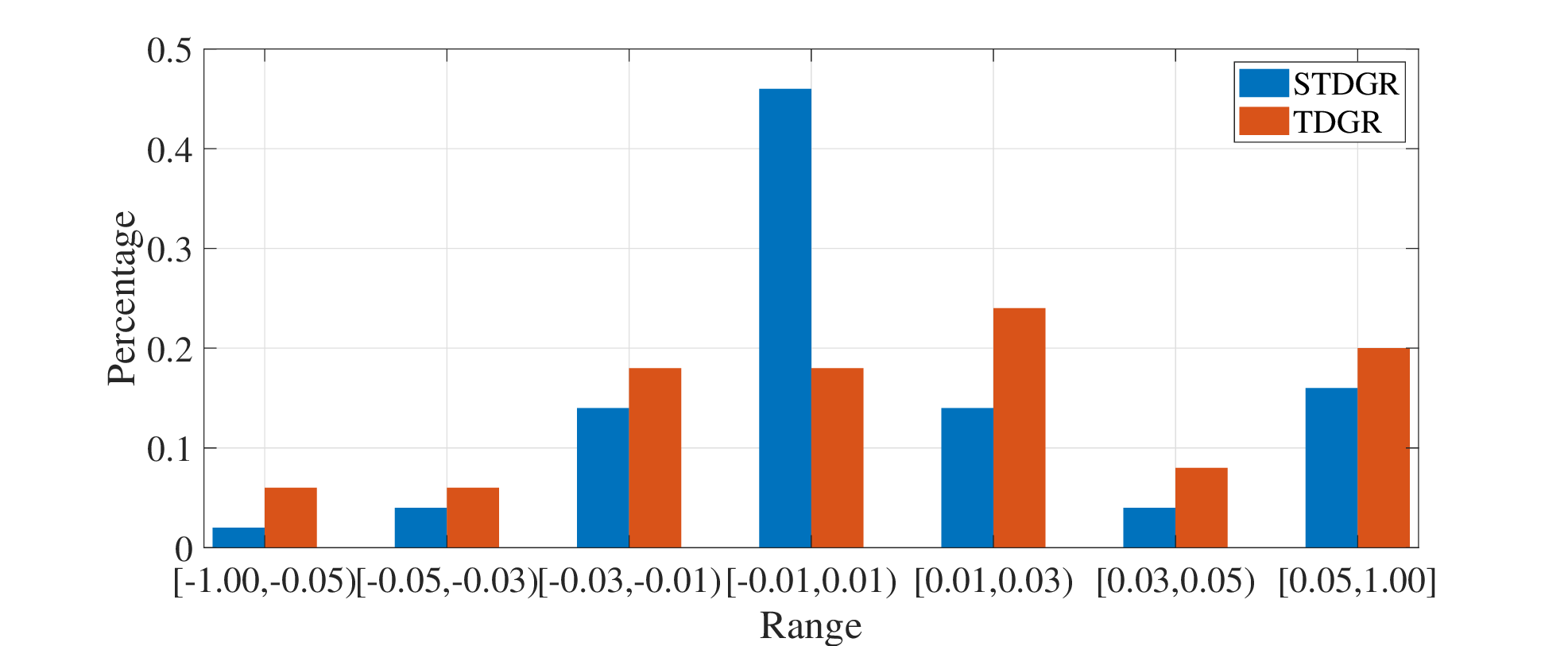}%
	}%
	\caption{Distributions of the elements of $\mathcal{G}$ obtained by non-sparse model (TDGR) and sparse model (STDGR) for the CCD and PeMS04 datasets, respectively.}
	\label{histogr}
\end{figure}

In the Tucker decomposition,
 the elements of the core tensor $\mathcal{G}$  reflect
 the interactions and connections between the components
 (columns) in different mode matrices, which
 only keeps the most significant
 connections between the components in different modes (response, predictor and
temporal).
Therefore, we impose the sparsity constraint on the core tensor. 
In Figure \ref{histogr}, we demonstrate the effectiveness of the sparsity constraint on the core tensor by numerical examples on the CCD and PeMS04 datasets (see 
the details in Sections \ref{NxMCCD} and \ref{Traffiddatas}).
In the figure, we show 
the distributions of the elements of $\mathcal{G}$ obtained by sparse model (STDGR)
and non-sparse model (TDGR). It is clear that many elements 
of $\mathcal{G}$ are close to zero (in $[-0.01,0.01)$), and significant 
 connections between the components in different modes (response, predictor and
temporal) are identified.  

Notice that the coefficient matrices
$\mathbf{A}_i$ are the low-dimensional representation of the sampled data in 
the response, predictor or temporal mode, which can reflect the similarity of sampled data in the low-dimensional subspace.
Hence, we utilize the graph prior to each factor matrix of Tucker decomposition of 
$\mathbf{A}_i$ to characterize the similarity of the sampled data along the response, predictor and temporal modes for VAR.


In this paper, we propose a sparse Tucker decomposition model combined with graph regularization for the factor matrices (called STDGR for short) as follows:
\begin{equation}\label{test3}
\begin{aligned}
\underset{\mathcal{G},\mathbf{A}_1,\mathbf{A}_2,\mathbf{A}_3}{\operatorname{min}}& \frac{1}{2T} \sum_{t=1}^{T}\left\|\mathbf{y}_t-\left(\mathcal{G} \times_1 \mathbf{A}_1 \times_2 \mathbf{A}_2 \times_3 \mathbf{A}_3\right)_{(1)}\mathbf{x}_t\right\|_2^2+\beta\|\mathcal{G}\|_1
+ \sum_{i=1}^{3}\alpha_i \operatorname{tr}\left(\mathbf{A}_i^{H}\mathbf{L}_i \mathbf{A}_i\right)\\
\text { s.t. }& \|\mathcal{G}\|_{\infty}\leq c, \mathbf{A}_i^{H}\mathbf{A}_i=\mathbf{I}_{r_i}, i=1,2,3,
\end{aligned}
\end{equation}
where $\beta>0, \alpha_i>0, c>0$ are given parameters
and $\mathbf{L}_i$ is the graph Laplacian matrix, $i=1,2,3$.
Here we note that for any matrix $\mathbf{A}\in\mathbb{R}^{m \times n}$,
its graph Laplacian regularization is defined as \cite{Qiu2022}
\begin{equation}\label{test1}
	\sum_{i,j=1}^m  \frac{1}{2} {z}_{i j}\left\|\mathbf{a}_i-\mathbf{a}_j\right\|_2^2= \operatorname{tr}\left(\mathbf{A}^{H} \mathbf{L}\mathbf{A}\right),
\end{equation}
where $\mathbf{a}_i$ denotes  the $i$-th row of $\mathbf{A}$, 
\begin{equation}\label{DefIoL}
\mathbf{L}=\mathbf{D}-{\mathbf{Z}}	
\end{equation}
represents the Laplacian matrix with each element of ${\mathbf{Z}}$ be ${z}_{i j}$,
and $\mathbf{D} \in$ $\mathbb{R}^{m\times m}$ is a diagonal
matrix with the $i$-th  diagonal element being  $\sum_{j=1}^m 
{z}_{i j}, i=1, \ldots, m$ \cite{chung1997spectral}.
In practice, $\mathbf{L}_i$ are given or constructed in advance. In Section \ref{ExperM} of numerical example, we discuss how to construct $\mathbf{L}_i$ in the implementation.


To prevent nonstationarity in the VAR process,
we assume that the core tensor $\mathcal{G}$ is upper bounded,
that is, $\|\mathcal{G}\|_{\infty}\leq c$.
The graph Laplacian regularization on the factor matrices in model (\ref{test3}) is used to explore
the local consistency of the low-dimensional subspace of the  coefficient tensor.
The factor matrices $\mathbf{A}_1, \mathbf{A}_2,$ and $\mathbf{A}_3$
represent the low-dimensional subspace loadings of
$\mathcal{W}$ along the response, predictor, and temporal modes, respectively \cite{Wang2022}.
Conversely, it can reflect the similarity of the sampled data for VAR.
For  traffic data as an example, adjacent sensors
show similar patterns and present temporal correlation properties \cite{Roughan2012},
and can exhibit similar low-dimensional loadings \cite{Cai2011}, which can be effectively captured by imposing graph Laplacian regularization on the response loading $\mathbf{A}_1$ and the predictor loading $\mathbf{A}_2$.
The Laplacian regularization on the temporal loading matrix $\mathbf{A}_3$ can capture temporal similarity across $p$ time lags, which facilitates adjacent lags to have similar subspace loadings.
Moreover, the sparse Tucker decomposition in model (\ref{test3})
is used to reduce the dimension and can address the over-parameterization issue in VAR. In fact, the parameters in VAR with sparse Tucker decomposition are only $mr_1+mr_2+pr_3+\tilde{s}r_1r_2r_3$, 
where $0<\tilde{s}<1$ denotes the sparse ratio of the core tensor. This is less than that of model (\ref{test0}) (i.e., $m^2p$) or the reduced-rank model \cite{Negahban2011} (i.e., $(mp+m-r_1)r_1$).


\begin{remark}
Compared with the model in \cite{Wang2022},
	the model in (\ref{test3}) used the sparse core tensor in Tucker decomposition, which can reduce the parameters of the transition tensor further.
	Moreover,
 the graph regularization is employed in model (\ref{test3}) to characterize the local consistency of the factor matrices,
which is capable of reflecting the similarity of the sampled data in response, 
predictor and temporal modes.
\end{remark}

\subsection{Non-asymptotic Error Bound}\label{TheErBR}

In this subsection, we establish the non-asymptotic error bound of the estimator of model (\ref{test3}).
Denote  the optimal solution of (\ref{test3}) by $(\widehat{\mathcal{G}},\widehat{\mathbf{A}}_1,\widehat{\mathbf{A}}_2,
\widehat{\mathbf{A}}_3)$.
Let $\widehat{\mathcal{W}}
=\widehat{\mathcal{G}}\times_1 \widehat{\mathbf{A}}_1 \times_2 \widehat{\mathbf{A}}_2 \times_3
\widehat{\mathbf{A}}_3$.
To derive the non-asymptotic error bound of $\widehat{\mathcal{W}}$, we first make the following assumptions.
\begin{assumption}\label{assum_stable}
The matrix polynomial $\mathbf{W}(z) = \mathbf{I}_m-
 \sum_{i=1}^p\mathbf{W}_iz^i, z\in \mathbb{C}$ satisfies
 $\det(\mathbf{W}(z))\neq 0$ for all $|z|=1$,
where $\mathbb{C}$ denotes the complex field and $\mathbf{W}_i$ is defined in (\ref{test0}).
\end{assumption}

\begin{assumption}\label{assum_E}
The error vector  $\{\boldsymbol{\varepsilon}_{t}\}$ in (\ref{test0}) are i.i.d. Gaussian random vectors
with mean zero and positive definite covariance matrix $\mathbf{\Sigma}_{\boldsymbol{\varepsilon}}$.
\end{assumption}

\begin{assumption}\label{assum_G}
The core tensor $\mathcal{G}\in \mathbb{R}^{r_1 \times r_2\times r_3}$ has at most $s$ nonzero entries. And
define $\bar{s}=1-s /(r_1 r_2r_3)\in[0,1]$,
which represents the sparse level of  $\mathcal{G}$, i.e., the proportion of zero entries.
\end{assumption}

\begin{assumption}\label{assum_spectral_A}
For $i=1,2,3,$
the nonzero singular values of $\mathcal{W}_{(i) }$ in (\ref{test2}) satisfy
$\sigma_{j-1}^2(\mathcal{W}_{(i)})-\sigma _j^2(\mathcal{W}_{(i)})
\geq$
$\rho \sigma _{j- 1}^2( \mathcal{W}_{(i)})$, $j=2,\ldots,r_i$, where $\rho>0$ is a constant.
\end{assumption}

\begin{assumption}\label{assum_delta_A}
Denote
$
\Delta_{\mathbf{A}_{i}}:=\widehat{\mathbf{A}}_{i}-\mathbf{A}_{i},
$
for $i=1,2,3$.
Suppose that $\Delta_{\mathbf{A}_{i}}$ satisfies
$
\sum_{l=1}^{r_i}\sum_{j=1}^{n_i}{(\Delta_{\mathbf{A}_{i}})}_{jl}^2$
$\gg
\sum_{l=1}^{r_i}\frac{1}{n_i}(\sum_{j=1}^{n_i}{(\Delta_{\mathbf{A}_{i}})}_{jl})^2
$ with $n_1=n_2=m$ and $n_3=p$.
\end{assumption}

Assumption \ref{assum_stable} guarantees the stability
of the VAR process \cite{Basu2015}.
Assumption \ref{assum_E} is used to employ the concentration inequalities for VAR models \cite{Basu2015}.
And Assumption \ref{assum_G} reflects the sparse level of core tensor $\mathcal{G}$.
Assumption \ref{assum_spectral_A} further ensures that the singular values of each unfolding matrices of $\mathcal{W}$
are sufficiently separated, which 
can  eliminate non-identifiability
and enable us to determine the upper bound for perturbation errors from Lemma \ref{HOSVD perturbation bound}.
Finally, Assumption \ref{assum_delta_A} is essential for the analysis of the graph regularization terms.

Let $\mathbf{W}^{H}(z)$ denote the conjugate transpose of $\mathbf{W}(z)$.
Following \cite{Basu2015,Wang2022}, we define
$$
\mu_{\min}(\mathbf{W}):=\min_{|z|=1}\lambda_{\min}(\mathbf{W}^{H}(z)\mathbf{W}(z)),
\ \mu_{\max}(\mathbf{W}):=\max_{|z|=1}\lambda_{\max}(\mathbf{W}^{H}(z)\mathbf{W}(z)),
$$
where $\lambda_{\min}(\cdot)$ and $\lambda_{\max}(\cdot)$ denote the minimum
and maximum eigenvalues of a matrix, respectively.
In addition,
we assume that the core tensor in model (\ref{test3}) also satisfies the constraint $\mathcal{U}:=\{\mathcal{G}:\|\mathcal{G}\|_0\leq s\}$.
Now we give the error bound about the estimator of  (\ref{test3}) in the following theorem.

%

\begin{theorem}\label{statistical}
Suppose that Assumptions \ref{assum_stable}-\ref{assum_delta_A} hold.
If the sample size of the VAR model satisfies $T\geq\max\{a_1\log (m^2p), \frac{
	8\varsigma^2}
{a_4}
\bar{p}\min(\log (mp), \log(21emp/\bar{p}))\}$,
and the regularization parameters satisfy $\alpha_{i}\geq{\frac{6c^2}{\lambda_{2}(\mathbf{L}_{i})}\sqrt{\frac{s\log (m^2p)}{T}}}$,
$\beta \geq 2
\pi\vartheta a_2\sqrt{\log(m^2p)/T}$, then
with probability at least $1-6\exp(-a$
$\log (m^2 p))-2 \exp (-\bar{p}\min(\log mp, \log(21emp/\bar{p})))$,
the estimator of (\ref{test3}) satisfies
\begin{equation}\label{WEstErB}
\|\widehat{\mathcal{W}}-\mathcal{W}\|_{F}\leq
\frac{4\beta(\sqrt{\tilde{s}}+b\sqrt{s})
+\hat{a}{\kappa}
+\bar{c}\sum_{i=1}^{3}\alpha_il_i \eta_i}{\varpi+{2\kappa\sqrt{s/(r_1 r_2 r_3)^2}}},
\end{equation}
and
\begin{equation}\label{Preerrorb}
\begin{aligned}
\frac{1}{T}\sum_{t=1}^{T}\|\widehat{\mathcal{W}}_{(1)}\mathbf{x}_{t}-\mathcal{W}_{(1)}\mathbf{x}_{t}\|_2^2
\leq
\frac{\lambda_{\max}
\left(\boldsymbol{\Sigma}_{\boldsymbol{\varepsilon}}\right) (4\beta(\sqrt{\tilde{s}}+b\sqrt{s})
+\hat{a}\kappa
+\bar{c}\sum_{i=1}^{3}\alpha_il_i \eta_i)^2}{\mu_{\min }(\mathbf{W})\left(\varpi+2\kappa\sqrt{s/(r_1 r_2 r_3)^2}\right)^2},
\end{aligned}
\end{equation}
where  $a_1, a_2\geq1, a,\hat{a}, a_4,\bar{c},b>0$ are absolute constants,
$\tilde{s}, \bar{p}$ depend on $\bar{s}$,
$\kappa=\sqrt{(\log(m^2p))/T}$,
$l_i=\|\mathbf{L}_{i}\|_F$,
$\varpi = \lambda _{\min }( \mathbf{\Sigma }_{\boldsymbol{\varepsilon}} )
/ \mu _{\max }(\mathbf{W})$,
$\eta_i=\left(3\sigma_1({\mathcal {W}}_{(i)})+
c\sqrt{r_1r_2r_3}\right)
\sqrt{\sum_{j=1}^{r_i}\frac{1}
	{\sigma_j^4
		(\mathcal{W}_{(i)})}}$,
	 $\varsigma=\frac{\lambda_{\min}
	\left(\boldsymbol{\Sigma}_{\boldsymbol{\varepsilon}}\right) / \mu_{\max }(\mathbf{W})}{\lambda_{\max}
	\left(\boldsymbol{\Sigma}_{\boldsymbol{\varepsilon}}\right) / \mu_{\min }(\mathbf{W})}$,
$\vartheta=
\lambda _{\max }(\mathbf{\Sigma}_{\boldsymbol{\varepsilon}})
( 1+ \mu _{\max }(\mathbf{W}) / \mu _{\min }(\mathbf{W} ) )$,
and $\lambda_{2}(\mathbf{L}_{i})>0$
denotes the second smallest eigenvalue of the Laplacian matrix $\mathbf{L}_{i}$, $i=1,2,3$.
\end{theorem}

In Theorem \ref{statistical}, we
establish a non-asymptotic error bound for model (\ref{test3}).
When the Tucker rank $\left(r_1, r_2, r_3\right)$ is fixed,
the inequality in (\ref{WEstErB}) shows that $\widehat{\mathcal{W}}$ is a consistent estimator
as
$T\rightarrow\infty$.
In fact, the estimation and prediction error bounds in (\ref{WEstErB}) and (\ref{Preerrorb})
are given by $O_p((\sqrt{s}+\sqrt{\tilde{s}})\sqrt{\log(m^2 p) / T})$ and
$O_p((\sqrt{s}+\sqrt{\tilde{s}})^2\log(m^2 p)/T)$, respectively.

\begin{remark}
	If $\tilde{s} \leq s$, where $\tilde{s}$ is defined  in  (\ref{S1S2St}),
	 the error bound in (\ref{WEstErB}) reduces to
	\begin{equation}\label{bound2}
		O_p(\sqrt{{s}\log(m^2 p) / T}).
	\end{equation}
Now we compare the estimation error bound in Theorem \ref{statistical} with the existing methods for VAR modeling.
	\begin{itemize}
		\item Basu et al. \cite{Basu2015} studied the estimation of stationary Gaussian
		VAR$(p)$ model with the sparse transition matrix.
The corresponding Lasso estimator  is given by
		$$
		\widehat{\mathcal{W}}_{\operatorname{Lasso}}=\arg\min_{\mathcal{W}}
		\frac{1}{T}
		\sum_{t=1}^T\|\mathbf{y}_t-\mathcal{W}_{(1)} \mathbf{x}_t\|_2^2
		+\lambda_{\operatorname{Lasso}}\|\mathcal{W}_{(1)}\|_1.
		$$
		It has been shown that the estimator $\widehat{\mathcal{W}}_{\operatorname{Lasso}}$ satisfies
		$\|\widehat{\mathcal{W}}_{\operatorname{Lasso}}-\mathcal{W}\|_{F}
		=O_p(\sqrt{k_0 \log (m^2 p) / T})$ \cite[Proposition 4.1]{Basu2015}, where
		$k_0=\|\mathcal{W}\|_0$.
In this case,
		it can be seen   that the error bound in Theorem \ref{statistical} is smaller than that of the Lasso estimator in  \cite{Basu2015} when $s\leq k_0$.
In general, the number of nonzero entries of the core tensor is much smaller than that of the full tensor in practice.
		\item Negahban et al. \cite{Negahban2011} proposed to use the nuclear norm (NN) to characterize the
		low-rankness of the transition matrix in the VAR model, but their analysis was restricted to the VAR$(1)$ case.
We extend the VAR$(1)$ model to the VAR$(p)$ case, yielding
	\begin{equation}\label{NN}
		\widehat{\mathcal{W}}_{\operatorname{NN}}=\arg \min \frac{1}{T} \sum_{t=1}^T\|\mathbf{y}_t-\mathcal{W}_{(1)} \mathbf{x}_t\|_2^2+\lambda_{\operatorname{NN}}\|\mathcal{W}_{(1)}\|_*.
		\end{equation}
		Furthermore, it has been shown that,
		$\|\widehat{\mathcal{W}}_{\operatorname{NN}}-\mathcal{W}\|_F=O_p(\sqrt{r_1 mp/T})$ \cite[Remark 10]{Wang2022}, where $r_1$ is the rank of mode-1 unfolding matrix of $\mathcal{W}$.
		Therefore,
		if $\frac{s}{r_1}\leq \frac{m p}{\log(m^2 p)}$,
		the estimation error bound of the nuclear norm method
		is larger than the error bound in (\ref{bound2}).
		This condition can be satisfied easily in real-world scenarios since $s$ is much smaller than $ r_1r_2r_3$,  $r_1, r_2$ are much smaller than $m$, and $r_3$ is much smaller than $p$.
		\item Recently, Wang et al. \cite{Wang2022} proposed an $\ell_1$-penalized sparse
		higher-order reduced-rank  estimator  for VAR, where the coefficient tensor $\mathcal{W}$ has the formulation of
		high-order singular value decomposition (HOSVD) and the factor matrices are sparse.
		When each column of  $\mathbf{A}_i$ has at most  $s_i$ nonzero entries,
		$i=1,2,3$,
		the resulting estimator, denoted by  $\widehat{\mathcal{W}}_{\operatorname{SHORR}}$, satisfies
		$\|\widehat{\mathcal{W}}_{\operatorname{SHORR}}-\mathcal{W}\|_{F}=O_p\big(\sqrt{(s_1 s_2 s_3)\log(m^2 p) / T}\big)$ \cite[Theorem 2]{Wang2022}.
		Therefore,
		if
		$s\leq s_1 s_2 s_3,
		$
		 the  error bound of our estimator in (\ref{bound2}) is smaller than that in \cite{Wang2022}.
		 Here the condition $s\leq s_1 s_2 s_3
		 $ can be easily satisfied in real-world applications since $s$ is much smaller than $ r_1r_2r_3$, $1\leq s_i\leq m,i=1,2, 1\leq s_3\leq p$,  $r_i$  is much smaller than $m,i=1,2$, $r_3$ is much smaller than p, and $s_i$ is not too small in general, $i=1,2,3$.
	\end{itemize}
\end{remark}

\section{Proximal Alternating Linearized Minimization Algorithm}\label{PALMG}

By letting $\mathbf{A}_i=\mathbf{U}_i$ and taking these equality constraints into the objective function, model (\ref{test3}) can be rewritten as follows:
\begin{equation}\label{PenOnbH}
\begin{aligned}
 \underset{\mathcal{G},\mathbf{A}_1,\mathbf{A}_2,\mathbf{A}_3, \atop \mathbf{U}_1,\mathbf{U}_2,\mathbf{U}_3}{\operatorname{min}}  &\frac{1}{2T} \sum_{t=1}^{T}\left\|\mathbf{y}_t-\left(\mathcal{G} \times_1 \mathbf{A}_1 \times_2 \mathbf{A}_2 \times_3 \mathbf{A}_3\right)_{(1)}  \mathbf{x}_t\right\|_2^2+\beta\|\mathcal{G}\|_1
 \\
&~~~ +\sum_{i=1}^{3}\alpha_i \operatorname{tr}\left(\mathbf{U}_i^{H} \mathbf{L}_i \mathbf{U}_i\right)
+\sum_{i=1}^{3}\frac{\gamma_i}{2}\left\|\mathbf{U}_i-\mathbf{A}_i\right\|_{F}^2\\
\text { s.t.} &  \ \|\mathcal{G}\|_{\infty}\leq c, \mathbf{A}_i^{H}\mathbf{A}_i=\mathbf{I}_{r_i}, i=1,2,3,
\end{aligned}
\end{equation}
where $\gamma_i>0$ are given constants.
By using the definition of the indicator function, model (\ref{PenOnbH}) can be reformulated
into the following unconstrained optimization problem:
\begin{equation}\label{model}
\begin{aligned}
& \underset{\mathcal{G},\mathbf{A}_1,\mathbf{A}_2,\mathbf{A}_3,\atop \mathbf{U}_1,\mathbf{U}_2,\mathbf{U}_3}{\operatorname{min}}  \frac{1}{2T} \sum_{t=1}^{T}\left\|\mathbf{y}_t-\left(\mathcal{G} \times_1 \mathbf{A}_1 \times_2 \mathbf{A}_2 \times_3 \mathbf{A}_3\right)_{(1)}  \mathbf{x}_t\right\|_2^2+\beta\|\mathcal{G}\|_1
 \\
&~~~~~~~~~~~~~~ +\sum_{i=1}^{3}\alpha_i \operatorname{tr}\left(\mathbf{U}_i^{H} \mathbf{L}_i \mathbf{U}_i\right)
+\sum_{i=1}^{3}\frac{\gamma_i}{2}\left\|\mathbf{U}_i-\mathbf{A}_i\right\|_{F}^2
+\sum_{i=1}^3 \delta_{\mathfrak{B}_i}(\mathbf{A}_i)
+\delta_{\mathcal{D}}(\mathcal{G}),
\end{aligned}
\end{equation}
where $\mathfrak{B}_i:=\{\mathbf{A}_i:\mathbf{A}_i^{H}\mathbf{A}_i=\mathbf{I}_{r_i}\}, i=1,2,3,
\mathcal{D}:=\{\mathcal{G}:\|\mathcal{G}\|_{\infty}\leq c\}$.

Now we design a
proximal alternating linearized minimization (PALM) algorithm \cite{Bolte2014}
to solve model (\ref{model}).
Before proceeding,
for simplicity, we denote
\begin{equation}\label{objecF}
	\begin{aligned}
		& F(\mathcal{G},\mathbf{A}_1,\mathbf{A}_2,\mathbf{A}_3,\mathbf{U}_1,
		\mathbf{U}_2,\mathbf{U}_3):=\frac{1}{2T} \sum_{t=1}^{T}
		\left\|\mathbf{y}_t-\left(\mathcal{G} \times_1 \mathbf{A}_1
		\times_2 \mathbf{A}_2 \times_3 \mathbf{A}_3\right)_{(1)}  \mathbf{x}_t\right\|_2^2 \\
		&~~~~+\beta\|\mathcal{G}\|_1
		+\sum_{i=1}^{3}\alpha_i \operatorname{tr}\left(\mathbf{U}_i^{H} \mathbf{L}_i
		\mathbf{U}_i\right)
		+\sum_{i=1}^{3}\frac{\gamma_i}{2}\left\|\mathbf{U}_i-\mathbf{A}_i\right\|_{F}^2
		+\sum_{i=1}^3 \delta_{\mathfrak{B}_i}(\mathbf{A}_i)
		+\delta_{\mathcal{D}}(\mathcal{G}),
	\end{aligned}
\end{equation}
\begin{equation}\label{DefQ}
	\begin{aligned}
		Q(\llbracket \mathcal{G}; \mathbf{A}_1, \mathbf{A}_2, \mathbf{A}_3\rrbracket)
		&:=\frac{1}{2T} \sum_{t=1}^{T}\left\|\mathbf{y}_t-\left(\mathcal{G} \times_1 \mathbf{A}_1 \times_2 \mathbf{A}_2 \times_3 \mathbf{A}_3\right)_{(1)}\mathbf{x}_t\right\|_2^2,
	\end{aligned}
\end{equation}
and
\begin{equation}\label{DefPsi}
\begin{aligned}
&\Psi(\mathcal{G},\mathbf{A}_1,\mathbf{A}_2,\mathbf{A}_3,\mathbf{U}_1,
\mathbf{U}_2,\mathbf{U}_3)
:=Q(\llbracket \mathcal{G}; \mathbf{A}_1, \mathbf{A}_2, \mathbf{A}_3\rrbracket)
+\sum_{i=1}^{3}\frac{\gamma_i}{2}\left\|\mathbf{U}_i-\mathbf{A}_i\right\|_{F}^2.
\end{aligned}
\end{equation}
Therefore, we get that
$$
\begin{aligned}
F(\mathcal{G},\mathbf{A}_1,\mathbf{A}_2,\mathbf{A}_3,\mathbf{U}_1,
\mathbf{U}_2,\mathbf{U}_3)& =\Psi(\mathcal{G},\mathbf{A}_1,\mathbf{A}_2,\mathbf{A}_3,\mathbf{U}_1,
\mathbf{U}_2,\mathbf{U}_3)+\beta\|\mathcal{G}\|_1
 \\
&~~~~~~ +\sum_{i=1}^{3}\alpha_i \operatorname{tr}\left(\mathbf{U}_i^{H} \mathbf{L}_i
\mathbf{U}_i\right)	+\sum_{i=1}^3 \delta_{\mathfrak{B}_i}\left(\mathbf{A}_i\right)
+\delta_{\mathcal{D}}(\mathcal{G}).
\end{aligned}
$$
Note that the gradient of the loss function $Q$ with respect to the coefficient tensor $\mathcal{W}$ is given by \cite[Section D.2]{Wang2022a}
\begin{equation}\label{gradient_w}
\nabla Q(\llbracket \mathcal{G}; \mathbf{A}_1, \mathbf{A}_2, \mathbf{A}_3\rrbracket) =\frac{1}{T}\sum_{t=1}^{T}\left(\mathcal{W}_{(1)} \mathbf{x}_t-\mathbf{y}_t\right) \circ \mathbf{X}_t = \frac{1}{T}\sum_{t=1}^{T}\left(-\boldsymbol{\varepsilon}_t\right) \circ \mathbf{X}_t,
\end{equation}
where 
$\mathbf{X}_t=(\mathbf{y}_{t-1},\mathbf{y}_{t-2},
\ldots,\mathbf{y}_{t-p})\in\mathbb{R}^{m \times p}$.
It follows from \cite[Lemma 2.1]{Han2022} that the partial gradients of  $Q$ are given as follows:
\begin{equation}\label{gradient_Q}
\left\{
\begin{aligned}
&\nabla_\mathcal{G} Q(\llbracket \mathcal{G}; \mathbf{A}_1, \mathbf{A}_2, \mathbf{A}_3\rrbracket)  =  \nabla Q(\llbracket \mathcal{G}; \mathbf{A}_1, \mathbf{A}_2, \mathbf{A}_3\rrbracket) \times_1 \mathbf{A}_1^{H} \times_2 \mathbf{A}_2^{H} \times_3 \mathbf{A}_3^{H},\\
&\nabla_{\mathbf{A}_1} Q(\llbracket \mathcal{G}; \mathbf{A}_1, \mathbf{A}_2, \mathbf{A}_3\rrbracket)  =  (\nabla Q(\llbracket \mathcal{G}; \mathbf{A}_1, \mathbf{A}_2, \mathbf{A}_3\rrbracket))_{(1)}\left(\mathbf{A}_3 \otimes \mathbf{A}_2\right) \mathcal{G}_{(1)}^{H}, \\
&\nabla_{\mathbf{A}_2} Q(\llbracket \mathcal{G}; \mathbf{A}_1, \mathbf{A}_2, \mathbf{A}_3\rrbracket)  =  (\nabla Q(\llbracket \mathcal{G}; \mathbf{A}_1, \mathbf{A}_2, \mathbf{A}_3\rrbracket))_{(2)}\left(\mathbf{A}_1 \otimes \mathbf{A}_3\right) \mathcal{G}_{(2)}^{H}, \\
&\nabla_{\mathbf{A}_3} Q(\llbracket \mathcal{G}; \mathbf{A}_1, \mathbf{A}_2, \mathbf{A}_3\rrbracket)  =  (\nabla Q(\llbracket \mathcal{G}; \mathbf{A}_1, \mathbf{A}_2, \mathbf{A}_3\rrbracket))_{(3)}\left(\mathbf{A}_2 \otimes \mathbf{A}_1\right) \mathcal{G}_{(3)}^{H}.
\end{aligned}
\right.
\end{equation}
Next, we can deduce the gradient of  $\Psi$ (defined in (\ref{DefPsi})) with respect to $\mathcal{G},\mathbf{A}_i,\mathbf{U}_i$ as follows:
\begin{equation}\label{gradient_H}
\left\{
\begin{aligned}
& \nabla_\mathcal{G} \Psi(\mathcal{G},\mathbf{A}_1,\mathbf{A}_2,\mathbf{A}_3,\mathbf{U}_1,
\mathbf{U}_2,\mathbf{U}_3) =  \nabla_\mathcal{G} Q(\llbracket \mathcal{G}; \mathbf{A}_1, \mathbf{A}_2, \mathbf{A}_3\rrbracket),\\
&\nabla_{\mathbf{A}_1} \Psi(\mathcal{G},\mathbf{A}_1,\mathbf{A}_2,\mathbf{A}_3,\mathbf{U}_1,
\mathbf{U}_2,\mathbf{U}_3)  =  \nabla_{\mathbf{A}_1} Q(\llbracket \mathcal{G}; \mathbf{A}_1, \mathbf{A}_2, \mathbf{A}_3\rrbracket) - \gamma_1(\mathbf{U}_1-\mathbf{A}_1), \\
&\nabla_{\mathbf{A}_2} \Psi(\mathcal{G},\mathbf{A}_1,\mathbf{A}_2,\mathbf{A}_3,\mathbf{U}_1,
\mathbf{U}_2,\mathbf{U}_3)  =  \nabla_{\mathbf{A}_2} Q(\llbracket \mathcal{G}; \mathbf{A}_1, \mathbf{A}_2, \mathbf{A}_3\rrbracket)-\gamma_2(\mathbf{U}_2-\mathbf{A}_2), \\
&\nabla_{\mathbf{A}_3} \Psi(\mathcal{G},\mathbf{A}_1,\mathbf{A}_2,\mathbf{A}_3,\mathbf{U}_1,
\mathbf{U}_2,\mathbf{U}_3)  =  \nabla_{\mathbf{A}_3} Q(\llbracket \mathcal{G}; \mathbf{A}_1, \mathbf{A}_2, \mathbf{A}_3\rrbracket)-\gamma_3(\mathbf{U}_3-\mathbf{A}_3),\\
&\nabla_{\mathbf{U}_1} \Psi(\mathcal{G},\mathbf{A}_1,\mathbf{A}_2,\mathbf{A}_3,\mathbf{U}_1,
\mathbf{U}_2,\mathbf{U}_3)  = \gamma_1(\mathbf{U}_1-\mathbf{A}_1),\\
&\nabla_{\mathbf{U}_2} \Psi(\mathcal{G},\mathbf{A}_1,\mathbf{A}_2,\mathbf{A}_3,\mathbf{U}_1,
\mathbf{U}_2,\mathbf{U}_3)  = \gamma_2(\mathbf{U}_2-\mathbf{A}_2),\\
& \nabla_{\mathbf{U}_3} \Psi(\mathcal{G},\mathbf{A}_1,\mathbf{A}_2,\mathbf{A}_3,\mathbf{U}_1,
\mathbf{U}_2,\mathbf{U}_3)  = \gamma_3(\mathbf{U}_3-\mathbf{A}_3).
\end{aligned}
\right.
\end{equation}
Then the iterative framework of the PALM algorithm is given as follows:
\begin{equation}\label{palm}
\left\{
\begin{aligned}
\mathcal{G}^{k+1}\in \underset{\mathcal{G}}{\arg \min }\ & \left\langle \nabla_{\mathcal{G}} \Psi(\mathcal{G}^{k}, \mathbf{A}_1^{k}, \mathbf{A}_2^{k}, \mathbf{A}_3^{k},\mathbf{U}_1^{k},
\mathbf{U}_2^{k},\mathbf{U}_3^{k}),\mathcal{G}- \mathcal{G}^k\right\rangle+ \beta\|\mathcal{G}\|_1\\
&+\delta_{\mathcal{D}}(\mathcal{G}) +
\frac{\rho_1}{2}\|\mathcal{G}- \mathcal{G}^k\|_F^2,\\
\mathbf{A}_1^{k+1} \in \underset{\mathbf{A}_1}{\arg\min }\ &\left\langle \nabla_{\mathbf{A}_1} \Psi(\mathcal{G}^{k+1}, \mathbf{A}_1^{k}, \mathbf{A}_2^{k}, \mathbf{A}_3^{k},\mathbf{U}_1^{k},
\mathbf{U}_2^{k},\mathbf{U}_3^{k}),\mathbf{A}_1-\mathbf{A}_1^{k}\right\rangle\\
& +\delta_{\mathfrak{B}_1}(\mathbf{A}_1)+\frac{\rho_2}{2}\|\mathbf{A}_1-\mathbf{A}_1^{k}\|_F^2, \\
\mathbf{A}_2^{k+1} \in \underset{\mathbf{A}_2}{\arg\min }\ &\left\langle \nabla_{\mathbf{A}_2} \Psi(\mathcal{G}^{k+1}, \mathbf{A}_1^{k+1}, \mathbf{A}_2^{k}, \mathbf{A}_3^{k},\mathbf{U}_1^{k},
\mathbf{U}_2^{k},\mathbf{U}_3^{k}),\mathbf{A}_2-\mathbf{A}_2^{k}\right\rangle \\
&+\delta_{\mathfrak{B}_2}(\mathbf{A}_2)+  \frac{\rho_3}{2}\|\mathbf{A}_2-\mathbf{A}_2^{k}\|_F^2, \\
\mathbf{A}_3^{k+1} \in \underset{\mathbf{A}_3}{\arg\min }\
&\left\langle \nabla_{\mathbf{A}_3} \Psi(\mathcal{G}^{k+1}, \mathbf{A}_1^{k+1}, \mathbf{A}_2^{k+1}, \mathbf{A}_3^{k},\mathbf{U}_1^{k},
\mathbf{U}_2^{k},\mathbf{U}_3^{k}),\mathbf{A}_3-\mathbf{A}_3^{k}\right\rangle
\\
&+\delta_{\mathfrak{B}_3}(\mathbf{A}_3) +\frac{\rho_4}{2}\|\mathbf{A}_3-\mathbf{A}_3^{k}\|_F^2, \\
\mathbf{U}_1^{k+1}\in \underset{\mathbf{U}_1}{\arg\min }\ &\left\langle \nabla_{\mathbf{U}_1} \Psi(\mathcal{G}^{k+1}, \mathbf{A}_1^{k+1}, \mathbf{A}_2^{k+1}, \mathbf{A}_3^{k+1},\mathbf{U}_1^{k},
\mathbf{U}_2^{k},\mathbf{U}_3^{k}),\mathbf{U}_1-\mathbf{U}_1^{k}\right\rangle\\
&+\alpha_1 \operatorname{tr}(\mathbf{U}_1^{H} \mathbf{L}_1 \mathbf{U}_1) + \frac{\rho_5}{2}\|\mathbf{U}_1-\mathbf{U}_1^{k}\|_F^2 , \\
\mathbf{U}_2^{k+1}  \in \underset{\mathbf{U}_2}{\arg\min }\
&\left\langle \nabla_{\mathbf{U}_2} \Psi(\mathcal{G}^{k+1}, \mathbf{A}_1^{k+1}, \mathbf{A}_2^{k+1}, \mathbf{A}_3^{k+1},\mathbf{U}_1^{k+1},
\mathbf{U}_2^{k},\mathbf{U}_3^{k}),\mathbf{U}_2-\mathbf{U}_2^{k}\right\rangle
\\
&\alpha_2 \operatorname{tr}(\mathbf{U}_2^{H} \mathbf{L}_2 \mathbf{U}_2) +
\frac{\rho_6}{2}\|\mathbf{U}_2-\mathbf{U}_2^{k}\|_F^2 , \\
\mathbf{U}_3^{k+1}  \in \underset{\mathbf{U}_3}{\arg\min }\
&\left\langle \nabla_{\mathbf{U}_3} \Psi(\mathcal{G}^{k+1}, \mathbf{A}_1^{k+1}, \mathbf{A}_2^{k+1}, \mathbf{A}_3^{k+1},\mathbf{U}_1^{k+1},
\mathbf{U}_2^{k+1},\mathbf{U}_3^{k}),\mathbf{U}_3-\mathbf{U}_3^{k}\right\rangle\\
&+\alpha_3 \operatorname{tr}(\mathbf{U}_3^{H} \mathbf{L}_3 \mathbf{U}_3)
+\frac{\rho_7}{2}\|\mathbf{U}_3-\mathbf{U}_3^{k}\|_F^2,
\end{aligned}
\right.
\end{equation}
where $\rho_i>0$ are given constants, $i=1,\ldots,7$.
Now we give the detailed  solution for each subproblem.
The subproblem about $\mathcal{G}$ in  (\ref{palm}) can be expressed as
$$
\begin{aligned}
\mathcal{G}^{k+1}  = &\ \mathop{\arg\min}\limits_{\|\mathcal{G}\|_{\infty}\leq c}
\beta\|\mathcal{G}\|_1
 +\left\langle\nabla_{\mathcal{G}}\Psi(\mathcal{G}^{k}, \mathbf{A}_1^{k}, \mathbf{A}_2^{k}, \mathbf{A}_3^{k},\mathbf{U}_1^{k},
\mathbf{U}_2^{k},\mathbf{U}_3^{k}) ,  \mathcal{G}-\mathcal{G}^k\right\rangle+\frac{\rho_1}{2}
 \|\mathcal{G}-\mathcal{G}^k\|_F^2 \\
=&\ \mathop{\arg\min}\limits_{\|\mathcal{G}\|_{\infty}\leq c}
\beta\|\mathcal{G}\|_1+\frac{\rho_1}{2}\left\|\mathcal{G}-\mathcal{G}^k+\frac{1}{\rho_1} \nabla_{\mathcal{G}} \Psi(\mathcal{G}^{k}, \mathbf{A}_1^{k}, \mathbf{A}_2^{k}, \mathbf{A}_3^{k},\mathbf{U}_1^{k},
\mathbf{U}_2^{k},\mathbf{U}_3^{k})\right\|_F^2 \\
=&\ \mathop{\arg\min}\limits_{\|\mathcal{G}\|_{\infty}\leq c} \beta\|\mathcal{G}\|_1+\frac{\rho_1}{2}\|\mathcal{G}-\mathcal{L}\|_F^2,
\end{aligned}
$$
where $\mathcal{L}=\mathcal{G}^k-\frac{1}{\rho_1} \nabla_{\mathcal{G}} \Psi(\mathcal{G}^{k}, \mathbf{A}_1^{k}, \mathbf{A}_2^{k}, \mathbf{A}_3^{k},\mathbf{U}_1^{k},
\mathbf{U}_2^{k},\mathbf{U}_3^{k})$. Simple calculations show that the closed form solution with respect to $\mathcal{G}$ can be given by \cite[Example 6.34]{Beck2017}
\begin{equation}\label{solutionG}
\mathcal{G}_{i j l}^{k+1}=\left\{\begin{array}{cl}
\operatorname{sign}\left(\mathcal{L}_{i j l}\right) \max \left\{\left|\mathcal{L}_{i j l}\right|-\beta/\rho_1, 0\right\}, & \ \textup{if} \  \left|\mathcal{L}_{i j l}\right| \leq c+\beta/\rho_1, \\
\operatorname{sign}\left(\mathcal{L}_{i j l}\right) c, & \ \textup{if} \ \left|\mathcal{L}_{i j l}\right|>c+\beta/\rho_1,
\end{array}\right.
\end{equation}
where $\operatorname{sign}(\cdot)$ denotes the signum function.

The subproblem with respect to $\mathbf{A}_1$ can be equivalently expressed as
$$
\begin{aligned}
\mathbf{A}_1^{k+1} = &\ \mathop{\arg\min}\limits_{\mathbf{A}_1}
\delta_{\mathfrak{B}_1}(\mathbf{A}_1)+\left\langle\nabla_{\mathbf{A}_1} \Psi(\mathcal{G}^{k+1}, \mathbf{A}_1^{k}, \mathbf{A}_2^{k}, \mathbf{A}_3^{k},\mathbf{U}_1^{k},
\mathbf{U}_2^{k},\mathbf{U}_3^{k}), \mathbf{A}_1-\mathbf{A}_1^k\right\rangle\\
&~~~~~~~+\frac{\rho_2}{2}\|\mathbf{A}_1-\mathbf{A}_1^k\|_F^2 \\
= &\ \mathop{\arg\min}\limits_{\mathbf{A}_1} \delta_{\mathfrak{B}_1}(\mathbf{A}_1)
+\frac{\rho_2}{2}\left\|\mathbf{A}_1-\left(\mathbf{A}_1^k-\frac{1}{\rho_2} \nabla_{\mathbf{A}_1} \Psi(\mathcal{G}^{k+1}, \mathbf{A}_1^{k}, \mathbf{A}_2^{k}, \mathbf{A}_3^{k},\mathbf{U}_1^{k},
\mathbf{U}_2^{k},\mathbf{U}_3^{k})\right)\right\|_F^2.
\end{aligned}
$$
It follows from \cite[Theorem 4]{Zou2006} that the optimal solution of the above problem is easily
computed as
\begin{equation}\label{solutionA1}
\mathbf{A}_1^{k+1}=\bar{\mathbf{U}}_1 \bar{\mathbf{V}}_1^H.
\end{equation}
Here $\mathbf{A}_1^k-\frac{1}{\rho_2} \nabla_{\mathbf{A}_1} \Psi(\mathcal{G}^{k+1}, \mathbf{A}_1^{k}, \mathbf{A}_2^{k}, \mathbf{A}_3^{k},\mathbf{U}_1^{k},
\mathbf{U}_2^{k},\mathbf{U}_3^{k})=\bar{\mathbf{U}}_1 \Sigma_1\bar{\mathbf{V}}_1^H$.

Similar techniques are employed to solve the subproblems with respect to $\mathbf{A}_2$
and $\mathbf{A}_3$.
By performing SVDs on the two matrices $\mathbf{A}_2^k-
\frac{1}{\rho_3} \nabla_{\mathbf{A}_2} \Psi(\mathcal{G}^{k+1}, \mathbf{A}_1^{k+1}, \mathbf{A}_2^{k}, \mathbf{A}_3^{k},\mathbf{U}_1^{k},
\mathbf{U}_2^{k},\mathbf{U}_3^{k})$ and $\mathbf{A}_3^k-
\frac{1}{\rho_4} \nabla_{\mathbf{A}_3} \Psi(\mathcal{G}^{k+1}, \mathbf{A}_1^{k+1}, \mathbf{A}_2^{k+1}, \mathbf{A}_3^{k},\mathbf{U}_1^{k},
\mathbf{U}_2^{k},\mathbf{U}_3^{k})$, i.e.,
 $$\mathbf{A}_2^k-\frac{1}{\rho_3}
\nabla_{\mathbf{A}_2} \Psi(\mathcal{G}^{k+1}, \mathbf{A}_1^{k+1}, \mathbf{A}_2^{k}, \mathbf{A}_3^{k},\mathbf{U}_1^{k},
\mathbf{U}_2^{k},\mathbf{U}_3^{k})=\bar{\mathbf{U}}_2 \Sigma_2\bar{\mathbf{V}}_2^H
$$
and
$$
\mathbf{A}_3^k-\frac{1}{\rho_4} \nabla_{\mathbf{A}_3} \Psi(\mathcal{G}^{k+1},\mathbf{A}_1^{k+1}, \mathbf{A}_2^{k+1}, \mathbf{A}_3^{k},\mathbf{U}_1^{k},
\mathbf{U}_2^{k},\mathbf{U}_3^{k})=
\bar{\mathbf{U}}_3 \Sigma_3\bar{\mathbf{V}}_3^H,
$$
we get that the closed form solutions with respect to $\mathbf{A}_2$ and $\mathbf{A}_3$ in  (\ref{palm}) are given by
\begin{equation}\label{solutionA2}
\mathbf{A}_2^{k+1}=\bar{\mathbf{U}}_2 \bar{\mathbf{V}}_2^H, \ \
\mathbf{A}_3^{k+1}=\bar{\mathbf{U}}_3 \bar{\mathbf{V}}_3^H.
\end{equation}



The subproblem about $\mathbf{U}_1$ in  (\ref{palm})  can be rewritten equivalently as
\begin{equation}\label{Uk1Op}
\begin{aligned}
\mathbf{U}_1^{k+1}
=&\ \mathop{\arg\min}\limits_{\mathbf{U}_1}\alpha_1
\operatorname{tr}(\mathbf{U}_1^{H} \mathbf{L}_1 \mathbf{U}_1)\\
&+\frac{\rho_5}{2}\left\|\mathbf{U}_1-\mathbf{U}_1^{k}
+\frac{1}{\rho_5}\nabla_{\mathbf{U}_1} \Psi(\mathcal{G}^{k+1}, \mathbf{A}_1^{k+1}, \mathbf{A}_2^{k+1}, \mathbf{A}_3^{k+1},\mathbf{U}_1^{k},
\mathbf{U}_2^{k},\mathbf{U}_3^{k})\right\|_F^2.
\end{aligned}
\end{equation}
By the optimality condition of (\ref{Uk1Op}),
 the closed form solution with respect to $\mathbf{U}_1$ can be given by
\begin{equation}\label{solutionU1}
\begin{aligned}
\mathbf{U}_1^{k+1} =(2 \alpha_1\mathbf{L}_1 +\rho_5 \mathbf{I}_m)^{-1}(\rho_5 \mathbf{U}_1^{k}-\nabla_{\mathbf{U}_1} \Psi(\mathcal{G}^{k+1}, \mathbf{A}_1^{k+1}, \mathbf{A}_2^{k+1}, \mathbf{A}_3^{k+1},\mathbf{U}_1^{k},
\mathbf{U}_2^{k},\mathbf{U}_3^{k})).
\end{aligned}
\end{equation}
Similar to the calculation of $\mathbf{U}_1^{k+1}$,
the optimal solution with respect to $\mathbf{U}_2$ and $\mathbf{U}_3$ are given explicitly by
\begin{equation}\label{solutionU2}
\begin{aligned}
&\mathbf{U}_2^{k+1} =(2 \alpha_2\mathbf{L}_2 +\rho_6 \mathbf{I}_m)^{-1}
(\rho_6 \mathbf{U}_2^{k}-\nabla_{\mathbf{U}_2} \Psi(\mathcal{G}^{k+1}, \mathbf{A}_1^{k+1}, \mathbf{A}_2^{k+1}, \mathbf{A}_3^{k+1},\mathbf{U}_1^{k+1},
\mathbf{U}_2^{k},\mathbf{U}_3^{k})),
\end{aligned}
\end{equation}
\begin{equation}\label{solutionU3}
\begin{aligned}
&\mathbf{U}_3^{k+1} =(2 \alpha_3\mathbf{L}_3 +\rho_7
\mathbf{I}_p)^{-1}(\rho_7 \mathbf{U}_3^{k}-\nabla_{\mathbf{U}_3} \Psi(\mathcal{G}^{k+1}, \mathbf{A}_1^{k+1}, \mathbf{A}_2^{k+1}, \mathbf{A}_3^{k+1},\mathbf{U}_1^{k+1},
\mathbf{U}_2^{k+1},\mathbf{U}_3^{k})).
\end{aligned}
\end{equation}

Now we summarize the PALM algorithm  for solving  problem (\ref{model}) in Algorithm \ref{alg:algorithm1}.

\begin{algorithm}[htbp]
    \caption{A PALM Algorithm for Solving Problem (\ref{model})}
    \label{alg:algorithm1}
    \begin{algorithmic}[1]
        \State \textbf{Input:} $\{\mathbf{y}_t\},\mathcal{G}^0,\mathbf{A}_1^0,\mathbf{A}_2^0,
        \mathbf{A}_3^0,\mathbf{U}_1^0,\mathbf{U}_2^0,\mathbf{U}_3^0$,
         parameters
         $\beta, \alpha_1,\alpha_2,\alpha_3,\gamma_1,\gamma_2,\gamma_3, \rho_i, i=1,\ldots, 7$.
        \Repeat
  		\State \textbf{Step 1}. Compute ${\mathcal{G}}^{k+1}$ by (\ref{solutionG}).
  		\State \textbf{Step 2}. Compute $\mathbf{A}_i^{k+1}$ by (\ref{solutionA1}) and (\ref{solutionA2}).
  		\State \textbf{Step 3}. Compute $\mathbf{U}_i^{k+1}$ by (\ref{solutionU1}), (\ref{solutionU2}), and (\ref{solutionU3}).
        \Until A stopping condition is satisfied.
        \State \textbf{Output:} ${\mathcal{G}}^{k+1}\times_1\mathbf{A}_1^{k+1}
        \times_2\mathbf{A}_2^{k+1}\times_3\mathbf{A}_3^{k+1}$.
	\end{algorithmic}
\end{algorithm}

The computational complexity of Algorithm \ref{alg:algorithm1} is given as follows.
The main computational complexity of
${\mathcal{G}}^{k+1}$ is that of computing $\nabla_{\mathcal{G}} \Psi$,
which is on the order of $O(m^2pT)$.
The computational cost of $\mathbf{A}_1^{k+1}$ is $O(m^2pT+r_1m^2p+r_1r_2r_3mp)$.
And the computational costs of $\mathbf{A}_2^{k+1}$ and $\mathbf{A}_3^{k+1}$ are $O(m^2pT+r_2m^2p+r_1r_2r_3mp)$
and $O(m^2pT+r_3m^2p+r_1r_2r_3m^2)$, respectively.
The computational complexities of $\mathbf{U}_1^{k+1}$ and $\mathbf{U}_2^{k+1}$ are both $O(m^3)$, and
the computational complexity of $\mathbf{U}_3^{k+1}$ is $O(p^3)$.
Therefore, the computational complexity of Algorithm \ref{alg:algorithm1} is $O((r_1+r_2+r_3+T)m^2p+r_1r_2r_3m(m+p)+m^3+p^3)$  at each iteration.

\subsection{Convergence Analysis}
In this subsection, we establish the global convergence  of Algorithm \ref{alg:algorithm1}.
First we give the sufficient decrease property of $F$ in the following lemma.

\begin{lemma}\label{suffcond}(Sufficient decrease)
Let the sequence $\{(\mathcal{G}^k,\mathbf{A}_1^k,\mathbf{A}_2^k,\mathbf{A}_3^k,\mathbf{U}_1^k,
\mathbf{U}_2^k,\mathbf{U}_3^k)\}$ be generated by Algorithm \ref{alg:algorithm1}.
Suppose that $\rho_1> L_1,\rho_2> L_2,\rho_3> L_3,
\rho_4>L_4,\rho_5> \gamma_1,\rho_6>\gamma_2,\rho_7>\gamma_3$.
Then there exists a constant $\bar{\rho}>0$ such that the following inequality holds
\begin{equation}\label{suffDEC}
\begin{aligned}
&F\left(\mathcal{G}^{k+1}, \mathbf{A}_1^{k+1}, \mathbf{A}_2^{k+1}, \mathbf{A}_3^{k+1},\mathbf{U}_1^{k+1},
\mathbf{U}_2^{k+1},\mathbf{U}_3^{k+1}\right)
+\frac{\bar{\rho}}{2}\left(\|\mathcal{G}^{k+1}- \mathcal{G}^k\|_F^2+
\|\mathbf{A}_1^{k+1}-\mathbf{A}_1^{k}\|_F^2
\right.\\
+&\left.
\|\mathbf{A}_2^{k+1}-\mathbf{A}_2^{k}\|_F^2
+\|\mathbf{A}_3^{k+1}-\mathbf{A}_3^{k}\|_F^2
+\|\mathbf{U}_1^{k+1}- \mathbf{U}_1^{k}\|_F^2
+\|\mathbf{U}_2^{k+1}- \mathbf{U}_2^{k}\|_F^2
+\|\mathbf{U}_3^{k+1}- \mathbf{U}_3^{k}\|_F^2\right)
\\
\leq \ &
F\left(\mathcal{G}^{k}, \mathbf{A}_1^{k}, \mathbf{A}_2^{k}, \mathbf{A}_3^{k},\mathbf{U}_1^{k},
\mathbf{U}_2^{k},\mathbf{U}_3^{k}\right),
\end{aligned}
\end{equation}
where $F$ is defined in (\ref{objecF}) and $\bar{\rho}:=\min\{\rho_1-L_1,\rho_2-L_2,\rho_3-L_3,
\rho_4-L_4,\rho_5-\gamma_1,\rho_6-\gamma_2,\rho_7-\gamma_3\}$.
\end{lemma}


Next we show that the relative error condition holds in the following lemma.

\begin{lemma}\label{Relacond}(Relative error condition)
Let the sequence $\{(\mathcal{G}^k,\mathbf{A}_1^k,\mathbf{A}_2^k,\mathbf{A}_3^k,\mathbf{U}_1^k,
\mathbf{U}_2^k,\mathbf{U}_3^k)\}$ be generated by Algorithm \ref{alg:algorithm1}. Then there exist $\mathcal{N}^{k+1}\in
\partial F(\mathcal{G}^{k+1},\mathbf{A}_1^{k+1},\mathbf{A}_2^{k+1},\mathbf{A}_3^{k+1},
\mathbf{U}_1^{k+1},
\mathbf{U}_2^{k+1},\mathbf{U}_3^{k+1})$  and
a constant $\bar{\delta}>0$ such that
\begin{equation}\label{Rela}
\begin{aligned}
\|\mathcal{N}^{k+1}\|_F
\leq \bar{\delta}\left(\|\mathcal{G}^{k+1}- \mathcal{G}^k\|_F+
\|\mathbf{A}_1^{k+1}-\mathbf{A}_1^{k}\|_F
+\|\mathbf{A}_2^{k+1}-\mathbf{A}_2^{k}\|_F
+\|\mathbf{A}_3^{k+1}-\mathbf{A}_3^{k}\|_F
\right.\\
+\left.\|\mathbf{U}_1^{k+1}- \mathbf{U}_1^{k}\|_F
+\|\mathbf{U}_2^{k+1}- \mathbf{U}_2^{k}\|_F
+\|\mathbf{U}_3^{k+1}- \mathbf{U}_3^{k}\|_F\right).
\end{aligned}
\end{equation}
\end{lemma}

Now, by combining Lemmas \ref{suffcond} and \ref{Relacond},
we establish the global convergence of Algorithm \ref{alg:algorithm1} in the following theorem.

\begin{theorem}\label{ConvRe}
Let the sequence $\{(\mathcal{G}^k,\mathbf{A}_1^k,\mathbf{A}_2^k,\mathbf{A}_3^k,\mathbf{U}_1^k,
\mathbf{U}_2^k,\mathbf{U}_3^k)\}$ be  generated by Algorithm \ref{alg:algorithm1}.
Suppose that $\rho_1> L_1,\rho_2> L_2,\rho_3> L_3,
\rho_4>L_4,\rho_5> \gamma_1,\rho_6>\gamma_2,\rho_7>\gamma_3$.
Then the sequence  $\{(\mathcal{G}^k,\mathbf{A}_1^k,\mathbf{A}_2^k,\mathbf{A}_3^k,\mathbf{U}_1^k,
\mathbf{U}_2^k,\mathbf{U}_3^k)\}$ converges to a critical point of (\ref{model}).
\end{theorem}
\textbf{Proof}.
It can be easily seen from (\ref{test1}) that
$\operatorname{tr}(\mathbf{U}_i^{H}\mathbf{L}_i \mathbf{U}_i)\geq0$
for $i=1,2,3$,
which in conjunction with the definition of $F(\mathcal{G},\mathbf{A}_1,\mathbf{A}_2,\mathbf{A}_3,$
$\mathbf{U}_1,
\mathbf{U}_2,\mathbf{U}_3)$ in (\ref{objecF})
immediately establishes
$F\geq0$.
In addition, note
that $F(\mathcal{G},\mathbf{A}_1,\mathbf{A}_2,\mathbf{A}_3,$
$\mathbf{U}_1,
\mathbf{U}_2,\mathbf{U}_3)$ tends to
infinity as $\|\mathcal{G}\|_F$, $\|\mathbf{A}_i\|_F$ or $\|\mathbf{U}_i\|_F$
tends to infinity, where $i=1,2,3$.
Therefore, we obtain that
$F$ is coercive. Note that $F$ is  the monotonically decreasing by Lemma \ref{suffcond}, we can get that
the sequence  $\{(\mathcal{G}^k,\mathbf{A}_1^k,\mathbf{A}_2^k,\mathbf{A}_3^k,
\mathbf{U}_1^k,\mathbf{U}_2^k,\mathbf{U}_3^k)\}$ is bounded.
Then there exists a subsequence $\{(\mathcal{G}^{k_j},\mathbf{A}_1^{k_j},\mathbf{A}_2^{k_j},\mathbf{A}_3^{k_j},
\mathbf{U}_1^{k_j},\mathbf{U}_2^{k_j},\mathbf{U}_3^{k_j})\}$ such that the sequence $\{(\mathcal{G}^{k_j},\mathbf{A}_1^{k_j},\mathbf{A}_2^{k_j},\mathbf{A}_3^{k_j},
\mathbf{U}_1^{k_j},\mathbf{U}_2^{k_j},\mathbf{U}_3^{k_j})\}$ converges to $\{\mathcal{G}^{*},\mathbf{A}_1^{*},\mathbf{A}_2^{*},\mathbf{A}_3^{*},
\mathbf{U}_1^{*},\mathbf{U}_2^{*},\mathbf{U}_3^{*}\}$ as $j$ tends to infinity.
It is worth noting that $\delta_{\mathfrak{B}_i}(\mathbf{A}_i), i=1,2,3$,
and $\delta_{\mathcal{D}}(\mathcal{G})$ in $F$ are discontinuous, whereas all other functions in $F$ are continuous. Moreover, $\mathbf{A}_i^k$ and $\mathcal{G}^k$ satisfy $(\mathbf{A}_i^k)^{H}{\mathbf{A}_i^k}=\mathbf{I}_{r_i}, i=1,2,3,$ and $\|\mathcal{G}^k\|_{\infty}\leq c$ in each iteration,
which result in $(\mathbf{A}_i^*)^{H}{\mathbf{A}_i^*}=\mathbf{I}_{r_i}, i=1,2,3, \|\mathcal{G}^*\|_{\infty}\leq c$.
Hence,
$\delta_{\mathfrak{B}_i}(\mathbf{A}_i^{k_j})$ tends to
$\delta_{\mathfrak{B}_i}(\mathbf{A}_i^{*})$
and $\delta_{\mathcal{D}}(\mathcal{G}^{k_j})$
tends to
$\delta_{\mathcal{D}}(\mathcal{G}^{*})$
as $j$ tends to infinity.
As a result,
$F(\mathcal{G}^{k_j},\mathbf{A}_1^{k_j},\mathbf{A}_2^{k_j},
\mathbf{A}_3^{k_j},\mathbf{U}_1^{k_j},
\mathbf{U}_2^{k_j},\mathbf{U}_3^{k_j})$
tends to
$F(\mathcal{G}^{*},\mathbf{A}_1^{*},\mathbf{A}_2^{*},
\mathbf{A}_3^{*},\mathbf{U}_1^{*},
\mathbf{U}_2^{*},\mathbf{U}_3^{*})$
as $j$ tends to  infinity.

Notice that
$\Psi(\mathcal{G},\mathbf{A}_1,\mathbf{A}_2,\mathbf{A}_3,\mathbf{U}_1,
\mathbf{U}_2,\mathbf{U}_3)$ defined in (\ref{DefPsi})
and $\delta_{\mathcal{D}}(\mathcal{G})$ are Kurdyka-{\L}ojasiewicz (KL) function \cite[Example 2]{Bolte2014}.
For any $i=1,2,3$,
$\delta_{\mathfrak{B}_i}(\mathbf{A}_i)$
is a KL function since $\mathfrak{B}_i$ is a semi-algebra set \cite{Qiu2024}.
In addition,
the function $\operatorname{tr}(\mathbf{U}_i^{H} \mathbf{L}_i \mathbf{U}_i)$
is also a KL function since it is a real analytic function \cite[Example 2]{Bolte2014}.
Then one can deduce that
$F(\mathcal{G},\mathbf{A}_1,\mathbf{A}_2,\mathbf{A}_3,\mathbf{U}_1,
\mathbf{U}_2,\mathbf{U}_3)$
is a KL function.
Combining with Lemma \ref{suffcond}, Lemma \ref{Relacond}
and \cite[Theorem 3.9]{Attouch2013}, we can establish
that
$\{\mathcal{G}^k,\mathbf{A}_1^k,\mathbf{A}_2^k,\mathbf{A}_3^k,\mathbf{U}_1^k,
\mathbf{U}_2^k,\mathbf{U}_3^k\}$
converges to a critical point
of $F(\mathcal{G},\mathbf{A}_1,\mathbf{A}_2,\mathbf{A}_3,\mathbf{U}_1,
\mathbf{U}_2,\mathbf{U}_3)$.
This completes the proof.
\qed

\begin{remark}\label{Condirho}
The assumptions in Theorem 	\ref{ConvRe} are very weak. We only need $\rho_1> L_1,\rho_2> L_2,\rho_3> L_3,
\rho_4>L_4,\rho_5> \gamma_1,\rho_6>\gamma_2,\rho_7>\gamma_3$. In the experiments, we will give  $\rho_i$ in detail, $i=1,\ldots, 7$.
\end{remark}

\section{Numerical Experiments}\label{ExperM}

In this section, some numerical experiments are conducted to demonstrate the effectiveness of the proposed STDGR model for high-dimensional time series forecasting. 
We compare the STDGR with the following five models: the Lasso type method for time series models (Lasso) \cite{Basu2015},
nuclear norm minimization of the transition matrix  (NNM) \cite{Negahban2011},  low-rank tensor learning via Tucker rank constraint (LRTLT) \cite{bahadori2014fast}, time-varying autoregression with
low-rank tensors (TVART) \cite{harris2021time}, sparse higher-order reduced-rank VAR (SHORR) \cite{Wang2022}.
All  experiments are performed in MATLAB
R2020b with a 12-core, 2.6 GHz Intel Core CPU and 16 GB of RAM.

The mean squared error (MSE) is used to evaluate the forecasting accuracy for high-dimensional time series data,
which is
defined as
$$
\text{MSE}:=\frac{1}{mT_0} \sum_{t=1}^{T_0}\|\tilde{\mathbf{y}}_t-\hat{\mathbf{y}}_t\|_2^2,
$$
where $T_0$ is the number of testing samples, $\tilde{\mathbf{y}}_t\in\mathbb{R}^m$ denotes the ground truth value and $\hat{\mathbf{y}}_t\in\mathbb{R}^m$ denotes the forecasting value.

For any $i\in\{1,2,3\}$, denote
$$
\Lambda_1=\frac{\|{\mathcal{G}}^{k+1}-{\mathcal{G}}^{k}\|_{F}}{\|{\mathcal{G}}^{k}\|_{F}},
\enspace
\Lambda_{i+1}=\frac{\|\mathbf{A}_i^{k+1}-\mathbf{A}_i^{k}\|_{F}}{\|\mathbf{A}_i^{k}\|_{F}},
\enspace
\Lambda_{i+4}=\frac{\|\mathbf{U}_i^{k+1}-\mathbf{U}_i^{k}\|_{F}}{\|\mathbf{U}_i^{k}\|_{F}}.
$$
Algorithm \ref{alg:algorithm1} will be terminated if
$\max_{1\leq j\leq 7}\{\Lambda_j\}\leq 3\times 10^{-3}$
or the number of iterations reaches $200$.

\subsection{The Selection of Tucker Rank}

In this subsection, we give the choice of the Tucker rank in  Algorithm \ref{alg:algorithm1}.
The NNM in (\ref{NN}) is used to get an initial estimator of our model,
where the mode-1 unfolding of the parameter tensor is used in the NNM model.
The estimator in (\ref{NN}) is a consistent initial estimator \cite{Negahban2011, Wang2022}.
We apply the ridge-type ratio estimator \cite[Section 5]{Wang2022} to estimate the Tucker rank in Algorithm \ref{alg:algorithm1},
which is defined as follows:
\begin{equation}\label{rank-selec}
	\widehat{r}_i=\arg \min _{1 \leq j \leq n_i-1} \frac{\sigma_{j+1}((\widehat{\mathcal{W}}_{\operatorname{NN}})_{(i)})+\bar{c}}
	{\sigma_j((\widehat{\mathcal{W}}_{\operatorname{NN}})_{(i)})+\bar{c}},
\end{equation}
where $i=1,2,3$, $n_1=n_2=m, n_3=p$, and  $\bar{c}=o_p(\sqrt{r_1mp/T})$.
It is shown in \cite[Theorem 3]{Wang2022} that the rank selection in (\ref{rank-selec}) is consistent, i.e., $\mathbb{P}(	\widehat{r}_1=r_1,	\widehat{r}_2=r_2,	\widehat{r}_3=r_3)\rightarrow1$ as $T\rightarrow \infty$.
In the experiments, we set 
$\bar{c}=\sqrt{mp\log (T) / (50 T)}$.

\subsection{The Construction of Laplacian Matrices}

The Laplacian matrices $\mathbf{L}_i$ in model  (\ref{test3}), $i=1,2,3$, are constructed as follows: 
First, we use the NNM  model in (\ref{NN}) to get an initial estimator $\widehat{\mathcal{W}}_{\operatorname{NN}}$ and then take the HOSVD \cite{DeLathauwer2000} to obtain the factor matrices $\widehat{\mathbf{A}}_i$ with rank $\widehat{r}_i$, $i=1,2,3$, where $\widehat{r}_i$ is computed via (\ref{rank-selec}). 
Let $(\widehat{\mathbf{a}}_i)_l$ denote the $l$-th row of $\widehat{\mathbf{A}}_i$.
Then the weight $(z_i)_{lt}$  in (\ref{test1}) 
between $(\widehat{\mathbf{a}}_i)_l$ and $(\widehat{\mathbf{a}}_i)_t$ is given  by $(z_i)_{l t}=e^ {-{d_{l t}^2}/{(2 \epsilon^2)}}$
and
$
d_{l t}^2=\|(\widehat{\mathbf{a}}_i)_l-(\widehat{\mathbf{a}}_i)_t\|_2^2,
$ where we set $\epsilon=0.2$ in all experiments for simplicity.
Recall that $\mathbf{L}_i$ is constructed
by $\mathbf{L}_i=\mathbf{D}_i-{\mathbf{Z}}_i$,
where the $(l,t)$-th element of ${\mathbf{Z}}_i$ is $({z}_i)_{l t}$ and $\mathbf{D}_i$ is a diagonal matrix with  the $l$-th  diagonal element being  the sum of elements of the $l$-th row of ${\mathbf{Z}}_i$. 

\subsection{Synthetic Data}

\begin{figure}[!t]
	\centering
	\subfigure[]{%
		\includegraphics[width=0.4\textwidth]{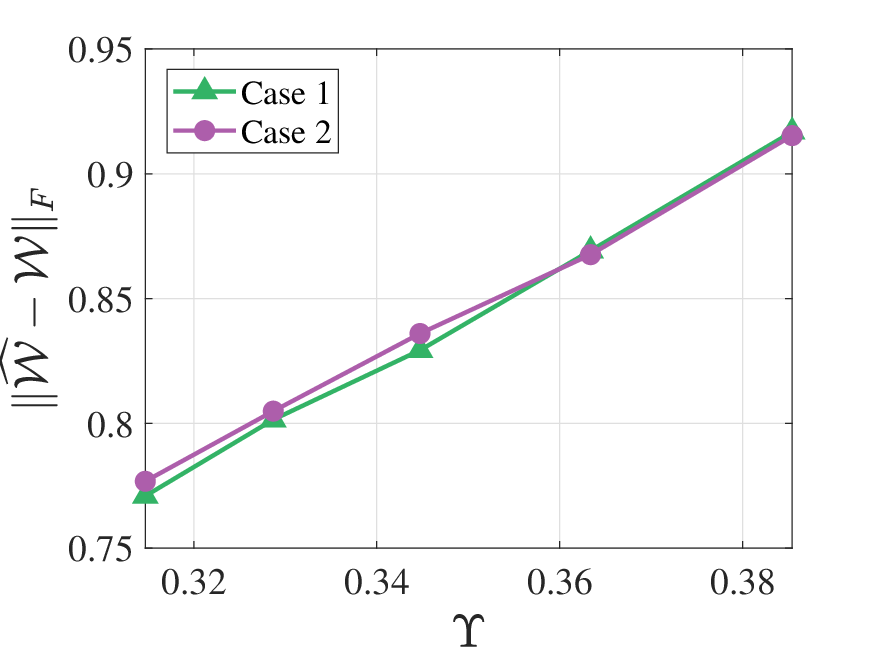}%
	}%
	\hspace{0.05mm}%
	\subfigure[]{%
		\includegraphics[width=0.4\textwidth]{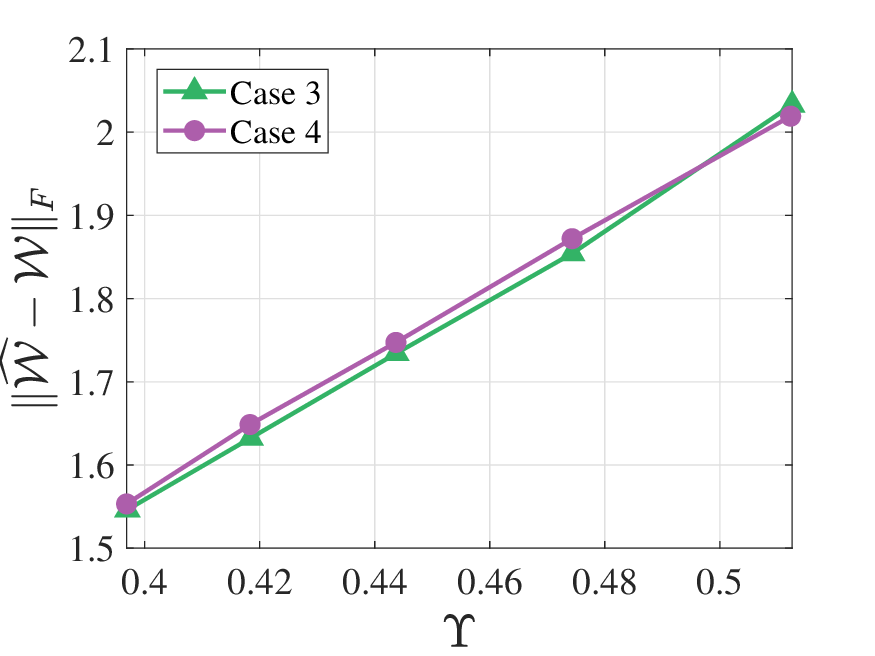}%
	}%
	\caption{The estimation error $\|\widehat{\mathcal{W}}-\mathcal{W}\|_{F}$  versus $\Upsilon=(\sqrt{s}+\sqrt{\tilde{s}})\sqrt{\log(m^2 p) / T}$ for different synthetic data. (a) $(m,p)=(50,8)$. (b)  $(m,p)=(80,12)$.}
	\label{syn_error}
\end{figure}

In this subsection,  we first verify the non-asymptotic error bound of the estimator of STDGR via synthetic data. 
We set $\gamma_1=\gamma_2=\gamma_3$ simply and 
choose them from the set $\{0.001,0.1,10,100\}$ to get the best performance.
By Theorem \ref{ConvRe},
we simply set $\rho_i=\bar{a}_1 L_i$ for $i=1,2,3,4$ and $\rho_{j+4}=\bar{a}_2 \gamma_j$ for $j=1,2,3$ 
to guarantee the convergence of Algorithm \ref{alg:algorithm1}, 
where $\bar{a}_1$ is set as $1.1$ and $\bar{a}_2$ is selected from the set $\{1.1,3,5,10\}$.
The parameter $c$ defined in (\ref{test3}) is chosen from $\{1,2\}$  and $\beta$ is set to be $0.001$.
For simplicity, we set $\alpha_1=\alpha_2=\alpha_3=\bar{\alpha}$ in all experiments, where $\bar{\alpha}$ is selected from the set $\{0.001,0.1,1\}$ to get the best performance.

Four cases are discussed for the synthetic data.
Specifically,
for the core tensor, we choose a   diagonal cube with superdiagonal elements being nonzero,
and $\mathbf{A}_i$ is generated as the $r_i$ left singular vectors of Gaussian random matrices.
Other parameters about the dimensions and Tucker rank are set as follows: 
Case 1.  $(m,p)=(50,8)$,  Tucker rank $(r_1,r_2,r_3)=(3,3,3)$, $\mathcal{G}_{111}= \mathcal{G}_{222}= \mathcal{G}_{333} =2$.
Case 2.  $(m,p)=(50,8)$, 
 Tucker rank $(r_1, r_2, r_3) = (6, 6, 6)$, $\mathcal{G}_{111}= \mathcal{G}_{222}= \mathcal{G}_{333}= \mathcal{G}_{444}=\mathcal{G}_{555}= \mathcal{G}_{666}=2$.
 Case 3.  $(m,p)=(80,12)$,  Tucker rank $(r_1, r_2, r_3) = (4, 4, 4)$,  $(\mathcal{G}_{111}, \mathcal{G}_{222}, \mathcal{G}_{333}, \mathcal{G}_{444}) = (2, 2, 2, 1)$.
 Case 4.  $(m,p)=(80,12)$,
 Tucker rank $(r_1, r_2, r_3) = (7, 7, 7)$, and  $(\mathcal{G}_{111}, \mathcal{G}_{222}, \mathcal{G}_{333}, \mathcal{G}_{444}, \mathcal{G}_{555}, \mathcal{G}_{666}, \mathcal{G}_{777})=(2,2,2,1,1,1,0.5)$.
 Moreover, $\boldsymbol{\varepsilon}_t \stackrel{\scriptsize{\mathrm{i.i.d.}}}{\sim} N(\mathbf{0}, \mathbf{I}_m)$.
 Due to the randomness, all experiments are performed 10 times in each case and we take the average results as the last result.
 It follows from Theorem \ref{statistical} that the error bound $\|\widehat{\mathcal{W}}-\mathcal{W}\|_{F}$ is on the order of  $O_p((\sqrt{s}+\sqrt{\tilde{s}})\sqrt{\log(m^2 p) / T})$ for fixed Tucker rank, where we set $\tilde{s}=s$.
 Let $\Upsilon=(\sqrt{s}+\sqrt{\tilde{s}})\sqrt{\log(m^2 p) / T}$ and we choose different sample sizes $T$ such as $\Upsilon= 0.314, 0.328, 0.344, 0.363, 0.385$ for Case 1 and Case 2, and $\Upsilon= 0.397, 0.418, 0.443, 0.474, 0.512$ for Case 3 and Case 4.
In Figure \ref{syn_error}, we plot $\|\widehat{\mathcal{W}}-\mathcal{W}\|_{F}$ 
versus $\Upsilon$ for the four cases. It can be seen that the estimation error $\|\widehat{\mathcal{W}}-\mathcal{W}\|_{F}$   generally increases linearly
in $\Upsilon$, and the two lines in each figure almost coincide,
which demonstrates the theoretical finding in Theorem \ref{statistical}.

Next we compare STDGR with Lasso, NNM, LRTLT, TVART, and SHORR for synthetic data.
For the noise term, we set  $\boldsymbol{\varepsilon}_t \stackrel{\scriptsize{\mathrm{i.i.d.}}}{\sim} N(\mathbf{0}, 0.25\mathbf{I}_m)$.
Additionally, the factor matrices $\mathbf{A}_i$ are constructed as follows \cite{Belkin2002}: we randomly generate weight matrices $\mathbf{Z}_i$ and 
get  the corresponding Laplacian matrices based on (\ref{DefIoL}). 
Then the first $r_i$ eigenvectors associated with the smallest eigenvalues are used to obtain $\mathbf{A}_i$. 
In spectral clustering, the entries of $\mathbf{A}_i$ measure the
similarities between data samples and clusters, as they indicate how likely data samples belong to specific clusters \cite{wang2014multi}.
Therefore, this construction of $\mathbf{A}_i$ can 
reflect the similarity of the rows  of $\mathbf{A}_i$.
Two cases are considered, where $\mathcal{G}$ is a   diagonal cube with superdiagonal elements being nonzero.
Case a: $(m,p)=(70,12)$,
Tucker rank $(r_1, r_2, r_3) = (5, 5, 5)$,
and $(\mathcal{G}_{111}, \mathcal{G}_{222}, \mathcal{G}_{333},\mathcal{G}_{444}, \mathcal{G}_{555}) = (0.2, 0.2, 0.2,0.1,0.1)$.
Case b: $(m,p)=(100,20)$,
Tucker rank $(r_1, r_2, r_3) = (7, 7, 7)$,
and $(\mathcal{G}_{111}, \mathcal{G}_{222}, \mathcal{G}_{333}, \mathcal{G}_{444}, \mathcal{G}_{555}, \mathcal{G}_{666}, \mathcal{G}_{777})$ $=(0.1,0.1,0.1,0.05,0.05,0.05,0.05)$.

\begin{figure}[!t]
	\centering
	\subfigure[]{%
		\includegraphics[width=0.4\textwidth,trim=0 40 0 30,clip]{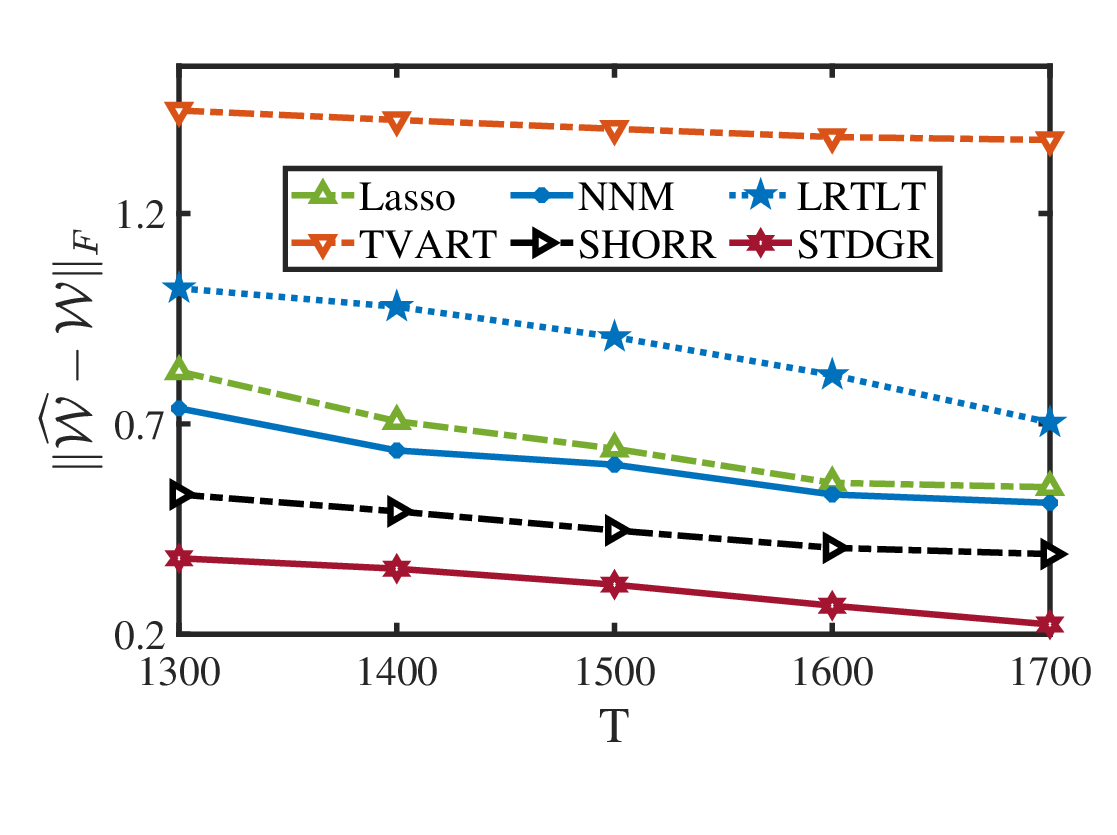}%
	}%
	\hspace{0.05mm}%
	\subfigure[]{%
		\includegraphics[width=0.4\textwidth,trim=0 40 0 30,clip]{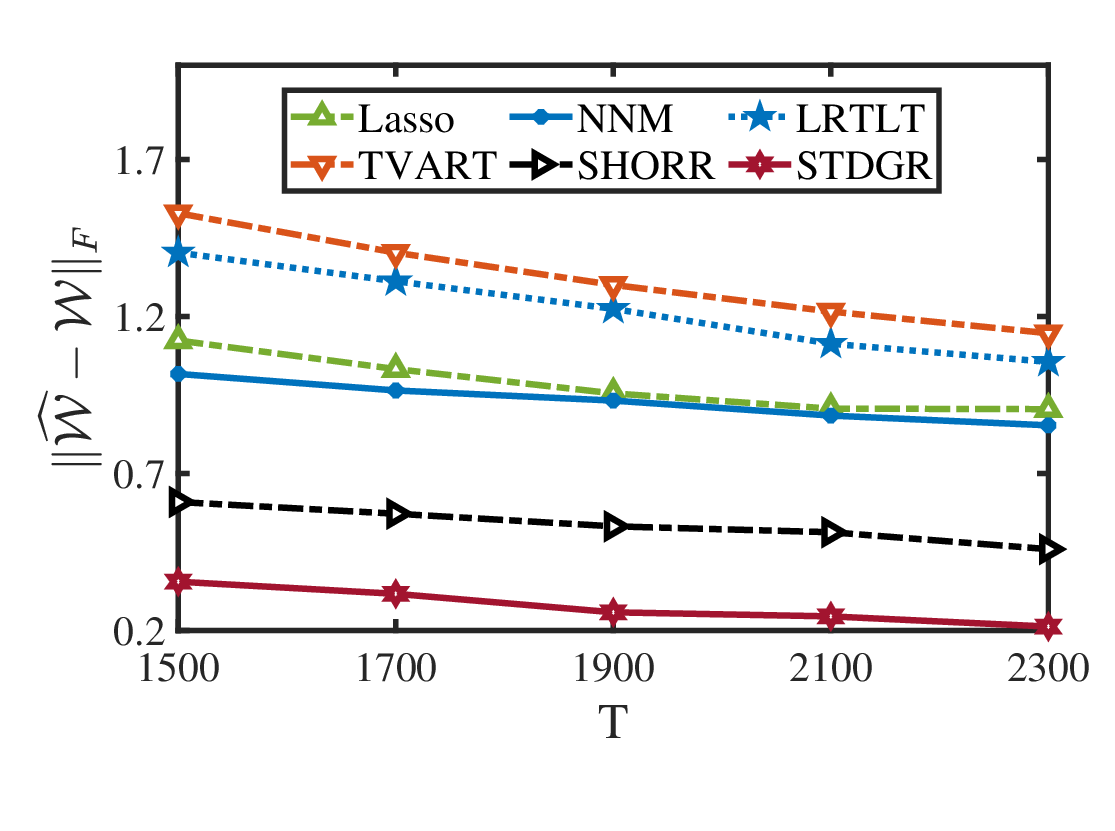}%
	}%
	\caption{The estimation error $\|\widehat{\mathcal{W}}-\mathcal{W}\|_{F}$ versus different sample size $T$ of different methods for synthetic data under two settings of $(m,p)$. (a) Case a. (b)  Case b.}\label{syn_compare}
\end{figure}

In Figure \ref{syn_compare}, we show the estimation error $\|\widehat{\mathcal{W}}-\mathcal{W}\|_{F}$ versus different sample sizes $T$ of different methods for synthetic data under the above two settings of $(m,p)$.
It can be observed that the estimation errors of different methods decrease as the size $T$ increases. 
Moreover, the errors obtained by STDGR are lower than those obtained by other methods, where the similarity of  $\mathbf{A}_i$ is exploited in STDGR.
Additionally, the SHORR outperforms Lasso, NNM, LRTLT, TVART in terms of estimation errors, which shows the Tucker decomposition based method performs better than the matrix or CP based methods.

\subsection{Comprehensive Climate Dataset}\label{NxMCCD}

The comprehensive climate dataset (CCD)\footnote{\url{https://viterbi-web.usc.edu/~liu32/data/NA-1990-2002-Monthly.csv}} contains
17 climate variables of North America from 1990 to 2002.
The data are collected and preprocessed by five federal agencies and interpolated on a \(2.5\times2.5\) degree grid, which results in 125 observation locations for each variable.
For each location, 156 time points (12 months per year $\times$ 13 years) are recorded.
All series are  performed de-seasonalization by removing seasonal averages  and then standardized to zero mean and unit variance.
For the CCD, we choose eight variables in the experiments, including Carbon dioxide (CO2),  methane (CH4), 
Carbon monoxide (CO), hydrogen (H2), wet day,  
temperature minimum  (TMN),  
temperature mean (TMP), and  global horizontal radiation (GLO).

For the VAR($p$) model, 
the order is simply set to $p=2$ for all methods, except for TVART, for which $p=1$ is employed.
We set $\gamma_1=\gamma_2=\gamma_3$ in each case and 
choose them from the set $\{50,80,90,100,140,150,180,300,350,1000,2000,3000\}$ to get the best forecasting performance.
By the condition in Theorem \ref{ConvRe},
we set $\rho_i=\bar{a}_1 L_i$ for $i=1,2,3,4$ and $\rho_{j+4}=\bar{a}_2 \gamma_j$ for $j=1,2,3$ to guarantee the convergence of Algorithm \ref{alg:algorithm1}, 
where $\bar{a}_1$ is chosen from the set $\{1.001, 1.1\}$ and $\bar{a}_2$ is chosen from the set $\{10, 50\}$.
The parameter $c$ defined in (\ref{test3}) is set to $1$ and $\beta$ is set to $0.001$.
For simplicity, we set $\alpha_1=\alpha_2=\alpha_3=\bar{\alpha}$ in all experiments, where $\bar{\alpha}$ is selected from the set $\{0.001,0.4,1,2,4,5,6,10,12,14,100\}$ to obtain the best forecasting performance in our experiments.

\begin{table}[!t]
	\centering
	\caption{MSE of different methods for CCD with eight variables.}
	\begin{tabular}{ccccccc}
		\toprule
		Variable & Lasso & NNM    & LRTLT  & TVART  & SHORR & STDGR \\
		\midrule
		CO2     & 0.9787 & 0.9775 & 1.0161 & 1.1477 & 0.9774 & \bf{0.9620} \\
		CH4     & 1.1942 & 1.1977 & 1.1619 & 1.2653 & 1.1655 & \bf{1.1107} \\
		CO     & 0.8832 & 1.0420 & 0.8972 & 0.9267 & 0.8374 & \bf{0.7986} \\
		H2     & 0.5687 & 0.5652 & 0.5575 & 0.6191 & 0.5299 & \bf{0.4730} \\
		Wet     & 0.8850 & 0.8788 & 1.1835 & 1.0781 & 0.8747 & \bf{0.8645} \\
		TMN    & 0.7620 & 0.7741 & 0.9386 & 0.8626 & 0.7416 & \bf{0.7332} \\
		TMP    & 0.7720 & 0.7711 & 0.9829 & 0.8553 & 0.7584 & \bf{0.7388} \\
		GLO    & 0.9603 & 0.9424 & 1.1006 & 0.9956 & 0.9236 & \bf{0.9101} \\
		\bottomrule
	\end{tabular}%
	\label{CCDS-60}%
\end{table}%

\begin{figure}[!t]
	\centering
	\subfigure[CO2]{%
		\includegraphics[width=5.22cm,height=3.5cm]{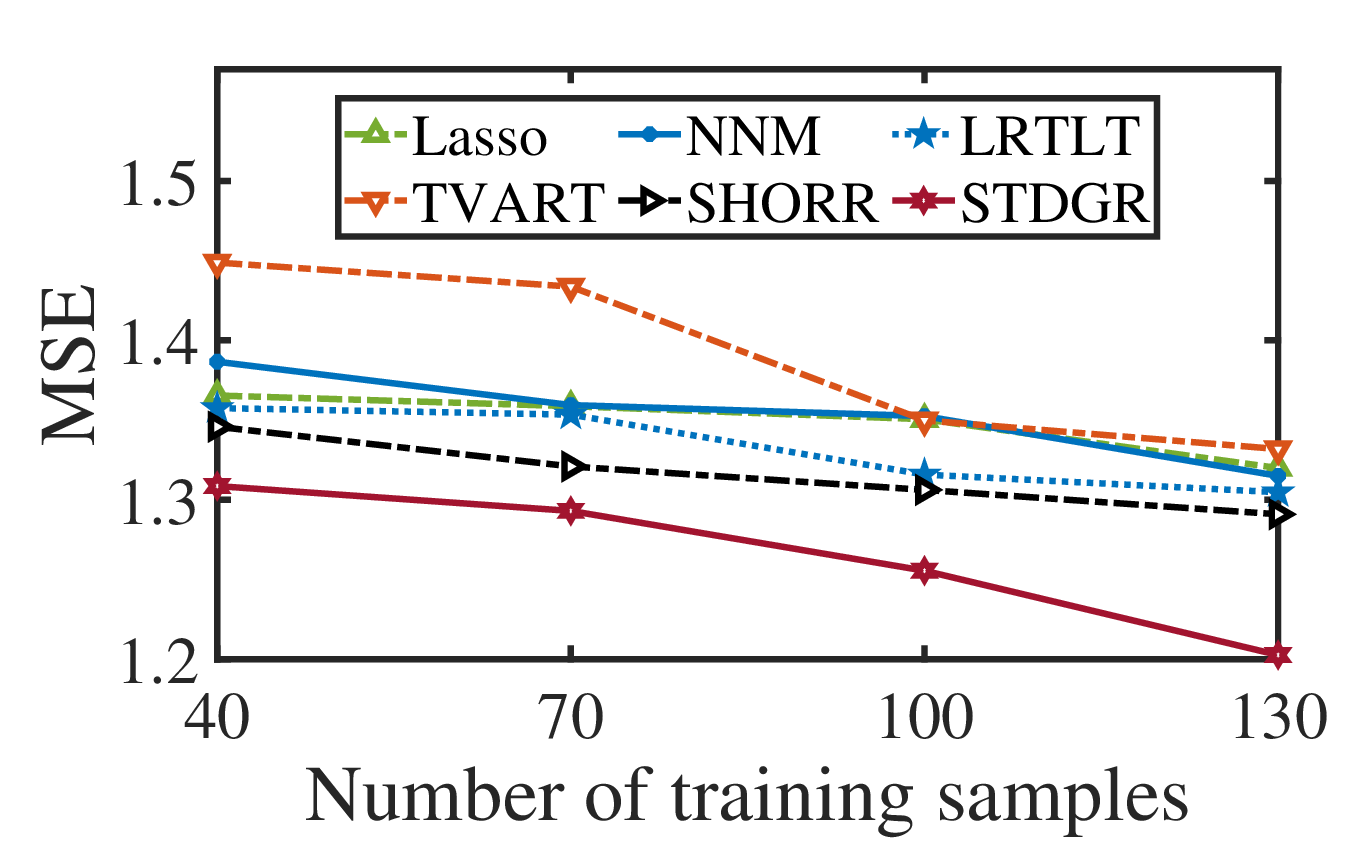}%
	}%
	\hspace{0.001mm}%
	\subfigure[CH4]{%
		\includegraphics[width=5.22cm,height=3.5cm]{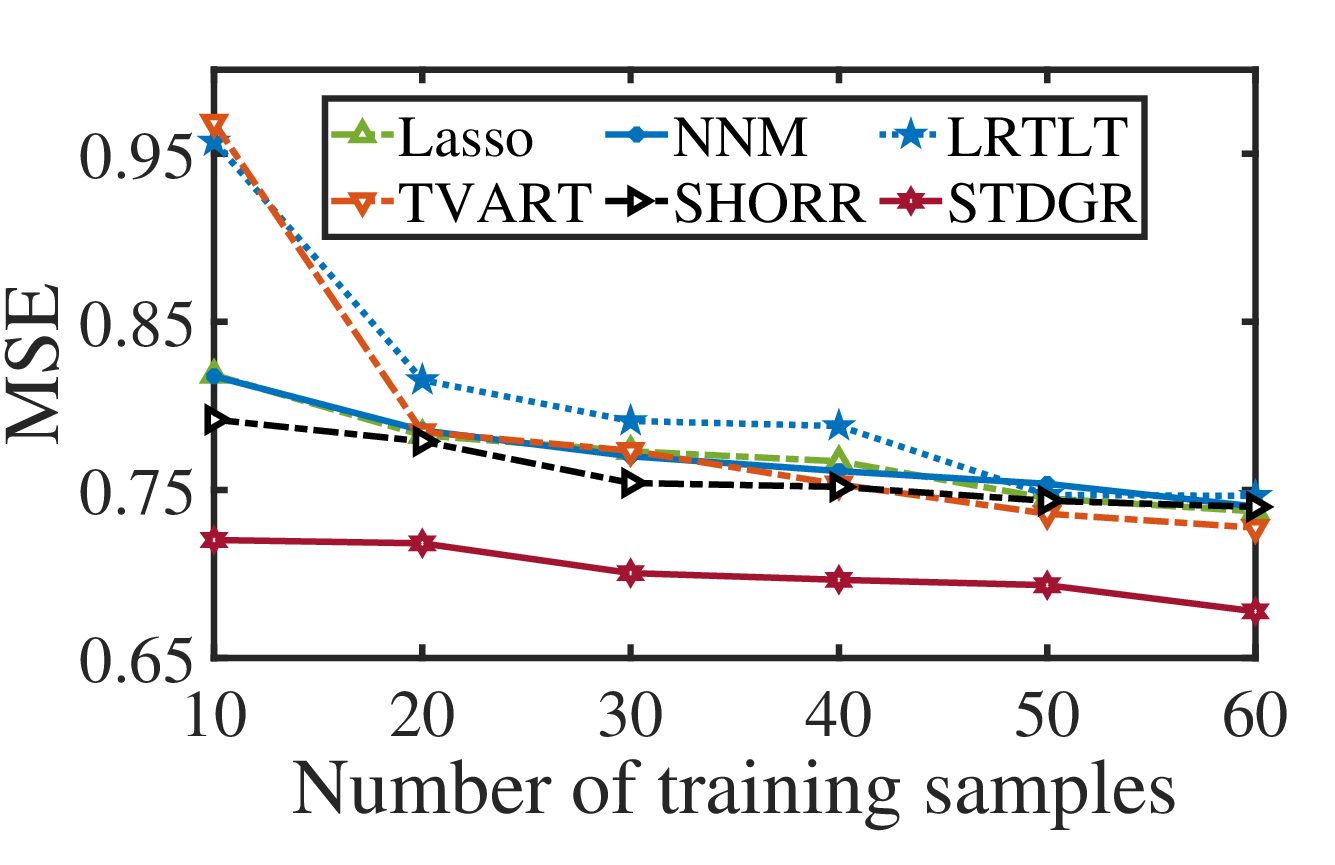}%
	}%
	\hspace{0.001mm}%
	\subfigure[H2]{%
		\includegraphics[width=5.22cm,height=3.5cm]{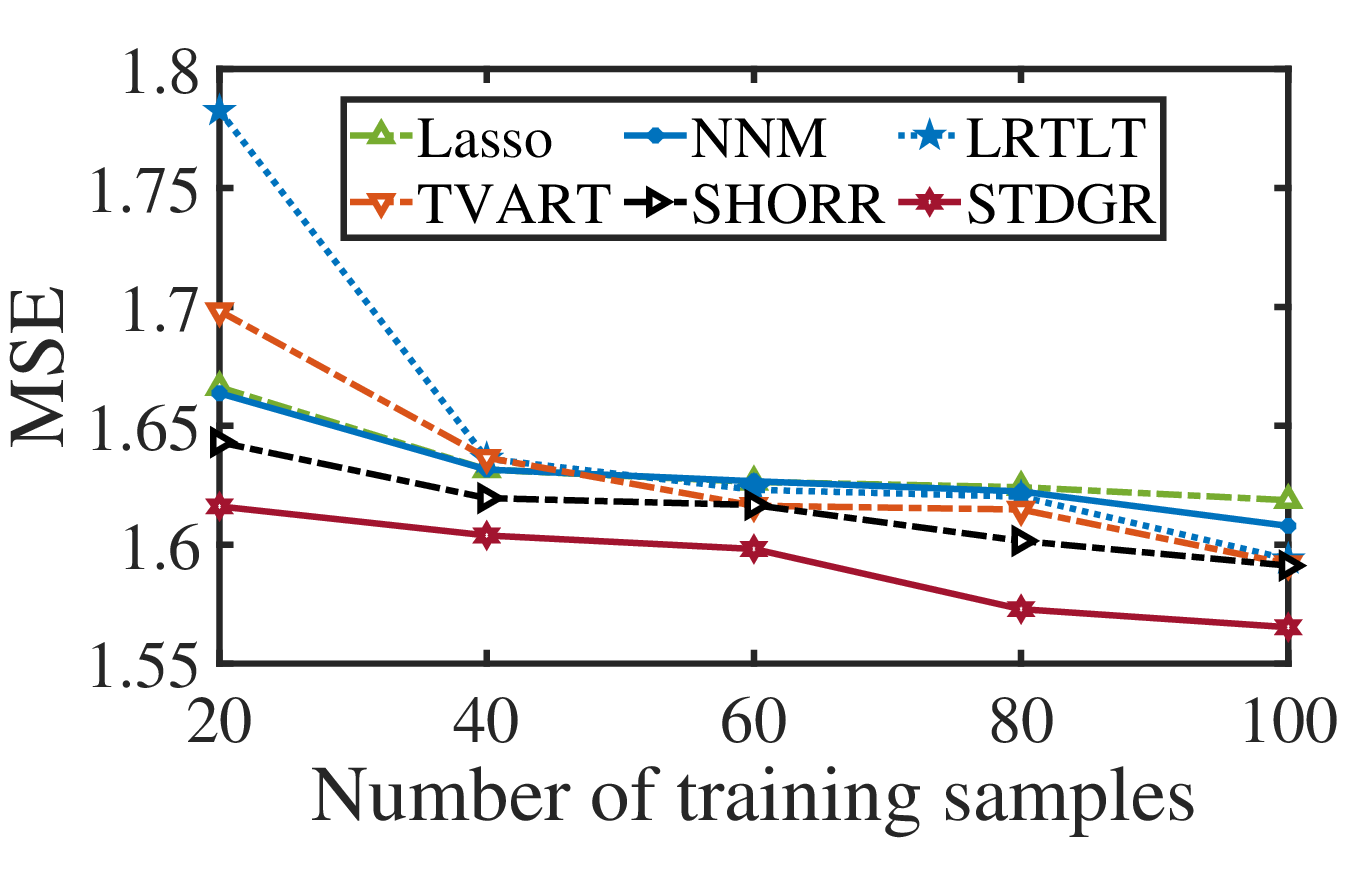}%
	}%
	\caption{MSE versus number of training samples of different methods for the CCD.}
	\label{CCDS}
\end{figure}

In Table \ref{CCDS-60}, we report the MSE of different methods  for the variables CO2, CH4, CO, H2, wet day, TMN, TMP, and GLO, where the first $60\%$ time series of each variable are used for training and the resting 40\% time series are used for testing.
The best results are
highlighted in bold.
It can be seen from this table that the STDGR achieves lower MSE than other methods. 
Furthermore, the SHORR outperforms Lasso, NNM, LRTLT, and TVART in terms of MSE.
The improvement of  STDGR is very impressive compared with SHORR for CH4, CO, H2.

In Figure \ref{CCDS}, we show the MSE of different methods versus the number of training samples for CO2, CH4, and H2, where the number of testing samples are 20 for CO2 and 10 for CH4, H2, respectively.
We can observe that the MSE of different methods decreases as the number of training samples increase for the three variables, 
which demonstrates that larger number of training samples can improve the forecasting accuracy for VAR.
Moreover, the MSE obtained by STDGR is lower than those obtained by other  methods, 
which indicates the superior  performance of STDGR for forecasting.

\subsection{Traffic Flow Dataset}\label{Traffiddatas}

In this subsection, we test two traffic flow datasets to demonstrate the effectiveness of STDGR for VAR.
The first one is the Hangzhou metro passenger flow dataset\footnote{\url{https://tianchi.aliyun.com/competition/entrance/231708/information}}, 
which  records incoming passenger flow of 80 metro stations over 25 days from January 1 to January 25, 2019, with a 10-minute resolution in Hangzhou, China. We discard the interval from 00:00 to 06:00 since metro service is unavailable at that time, 
retaining the remaining 108 time intervals per day.
We select 2500 time intervals in this dataset.
The second one is the PeMS04 dataset\footnote{\url{https://ieee-dataport.org/documents/pems03-and-pems04}},
which is a public traffic flow dataset collected by the California Transportation Agencies Performance Measurement System (PeMS) in the Bay Area of United States. 
The PeMS04 dataset contains flow volumes of 307 sensors from January 1 to February 28, 2018, which are recorded every five minutes.
As a result, there are 16992 time intervals and 
we select 7500 time intervals in the following experiments.

For the traffic flow dataset, the order of the VAR($p$) model is set as $p=7$ for all methods, except for TVART, in which  we set $p=1$.
The parameters $\gamma_1,\gamma_2,\gamma_3$ are set as the same in each case and 
chosen  from the set $\{0.001, 0.01, 0.1, 1, 100\}$ to obtain the best forecasting performance.
We fix $c=1$ defined in (\ref{test3}) and set $\beta=0.001$.
In all experiments, 
we set $\alpha_1=\alpha_2=\alpha_3=\bar{\alpha}$, 
where $\bar{\alpha}$ is chosen from the set $\{0.001, 0.1\}$. 
According to Remark \ref{Condirho}, 
we set $\rho_i=\bar{a}_1 L_i$ for $i=1,2,3,4$ and 
$\rho_{j+4}=\bar{a}_2 \gamma_j$ for $j=1,2,3$ to 
ensure the convergence of Algorithm \ref{alg:algorithm1}. 
Here, $\bar{a}_1$ and $\bar{a}_2$ are chosen from 
the set $\{1.1, 1.3, 10\}$ to get the best forecasting performance.

\begin{table}[!t]
	\centering
	\caption{MSE  ($\times 10^{-3}$) of forecasting results of different methods for the traffic flow dataset.}
	\begin{tabular}{ccccccc}
		\toprule
		Dataset & Lasso & NNM    & LRTLT  & TVART  & SHORR & STDGR \\
		\midrule
		Hangzhou     & 1.9612 & 1.7411 & 1.4653 & 4.2263
		& 1.2394 & \bf{0.5961} \\
		PeMS04     & 3.3012 & 2.7345 & 5.4871 & 5.2348 & 2.5538 & \bf{2.0632}\\
		\bottomrule
	\end{tabular}%
	\label{Traffic-70}%
\end{table}

\begin{figure}[!t]
	\centering
\subfigure[Hangzhou]{%
	\includegraphics[width=0.4\textwidth,trim=0 8mm 5mm 0, clip]{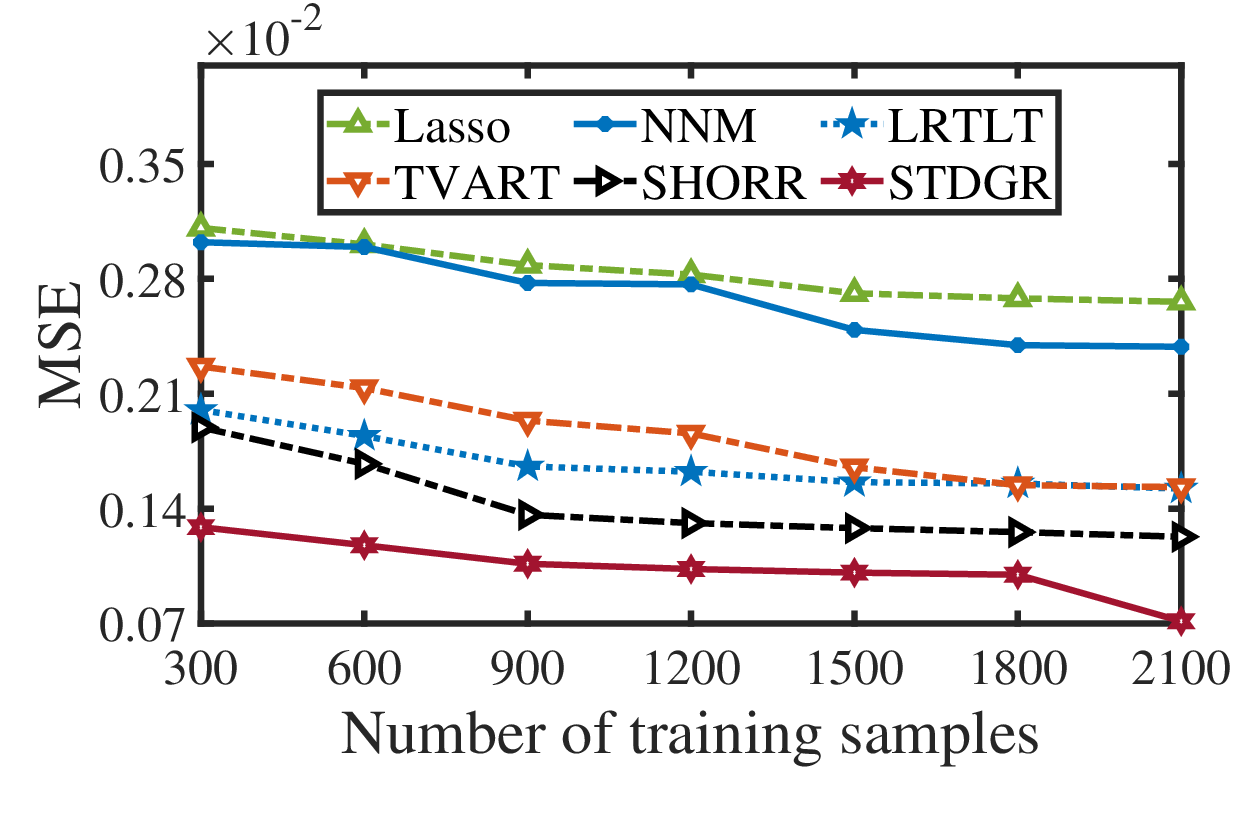}%
}%
\hspace{0.3mm}%
\subfigure[PeMS04]{%
	\includegraphics[width=0.4\textwidth,trim=0 8mm 5mm 0, clip]{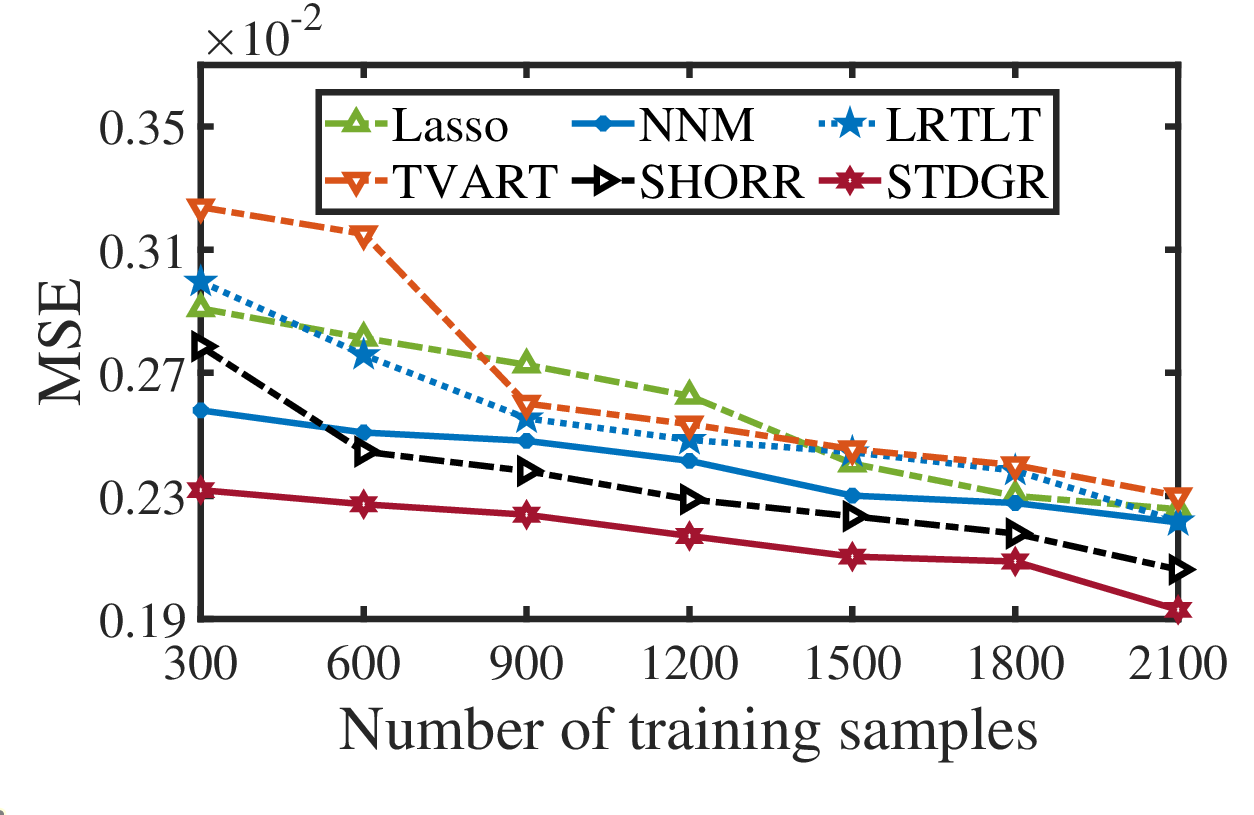}%
}%
	\caption{MSE versus number of training samples of different methods for the two traffic flow datasets.}
	\label{Traffic-300}
\end{figure}

\begin{figure}[!t]
	\centering
	\subfigure[Station 74]{%
		\includegraphics[width=0.49\textwidth]{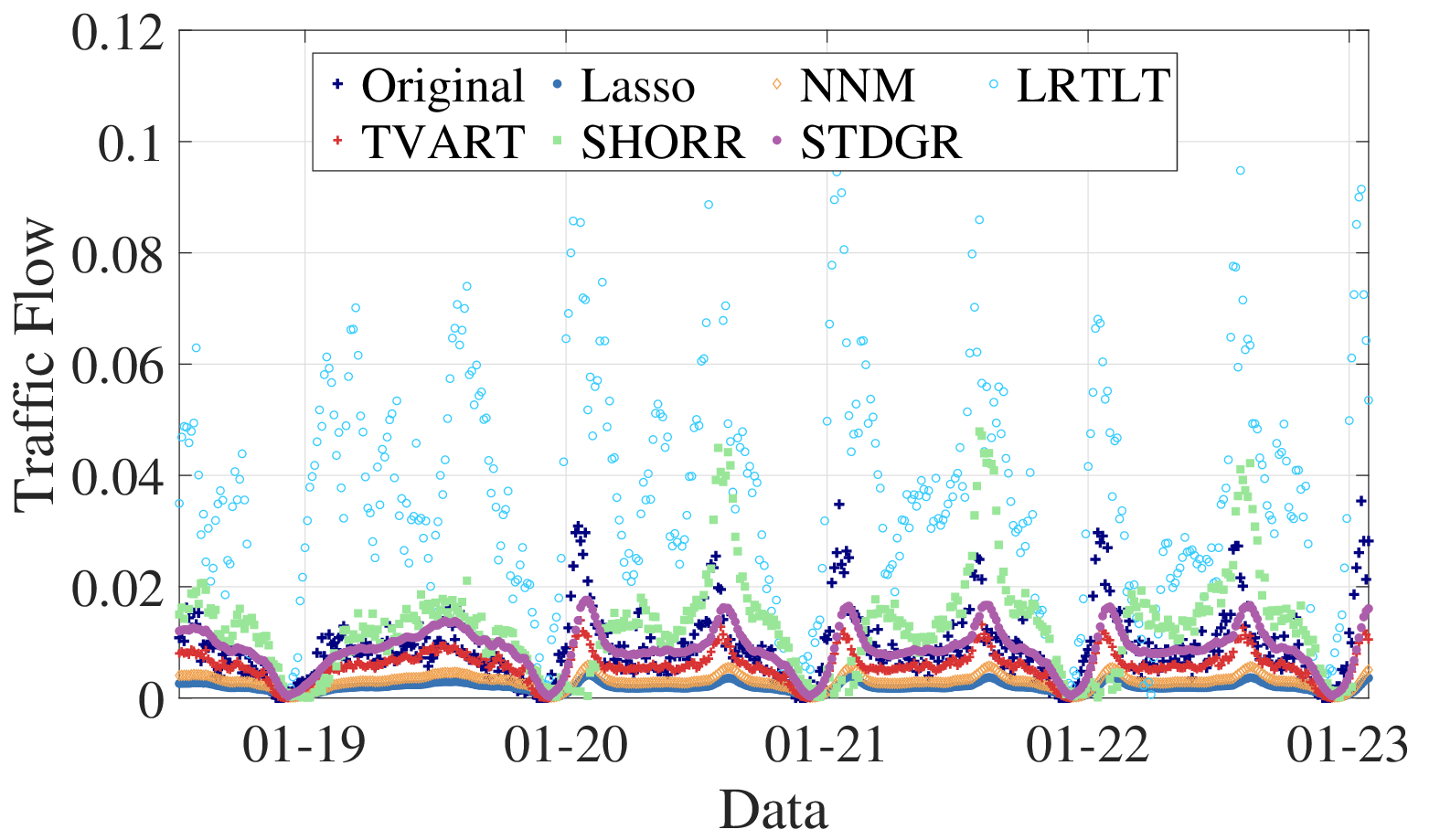}%
	}%
	\hspace{0.05mm}%
	\subfigure[Station 79]{%
		\includegraphics[width=0.49\textwidth]{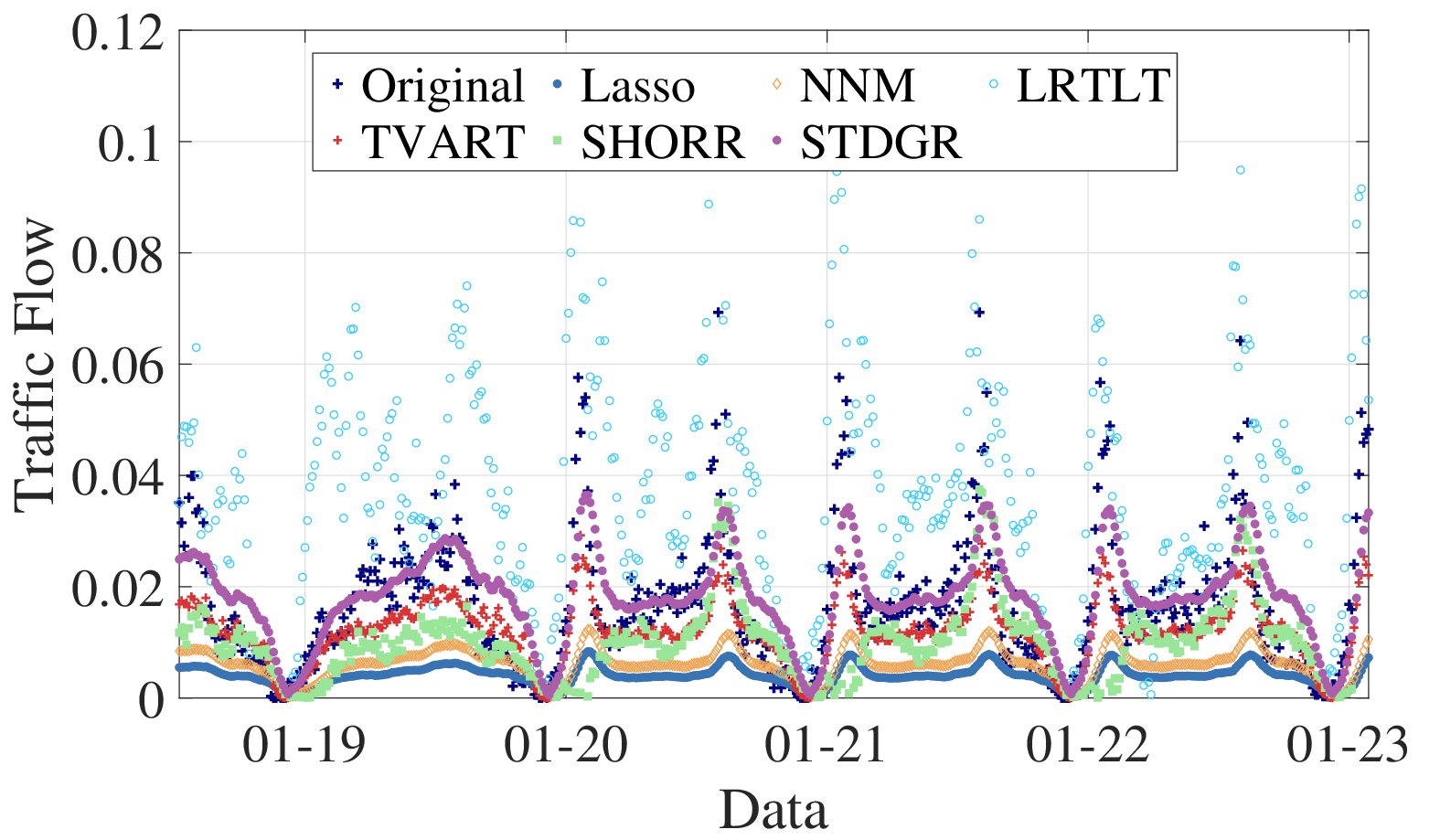}%
	}%
	\caption{Traffic flow prediction of different methods at different time points for the 74th and 79th metro stations of the Hangzhou dataset, respectively.}
	\label{hangzhou-80}
\end{figure}

\begin{figure}[!t]
	\centering
	\subfigure[Sensor 71]{%
		\includegraphics[width=0.49\textwidth]{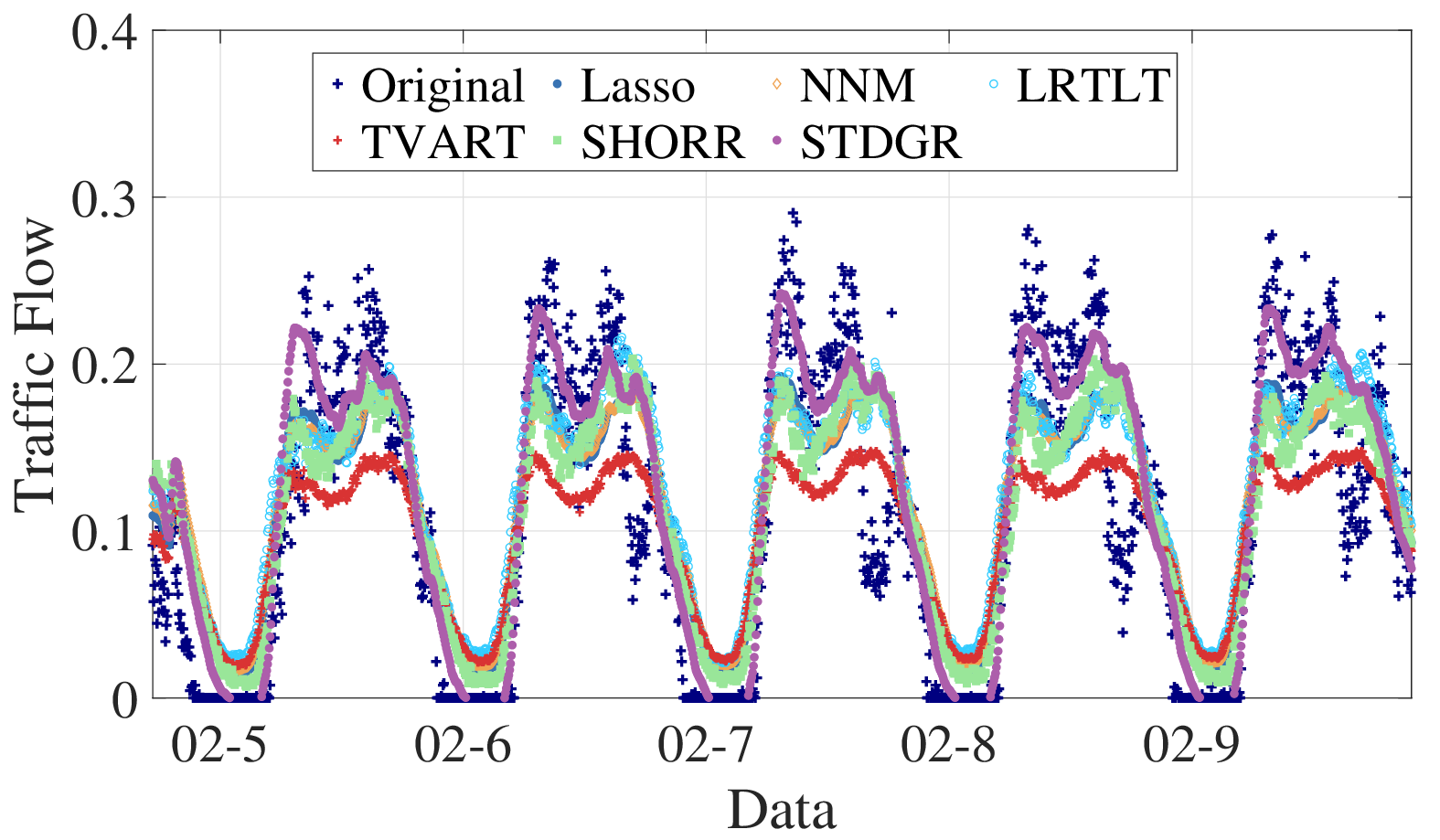}%
	}%
	\hspace{0.05mm}%
	\subfigure[Sensor 74]{%
		\includegraphics[width=0.49\textwidth]{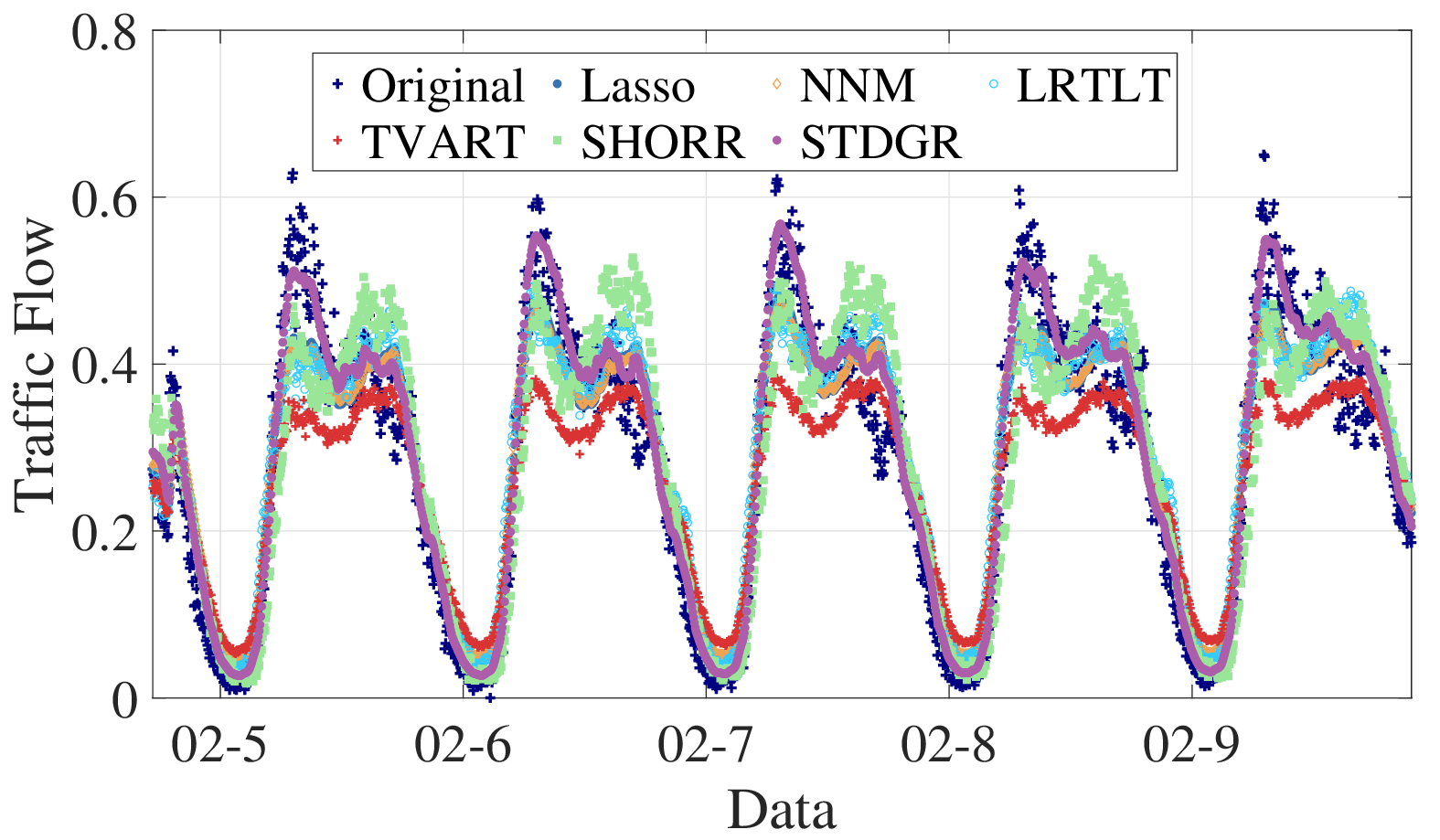}%
	}%
	\caption{Traffic flow prediction of different methods at different time points for the 71st and 74th sensors of the PeMS04 dataset, respectively.}
	\label{PEMS04-80}
\end{figure}

In Table \ref{Traffic-70}, we report the prediction errors of different methods for the Hangzhou and PeMS04 datasets, where we use 70\% samples for training and the resting 30\%  samples for testing in each dataset. 
The best results are highlighted in bold.
It can be seen from this table that the STDGR achieves the smallest prediction error compared with other methods for the  two traffic flow datasets, which exhibits the superior forecasting performance of STDGR.
Moreover, 
the SHORR performs better than Lasso, NNM, LRTLT, and TVART in terms of MSE.

In Figure \ref{Traffic-300}, we show the MSE  ($\times 10^{-2}$) versus number of training samples of different methods for the Hangzhou and PeMS04 datasets, 
where the number of training samples varies from 300 to 2100 with step size 300 and the number of testing samples is 300.
As it can be seen that the MSE of all methods generally decreases as the number of the training samples increases. 
And the STDGR consistently achieves the lowest MSE 
for all testing cases and exhibits great improvement when the number of training samples increases. 
Moreover, the SHORR performs better than Lasso, NNM, LRTLT, and TVART for most cases, especially for large number of training samples.
These results demonstrate that the STDGR obtains higher forecasting accuracy compared with 
the competition approaches, even when the number of training samples is relatively small.

In Figures \ref{hangzhou-80} and \ref{PEMS04-80}, we show the traffic flow prediction of different methods at different time points for the Hangzhou and PeMS04 datasets, respectively,
where the 74th and 79th metro stations of the Hangzhou dataset are shown in Figure \ref{hangzhou-80}  and 
the 71st and 74th sensors of the PeMS04 dataset are shown in Figure \ref{PEMS04-80}.
Here 80\% samples are used for training and the remaining 20\% samples are used for testing for the two datasets.
It can be observed that the STDGR can fit the original traffic flow  better compared with other methods for different sensors or stations.
Furthermore, the STDGR can predict the trends of traffic flows better for the two datasets, which exhibits the superior performance of STDGR for prediction.

\subsection{Macroeconomic Dataset}

In this subsection, we test the performance of STDGR for the macroeconomic dataset, 
which contains  40 quarterly  
macroeconomic variables of the United States from 1959 to 2007 with 194 time points 
for each  variable\footnote{\url{https://www.stlouisfed.org/research/economists/mccracken/fred-databases}} \cite{Koop2013}.
The variables are seasonally adjusted  except for the financial series,
and transformed to stationary with zero mean and unit variance \cite{Wang2022}.

We set  $p = 4$ in the VAR($p$) model and  $\gamma_1=\gamma_2=\gamma_3=0.001$ for simplicity.
As suggested in Remark \ref{Condirho}, the penalty parameters are simply set to $\rho_i=\bar{a}_1 L_i$ for $i=1,\ldots,4$ and $\rho_{j+4}=\bar{a}_2\gamma_j$ for $j=1,2,3$, where  $\bar{a}_1=1.001$ and $\bar{a}_2=250$. 
Besides, $c$ is set to $1$ for simplicity.
For the parameters $\beta$, $\alpha_i$ in the STDGR model, we set $\beta=200$, $\alpha_1=\alpha_2=\alpha_3=0.001$.
Moreover, the Tucker ranks of SHORR and STDGR are chosen as the same for this dataset.

Table \ref{macro_MSE} reports the MSE and training time (in seconds) of Lasso, NNM, LRTLT, TVART, SHORR, and STDGR, where we use 85\% and 90\% samples for training and the remaining samples for testing, respectively.
It can be seen from this table that the MSE achieved  by STDGR is lower than those achieved  by other methods for the testing two cases,
which shows that the forecasting accuracy of STDGR is better than other comparison methods.
Furthermore, the SHORR achieves the second best forecasting performance, while it requires more running time compared with other methods.
Although the STDGR takes more running time  than Lasso, NNM, LRTLT, TVART, the MSE obtained by STDGR is lower  than these methods.

\begin{table}[!t]
	\centering
	\caption{Comparisons of  MSE and running time (in seconds) of different methods for the macroeconomic dataset.}
	\begin{tabular}{cccccccc}
		\toprule
		Index & \multicolumn{1}{l}{Training size} & Lasso & NNM   & LRTLT & TVART & SHORR & STDGR \\
		\midrule
		\multirow{2}[0]{*}{MSE} & 85\% & 5.3213  & 5.4821  & 5.7142 & 6.1278 & 5.0664 & \bf{4.3838} \\
		& 90\% & 4.6321 & 4.6500 & 4.2347 & 4.8468 & 3.8021 & \bf{3.1845} \\
		\midrule
		\multirow{2}[0]{*}{Time} & 85\% & 0.02 & 3.00 & 0.05 & 0.26 & 537.14 & \bf{17.08} \\
		& 90\% & 0.02 & 2.88 & 0.07 & 0.25 & 188.51 & \bf{20.15} \\
		\bottomrule
	\end{tabular}%
	\label{macro_MSE}%
\end{table}%

\subsection{Ablation Analysis}

In this subsection, we show the effectiveness of STDGR for both the sparse Tucker decomposition and graph regularization.
Specifically, we compare with the Tucker decomposition model with sparse core tensor (called TDSC for short), which uses the sparsity constraint on the core tensor
rather than on the factor matrices compared with \cite{Wang2022} and does not impose the graph regularization.
Moreover, we also compare with the Tucker decomposition and graph regularization without sparse core tensor constraint (called TDGR for short).
And the SHORR is also compared to show the effectiveness of the sparsity constraint on the core tensor rather than on the factor matrices in Tucker decomposition.

\begin{table}[!t]
	\centering
	\caption{Comparisons of MSE ($\times 10^{-3}$) of SHORR, TDSC, TDGR, STDGR for the Hangzhou dataset.}
	\begin{tabular}{ccccc}
		\toprule
		Training size & SHORR &  TDSC   &  TDGR  &  STDGR \\
		\midrule
		75\%     & 0.4475 & 0.3941 & 0.4782 & 0.3532 \\
		85\%     & 0.3714 & 0.3370 & 0.3921 & 0.2743  \\
		\bottomrule
	\end{tabular}%
	\label{Compare_methods}%
\end{table}

In Table \ref{Compare_methods}, we report the  MSE ($\times 10^{-3}$) of SHORR, TDSC, TDGR, and STDGR for the Hangzhou dataset, where  75\% and 85\%  samples are used for training and the remaining samples are used for testing, respectively.
It can be seen from this table that the MSE of TDSC is lower than those of SHORR, which shows the superiority of imposing sparsity on the core tensor rather than on the factor matrices.
Moreover, the SHORR and TDSC perform better than TDGR for the testing  cases, which implies that the sparsity constraints on the core tensor or factor matrices are more effective than the graph regularization.
Besides, the STDGR substantially outperforms  other three methods. This demonstrates that the sparsity constraint on the core tensor and graph regularization on the factor matrices are both effective for VAR modeling.

\section{Concluding Remarks}\label{CondRem}

In this paper, we have proposed a sparse Tucker decomposition model with
graph Laplacian regularization for high-dimensional VAR time series forecasting.
By stacking the transition matrices into a third-order tensor, the Tucker decomposition is employed for the transition tensor
to address the over-parameterization issue, which can restrict the parameter space along three directions simultaneously.
Besides, the sparsity of the core tensor is enforced to further reduce the parameters and add the interpretability.
Moreover, the graph regularization is employed to characterize the local consistency of the factor matrices,
which can characterize the similarity of the sampled data  for VAR time series.
A PALM algorithm is designed to solve the resulting model  with global convergence guarantee.
Numerical experiments on simulated and real data show the superior forecasting performance of the proposed method over several baseline models.

In the future, we are going to explore the tensor value autoregressive time series based on other types of tensor decomposition, 
which may explore the correlation of data better (cf.  \cite{chen2022factor, wang2024high}).
Moreover, the second-order information about Tucker decomposition and the $\ell_1$ norm on the core tensor  will be exploited to design a more efficient algorithm, which is capable of dealing with higher-dimensional VAR modeling.

\section*{Appendix A. Auxiliary Lemmas}

\begin{Prop}\label{kron_pro}
	For any matrix $\mathbf{A}\in \mathbb{C}^{n_1 \times n_2}, \mathbf{B}\in
	\mathbb{C}^{n_3 \times n_4}$, the following assertions hold:

	(1) $(\mathbf{A} \otimes \mathbf{B})^{H}=\mathbf{A}^{H} \otimes \mathbf{B}^{H}$ \cite[Equation (1.3.1) 4.2.4]{golub2013matrix}.
	
	(2) Given matrices $\mathbf{D}\in \mathbb{C}^{n_2 \times n_5}, \mathbf{E}\in
	\mathbb{C}^{n_4 \times n_6}$, $(\mathbf{A} \otimes \mathbf{B})(\mathbf{D} \otimes \mathbf{E})=(\mathbf{A D}) \otimes (\mathbf{B E})$ \cite[Equation (1.3.2) 4.2.4]{golub2013matrix}.
	
%
	(3) $\|\mathbf{A} \otimes \mathbf{B}\|=\|\mathbf{A}\|\|\mathbf{B}\|$
	\cite[Equation (1.10)]{Gueridi2022}.
	
	(4) $\|\mathbf{A} \otimes \mathbf{B}\|_F=\|\mathbf{A}\|_F\|\mathbf{B}\|_F$ \cite[Equation (1.11)]{Gueridi2022}.
\end{Prop}

\begin{Prop}\label{kron_untari}
For any matrix $\mathbf{A}\in \mathbb{C}^{n_1 \times n_2}$ with $\mathbf{A}^H\mathbf{A}=\mathbf{I}_{n_2}$, $\mathbf{B}\in
\mathbb{C}^{n_3 \times n_4}$ with $\mathbf{B}^H\mathbf{B}=\mathbf{I}_{n_4}$,
it holds that
$$
\|\mathbf{A} \otimes \mathbf{B}\|\leq 1.
$$
\end{Prop}
\textbf{Proof}.\enspace
For any matrix $\mathbf{A}$ satisfying $\mathbf{A}^{H}\mathbf{A}=\mathbf{I}_{n_2}$, one has
$\|\mathbf{A}\|\leq1$ \cite[Theorem 3.9]{Stewart1990}.
Taking this together with Proposition \ref{kron_pro}(4) yields
$$
\|\mathbf{A} \otimes \mathbf{B}\|=\|\mathbf{A}\|\|\mathbf{B}\|\leq 1,
$$
which completes the proof.
\qed

\begin{Prop}\label{out_kron1}
For any vector $\mathbf{a}\in \mathbb{R}^{n_1}$ and
 matrix $\mathbf{X}\in \mathbb{R}^{n_2 \times n_3}$,
one has
$$
(\mathbf{a}\circ\mathbf{X})_{(2)}=\mathbf{X}\otimes \mathbf{a}^{H},\
(\mathbf{a}\circ\mathbf{X})_{(3)}=\mathbf{X}^{H}\otimes \mathbf{a}^{H}.
$$
\end{Prop}
\textbf{Proof}.\enspace
For vector $\mathbf{a}=(a_1,\ldots,a_{n_1})^H\in \mathbb{R}^{n_1}$ and
matrix $\mathbf{X}\in \mathbb{R}^{n_2 \times n_3}$ with the $(j,l)$-th element
 $x_{jl} ~(1 \leq j \leq n_2, 1 \leq l \leq n_3)$,
the outer product of $\mathbf{a}$ and $\mathbf{X}$, denoted by $\mathcal{T}=\mathbf{a}\circ\mathbf{X}$, is a
third-order tensor with size $n_1 \times n_2 \times n_3$, and its each entry is given by
$$
\mathcal{T}_{ijl}=a_ix_{jl}. 
$$
For the mode-2 unfolding of $\mathcal{T}$, denoted $\mathcal{T}_{(2)}
\in \mathbb{R}^{n_2 \times n_1 n_3}$,
the corresponding element-wise representation is given by
\begin{equation}\label{product1}
(\mathcal{T}_{(2)})_{j,(l-1) n_3+i}=\mathcal{T}_{i j  l}=a_ix_{j l}.
\end{equation}
Let $\mathbf{B}=\mathbf{X} \otimes \mathbf{a}^{H}\in \mathbb{R}^{n_2 \times n_1 n_3}$.
Then the $(j,(l-1)n_3+i)$-th element of $\mathbf{B}$, denoted by $\mathbf{B}_{j,(l-1) n_3+i}$, is given by
$$
\mathbf{B}_{j,(l-1) n_3+i}=x_{j l} a_i,
$$
which combined with (\ref{product1}) yields that
$
\mathcal{T}_{(2)}=(\mathbf{a}\circ\mathbf{X})_{(2)}=\mathbf{X}\otimes \mathbf{a}^{H}.
$
Similarly,
we can deduce that
$(\mathbf{a}\circ\mathbf{X})_{(3)}=\mathbf{X}^{H}\otimes \mathbf{a}^{H}$.
\qed

\begin{Prop}\label{out_kronFnorm}
For any vector $\mathbf{a}\in \mathbb{R}^{n_1}$ and matrix $\mathbf{X}\in \mathbb{R}^{n_2 \times n_3}$,
one has
$$
\|\mathbf{a}\circ\mathbf{X}\|_F=\|\mathbf{a}\|_2\|\mathbf{X}\|_F.
$$
\end{Prop}
\textbf{Proof}.\enspace
Combining Proposition \ref{out_kron1}  with
 Proposition \ref{kron_pro}, we have
$$
\|\mathbf{a}\circ\mathbf{X}\|_F=
\|(\mathbf{a}\circ\mathbf{X})_{(2)}\|_F=\|\mathbf{X}\otimes
\mathbf{a}^{H}\|_F=\|\mathbf{X}\|_F\|\mathbf{a}^{H}\|_2=\|\mathbf{X}\|_F\|\mathbf{a}\|_2.
$$
The proof is completed.
\qed

 Consider the weighted and undirected graph $G_i=(V_i,E_i,\widetilde{\mathbf{W}}_i), i=1,2,3$, where $V_i$ is the finite set
of vertices based on the row of $\mathbf{A}_{i}$ and $E_i$ the finite set of edges, and $\widetilde{\mathbf{W}}_i$ denotes the weighted adjacency
matrix.
Assume that the  graph $G_i$ is connected. Then we have the following inequality with respect to the graph regularization term.

\begin{lemma}\label{graph}
Denote $\Delta_{\mathbf{A}_{i}}:=\widehat{\mathbf{A}}_{i}-\mathbf{A}_{i}$.
Suppose that  Assumption \ref{assum_delta_A} holds and  the  graph $G_i$ is connected.
Then the matrix $\mathbf{A}_{i}\in\mathbb{R}^{m\times r_i},i=1,2,\mathbf{A}_{3}\in\mathbb{R}^{p\times r_3}$ satisfy the following inequality:
$$
\operatorname{tr}(\Delta_{\mathbf{A}_{i}}^{H}\mathbf{L}_{i}\Delta_{\mathbf{A}_{i}})
\geq \frac{\lambda_{2}(\mathbf{L}_{i})}{2}\|\widehat{\mathbf{A}}_{i}-\mathbf{A}_{i}\|_F^2,
$$
where $\lambda_{2}(\mathbf{L}_{i})>0$ denotes the second smallest eigenvalue of $\mathbf{L}_{i},i=1,2,3$.
\end{lemma}
\textbf{Proof}.
Let the eigenvalue decomposition of $\mathbf{L}_{1}\in\mathbb{R}^{m\times m}$ be
$\mathbf{L}_{1}=\mathbf{Q} \Lambda\mathbf{Q}^H$.
Define the matrix $\mathbf{M}:=\mathbf{Q}^H\Delta_{\mathbf{A}_{1}}\in\mathbb{R}^{m\times r_1}$.
Then  the $l$-th column of $\mathbf{M}$  is given by $\mathbf{M}_{:,l}=\mathbf{Q}^H (\Delta_{\mathbf{A}_{1}})_{:,l}$, where $(\Delta_{\mathbf{A}_{1}})_{:,l}$ denotes the $l$-th column of $\Delta_{\mathbf{A}_{1}}$.
Note that $\mathbf{L}_{1}$ is a symmetric semi-positive definite matrix with real and
non-negative eigenvalues.
Let $\lambda_j(\mathbf{L}_{1})$ be the $j$-th smallest eigenvalue of $\mathbf{L}_{1}$
with $\lambda_{1}(\mathbf{L}_{1}) \leq \lambda_{2}(\mathbf{L}_{1}) \leq \cdots \leq \lambda_m(\mathbf{L}_{1})$.
It can be easily demonstrated that
\begin{equation}\label{graph1}
\begin{split}
\operatorname{tr}\left(\Delta_{\mathbf{A}_{1}}^H \mathbf{L}_{1} \Delta_{\mathbf{A}_{1}}\right)
&=\sum_{l=1}^{r_1} (\Delta_{\mathbf{A}_{1}})_{:,l}^H \mathbf{L}_{1}  (\Delta_{\mathbf{A}_{1}})_{:,l} \\
& =\sum_{l=1}^{r_1} (\Delta_{\mathbf{A}_{1}})_{:,l}^H \mathbf{Q} \Lambda\mathbf{Q}^H  (\Delta_{\mathbf{A}_{1}})_{:,l}
=\sum_{l=1}^{r_1} \sum_{j=1}^m \lambda_j(\mathbf{L}_{1}) \mathbf{M}_{jl}^2,
\end{split}
\end{equation}
where $\mathbf{M}_{jl}$ denotes the $(j,l)$-th entry of $\mathbf{M}$.
By \cite[Section 3.1]{Yang2019} and \cite[Theorem 5]{Matthew2021},
we can conclude that
$\lambda_1(\mathbf{L}_{1})=0$ and $\lambda_{2}(\mathbf{L}_{1})>0$.
Then, we have
\begin{equation}\label{graph2}
\begin{aligned}
\sum_{j=1}^m \lambda_j(\mathbf{L}_{1}) \mathbf{M}_{jl}^2=\sum_{j=2}^m \lambda_j (\mathbf{L}_{1}) \mathbf{M}_{jl}^2
\geq \sum_{j=2}^m \lambda_2(\mathbf{L}_{1}) \mathbf{M}_{jl}^2
&=\sum_{j=1}^m \lambda_2(\mathbf{L}_{1}) \mathbf{M}_{jl}^2-\lambda_2 (\mathbf{L}_{1})\mathbf{M}_{1l}^2\\
&=\lambda_2(\mathbf{L}_{1})\left\|(\Delta_{\mathbf{A}_{1}})_{:,l}\right\|_2^2
-\lambda_2(\mathbf{L}_{1}) \mathbf{M}_{1l}^2.
\end{aligned}
\end{equation}
Substituting (\ref{graph2}) into (\ref{graph1}) yields
\begin{equation}\label{graph3}
\begin{aligned}
\operatorname{tr}\left(\Delta_{\mathbf{A}_{1}}^H \mathbf{L}_{1} \Delta_{\mathbf{A}_{1}}\right)
&\geq \sum_{l=1}^{r_1}\left(\lambda_2(\mathbf{L}_{1})\left\|(\Delta_{\mathbf{A}_{1}})_{:,l}\right\|_2^2
-\lambda_2(\mathbf{L}_{1}) \mathbf{M}_{1l}^2\right) \\
&=\lambda_2(\mathbf{L}_{1})\|\Delta_{\mathbf{A}_{1}}\|_F^2-\lambda_2(\mathbf{L}_{1})
\left\|(\mathbf{Q}_{:,1})^H \Delta_{\mathbf{A}_{1}}\right\|_2^2,
\end{aligned}
\end{equation}
where
the last equality follows from the fact
$\mathbf{M}=\mathbf{Q}^H\Delta_{\mathbf{A}_{1}}$.
Note that the first column of $\mathbf{Q}$ satisfies
$\mathbf{Q}_{:,1}=(1 / \sqrt{m}, 1 / \sqrt{m}, \ldots, 1 / \sqrt{m})^H$ \cite[Section 7]{Matthew2021}.
Consequently, we get that  $\left\|(\mathbf{Q}_{:,1})^H \Delta_{\mathbf{A}_{1}}\right\|_2^2=\sum_{l=1}^{r_1} \frac{1}{m}\left(\sum_{j=1}^m (\Delta_{\mathbf{A}_{1}})_{jl}\right)^2$.
Combining $\|\Delta_{\mathbf{A}_{1}}\|_F^2=\sum_{l=1}^{r_1} \sum_{j=1}^m (\Delta_{\mathbf{A}_{1}})_{jl}^2$
with (\ref{graph3}) yields
$$
\begin{aligned}
\operatorname{tr}\left(\Delta_{\mathbf{A}_{1}}^H \mathbf{L}_{1} \Delta_{\mathbf{A}_{1}}\right)
& \geq \lambda_2(\mathbf{L}_{1}) \sum_{l=1}^{r_1} \sum_{j=1}^m (\Delta_{\mathbf{A}_{1}})_{jl}^2-\lambda_2(\mathbf{L}_{1}) \sum_{l=1}^{r_1} \frac{1}{m}\big(\sum_{j=1}^m (\Delta_{\mathbf{A}_{1}})_{jl}\big)^2.
\end{aligned}
$$
By Assumption \ref{assum_delta_A},
we can conclude that
$\sum_{l=1}^{r_1}
\frac{1}{m}
(\sum_{j=1}^m (\Delta_{\mathbf{A}_{1}})_{j l})^2\leq\frac{1}{2} \sum_{l=1}^{r_1} \sum_{j=1}^m (\Delta_{\mathbf{A}_{1}})_{jl}^2$.
Hence, we obtain
$$
\begin{aligned}
\operatorname{tr}\left(\Delta_{\mathbf{A}_{1}}^H \mathbf{L}_{1} \Delta_{\mathbf{A}_{1}}\right) &\geq \lambda_{2}(\mathbf{L}_{1}) \sum_{l=1}^{r_1} \sum_{j=1}^m (\Delta_{\mathbf{A}_{1}})_{jl}^2
-\frac{\lambda_{2}(\mathbf{L}_{1})}{2} \sum_{l=1}^{r_1} \sum_{j=1}^m (\Delta_{\mathbf{A}_{1}})_{jl}^2\\
&=\frac{\lambda_{2}(\mathbf{L}_{1})}{2}\|\Delta_{\mathbf{A}_{1}}\|_F^2
=\frac{\lambda_{2}(\mathbf{L}_{1})}{2}\|\widehat{\mathbf{A}}_{1}-\mathbf{A}_{1}\|_F^2.
\end{aligned}
$$

By Assumption \ref{assum_delta_A},
using similar arguments leads to
\begin{equation}\label{trace_A3}
\begin{aligned}
\operatorname{tr}(\Delta_{\mathbf{A}_{2}}^{H}\mathbf{L}_{2}\Delta_{\mathbf{A}_{2}})
&\geq \frac{\lambda_{2}(\mathbf{L}_{2})}{2}\|\widehat{\mathbf{A}}_{2}-\mathbf{A}_{2}\|_F^2,\\
\operatorname{tr}(\Delta_{\mathbf{A}_{3}}^{H}\mathbf{L}_{3}\Delta_{\mathbf{A}_{3}})
&\geq \frac{\lambda_{2}(\mathbf{L}_{3})}{2}\|\widehat{\mathbf{A}}_{3}-\mathbf{A}_{3}\|_F^2,
\end{aligned}
\end{equation}
where $\lambda_{2}(\mathbf{L}_{2}), \lambda_{2}(\mathbf{L}_{3})>0$
represent the second smallest eigenvalues of $\mathbf{L}_{2}$
and $\mathbf{L}_{3}$, respectively.
This completes the proof. \qed

\begin{lemma}\label{inner_product}
Let $\widehat{\Delta}$ be defined in (\ref{DelDf}), $\mathbf{Z}$ and $\mathbf{e}$ be defined in (\ref{VAR_vec}).
Suppose that Assumptions \ref{assum_stable} and \ref{assum_E} hold, and the sample size $T\geq a_1\log(m^2p)$.
Then the following inequality holds with probability at least $1-6\exp(-a\log (m^2 p))$:
$$
\frac{1}{T}\langle\mathbf{Z}\widehat{\Delta},\mathbf{e}\rangle
\leq 2
\pi\vartheta a_2\sqrt{\log(m^2p)/T}\|\widehat{{\Delta}}\|_1,
$$
where $\vartheta:=
\lambda _{\max }(\mathbf{\Sigma}_{\boldsymbol{\varepsilon}})
( 1+ \mu _{\max }( \mathbf{W} )/\mu_{\min}(\mathbf{W}))$ and the constants $a>0, a_1, a_2\geq1$.
\end{lemma}
\textbf{Proof}.
By H\"older's inequality, we have that
\begin{equation}\label{sparse_1}
\frac{1}{T}\langle\mathbf{Z}\widehat{\Delta},\mathbf{e}\rangle
=\frac{1}{T}\langle \mathbf{Z}^H\mathbf{e}, \widehat{\Delta}\rangle
\leq\frac{1}{T}\|\mathbf{Z}^H\mathbf{e}\|_{\infty}\|\widehat{\Delta}\|_1.
\end{equation}
By a simple calculation, we can deduce
$\|\mathbf{Z}^H\mathbf{e}\|_{\infty}=\|\mathbf{X}^{H} \mathbf{E}\|_{\infty}=
\max _{1 \leq i \leq mp, 1 \leq j \leq m}\left|\mathbf{e}_i^{H} \mathbf{X}^{H}
\mathbf{E} \mathbf{e}_j\right|$,
where $\mathbf{e}_i$ denotes the standard basis vector with the $i$-th entry being $1$
and other  entries being $0$.
According to \cite[Proposition 2.4]{Basu2015} with Assumptions  \ref{assum_stable} and \ref{assum_E}, for any vectors $\mathbf{u}$ and $\mathbf{v}$ with
$\|\mathbf{u}\|_2\leq 1, \|\mathbf{v}\|_2\leq 1$ and any $\eta>0$, there exists a constant $a_3>0$ such that
$$
\begin{aligned}
\ & \mathbb{P}\left[\left|\mathbf{u}^{H}\left(\mathbf{X}^{H} \mathbf{E} / T\right) \mathbf{v}\right|>
2 \pi\left(\lambda_{\max }\left(\mathbf{\Sigma}_{\boldsymbol{\varepsilon}}\right)\left(1+\frac{\mu_{\max }(\mathbf{W})}{\mu_{\min }(\mathbf{W})}\right)\right) \eta\right] \\
\leq \ &   6 \exp (-a_3T \min (\eta, \eta^2)).
\end{aligned}
$$
Denote
$\vartheta:=\lambda_{\max }\left(\mathbf{\Sigma}_{\boldsymbol{\varepsilon}}\right)
\left[1+\mu_{\max }(\mathbf{W}) / \mu_{\min }(\mathbf{W})\right]$. By using a union bound, we get that
\begin{equation}\label{sparse_4}
\begin{aligned}
\ &\mathbb{P}\left[\max _{1 \leq i \leq mp, 1 \leq j \leq m}
\left|\mathbf{e}_i^{H} \mathbf{X}^{H} \mathbf{E} \mathbf{e}_j / T
\right|>2 \pi \vartheta \eta\right] \\
\leq \ & \sum _{1 \leq i \leq mp, 1 \leq j \leq m}\mathbb{P}\left[
\left|\mathbf{e}_i^{H} \mathbf{X}^{H} \mathbf{E} \mathbf{e}_j / T
\right|>2 \pi \vartheta \eta\right]\\
\leq  \ & 6 m^2 p
\exp (-a_3T \min (\eta, \eta^2)).
\end{aligned}
\end{equation}

Suppose that $a_3>1$,
and let $\eta=\sqrt{(\log (m^2 p)) / T}$,
where $T\geq\log(m^2p)$.
Then,
it is straightforwardly obtain that $\min (\eta, \eta^2)=(\log (m^2 p))/ T$.
As a consequence, we get
\begin{equation}\label{sparse_2}
\begin{aligned}
6 m^2 p
\exp (-a_3T \min (\eta, \eta^2))
&=6 \exp(\log(m^2p))\exp (-a_3\log(m^2p))\\
&=6 \exp (-(a_3-1)\log(m^2p)).
\end{aligned}
\end{equation}

Suppose that $a_3\in(0,1]$
and $T\geq\frac{2}{a_3}\log(m^2p)$.
Let $\eta=\sqrt{\frac{2}{a_3}}\sqrt{\frac{\log(m^2 p)}{T}}$, where $\eta\leq 1$.
Then, we have
\begin{equation}\label{sparse_3}
\begin{aligned}
6m^2 p
\exp (-a_3T \min (\eta, \eta^2))
&=6m^2 p \exp (-2\frac{a_3}{2}T \min (\eta, \eta^2))\\
&=6m^2 p \exp (-2T \frac{a_3}{2}\eta^2)\\
&=6m^2 p \exp (-2\log \left(m^2 p\right))\\
&= 6\exp (-\log \left(m^2 p\right)).
\end{aligned}
\end{equation}
By combining
(\ref{sparse_2}), (\ref{sparse_3}),
and the condition $T\geq a_1\log(m^2p)$,
the following inequality holds:
$$
\mathbb{P}\left[\max _{1 \leq i \leq mp, 1 \leq j \leq m}
\left|\mathbf{e}_i^{H}
\mathbf{X}^{H} \mathbf{E} \mathbf{e}_j / T\right|>2
\pi\vartheta a_2\sqrt{\log \left(m^2 p\right) / T}\right]
\leq 6\exp (-a\log (m^2 p)),
$$
where $a_1, a_2 \geq1$ and $a>0$.
This taken together with (\ref{sparse_1}) indicates that
$$
\frac{1}{T}\langle\mathbf{Z}\widehat{\Delta},\mathbf{e}\rangle
\leq\frac{1}{T}\|\mathbf{Z}^H\mathbf{e}\|_{\infty}\|\widehat{\Delta}\|_1
\leq 2
\pi\vartheta a_2\sqrt{\log(m^2p)/T}\|\widehat{{\Delta}}\|_1
$$
with probability at least  $1-6\exp(-a\log (m^2 p))$.
\qed

\begin{lemma}\label{Restricted eigenvalue}
	Let $\widehat{\Delta}$ and $\mathbf{X}$ be defined in (\ref{DelDf}) and (\ref{VAR_matrix}), respectively.
	Suppose that Assumptions
 \ref{assum_stable}, \ref{assum_E} and \ref{assum_G} hold, and
the sample size $T\geq\frac{
8\varsigma^2}
{a_4}
\bar{p}\min(\log (mp), \log(21emp/\bar{p}))$, where $\varsigma:=\frac{\lambda_{\min}
	\left(\boldsymbol{\Sigma}_{\boldsymbol{\varepsilon}}\right) / \mu_{\max }(\mathbf{W})}{\lambda_{\max}
	\left(\boldsymbol{\Sigma}_{\boldsymbol{\varepsilon}}\right) / \mu_{\min }(\mathbf{W})}$,
$a_4>0$ is a constant and  $\bar{p}$ depends on $\bar{s}$.
Then
$$
T^{-1}\left\|\left(\mathbf{I}_m \otimes \mathbf{X}\right)
\widehat{\Delta}\right\|_2^2 \geq \varpi\|\widehat{\Delta}\|_2^2 / 2
$$
holds with probability at least $1-2 \exp (-\bar{p}\min(\log mp, \log(21emp/\bar{p})))$,
where
$\varpi:=\lambda_{\min }\left(\mathbf{\Sigma}_{\boldsymbol{\varepsilon}}
\right)/ \mu_{\max }(\mathbf{W})$.
\end{lemma}
\textbf{Proof}.
For $\widehat{{\Delta}}$ in (\ref{DelDf}), we partition   it  into $m$ parts,
namely $\widehat{{\Delta}}=(\widehat{{\Delta}}_1^{H},
\ldots, \widehat{{\Delta}}_m^{H})^{H}$,
where $\widehat{{\Delta}}_i \in$ $\mathbb{R}^{mp}$.
Then we get
$$
\left\|\left(\mathbf{I}_m \otimes \mathbf{X}\right)
\widehat{{\Delta}}\right\|_2^2=\sum_{i=1}^m\left\|\mathbf{X}
\widehat{{\Delta}}_i\right\|_2^2.
$$
Denote $\widehat{\mathbf{\Gamma}}=\mathbf{X}^{H} \mathbf{X} / T$ and $\mathbf{\Gamma}=\mathbb{E} \widehat{\mathbf{\Gamma}}$.
It can be seen that
\begin{equation}\label{RE6}
	\begin{aligned}
		\widehat{{\Delta}}^{H}\left(\mathbf{I}_m \otimes \mathbf{\Gamma}\right) \widehat{{\Delta}}
		=\widehat{{\Delta}}^{H}\big(\mathbf{I}_m \otimes \mathbb{E} \widehat{\mathbf{\Gamma}}\big) \widehat{{\Delta}}
		&=\widehat{{\Delta}}^{H}\Big(\mathbf{I}_m \otimes \mathbb{E} \big(\mathbf{X}^{H} \mathbf{X}/T\big)\Big) \widehat{{\Delta}}\\
		&=\widehat{{\Delta}}^{H}\mathbb{E} \Big(\mathbf{I}_m \otimes \big(\mathbf{X}^{H} \mathbf{X}/T\big)\Big) \widehat{{\Delta}}\\
		&=\mathbb{E} \Big(\widehat{{\Delta}}^{H}\Big(\mathbf{I}_m \otimes \big(\mathbf{X}^{H} \mathbf{X}/T\big)\Big) \widehat{{\Delta}}\Big)\\
		&=T^{-1}\mathbb{E} \Big(\widehat{{\Delta}}^{H}\Big(\mathbf{I}_m \otimes \big(\mathbf{X}^{H} \mathbf{X}\big)\Big) \widehat{{\Delta}}\Big)\\
		&=T^{-1}\mathbb{E} \left(\|(\mathbf{I}_m
		\otimes\mathbf{X}) \widehat{{\Delta}} \|_2^2\right).
	\end{aligned}
\end{equation}
Additionally,
%
%
\begin{equation}\label{RE1}
\begin{aligned}
T^{-1}\|\left(\mathbf{I}_m \otimes \mathbf{X}\right) \widehat{{\Delta}}\|_2^2
& =\widehat{{\Delta}}^{H}(\mathbf{I}_m \otimes \widehat{\mathbf{\Gamma}}) \widehat{{\Delta}} \\
& =\widehat{{\Delta}}^{H}\left(\mathbf{I}_m \otimes \mathbf{\Gamma}\right) \widehat{{\Delta}}
+\widehat{{\Delta}}^{H}\left(\mathbf{I}_m \otimes(\widehat{\mathbf{\Gamma}}-\mathbf{\Gamma})\right) \widehat{{\Delta}} \\
& =\widehat{{\Delta}}^{H}\left(\mathbf{I}_m \otimes \mathbf{\Gamma}\right) \widehat{{\Delta}}
+\sum_{i=1}^m \widehat{{\Delta}}_i^{H}(\widehat{\mathbf{\Gamma}}-\mathbf{\Gamma}) \widehat{{\Delta}}_i\\
& =T^{-1} \mathbb{E}\left(\|(\mathbf{I}_m
\otimes\mathbf{X}) \widehat{{\Delta}} \|_2^2\right)
+\sum_{i=1}^m \widehat{{\Delta}}_i^{H}(\widehat{\mathbf{\Gamma}}-\mathbf{\Gamma}) \widehat{{\Delta}}_i,
\end{aligned}
\end{equation}
where the last equality follows from (\ref{RE6}).
Specifically,
it follows from the proof of Proposition 4.2 in \cite{Basu2015} and Assumptions \ref{assum_stable}, \ref{assum_E}  that
\begin{equation}\label{cov-spectral}
\lambda_{\min }(\mathbf{\Gamma})
\geq \lambda_{\min }\left(\boldsymbol{\Sigma}_{\boldsymbol{\varepsilon}}\right) /
\mu_{\max }(\mathbf{W}).
\end{equation}
Then, we have
\begin{equation}\label{RE2}
\begin{aligned}
T^{-1} \mathbb{E}\left(\|(\mathbf{I}_m
\otimes\mathbf{X}) \widehat{{\Delta}} \|_2^2\right)=
\widehat{{\Delta}}^{H}\left(\mathbf{I}_m \otimes \mathbf{\Gamma}\right)
\widehat{{\Delta}}
&\geq \lambda_{\min }(\mathbf{\Gamma})\|\widehat{{\Delta}}\|_2^2\\
&\geq \lambda_{\min}
\left(\boldsymbol{\Sigma}_{\boldsymbol{\varepsilon}}\right) / \mu_{\max }(\mathbf{W})
\|\widehat{{\Delta}}\|_2^2\\
&=\varpi\|\widehat{{\Delta}}\|_2^2,
\end{aligned}
\end{equation}
where
the first inequality follows from the Rayleigh\text{-}Ritz theorem \cite[Theorem 4.2.2]{Horn2013} and
$\lambda_{\min }(\mathbf{\Gamma})=\lambda_{\min }(\mathbf{I}_m \otimes \mathbf{\Gamma})$,
the equality follows from $\varpi:=\frac{\lambda_{\min}
\left(\boldsymbol{\Sigma}_{\boldsymbol{\varepsilon}}\right)}{\mu_{\max }(\mathbf{W})}$.
By Lemma \ref{G_sparsity},
for any $i=1,2,\ldots,m$,
we have $\|\widehat{{\Delta}}_i\|_0\leq q_imp$, where $q_i\in[0,1]$ is related to $\bar{s}$.
Let $\bar{p}=\max\{q_1mp,q_2mp,\ldots,q_m mp\}$.
We define $\mathcal{K}(\bar{p}):=\left\{\boldsymbol{v} \in \mathbb{R}^{mp}:
\|\boldsymbol{v}\|_0 \leq \bar{p},\|\boldsymbol{v}\|_2 \leq 1\right\}$ as the
set of $\bar{p}$-sparse vectors.
Note that
$$
\begin{aligned}
\ &T^{-1}\|\left(\mathbf{I}_m \otimes \mathbf{X}\right) \widehat{{\Delta}}\|_2^2\\
\geq \ & T^{-1} \mathbb{E}\left(\|(\mathbf{I}_m
\otimes\mathbf{X}) \widehat{{\Delta}} \|_2^2\right)
-\sup _{\widehat{{\Delta}}_i \in \mathcal{K}(\bar{p})}\left|T^{-1}\|\left(\mathbf{I}_m \otimes \mathbf{X}\right) \widehat{{\Delta}}\|_2^2
-T^{-1} \mathbb{E}\left(\|(\mathbf{I}_m
\otimes\mathbf{X}) \widehat{{\Delta}} \|_2^2\right)\right|\\
= \ &T^{-1} \mathbb{E}\left(\|(\mathbf{I}_m
\otimes\mathbf{X}) \widehat{{\Delta}} \|_2^2\right)
-\sup _{\widehat{{\Delta}}_i \in \mathcal{K}(\bar{p})}\left|\sum_{i=1}^m \widehat{{\Delta}}_i^{H}(\widehat{\mathbf{\Gamma}}-\mathbf{\Gamma}) \widehat{{\Delta}}_i\right|\\
\geq \  &\varpi\|\widehat{{\Delta}}\|_2^2
-\sup _{\widehat{{\Delta}}_i \in \mathcal{K}(\bar{p})}\left|\sum_{i=1}^m \widehat{{\Delta}}_i^{H}(\widehat{\mathbf{\Gamma}}-\mathbf{\Gamma}) \widehat{{\Delta}}_i\right|,
\end{aligned}
$$
where the equality comes from (\ref{RE1})
and the second inequality comes from (\ref{RE2}).
Then, we have
$$
T^{-1}\|\left(\mathbf{I}_m \otimes \mathbf{X}\right) \widehat{{\Delta}}\|_2^2-
\frac{\varpi}{2}\|\widehat{{\Delta}}\|_2^2\geq
\frac{\varpi}{2}\|\widehat{{\Delta}}\|_2^2
-\sup _{\widehat{{\Delta}}_i \in \mathcal{K}(\bar{p})}\left|\sum_{i=1}^m \widehat{{\Delta}}_i^{H}(\widehat{\mathbf{\Gamma}}-\mathbf{\Gamma}) \widehat{{\Delta}}_i\right|,
$$
which readily implies that
\begin{equation}\label{RE3}
\begin{aligned}
 \ &\mathbb{P}\left[T^{-1}\|\left(\mathbf{I}_m \otimes \mathbf{X}\right) \widehat{{\Delta}}\|_2^2-
\frac{\varpi}{2}\|\widehat{{\Delta}}\|_2^2\geq0\right]\\
\geq \ &
\mathbb{P}\left[
\frac{\varpi}{2}\|\widehat{{\Delta}}\|_2^2
-\sup _{\widehat{{\Delta}}_i \in \mathcal{K}(\bar{p})}\left|\sum_{i=1}^m \widehat{{\Delta}}_i^{H}(\widehat{\mathbf{\Gamma}}-\mathbf{\Gamma}) \widehat{{\Delta}}_i\right|\geq0\right]\\
= \ &
\mathbb{P}\left[
\frac{\varpi}{2}\sum_{i=1}^m\|\widehat{{\Delta}}_i\|_2^2
-\sup _{\widehat{{\Delta}}_i \in \mathcal{K}(\bar{p})}\left|\sum_{i=1}^m \widehat{{\Delta}}_i^{H}(\widehat{\mathbf{\Gamma}}-\mathbf{\Gamma}) \widehat{{\Delta}}_i\right|\geq0\right].
\end{aligned}
\end{equation}
For stable VAR($p$) processes,
as stated in \cite[Section 4]{Basu2015},
the spectral density function $f_{\boldsymbol{X}}(\theta)$
is defined as $f_{\boldsymbol{X}}(\theta):=\frac{1}{2 \pi}\left(\mathbf{W}^{-1}\left(e^{-i \theta}\right)\right)\boldsymbol{\Sigma}_{\boldsymbol{\varepsilon}}\left(\mathbf{W}^{-1}\left(e^{-i \theta}\right)\right)^H$ for $\theta \in[-\pi, \pi]$.
We then define $\mathcal{M}\left(f_{\boldsymbol{X}}\right):=\underset{\theta \in[-\pi, \pi]}{\operatorname{ess} \sup } \lambda_{\max }\left(f_{\boldsymbol{X}}(\theta)\right)$ \cite[Assumption 2.1]{Basu2015}.
According to \cite[Proposition 2.4]{Basu2015},
for any $\eta>0$, $\|\mathbf{u}\|_2\leq 1$, there exists a constant $a_4>0$
such that
\begin{equation}\label{RE4}
\begin{aligned}
\mathbb{P}\left[\left|\mathbf{u}^{H}
(\widehat{\mathbf{\Gamma}}-\mathbf{\Gamma}) \mathbf{u}\right|>
2 \pi \mathcal{M}\left(f_{\boldsymbol{X}}\right) \eta\right]
\leq 2 \exp \left(- a_4 T \min (\eta, \eta^2)\right).
\end{aligned}
\end{equation}
Let $\eta:=\frac{\lambda_{\min}
\left(\boldsymbol{\Sigma}_{\boldsymbol{\varepsilon}}\right) / \mu_{\max }(\mathbf{W})}{4 \pi \mathcal{M}\left(f_{\boldsymbol{X}}\right) }$.
Therefore, we get that
\begin{equation}
\begin{aligned}
\mathbb{P}\left[\left|\mathbf{u}^{H}
(\widehat{\mathbf{\Gamma}}-\mathbf{\Gamma}) \mathbf{u}\right|>
\frac{\varpi}{2}\right]
\leq 2 \exp \left(- a_4 T \min (\eta, \eta^2)\right).
\end{aligned}
\end{equation}
By \cite[Proposition 2.3]{Basu2015},
we have
$2 \pi \mathcal{M}\left(f_{\boldsymbol{X}}\right)\geq\lambda_{\min }(\mathbf{\Gamma})$.
This combined with (\ref{cov-spectral}) yields
$
2 \pi \mathcal{M}\left(f_{\boldsymbol{X}}\right)\geq\lambda_{\min }(\mathbf{\Gamma})
\geq \lambda_{\min }\left(\boldsymbol{\Sigma}_{\boldsymbol{\varepsilon}}\right) /
\mu_{\max }(\mathbf{W}).
$
As a consequence, $\eta\leq 1/2$.
It can be shown that (\ref{RE3}) together with (\ref{RE4}) implies
\begin{equation}\label{RE5}
\begin{aligned}
\ & \mathbb{P}\left[T^{-1}\left\|\left(\mathbf{I}_m\otimes \mathbf{X}\right)
\widehat{{\Delta}}\right\|_2^2\geq \frac{\varpi}{2}\|\widehat{{\Delta}}\|_2^2
\right] \\
\geq  \ & \mathbb{P}\left[
\sup _{\widehat{{\Delta}}_i \in \mathcal{K}(\bar{p})}\left|\sum_{i=1}^m \widehat{{\Delta}}_i^{H}(\widehat{\mathbf{\Gamma}}-\mathbf{\Gamma}) \widehat{{\Delta}}_i\right|\leq
\frac{\varpi}{2}\sum_{i=1}^m\|\widehat{{\Delta}}_i\|_2^2
\right] \\
\geq  \ & \mathbb{P}\left[
\sup _{\widehat{{\Delta}}_i \in \mathcal{K}(\bar{p})}\left|\widehat{{\Delta}}_i^{H}(\widehat{\mathbf{\Gamma}}-\mathbf{\Gamma})
\widehat{{\Delta}}_i\right|\leq\frac{\varpi}{2}\|\widehat{{\Delta}}_i\|_2^2
\right] \\
= \ & \mathbb{P}\left[
\sup _{\widehat{{\Delta}}_i \in \mathcal{K}(\bar{p})}\left|\frac{\widehat{{\Delta}}_i^{H}}{\|\widehat{{\Delta}}_i\|_2}(\widehat{\mathbf{\Gamma}}
-\mathbf{\Gamma})
\frac{\widehat{{\Delta}}_i}{\|\widehat{{\Delta}}_i\|_2}\right|
\leq
\frac{\varpi}{2}\right] \\
\geq \ & 1-2 \exp \left(- a_4 T \frac{\varpi^2}{
(4 \pi \mathcal{M}\left(f_{\boldsymbol{X}}\right) )^2}
+\bar{p}\min(\log (mp), \log(21emp/\bar{p}))
\right)\\
\geq \ & 1-2 \exp \left(- a_4 T \frac{\varpi^2}{
4(\lambda_{\max}
\left(\boldsymbol{\Sigma}_{\boldsymbol{\varepsilon}}\right) / \mu_{\min }(\mathbf{W}))^2}
+\bar{p}\min(\log (mp), \log(21emp/\bar{p}))
\right),
\end{aligned}
\end{equation}
where the third inequality follows from \cite[Lemma F.2]{Basu2015}
with $\eta=\frac{\lambda_{\min}
\left(\boldsymbol{\Sigma}_{\boldsymbol{\varepsilon}}\right) / \mu_{\max }(\mathbf{W})}{4 \pi \mathcal{M}\left(f_{\boldsymbol{X}}\right) }=\frac{\varpi}{4 \pi \mathcal{M}\left(f_{\boldsymbol{X}}\right) }$
and the last inequality follows from $2\pi \mathcal{M}\left(f_{\boldsymbol{X}}\right)\leq  \lambda_{\max}
\left(\boldsymbol{\Sigma}_{\boldsymbol{\varepsilon}}\right) / \mu_{\min }(\mathbf{W})$ \cite[equation (4.1)]{Basu2015}.

Suppose that $T\geq  \frac{
8(\lambda_{\max}
\left(\boldsymbol{\Sigma}_{\boldsymbol{\varepsilon}}\right) / \mu_{\min }(\mathbf{W}))^2}
{a_4\varpi^2}
\bar{p}\min\left(\log (mp), \log(21emp/\bar{p})\right)$,
then
$$
\begin{aligned}
\ &1-2 \exp \left(- a_4 T \frac{\varpi^2}{
4(\lambda_{\max}
\left(\boldsymbol{\Sigma}_{\boldsymbol{\varepsilon}}\right) / \mu_{\min }(\mathbf{W}))^2}
+\bar{p}\min(\log (mp), \log(21emp/\bar{p}))
\right)\\
\geq \  &1-2 \exp \left(-\bar{p}\min(\log (mp), \log(21emp/\bar{p}))
\right).
\end{aligned}
$$
This combined with (\ref{RE5}) yields
$$
\mathbb{P}\left[T^{-1}\left\|\left(\mathbf{I}_m\otimes \mathbf{X}\right)
\widehat{{\Delta}}\right\|_2^2\geq \frac{\varpi}{2}\|\widehat{{\Delta}}\|_2^2
\right]\geq 1-2 \exp (-\bar{p}\min(\log mp, \log(21emp/\bar{p}))).
$$
This completes the proof.
\qed

\begin{lemma}\label{delta_A}
Let
$\widehat{\Delta}$ be defined in (\ref{DelDf}).
Then
$$
\|\widehat{\Delta}\|_2^2\leq2 \|\operatorname{vec}
(\widehat{\mathcal{G}}_{(1)}^{H})-\operatorname{vec}
(\mathcal{G}_{(1)}^{H})\|_2^2 \\
+6 c^2 r_1 r_2 r_3\left(\|\widehat{\mathbf{A}}_3-\mathbf{A}_3\|_F^2+
\|\widehat{\mathbf{A}}_2-\mathbf{A}_2\|_F^2
+\|\widehat{\mathbf{A}}_1-\mathbf{A}_1\|_F^2\right),
$$
where
$\widehat{\mathcal{G}}, \mathcal{G} \in \mathbb{R}^{r_1 \times r_2 \times r_3}$ with
$\|\widehat{\mathcal{G}}\|_\infty \leq c, \|\mathcal{G}\|_\infty \leq c$,  $c$ is defined in (\ref{test3}),
$\widehat{\mathbf{A}}_i, \mathbf{A}_i \in \mathbb{R}^{n_i \times r_i},$
and $\widehat{\mathbf{A}}_i^H\widehat{\mathbf{A}}_i=\mathbf{I}_{r_i},
\mathbf{A}_i^H\mathbf{A}_i=\mathbf{I}_{r_i},i=1,2,3$.
\end{lemma}
\textbf{Proof}.
From the definitions of $\widehat{\mathbf{A}}, \mathbf{A}, \widehat{\mathcal{G}}, \mathcal{G}$,
it can be easily seen that
$$
\begin{aligned}
0\leq\|\widehat{\Delta}\|_2=& \left\|(\widehat{\mathbf{A}}_1\otimes\widehat{\mathbf{A}}_3\otimes
\widehat{\mathbf{A}}_2)\operatorname{vec}
(\widehat{\mathcal{G}}_{(1)}^{H})-(\mathbf{A}_1\otimes\mathbf{A}_3\otimes
\mathbf{A}_2)\operatorname{vec}
(\mathcal{G}_{(1)}^{H})\right\|_2 \\
=& \left\|(\widehat{\mathbf{A}}_1\otimes\widehat{\mathbf{A}}_3\otimes
\widehat{\mathbf{A}}_2)\operatorname{vec}
(\widehat{\mathcal{G}}_{(1)}^{H})
-(\mathbf{A}_1\otimes\mathbf{A}_3\otimes
\mathbf{A}_2)\operatorname{vec}(\widehat{\mathcal{G}}_{(1)}^{H})
\right.\\
&+\left.
(\mathbf{A}_1\otimes\mathbf{A}_3\otimes
\mathbf{A}_2)\operatorname{vec}
(\widehat{\mathcal{G}}_{(1)}^{H})-
(\mathbf{A}_1\otimes\mathbf{A}_3\otimes\mathbf{A}_2)\operatorname{vec}
(\mathcal{G}_{(1)}^{H})\right\|_2 \\
\leq& \left\|(\widehat{\mathbf{A}}_1\otimes\widehat{\mathbf{A}}_3\otimes
\widehat{\mathbf{A}}_2)\operatorname{vec}(\widehat{\mathcal{G}}_{(1)}^{H})
-(\mathbf{A}_1\otimes\mathbf{A}_3\otimes
\mathbf{A}_2)\operatorname{vec}
(\widehat{\mathcal{G}}_{(1)}^{H})\right\|_2 \\
& +\left\|(\mathbf{A}_1\otimes\mathbf{A}_3\otimes
\mathbf{A}_2)\operatorname{vec}(\widehat{\mathcal{G}}_{(1)}^{H})
-(\mathbf{A}_1\otimes\mathbf{A}_3\otimes
\mathbf{A}_2)\operatorname{vec}(\mathcal{G}_{(1)}^{H})\right\|_2.
\end{aligned}
$$
Then, we deduce
\begin{equation}\label{delta_A_1}
\begin{aligned}
 \|\widehat{\Delta}\|_2^2
\leq  & \left(\left\|(\widehat{\mathbf{A}}_1\otimes\widehat{\mathbf{A}}_3\otimes
\widehat{\mathbf{A}}_2)\operatorname{vec}(\widehat{\mathcal{G}}_{(1)}^{H})
-(\mathbf{A}_1\otimes\mathbf{A}_3\otimes
\mathbf{A}_2)\operatorname{vec}
(\widehat{\mathcal{G}}_{(1)}^{H})\right\|_2
\right.\\
& \ +\left.
\left\|(\mathbf{A}_1\otimes\mathbf{A}_3\otimes
\mathbf{A}_2)\operatorname{vec}(\widehat{\mathcal{G}}_{(1)}^{H})
-(\mathbf{A}_1\otimes\mathbf{A}_3\otimes
\mathbf{A}_2)\operatorname{vec}(\mathcal{G}_{(1)}^{H})\right\|_2\right)^2\\
\leq \ & 2\left\|(\widehat{\mathbf{A}}_1\otimes\widehat{\mathbf{A}}_3\otimes
\widehat{\mathbf{A}}_2)\operatorname{vec}(\widehat{\mathcal{G}}_{(1)}^{H})
-(\mathbf{A}_1\otimes\mathbf{A}_3\otimes
\mathbf{A}_2)\operatorname{vec}
(\widehat{\mathcal{G}}_{(1)}^{H})\right\|_2^2 \\
& \  +2\left\|(\mathbf{A}_1\otimes\mathbf{A}_3\otimes
\mathbf{A}_2)\operatorname{vec}(\widehat{\mathcal{G}}_{(1)}^{H})
-(\mathbf{A}_1\otimes\mathbf{A}_3\otimes
\mathbf{A}_2)\operatorname{vec}(\mathcal{G}_{(1)}^{H})\right\|_2^2,
\end{aligned}
\end{equation}
where the last inequality follows from the fact that $\frac{a+b}{2} \leq \sqrt{\frac{a^2+b^2}{2}}$
for any $a,b\geq0$.
Note that
\begin{equation}\label{delta_A_2}
\begin{aligned}
&\left\|(\mathbf{A}_1\otimes\mathbf{A}_3\otimes\mathbf{A}_2)\operatorname{vec}
(\widehat{\mathcal{G}}_{(1)}^{H})
-(\mathbf{A}_1\otimes\mathbf{A}_3\otimes
\mathbf{A}_2)\operatorname{vec}(\mathcal{G}_{(1)}^{H})\right\|_2\\
= &\left\|(\mathbf{A}_1\otimes\mathbf{A}_3\otimes
\mathbf{A}_2)(\operatorname{vec}(\widehat{\mathcal{G}}_{(1)}^{H})
-\operatorname{vec}
(\mathcal{G}_{(1)}^{H}))\right\|_2\\
\leq &\left\|(\mathbf{A}_1\otimes\mathbf{A}_3\otimes
\mathbf{A}_2)\right\|
\left\|\operatorname{vec}
(\widehat{\mathcal{G}}_{(1)}^{H})
-\operatorname{vec}
(\mathcal{G}_{(1)}^{H})\right\|_2\\
= &\left\|\mathbf{A}_3\right\|
\left\|\mathbf{A}_2\right\|
\left\|\mathbf{A}_1\right\|
\left\|\operatorname{vec}
(\widehat{\mathcal{G}}_{(1)}^{H})-\operatorname{vec}
(\mathcal{G}_{(1)}^{H})\right\|_2\\
\leq &\left\|\operatorname{vec}
(\widehat{\mathcal{G}}_{(1)}^{H})-\operatorname{vec}
(\mathcal{G}_{(1)}^{H})\right\|_2,
\end{aligned}
\end{equation}
where the first inequality follows from \cite[Lemma 2.1]{Chen2021},
the second equality follows from Proposition \ref{kron_pro}
and the second inequality follows from Proposition \ref{kron_untari}.

On the other hand, observe that
\begin{equation}\label{delta_A_3}
\begin{aligned}
& \left\|(\widehat{\mathbf{A}}_1\otimes\widehat{\mathbf{A}}_3\otimes
\widehat{\mathbf{A}}_2)\operatorname{vec}(\widehat{\mathcal{G}}_{(1)}^{H})
-(\mathbf{A}_1\otimes\mathbf{A}_3\otimes
\mathbf{A}_2)\operatorname{vec}(\widehat{\mathcal{G}}_{(1)}^{H})\right\|_2\\
=& \left\|(\widehat{\mathbf{A}}_1\otimes\widehat{\mathbf{A}}_3\otimes
\widehat{\mathbf{A}}_2)\operatorname{vec}(\widehat{\mathcal{G}}_{(1)}^{H})
-(\widehat{\mathbf{A}}_1\otimes\widehat{\mathbf{A}}_3\otimes
\mathbf{A}_2)\operatorname{vec}(\widehat{\mathcal{G}}_{(1)}^{H})
\right.\\
&+\left.(\widehat{\mathbf{A}}_1\otimes\widehat{\mathbf{A}}_3\otimes
\mathbf{A}_2)\operatorname{vec}
(\widehat{\mathcal{G}}_{(1)}^{H})
- (\widehat{\mathbf{A}}_1\otimes\mathbf{A}_3\otimes
\mathbf{A}_2)\operatorname{vec}(\widehat{\mathcal{G}}_{(1)}^{H})
\right.\\
&+\left.
(\widehat{\mathbf{A}}_1\otimes\mathbf{A}_3\otimes
\mathbf{A}_2)\operatorname{vec}
(\widehat{\mathcal{G}}_{(1)}^{H})
-(\mathbf{A}_1\otimes\mathbf{A}_3\otimes
\mathbf{A}_2)\operatorname{vec}
(\widehat{\mathcal{G}}_{(1)}^{H})\right\|_2 \\
\leq & \left\|(\widehat{\mathbf{A}}_1\otimes\widehat{\mathbf{A}}_3\otimes
\widehat{\mathbf{A}}_2)\operatorname{vec}
(\widehat{\mathcal{G}}_{(1)}^{H})-(\widehat{\mathbf{A}}_1\otimes\widehat{\mathbf{A}}_3\otimes
\mathbf{A}_2)\operatorname{vec}
(\widehat{\mathcal{G}}_{(1)}^{H})\right\|_2\\
& +\left\|(\widehat{\mathbf{A}}_1\otimes\widehat{\mathbf{A}}_3\otimes
\mathbf{A}_2)\operatorname{vec}
(\widehat{\mathcal{G}}_{(1)}^{H}) - (\widehat{\mathbf{A}}_1\otimes\mathbf{A}_3\otimes
\mathbf{A}_2)\operatorname{vec}
(\widehat{\mathcal{G}}_{(1)}^{H})\right\|_2 \\
&+\left\|(\widehat{\mathbf{A}}_1\otimes\mathbf{A}_3\otimes
\mathbf{A}_2)\operatorname{vec}
(\widehat{\mathcal{G}}_{(1)}^{H})-(\mathbf{A}_1\otimes\mathbf{A}_3\otimes
\mathbf{A}_2)\operatorname{vec}
(\widehat{\mathcal{G}}_{(1)}^{H})\right\|_2.
\end{aligned}
\end{equation}
For the first term on the right-hand side of the above inequality,
it can be seen that
\begin{equation}\label{delta_A_4}
\begin{aligned}
&\left\|(\widehat{\mathbf{A}}_1\otimes\widehat{\mathbf{A}}_3\otimes
\widehat{\mathbf{A}}_2)\operatorname{vec}
(\widehat{\mathcal{G}}_{(1)}^{H})-(\widehat{\mathbf{A}}_1\otimes\widehat{\mathbf{A}}_3\otimes
\mathbf{A}_2)\operatorname{vec}
(\widehat{\mathcal{G}}_{(1)}^{H})\right\|_2\\
= & \left\|(\widehat{\mathbf{A}}_1\otimes\widehat{\mathbf{A}}_3\otimes
(\widehat{\mathbf{A}}_2-\mathbf{A}_2))\operatorname{vec}
(\widehat{\mathcal{G}}_{(1)}^{H})\right\|_2\\
\leq & \left\|\widehat{\mathbf{A}}_1\otimes\widehat{\mathbf{A}}_3\otimes
(\widehat{\mathbf{A}}_2-\mathbf{A}_2)\right\|\|\operatorname{vec}
(\widehat{\mathcal{G}}_{(1)}^{H})\|_2\\
= \ & \|\widehat{\mathbf{A}}_1\|
\|\widehat{\mathbf{A}}_3\|\|\widehat{\mathbf{A}}_2-\mathbf{A}_2\|
\|\operatorname{vec}
(\widehat{\mathcal{G}}_{(1)}^{H})\|_2\\
\leq  \ & \|\widehat{\mathbf{A}}_2-\mathbf{A}_2\|
\|\operatorname{vec}
(\widehat{\mathcal{G}}_{(1)}^{H})\|_2\\
=  \ & \|\widehat{\mathbf{A}}_2-\mathbf{A}_2\|
\|\widehat{\mathcal{G}}\|_F\\
\leq \  & c \sqrt{r_1 r_2 r_3}\|\widehat{\mathbf{A}}_2-\mathbf{A}_2\|\\
\leq \ & c \sqrt{r_1 r_2 r_3}\|\widehat{\mathbf{A}}_2-\mathbf{A}_2\|_F,
\end{aligned}
\end{equation}
where the first inequality comes from \cite[Lemma 2.1]{Chen2021},
the second equality comes from Proposition \ref{kron_pro},
the second inequality comes from Proposition \ref{kron_untari},
the third inequality comes from
$\|\widehat{\mathcal{G}}\|_F\leq \sqrt{r_1r_2r_3}\|\widehat{\mathcal{G}}\|_{\infty}\leq \sqrt{r_1r_2r_3}c$,
and the last inequality comes from
$\|\mathbf{A}\|\leq \|\mathbf{A}\|_F$ for any matrix $\mathbf{A}$.
Similarly,
\begin{equation}\label{delta_A_5}
\begin{aligned}
\left\|(\widehat{\mathbf{A}}_1\otimes\widehat{\mathbf{A}}_3\otimes
\mathbf{A}_2)\operatorname{vec}
(\widehat{\mathcal{G}}_{(1)}^{H}) - (\widehat{\mathbf{A}}_1\otimes\mathbf{A}_3\otimes
\mathbf{A}_2)\operatorname{vec}
(\widehat{\mathcal{G}}_{(1)}^{H})\right\|_2
&\leq c \sqrt{r_1 r_2 r_3}\|\widehat{\mathbf{A}}_3-\mathbf{A}_3\|_F,\\
\left\|(\widehat{\mathbf{A}}_1\otimes\mathbf{A}_3\otimes
\mathbf{A}_2)\operatorname{vec}
(\widehat{\mathcal{G}}_{(1)}^{H})-(\mathbf{A}_1\otimes\mathbf{A}_3\otimes
\mathbf{A}_2)\operatorname{vec}
(\widehat{\mathcal{G}}_{(1)}^{H})\right\|_2
&\leq c \sqrt{r_1 r_2 r_3}\|\widehat{\mathbf{A}}_1-\mathbf{A}_1\|_F.
\end{aligned}
\end{equation}
Substituting (\ref{delta_A_4}), (\ref{delta_A_5}) into (\ref{delta_A_3}),
we obtain
$$
\begin{aligned}
0&\leq \left\|(\widehat{\mathbf{A}}_1\otimes\widehat{\mathbf{A}}_3\otimes
\widehat{\mathbf{A}}_2)\operatorname{vec}
(\widehat{\mathcal{G}}_{(1)}^{H})-(\mathbf{A}_1\otimes\mathbf{A}_3\otimes
\mathbf{A}_2)\operatorname{vec}
(\widehat{\mathcal{G}}_{(1)}^{H})\right\|_2\\
& \leq c \sqrt{r_1 r_2 r_3}(\|\widehat{\mathbf{A}}_3-\mathbf{A}_3\|_F+\|\widehat{\mathbf{A}}_2-\mathbf{A}_2\|_F+ \|\widehat{\mathbf{A}}_1-\mathbf{A}_1\|_F).
\end{aligned}
$$
Further, we have
\begin{equation}\label{delta_A_6}
\begin{aligned}
\ & \left\|(\widehat{\mathbf{A}}_1\otimes\widehat{\mathbf{A}}_3\otimes
\widehat{\mathbf{A}}_2)\operatorname{vec}
(\widehat{\mathcal{G}}_{(1)}^{H})-(\mathbf{A}_1\otimes\mathbf{A}_3\otimes
\mathbf{A}_2)\operatorname{vec}
(\widehat{\mathcal{G}}_{(1)}^{H})\right\|_2^2\\
\leq \ &c^2r_1 r_2 r_3\left(\|\widehat{\mathbf{A}}_3-\mathbf{A}_3\|_F+
\|\widehat{\mathbf{A}}_2-\mathbf{A}_2\|_F+
\|\widehat{\mathbf{A}}_1-\mathbf{A}_1\|_F\right)^2\\
= \ &c^2r_1 r_2 r_3\left(\|\widehat{\mathbf{A}}_3-\mathbf{A}_3\|_F^2+
\|\widehat{\mathbf{A}}_2-\mathbf{A}_2\|_F^2+
\|\widehat{\mathbf{A}}_1-\mathbf{A}_1\|_F^2
+2\|\widehat{\mathbf{A}}_3-\mathbf{A}_3\|_F\|\widehat{\mathbf{A}}_2-\mathbf{A}_2\|_F
\right.\\
 \  &+\left.
2\|\widehat{\mathbf{A}}_3-\mathbf{A}_3\|_F\|\widehat{\mathbf{A}}_1-\mathbf{A}_1\|_F
+2\|\widehat{\mathbf{A}}_1-\mathbf{A}_1\|_F\|\widehat{\mathbf{A}}_2-\mathbf{A}_2\|_F\right).
\end{aligned}
\end{equation}
For any non-negative numbers $a,b,c$, it is known that
$a^2+b^2+c^2\geq ab+ac+bc$.
Therefore, we can conclude that
$$
\begin{aligned}
\ &\|\widehat{\mathbf{A}}_3-\mathbf{A}_3\|_F\|\widehat{\mathbf{A}}_2-\mathbf{A}_2\|_F
+\|\widehat{\mathbf{A}}_3-\mathbf{A}_3\|_F\|\widehat{\mathbf{A}}_1-\mathbf{A}_1\|_F
+\|\widehat{\mathbf{A}}_1-\mathbf{A}_1\|_F\|\widehat{\mathbf{A}}_2-\mathbf{A}_2\|_F\\
\leq \ & \|\widehat{\mathbf{A}}_3-\mathbf{A}_3\|_F^2+
\|\widehat{\mathbf{A}}_2-\mathbf{A}_2\|_F^2+
\|\widehat{\mathbf{A}}_1-\mathbf{A}_1\|_F^2,
\end{aligned}
$$
which taken collectively with (\ref{delta_A_6}) gives
\begin{equation}\label{delta_A_7}
\begin{aligned}
 \ & \left\|(\widehat{\mathbf{A}}_3\otimes\widehat{\mathbf{A}}_2\otimes
\widehat{\mathbf{A}}_1)\operatorname{vec}
(\widehat{\mathcal{G}}_{(1)}^{H})-(\mathbf{A}_3\otimes\mathbf{A}_2\otimes
\mathbf{A}_1)\operatorname{vec}
(\widehat{\mathcal{G}}_{(1)}^{H})\right\|_2^2\\
\leq \ &c^2r_1 r_2 r_3\left(3\|\widehat{\mathbf{A}}_3-\mathbf{A}_3\|_F^2+
3\|\widehat{\mathbf{A}}_2-\mathbf{A}_2\|_F^2+
3\|\widehat{\mathbf{A}}_1-\mathbf{A}_1\|_F^2\right).
\end{aligned}
\end{equation}
Substituting (\ref{delta_A_7}), (\ref{delta_A_2}) into (\ref{delta_A_1}) yields
$$
\|\widehat{\Delta}\|_2^2\leq 2\|\operatorname{vec}
(\widehat{\mathcal{G}}_{(1)}^{H})-\operatorname{vec}
(\mathcal{G}_{(1)}^{H})\|_2^2
+6 c^2 r_1 r_2 r_3\left(\|\widehat{\mathbf{A}}_3-\mathbf{A}_3\|_F^2+
\|\widehat{\mathbf{A}}_2-\mathbf{A}_2\|_F^2
+\|\widehat{\mathbf{A}}_1-\mathbf{A}_1\|_F^2\right).
$$
This completes the proof.
\qed

\begin{lemma}\label{HOSVD perturbation bound}
Let $\mathcal{W}=\llbracket \mathcal{G} ;
\mathbf{A}_1, \mathbf{A}_2, \mathbf{A}_3 \rrbracket$ and
$\widetilde{\mathcal{W}}=\llbracket \widetilde{\mathcal{G}} ; \widetilde{\mathbf{A}}_1,
\widetilde{\mathbf{A}}_2, \widetilde{\mathbf{A}}_3 \rrbracket$ be the HOSVD of
$\mathcal{W}\in\mathbb{R}^{m\times m\times p}$ and $\widetilde{\mathcal{W}}\in\mathbb{R}^{m\times m\times p}$ with the same multilinear ranks $\left(r_1, r_2, r_3\right)$, respectively.
Suppose that Assumption \ref{assum_spectral_A} holds for $\mathcal{W}$.
Furthermore, assume that
$\|\widetilde{\mathcal{G}}\|_\infty\leq c$ and $\|\mathcal{G}\|_\infty\leq c$.
Then, we have
$$
\|\widetilde{\mathcal{G}}-\mathcal{G}\|_F \leq b\|\widetilde{\mathcal{W}}-\mathcal{W}\|_F,
$$
and
$$
\|\widetilde{\mathbf{A}}_i-\mathbf{A}_i\|_{F}
\leq 2^{\frac{3}{2}}\left(\eta_i / \rho\right)\|\widetilde{\mathcal{W}}-\mathcal{W}\|_{F},
$$
where $\eta_i=\left(3\sigma_1({\mathcal {W}}_{(i)})+
c\sqrt{r_1r_2r_3}\right)
\sqrt{\sum_{j=1}^{r_i}\frac{1}
{\sigma_j^4
(\mathcal{W}_{(i)})}}$,
$b=1+
\frac{2^{\frac{3}{2}}\sigma_1({\mathcal{W}}_{(1)})\sum_{i=1}^{3}\eta_i}{\rho}$
and $\rho$ is defined in Assumption \ref{assum_spectral_A}.
\end{lemma}
\textbf{Proof.}
By \cite[Theorem 3]{Yu2015}, under Assumption \ref{assum_spectral_A},
for the $j$-th singular vector of $\mathcal{W}_{(i)}$,
\begin{equation}\label{HOSVD5}
\begin{aligned}
\left\|(\widetilde{\mathbf{A}}_{i})_{:,j}-(\mathbf{A}_{i})_{:,j}\right\|_{2}
&\leq  \frac{2^{\frac{3}{2}}\left(2 \sigma_1({\mathcal {W}}_{(i)})+
\|\widetilde{\mathcal{W}}_{(i)}-\mathcal{W}_{(i)}\|\right)
\|\widetilde{\mathcal{W}}_{(i)}-\mathcal{W}_{(i)}\|_{F}}
{\min \left[\sigma_{j-1}^2(\mathcal{W}_{(i)})-\sigma_j^2(\mathcal{W}_{(i)}),
\sigma_j^2(\mathcal{W}_{(i)})-\sigma_{j+1}^2(\mathcal{W}_{(i)})\right]}\\
&\leq \frac{2^{\frac{3}{2}}\left(2 \sigma_1({\mathcal {W}}_{(i)})+
\|\widetilde{\mathcal{W}}_{(i)}-\mathcal{W}_{(i)}\|\right)
\|\widetilde{\mathcal{W}}_{(i)}-\mathcal{W}_{(i)}\|_{F}}
{\rho \sigma_j^2
(\mathcal{W}_{(i)})},
\end{aligned}
\end{equation}
where $(\mathbf{A}_{i})_{:,j}$ denotes the
$j$-th left singular vector of the singular value decomposition of $\mathcal{W}_{(i)}$
,
$(\widetilde{\mathbf{A}}_{i})_{:,j}$ denotes the
$j$-th left singular vector of the singular value decomposition of $\widetilde{\mathcal{W}}_{(i)}$
and
$\sigma_1({\mathcal {W}}_{(i)})=\|{\mathcal {W}}_{(i)}\|$.
Additionally, let
the SVD of $\widetilde{\mathcal{G}}_{(1)}$ be
$\widetilde{\mathcal{G}}_{(1)}=\mathbf{\bar{U}}_1\mathbf{\bar{S}}_1\mathbf{\bar{V}}_1^H$, where  $\mathbf{\bar{U}}_1\in\mathbb{R}^{r_1\times r_1}$ and $\mathbf{\bar{V}}_1\in\mathbb{R}^{r_2r_3\times r_2r_3}$ are orthogonal matrices,
$\mathbf{\bar{S}}_1\in\mathbb{R}^{r_1\times r_2r_3}$ is a diagonal matrix.
It follows from Proposition \ref{kron_pro} that
$$
(\widetilde{\mathbf{A}}_3\otimes
\widetilde{\mathbf{A}}_2)^H(\widetilde{\mathbf{A}}_3\otimes
\widetilde{\mathbf{A}}_2)=(\widetilde{\mathbf{A}}_3^H\otimes
\widetilde{\mathbf{A}}_2^H)(\widetilde{\mathbf{A}}_3\otimes
\widetilde{\mathbf{A}}_2)=(\widetilde{\mathbf{A}}_3^H\widetilde{\mathbf{A}}_3)\otimes
(\widetilde{\mathbf{A}}_2^H\widetilde{\mathbf{A}}_2)=\mathbf{I}_{r_3}\otimes\mathbf{I}_{r_2}=\mathbf{I}_{r_2r_3},
$$
which implies that $\widetilde{\mathbf{A}}_3\otimes
\widetilde{\mathbf{A}}_2$ has orthogonal columns.
Further, note that
\begin{equation}\label{SVD}
\begin{aligned}
\widetilde{\mathcal{W}}_{(1)}=\widetilde{\mathbf{A}}_1\widetilde{\mathcal{G}}_{(1)}
(\widetilde{\mathbf{A}}_3\otimes
\widetilde{\mathbf{A}}_2)^H=&\widetilde{\mathbf{A}}_1\mathbf{\bar{U}}_1\mathbf{\bar{S}}_1
\mathbf{\bar{V}}_1^H
(\widetilde{\mathbf{A}}_3\otimes
\widetilde{\mathbf{A}}_2)^H\\
=&(\widetilde{\mathbf{A}}_1\mathbf{\bar{U}}_1)\mathbf{\bar{S}}_1
((\widetilde{\mathbf{A}}_3\otimes
\widetilde{\mathbf{A}}_2)\mathbf{\bar{V}}_1)^H.
\end{aligned}
\end{equation}
Here,
$\widetilde{\mathbf{A}}_1\mathbf{\bar{U}}_1$ and $(\widetilde{\mathbf{A}}_3\otimes
\widetilde{\mathbf{A}}_2)\mathbf{\bar{V}}_1$ have orthogonal columns,
and   (\ref{SVD}) is the SVD of $\widetilde{\mathcal{W}}_{(1)}$.
This consequently implies that $\|\widetilde{\mathcal{W}}_{(1)}\|=\|\widetilde{\mathcal{G}}_{(1)}\|$.
Similarly, we can obtain that $\|\widetilde{\mathcal{W}}_{(2)}\|=\|\widetilde{\mathcal{G}}_{(2)}\|$
and $\|\widetilde{\mathcal{W}}_{(3)}\|=\|\widetilde{\mathcal{G}}_{(3)}\|$.
Then, for any $i=1,2,3$, we have
$$
\|\widetilde{\mathcal{W}}_{(i)}-\mathcal{W}_{(i)}\|
\leq \|\widetilde{\mathcal{W}}_{(i)}\|+\|\mathcal{W}_{(i)}\|
=\|\widetilde{\mathcal{G}}_{(i)}\|+\sigma_1({\mathcal {W}}_{(i)})
\leq c\sqrt{r_1r_2r_3}+\sigma_1({\mathcal {W}}_{(i)}),
$$
where the last inequality follows from
$\|\widetilde{\mathcal{G}}_{(i)}\|\leq\|\widetilde{\mathcal{G}}_{(i)}\|_F
=\|\widetilde{\mathcal{G}}\|_F\leq c\sqrt{r_1r_2r_3}$.
Combining the above inequality with (\ref{HOSVD5}) immediately yields
$$
\begin{aligned}
\left\|(\widetilde{\mathbf{A}}_{i})_{:,j}-(\mathbf{A}_{i})_{:,j}\right\|_{2}
&\leq \frac{2^{\frac{3}{2}}\left(3\sigma_1({\mathcal {W}}_{(i)})+
c\sqrt{r_1r_2r_3}\right)
\|\widetilde{\mathcal{W}}_{(i)}-\mathcal{W}_{(i)}\|_{F}}
{\rho \sigma_j^2
(\mathcal{W}_{(i)})}.
\end{aligned}
$$
Consequently,
$$
\left\|\widetilde{\mathbf{A}}_i-\mathbf{A}_i\right\|_{F}^2=
\sum_{j=1}^{r_i}\left\|(\widetilde{\mathbf{A}}_{i})_{:,j}-(\mathbf{A}_{i})_{:,j}\right\|_{2}^2
\leq
\sum_{j=1}^{r_i}\frac{2^{3}\left(3\sigma_1\left({\mathcal {W}}_{(i)}\right)+
c\sqrt{r_1r_2r_3}\right)^2
\|\widetilde{\mathcal{W}}_{(i)}-\mathcal{W}_{(i)}\|_{F}^2}
{\rho^2 \sigma_j^4
\left(\mathcal{W}_{(i)}\right)}.
$$
Then, we can further demonstrate that
\begin{equation}\label{HOSVD1}
\left\|\widetilde{\mathbf{A}}_i-\mathbf{A}_i\right\|_{F}
\leq
2^{\frac{3}{2}}\left(\eta_i / \rho\right)\|\widetilde{\mathcal{W}}-\mathcal{W}\|_{F},
\end{equation}
where
$\eta_i=\left(3 \sigma_1\left({\mathcal {W}}_{(i)}\right)+
c\sqrt{r_1r_2r_3}\right)
\sqrt{\sum_{j=1}^{r_i}\frac{1}
{\sigma_j^4
\left(\mathcal{W}_{(i)}\right)}}$.

On the other hand, we make the
observation that
\begin{equation}\label{HOSVD2}
\begin{aligned}
\|\widetilde{\mathcal{G}}-\mathcal{G}\|_F
=&\left\|\widetilde{\mathcal{W}}\times_1 \widetilde{\mathbf{A}}_1^H \times_2 \widetilde{\mathbf{A}}_2^H \times_3
\widetilde{\mathbf{A}}_3^H-\mathcal{W}\times_1 \mathbf{A}_1^H \times_2 \mathbf{A}_2^H \times_3
\mathbf{A}_3^H\right\|_F\\
=&\left\|\widetilde{\mathcal{W}}\times_1 \widetilde{\mathbf{A}}_1^H \times_2 \widetilde{\mathbf{A}}_2^H \times_3
\widetilde{\mathbf{A}}_3^H-
\mathcal{W}\times_1 \widetilde{\mathbf{A}}_1^H \times_2 \widetilde{\mathbf{A}}_2^H \times_3
\widetilde{\mathbf{A}}_3^H
\right.\\
&+\left.
\mathcal{W}\times_1 \widetilde{\mathbf{A}}_1^H \times_2 \widetilde{\mathbf{A}}_2^H \times_3
\widetilde{\mathbf{A}}_3^H
-\mathcal{W}\times_1 \mathbf{A}_1^H \times_2 \widetilde{\mathbf{A}}_2^H \times_3
\widetilde{\mathbf{A}}_3^H
\right.\\
&+\left.
\mathcal{W}\times_1 \mathbf{A}_1^H \times_2 \widetilde{\mathbf{A}}_2^H \times_3
\widetilde{\mathbf{A}}_3^H-
\mathcal{W}\times_1 \mathbf{A}_1^H \times_2 \mathbf{A}_2^H \times_3
\widetilde{\mathbf{A}}_3^H
\right.\\
&+\left.
\mathcal{W}\times_1 \mathbf{A}_1^H \times_2 \mathbf{A}_2^H \times_3
\widetilde{\mathbf{A}}_3^H
-\mathcal{W}\times_1 \mathbf{A}_1^H \times_2 \mathbf{A}_2^H \times_3
\mathbf{A}_3^H\right\|_F\\
\leq&\left\|\widetilde{\mathcal{W}}\times_1 \widetilde{\mathbf{A}}_1^H \times_2 \widetilde{\mathbf{A}}_2^H \times_3
\widetilde{\mathbf{A}}_3^H-
\mathcal{W}\times_1 \widetilde{\mathbf{A}}_1^H \times_2 \widetilde{\mathbf{A}}_2^H \times_3
\widetilde{\mathbf{A}}_3^H
\right\|_F\\
&+
\left\|
\mathcal{W}\times_1 \widetilde{\mathbf{A}}_1^H \times_2 \widetilde{\mathbf{A}}_2^H \times_3
\widetilde{\mathbf{A}}_3^H
-\mathcal{W}\times_1 \mathbf{A}_1^H \times_2 \widetilde{\mathbf{A}}_2^H \times_3
\widetilde{\mathbf{A}}_3^H
\right\|_F\\
&+
\left\|
\mathcal{W}\times_1 \mathbf{A}_1^H \times_2 \widetilde{\mathbf{A}}_2^H \times_3
\widetilde{\mathbf{A}}_3^H-
\mathcal{W}\times_1 \mathbf{A}_1^H \times_2 \mathbf{A}_2^H \times_3
\widetilde{\mathbf{A}}_3^H
\right\|_F\\
&+
\left\|
\mathcal{W}\times_1 \mathbf{A}_1^H \times_2 \mathbf{A}_2^H \times_3
\widetilde{\mathbf{A}}_3^H
-\mathcal{W}\times_1 \mathbf{A}_1^H \times_2 \mathbf{A}_2^H \times_3
\mathbf{A}_3^H\right\|_F\\
\leq \  &\|\widetilde{\mathcal{W}}-\mathcal{W}\|_F
+
\left\|
\mathcal{W}\times_1 \widetilde{\mathbf{A}}_1^H \times_2 \widetilde{\mathbf{A}}_2^H \times_3
\widetilde{\mathbf{A}}_3^H
-\mathcal{W}\times_1 \mathbf{A}_1^H \times_2\widetilde{\mathbf{A}}_2^H \times_3
\widetilde{\mathbf{A}}_3^H
\right\|_F\\
&+
\left\|
\mathcal{W}\times_1 \mathbf{A}_1^H \times_2 \widetilde{\mathbf{A}}_2^H \times_3
\widetilde{\mathbf{A}}_3^H-
\mathcal{W}\times_1 \mathbf{A}_1^H \times_2 \mathbf{A}_2^H \times_3
\widetilde{\mathbf{A}}_3^H
\right\|_F\\
&+
\left\|
\mathcal{W}\times_1 \mathbf{A}_1^H \times_2 \mathbf{A}_2^H \times_3
\widetilde{\mathbf{A}}_3^H
-\mathcal{W}\times_1 \mathbf{A}_1^H \times_2 \mathbf{A}_2^H \times_3
\mathbf{A}_3^H\right\|_F.
\end{aligned}
\end{equation}
Similar to (\ref{delta_A_4}), we have
$$
\begin{aligned}
\ &\left\|
\mathcal{W}\times_1 \widetilde{\mathbf{A}}_1^H \times_2 \widetilde{\mathbf{A}}_2^H \times_3
\widetilde{\mathbf{A}}_3^H
-\mathcal{W}\times_1 \mathbf{A}_1^H \times_2 \widetilde{\mathbf{A}}_2^H \times_3
\widetilde{\mathbf{A}}_3^H
\right\|_F\\
= \ &\left\|
\widetilde{\mathbf{A}}_1^H
\mathcal{W}_{(1)}(\widetilde{\mathbf{A}}_3\otimes\widetilde{\mathbf{A}}_2)
-\mathbf{A}_1^H
\mathcal{W}_{(1)}(\widetilde{\mathbf{A}}_3\otimes\widetilde{\mathbf{A}}_2)
\right\|_F\\
= \ &\left\|
(\widetilde{\mathbf{A}}_1^H-\mathbf{A}_1^H)
\mathcal{W}_{(1)}(\widetilde{\mathbf{A}}_3\otimes\widetilde{\mathbf{A}}_2)
\right\|_F\\
\leq \  &\|
\widetilde{\mathbf{A}}_1^H-\mathbf{A}_1^H\|_F
\|
\mathcal{W}_{(1)}(\widetilde{\mathbf{A}}_3\otimes\widetilde{\mathbf{A}}_2)
\|\\
\leq \ &\|
\widetilde{\mathbf{A}}_1^H-\mathbf{A}_1^H\|_F
\|
\mathcal{W}_{(1)}\|
\|
\widetilde{\mathbf{A}}_3\otimes\widetilde{\mathbf{A}}_2
\|\\
\leq \ &\|
\widetilde{\mathbf{A}}_1^H-\mathbf{A}_1^H\|_F
\|\mathcal{W}_{(1)}\|\\
= \  &\sigma_1({\mathcal {W}}_{(1)})\|
\widetilde{\mathbf{A}}_1^H-\mathbf{A}_1^H\|_F.
\end{aligned}
$$
The above inequality combined with (\ref{HOSVD1}) yields
\begin{equation}\label{HOSVD3}
\begin{aligned}
\left\|
\mathcal{W}\times_1 \widetilde{\mathbf{A}}_1^H \times_2 \widetilde{\mathbf{A}}_2^H \times_3
\widetilde{\mathbf{A}}_3^H
-\mathcal{W}\times_1 \mathbf{A}_1^H \times_2 \widetilde{\mathbf{A}}_2^H \times_3
\widetilde{\mathbf{A}}_3^H
\right\|_F
\leq
\frac{2^{\frac{3}{2}}\eta_1\sigma_1({\mathcal {W}}_{(1)})}{\rho}\|\widetilde{\mathcal{W}}-\mathcal{W}\|_{F}.
\end{aligned}
\end{equation}
Similarly,
\begin{equation}\label{HOSVD4}
\begin{aligned}
\left\|
\mathcal{W}\times_1 \mathbf{A}_1^H \times_2 \widetilde{\mathbf{A}}_2^H \times_3
\widetilde{\mathbf{A}}_3^H
-\mathcal{W}\times_1 \mathbf{A}_1^H \times_2 \mathbf{A}_2^H \times_3
\widetilde{\mathbf{A}}_3^H
\right\|_F
&\leq
\frac{2^{\frac{3}{2}}\eta_2\sigma_1({\mathcal{W}}_{(1)})}{\rho}\|\widetilde{\mathcal{W}}-\mathcal{W}\|_{F},\\
\left\|
\mathcal{W}\times_1 \mathbf{A}_1^H \times_2 \mathbf{A}_2^H \times_3
\widetilde{\mathbf{A}}_3^H
-\mathcal{W}\times_1 \mathbf{A}_1^H \times_2 \mathbf{A}_2^H \times_3
\mathbf{A}_3^H
\right\|_F
&\leq
\frac{2^{\frac{3}{2}}\eta_3\sigma_1({\mathcal{W}}_{(1)})}{\rho}\|\widetilde{\mathcal{W}}-\mathcal{W}\|_{F}.
\end{aligned}
\end{equation}
Substituting (\ref{HOSVD3}), (\ref{HOSVD4}) into (\ref{HOSVD2}) then gives
$$
\begin{aligned}
\|\widetilde{\mathcal{G}}-\mathcal{G}\|_F
\leq \ &\|\widetilde{\mathcal{W}}-\mathcal{W}\|_{F}+
\sum_{i=1}^{3}\frac{2^{\frac{3}{2}}\eta_i\sigma_1({\mathcal{W}}_{(1)})}{\rho}\|\widetilde{\mathcal{W}}-\mathcal{W}\|_{F}\\
= \ &b\|\widetilde{\mathcal{W}}-\mathcal{W}\|_{F},
\end{aligned}
$$
where $b:=1+
\frac{2^{\frac{3}{2}}\sigma_1({\mathcal{W}}_{(1)})\sum_{i=1}^{3}\eta_i}{\rho}>0$.
\qed

\begin{lemma}\label{G_sparsity}
Suppose that Assumption \ref{assum_G} holds. Then
	$\|\mathbf{A}\mathbf{g}\|_0\leq s_1m^2p$, where $\mathbf{A}$ and $\mathbf{g}$ are defined in (\ref{VAR_vec}), and  $s_1\in[0,1]$ is the sparsity threshold related to $\bar{s}$.
%
\end{lemma}
\textbf{Proof.}
Recall that
$\mathbf{g}=\operatorname{vec}
(\mathcal{G}_{(1)}^{H})\in\mathbb{R}^{r_1r_2r_3}$.
The ratio of zero entries in a vector or matrix is denoted as $\bar{z}(\cdot)$ \cite{Wang2024}.
It follows from the fact that $\mathcal{G}\in\mathbb{R}^{r_1\times r_2 \times r_3}$ belongs to the set $\mathcal{U}$ that
the sparsity of the vector $\mathbf{g}$ is equal to $\bar{s}$, i.e., $\bar{z}(\mathbf{g})=\bar{s}$.
The matrix $\mathbf{G}\in\mathbb{R}^{r_1r_2r_3\times m^2p}$ is constructed such that
$\mathbf{G}_{:,j}=\mathbf{g}$ for each $j=1,2,\ldots,m^2p$,
where $\mathbf{G}_{:,j}$ denotes the $j$-th column of $\mathbf{G}$.
The sparsity of each row of $\mathbf{G}$ is denoted as $\bar{z}(\mathbf{G}_{i,:})$,
where $\mathbf{G}_{i,:}$ represents the $i$-th row of $\mathbf{G}$ and $\bar{z}(\mathbf{G}_{i,:})\in\{0,1\}$ for $i=1,2,\ldots,r_1r_2r_3$.
Similarly,
the sparsity of the $i$-th column of $\mathbf{A}\in\mathbb{R}^{m^2p\times r_1r_2r_3}$ is defined as $\bar{z}(\mathbf{A}_{:,i})$.
According to \cite[Proposition 1]{Wang2024},
we can deduce that
\begin{equation}\label{Defz1}
\bar{z}(\mathbf{A}\mathbf{G})\geq\max \left\{0,1+\sum_{i=1}^{r_1r_2r_3}\left(\bar{z}(\mathbf{G}_{i,:})+\bar{z}\left(\mathbf{A}_{:,i }\right)-\bar{z}(\mathbf{G}_{i,:})\bar{z}\left(\mathbf{A}_{:,i}\right)\right)-r_1r_2r_3\right\}:=z_1.
\end{equation}
It can be easily verified that
$
\bar{z}(\mathbf{A}\mathbf{g})=\bar{z}(\mathbf{A}\mathbf{G})\geq z_1.
$
Consequently, we have $\|\mathbf{A}\mathbf{g}\|_0\leq (1-z_1)m^2p:=s_1m^2p$.
Note that
$\sum_{i=1}^{r_1r_2r_3}\bar{z}(\mathbf{G}_{i,:})/(r_1r_2r_3)=\bar{s}$, which implies that $s_1$ depends on $\bar{s}$.
\qed

By using similar techniques,
it follows from the fact that $\widehat{\mathcal{G}}\in\mathbb{R}^{r_1\times r_2 \times r_3}$ belongs to the set $\mathcal{U}$ that
$\|\widehat{\mathbf{A}}\widehat{\mathbf{g}}\|_0\leq s_2m^2p$,
where $s_2\in[0,1]$ related to $\bar{s}$.

\section*{Appendix B. Proof of Theorem \ref{statistical}}

The VAR model in (\ref{test0}) can be expressed equivalently as
\begin{equation}\label{VAR_matrix}
\begin{gathered}
\underbrace{\begin{bmatrix}
\mathbf{y}_1^{H} \\
\mathbf{y}_2^{H} \\
\vdots \\
\mathbf{y}_T^{H}
\end{bmatrix}}_{\mathbf{Y}}
=\underbrace{
\begin{bmatrix}
\mathbf{y}_0^{H} & \mathbf{y}_{-1}^{H} & \ldots & \mathbf{y}_{-p+1}^{H} \\
\mathbf{y}_1^{H} & \mathbf{y}_0^{H} & \ldots & \mathbf{y}_{-p+2}^{H} \\
\vdots & \vdots & \ddots & \vdots \\
\mathbf{y}_{T-1}^{H} & \mathbf{y}_{T-2}^{H} & \ldots & \mathbf{y}_{T-p}^{H}
\end{bmatrix}}_{\mathbf{X}}\underbrace{
\begin{bmatrix}
\mathbf{W}_1^{H} \\
\mathbf{W}_2^{H} \\
\vdots \\
\mathbf{W}_p^{H}
\end{bmatrix}}_{\mathcal{W}_{(1)}^{H}}+\underbrace{
\begin{bmatrix}
\boldsymbol{\varepsilon}_{1}^{H} \\
\boldsymbol{\varepsilon}_{2}^{H} \\
\vdots \\
\boldsymbol{\varepsilon}_T^{H}
\end{bmatrix}}_{\mathbf{E}}.
\end{gathered}
\end{equation}
Notice that  $\mathcal{W}=\mathcal{G} \times_1 \mathbf{A}_1 \times_2 \mathbf{A}_2 \times_3 \mathbf{A}_3$
represents the HOSVD of $\mathcal{W}$ and $\mathcal{W}_{(1)}=\mathbf{A}_1\mathcal{G}_{(1)}(\mathbf{A}_3\otimes\mathbf{A}_2)^H$ \cite{kolda2009tensor}.
Based on the formulation in (\ref{VAR_matrix}), we can rewrite (\ref{test0}) as
$$
\mathbf{Y}=\mathbf{X}(\mathbf{A}_3\otimes\mathbf{A}_2)
\mathcal{G}_{(1)}^{H}\mathbf{A}_1^{H}+\mathbf{E}.
$$
The above equality is equivalent to
\begin{equation}\label{VAR_vec}
\begin{aligned}
\mathbf{y}:=\operatorname{vec}(\mathbf{Y})
&=\operatorname{vec}\left(\mathbf{X}(\mathbf{A}_3\otimes\mathbf{A}_2)
\mathcal{G}_{(1)}^{H}\mathbf{A}_1^{H}\right)+\operatorname{vec}(\mathbf{E}) \\
& =\operatorname{vec}\left(\mathbf{X}(\mathbf{A}_3\otimes\mathbf{A}_2)
\mathcal{G}_{(1)}^{H}\mathbf{A}_1^{H}\mathbf{I}_m\right)+\operatorname{vec}(\mathbf{E}) \\
& =\left(\mathbf{I}_m^H \otimes \mathbf{X}\right) \operatorname{vec}\left((\mathbf{A}_3\otimes\mathbf{A}_2)
\mathcal{G}_{(1)}^{H}\mathbf{A}_1^{H}\right)+\operatorname{vec}(\mathbf{E}) \\
& =\left(\mathbf{I}_m \otimes \mathbf{X}\right)
\left(\mathbf{A}_1\otimes(\mathbf{A}_3\otimes\mathbf{A}_2)\right) \operatorname{vec}(\mathcal{G}_{(1)}^{H})+\operatorname{vec}(\mathbf{E})\\
& =\underbrace{(\mathbf{I}_m\otimes\mathbf{X})}_{\mathbf{Z}}
\underbrace{(\mathbf{A}_1\otimes\mathbf{A}_3\otimes
\mathbf{A}_2)}_{\mathbf{A}}\underbrace{\operatorname{vec}
(\mathcal{G}_{(1)}^{H})}_{\mathbf{g}}+\underbrace{\operatorname{vec}
(\mathbf{E})}_{\mathbf{e}},
\end{aligned}
\end{equation}
where the third and fourth equalities follow from \cite[Equation 12.3.9]{golub2013matrix}.
Note that $\widehat{\mathcal{W}}
=\widehat{\mathcal{G}}\times_1 \widehat{\mathbf{A}}_1 \times_2 \widehat{\mathbf{A}}_2 \times_3
\widehat{\mathbf{A}}_3$. Denote
\begin{equation}\label{hat_AG}
\widehat{\mathbf{A}}=\widehat{\mathbf{A}}_1\otimes\widehat{\mathbf{A}}_3\otimes
\widehat{\mathbf{A}}_2, \ \widehat{\mathbf{g}}=\operatorname{vec}
(\widehat{\mathcal{G}}_{(1)}^{H}).
\end{equation}
Then, it can be easily verified that
\begin{equation}\label{lossf}
\begin{aligned}
\frac{1}{2T} \sum_{t=1}^{T}\left\|\mathbf{y}_t-\left(\mathcal{G} \times_1 \mathbf{A}_1 \times_2
\mathbf{A}_2 \times_3 \mathbf{A}_3\right)_{(1)}\mathbf{x}_t\right\|_2^2
&=\frac{1}{2T}\|\mathbf{y}-\mathbf{Z}\mathbf{A}\mathbf{g}\|_{2}^{2},\\
\frac{1}{2T} \sum_{t=1}^{T}\left\|\mathbf{y}_t-(\widehat{\mathcal{G}}\times_1
\widehat{\mathbf{A}}_1 \times_2
\widehat{\mathbf{A}}_2 \times_3 \widehat{\mathbf{A}}_3)_{(1)}\mathbf{x}_t\right\|_2^2
&=\frac{1}{2T}\|\mathbf{y}-\mathbf{Z}\widehat{\mathbf{A}}\widehat{\mathbf{g}}\|_{2}^{2}.
\end{aligned}
\end{equation}
Due to the optimality of $\widehat{\mathcal{W}}$ in (\ref{test3}), we have
$$
\begin{aligned}
 \ &\frac{1}{2T} \sum_{t=1}^{T}\left\|\mathbf{y}_t-(\widehat{\mathcal{G}}\times_1
\widehat{\mathbf{A}}_1 \times_2
\widehat{\mathbf{A}}_2 \times_3 \widehat{\mathbf{A}}_3)_{(1)}\mathbf{x}_t\right\|_2^2
+
\beta\|\widehat{\mathbf{g}}\|_{1}+
\sum_{i=1}^{3}\alpha_{i}\operatorname{tr}(\widehat{\mathbf{A}}_{i}^{H}
\mathbf{L}_{i}\widehat{\mathbf{A}}_{i})\\
\leq \ & \frac{1}{2T} \sum_{t=1}^{T}\left\|\mathbf{y}_t-\left(\mathcal{G} \times_1 \mathbf{A}_1 \times_2
\mathbf{A}_2 \times_3 \mathbf{A}_3\right)_{(1)}\mathbf{x}_t\right\|_2^2+\beta\|\mathbf{g}\|_{1}
+\sum_{i=1}^{3}\alpha_{i}\operatorname{tr}(\mathbf{A}_{i}^{H}\mathbf{L}_{i}\mathbf{A}_{i}).
\end{aligned}
$$
Using (\ref{lossf}) further leads to
\begin{equation}\label{opti_model}
\begin{aligned}
 \ &\frac{1}{2T}\|\mathbf{y}-\mathbf{Z}\widehat{\mathbf{A}}\widehat{\mathbf{g}}\|_{2}^{2}+
\beta\|\widehat{\mathbf{g}}\|_{1}+
\sum_{i=1}^{3}\alpha_{i}\left(\operatorname{tr}(\widehat{\mathbf{A}}_{i}^{H}
\mathbf{L}_{i}\widehat{\mathbf{A}}_{i})-
\operatorname{tr}(\mathbf{A}_{i}^{H}\mathbf{L}_{i}\mathbf{A}_{i})\right) \\
\leq \ &\frac{1}{2T}\|\mathbf{y}-\mathbf{Z}\mathbf{A}\mathbf{g}\|_{2}^{2}+\beta\|\mathbf{g}\|_{1}.
\end{aligned}
\end{equation}
Denote
$
\Delta_{\mathbf{A}_{i}}=\widehat{\mathbf{A}}_{i}-\mathbf{A}_{i}, i=1,2,3.
$
Then we get
\begin{equation}\label{trace_sum}
\begin{aligned}
 \ &\operatorname{tr}(\widehat{\mathbf{A}}_{1}^{H}
\mathbf{L}_{1}\widehat{\mathbf{A}}_{1})-
\operatorname{tr}(\mathbf{A}_{1}^{H}\mathbf{L}_{1}\mathbf{A}_{1}) \\
= \ &\operatorname{tr}((\mathbf{A}_{1}+\Delta_{\mathbf{A}_{1}})^{H}\mathbf{L}_{1}
(\mathbf{A}_{1}+\Delta_{\mathbf{A}_{1}}))
-\operatorname{tr}(\mathbf{A}_{1}^{H}\mathbf{L}_{1}\mathbf{A}_{1})\\
= \ &\operatorname{tr}(\mathbf{A}_{1}^{H}\mathbf{L}_{1}\mathbf{A}_{1}
+\mathbf{A}_{1}^{H}\mathbf{L}_{1}\Delta_{\mathbf{A}_{1}}
+\Delta_{\mathbf{A}_{1}}^{H}\mathbf{L}_{1}\mathbf{A}_{1}+
\Delta_{\mathbf{A}_{1}}^{H}\mathbf{L}_{1}\Delta_{\mathbf{A}_{1}})
-\operatorname{tr}(\mathbf{A}_{1}^{H}\mathbf{L}_{1}\mathbf{A}_{1}) \\
= \ &2\operatorname{tr}(\mathbf{A}_{1}^{H}\mathbf{L}_{1}\Delta_{\mathbf{A}_{1}})
+\operatorname{tr}(\Delta_{\mathbf{A}_{1}}^{H}\mathbf{L}_{1}\Delta_{\mathbf{A}_{1}}),
\end{aligned}
\end{equation}
where the last equality follows from the fact that $\mathbf{L}_{1}$ is symmetric.
Note that
\begin{equation}\label{trace_inner}
\begin{aligned}
\operatorname{tr}(\mathbf{A}_{1}^{H}\mathbf{L}_{1}\Delta_{\mathbf{A}_{1}})
&=\langle\mathbf{L}_{1}\mathbf{A}_{1},\Delta_{\mathbf{A}_{1}}\rangle
\geq -\|\mathbf{L}_{1}\mathbf{A}_{1}\|_F\|\Delta_{\mathbf{A}_{1}}\|_F,
\end{aligned}
\end{equation}
where the inequality holds by Cauchy-Schwarz inequality.
Combining the above inequality with Lemma \ref{graph} and (\ref{trace_sum}) yields
\begin{equation}\label{trace_A1}
\begin{aligned}
&\operatorname{tr}(\widehat{\mathbf{A}}_{1}^{H}
\mathbf{L}_{1}\widehat{\mathbf{A}}_{1})-
\operatorname{tr}(\mathbf{A}_{1}^{H}\mathbf{L}_{1}\mathbf{A}_{1})
\geq \frac{\lambda_{2}(\mathbf{L}_{1})}{2}\|\widehat{\mathbf{A}}_{1}-\mathbf{A}_{1}\|_F^2
-2\|\mathbf{L}_{1}\mathbf{A}_{1}\|_F\|\Delta_{\mathbf{A}_{1}}\|_F.
\end{aligned}
\end{equation}
By using similar techniques,
we can easily obtain
\begin{equation}\label{trace_A2}
\begin{aligned}
&\operatorname{tr}(\widehat{\mathbf{A}}_{2}^{H}
\mathbf{L}_{2}\widehat{\mathbf{A}}_{2})-
\operatorname{tr}(\mathbf{A}_{2}^{H}\mathbf{L}_{2}\mathbf{A}_{2})
\geq \frac{\lambda_{2}(\mathbf{L}_{2})}{2}\|\widehat{\mathbf{A}}_{2}-\mathbf{A}_{2}\|_F^2
-2\|\mathbf{L}_{2}\mathbf{A}_{2}\|_F\|\Delta_{\mathbf{A}_{2}}\|_F,\\
&\operatorname{tr}(\widehat{\mathbf{A}}_{3}^{H}
\mathbf{L}_{3}\widehat{\mathbf{A}}_{3})-
\operatorname{tr}(\mathbf{A}_{3}^{H}\mathbf{L}_{3}\mathbf{A}_{3})
\geq \frac{\lambda_{2}(\mathbf{L}_{3})}{2}\|\widehat{\mathbf{A}}_{3}-\mathbf{A}_{3}\|_F^2
-2\|\mathbf{L}_{3}\mathbf{A}_{3}\|_F\|\Delta_{\mathbf{A}_{3}}\|_F.
\end{aligned}
\end{equation}
Substituting (\ref{trace_A1}), (\ref{trace_A2}) into (\ref{opti_model}) yields
\begin{equation}\label{model_ineq}
\begin{aligned}
 \ &\frac{1}{2T}\|\mathbf{y}-\mathbf{Z}\widehat{\mathbf{A}}\widehat{\mathbf{g}}\|_{2}^{2}+
\beta\|\widehat{\mathbf{g}}\|_{1}+
\sum_{i=1}^{3}\alpha_{i}\left(\frac{\lambda_{2}(\mathbf{L}_{i})}{2}\|\widehat{\mathbf{A}}_{i}-
\mathbf{A}_{i}\|_F^2
-2\|\mathbf{L}_{i}\mathbf{A}_{i}\|_F\|\Delta_{\mathbf{A}_{i}}\|_F\right) \\
\leq \ &\frac{1}{2T}\|\mathbf{y}-\mathbf{Z}\mathbf{A}\mathbf{g}\|_{2}^{2}+\beta\|\mathbf{g}\|_{1}.
\end{aligned}
\end{equation}

Denote
\begin{equation}\label{DelDf}
\widehat{\Delta}:=\widehat{\mathbf{A}}\widehat{\mathbf{g}}-\mathbf{A}\mathbf{g},
\end{equation}
where $\mathbf{A}, \mathbf{g}$ and $\widehat{\mathbf{A}}, \widehat{\mathbf{g}}$ are defined in  (\ref{VAR_vec}) and
(\ref{hat_AG}), respectively.
We can easily get that
$$
\begin{aligned}
\frac{1}{2T}\|\mathbf{y}-\mathbf{Z}\widehat{\mathbf{A}}\widehat{\mathbf{g}}\|_{2}^{2}=
\frac{1}{2T}\|\mathbf{y}-\mathbf{Z}(\widehat{\Delta}+\mathbf{A}\mathbf{g})\|_{2}^{2}
= \ &\frac{1}{2T}\|-\mathbf{Z}\widehat{\Delta}+\mathbf{e}\|_{2}^{2}\\
= \ &\frac{1}{2T}\|\mathbf{Z}\widehat{\Delta}\|_{2}^{2}+\frac{1}{2T}\|\mathbf{e}\|_{2}^{2}
-\frac{1}{T}\langle\mathbf{Z}\widehat{\Delta},\mathbf{e}\rangle,
\end{aligned}
$$
where the second equality follows from (\ref{VAR_vec}).
Combining the above result with (\ref{VAR_vec}) and (\ref{model_ineq}) immediately yields
$$
\begin{aligned}
& \frac{1}{2 T}\|\mathbf{Z}\widehat{\Delta}\|_2^2+\frac{1}{2 T}\|\mathbf{e}\|_2^2-\frac{1}{T}\langle\mathbf{Z}\widehat{\Delta},\mathbf{e}\rangle
+\beta\|\widehat{\mathbf{g}}\|_1 \\
& +
\sum_{i=1}^{3}\alpha_{i}\left(\frac{\lambda_{2}(\mathbf{L}_{i})}{2}\|\widehat{\mathbf{A}}_{i}-
\mathbf{A}_{i}\|_F^2
-2\|\mathbf{L}_{i}\mathbf{A}_{i}\|_F\|\Delta_{\mathbf{A}_{i}}\|_F\right)
\leq \frac{1}{2 T}\|\mathbf{e}\|_2^2+\beta\|\mathbf{g}\|_1,
\end{aligned}
$$
which is equivalent to
$$
\begin{aligned}
 \ &\frac{1}{2 T}\|\mathbf{Z}\widehat{\Delta}\|_2^2
+\beta\|\widehat{\mathbf{g}}\|_1
+
\sum_{i=1}^{3}\frac{\alpha_{i}\lambda_{2}(\mathbf{L}_{i})}{2}\|\widehat{\mathbf{A}}_{i}-
\mathbf{A}_{i}\|_F^2\\
\leq \ & \beta\|\mathbf{g}\|_1+\frac{1}{T}\langle\mathbf{Z}\widehat{\Delta},\mathbf{e}\rangle
+\sum_{i=1}^{3}2\alpha_{i}\|\mathbf{L}_{i}\mathbf{A}_{i}\|_F\|\Delta_{\mathbf{A}_{i}}\|_F.
\end{aligned}
$$
Then, by Lemma \ref{inner_product} and the assumption  $\beta \geq 2
\pi\vartheta a_2\sqrt{\log(m^2p)/T}$,
with probability at least  $1-6\exp(-a\log (m^2 p))$,
the following inequality holds:
\begin{equation}\label{model_ineq2}
\begin{aligned}
 \ &\frac{1}{2 T}\|\mathbf{Z}\widehat{\Delta}\|_2^2
+\beta\|\widehat{\mathbf{g}}\|_1
+\sum_{i=1}^{3}\frac{\alpha_{i}\lambda_{2}(\mathbf{L}_{i})}{2}\|\widehat{\mathbf{A}}_{i}-
\mathbf{A}_{i}\|_F^2 \\
\leq \ &  \beta\|\mathbf{g}\|_1+\beta\|\widehat{\Delta}\|_1
+\sum_{i=1}^{3}2\alpha_{i}\|\mathbf{L}_{i}\mathbf{A}_{i}\|_F\|\Delta_{\mathbf{A}_{i}}\|_F.
\end{aligned}
\end{equation}

Let $\Delta_{\mathbf{g}}:=\widehat{\mathbf{g}}-\mathbf{g}$.
Denote the nonzero index set of $\operatorname{vec}(\mathcal G^{H})$ by $\mathbb{S}_g$,
and let $\mathbb{S}_{g^c}$ represent its complementary set.
Then, we have
$$
\begin{aligned}
\|\widehat{\mathbf{g}}\|_1 =\left\|\widehat{\mathbf{g}}_{\mathbb{S}_g}\right\|_1
+\left\|\widehat{\mathbf{g}}_{\mathbb{S}_{g^c}}\right\|_1
& =\left\|\left(\mathbf{g}+\Delta_{\mathbf{g}}\right) _{\mathbb{S}_g}\right\|_1+\left\|\widehat{\mathbf{g}}_{\mathbb{S}_{g^c}}\right\|_1 \\
& \geq\left\|\mathbf{g}_{\mathbb{S}_{g}}\right\|_1-
\left\|(\Delta_{\mathbf{g}})_{\mathbb{S}_{g}}\right\|_1
+\left\|\widehat{\mathbf{g}}_{\mathbb{S}_{g^c}}\right\|_1,
\end{aligned}
$$
which taken collectively with (\ref{model_ineq2}) gives
$$
\begin{aligned}
\ & \frac{1}{2 T}\|\mathbf{Z}\widehat{\Delta}\|_2^2
+\beta\left\|\mathbf{g}_{\mathbb{S}_{g}}\right\|_1-
\beta\left\|(\Delta_{\mathbf{g}})_{\mathbb{S}_{g}}\right\|_1
+\beta\left\|\widehat{\mathbf{g}}_{\mathbb{S}_{g^c}}\right\|_1
+
\sum_{i=1}^{3}\frac{\alpha_{i}\lambda_{2}(\mathbf{L}_{i})}{2}\|\widehat{\mathbf{A}}_{i}-
\mathbf{A}_{i}\|_F^2\\
\leq \ & \beta\|\mathbf{g}\|_1+\beta\|\widehat{\Delta}\|_1
+\sum_{i=1}^{3}2\alpha_{i}\|\mathbf{L}_{i}\mathbf{A}_{i}\|_F\|\Delta_{\mathbf{A}_{i}}\|_F.
\end{aligned}
$$
The above inequality can be rewritten as
$$
\begin{aligned}
\ & \frac{1}{2 T}\|\mathbf{Z}\widehat{\Delta}\|_2^2
+\beta\left\|\widehat{\mathbf{g}}_{\mathbb{S}_{g^c}}\right\|_1
+
\sum_{i=1}^{3}\frac{\alpha_{i}\lambda_{2}(\mathbf{L}_{i})}{2}\|\widehat{\mathbf{A}}_{i}-
\mathbf{A}_{i}\|_F^2\\
\leq  \ &\beta\|\widehat{\Delta}\|_1
+\sum_{i=1}^{3}2\alpha_{i}\|\mathbf{L}_{i}\mathbf{A}_{i}\|_F\|\Delta_{\mathbf{A}_{i}}\|_F
+\beta\left\|(\Delta_{\mathbf{g}})_{\mathbb{S}_{g}}\right\|_1.
\end{aligned}
$$
It follows from $\beta\left\|\widehat{\mathbf{g}}_{\mathbb{S}_{g^c}}\right\|_1\geq 0$
that
$$
\begin{aligned}
& \frac{1}{2 T}\|\mathbf{Z}\widehat{\Delta}\|_2^2
+
\sum_{i=1}^{3}\frac{\alpha_{i}\lambda_{2}(\mathbf{L}_{i})}{2}\|\widehat{\mathbf{A}}_{i}-
\mathbf{A}_{i}\|_F^2
\leq \beta\|\widehat{\Delta}\|_1
+\sum_{i=1}^{3}2\alpha_{i}\|\mathbf{L}_{i}\mathbf{A}_{i}\|_F\|\Delta_{\mathbf{A}_{i}}\|_F
+\beta\left\|(\Delta_{\mathbf{g}})_{\mathbb{S}_{g}}\right\|_1.
\end{aligned}
$$
Note that $T\geq\frac{
8\varsigma^2}
{a_4}
\bar{p}\min(\log (mp), \log(21emp/\bar{p}))$.
Lemma \ref{Restricted eigenvalue} implies that
the following inequality holds
\begin{equation}\label{model_ineq3}
\begin{aligned}
\frac{\varpi}{4}\|\widehat{\Delta}\|_2^2
+
\sum_{i=1}^{3}\frac{\alpha_{i}\lambda_{2}(\mathbf{L}_{i})}{2}\|\widehat{\mathbf{A}}_{i}-
\mathbf{A}_{i}\|_F^2
\leq \beta\|\widehat{\Delta}\|_1
+\sum_{i=1}^{3}2\alpha_{i}\|\mathbf{L}_{i}\mathbf{A}_{i}\|_F\|\Delta_{\mathbf{A}_{i}}\|_F
+\beta\left\|(\Delta_{\mathbf{g}})_{\mathbb{S}_{g}}\right\|_1
\end{aligned}
\end{equation}
with probability at least $1-2 \exp (-\bar{p}\min(\log mp, \log(21emp/\bar{p})))$, where $\bar{p}$ is related to $\bar{s}$.

Note that
$\alpha_{i}\geq \frac{6c^2}{\lambda_{2}(\mathbf{L}_{i})}\sqrt{\frac{s\log (m^2p)}{T}}=\frac{6c^2 r_1 r_2 r_3}{\lambda_{2}(\mathbf{L}_{i})}\sqrt{\frac{s\log (m^2p)}{(r_1 r_2 r_3)^2T}}.
$
Then, we have
$$
\begin{aligned}
 \ &\sum_{i=1}^{3}\alpha_{i}\frac{\lambda_{2}(\mathbf{L}_{i})}{2}\|\widehat{\mathbf{A}}_{i}-
\mathbf{A}_{i}\|_F^2
+{\sqrt{\frac{s\log (m^2p)}{(r_1 r_2 r_3)^2T}}}\|\operatorname{vec}
(\widehat{\mathcal{G}}_{(1)}^{H})-\operatorname{vec}
(\mathcal{G}_{(1)}^{H})\|_2^2\\
\geq \ & 3 c^2 r_1 r_2 r_3{\sqrt{\frac{s\log (m^2p)}{(r_1 r_2 r_3)^2T}}}\left(\|\widehat{\mathbf{A}}_3-\mathbf{A}_3\|_F^2+
\|\widehat{\mathbf{A}}_2-\mathbf{A}_2\|_F^2
+\|\widehat{\mathbf{A}}_1-\mathbf{A}_1\|_F^2\right)\\
&+{\sqrt{\frac{s\log (m^2p)}{(r_1 r_2 r_3)^2T}}}\|\operatorname{vec}
(\widehat{\mathcal{G}}_{(1)}^{H})-\operatorname{vec}
(\mathcal{G}_{(1)}^{H})\|_2^2\\
\geq  \ & {\frac{1}{2}\sqrt{\frac{s\log (m^2p)}{(r_1 r_2 r_3)^2T}}}\|\widehat{\Delta}\|_2^2,
\end{aligned}
$$
where the last inequality follows from Lemma \ref{delta_A}.
This combined with (\ref{model_ineq3}) leads to
\begin{equation}\label{model_ineq4}
\begin{aligned}
 \ & \frac{\varpi}{4}\|\widehat{\Delta}\|_2^2
+\frac{1}{2}\sqrt{\frac{s\log (m^2p)}{(r_1 r_2 r_3)^2T}}\|\widehat{\Delta}\|_2^2\\
\leq \ & \beta\|\widehat{\Delta}\|_1
+\sum_{i=1}^{3}2\alpha_i\|\mathbf{L}_{i}\mathbf{A}_{i}\|_F\|\Delta_{\mathbf{A}_{i}}\|_F
+\beta\left\|(\Delta_{\mathbf{g}})_{\mathbb{S}_{g}}\right\|_1
+{\sqrt{\frac{s\log (m^2p)}{(r_1 r_2 r_3)^2T}}}\|\operatorname{vec}
(\widehat{\mathcal{G}}_{(1)}^{H})-\operatorname{vec}
(\mathcal{G}_{(1)}^{H})\|_2^2
\end{aligned}
\end{equation}
with probability at least $1-2 \exp (-\bar{p}\min(\log mp, \log(21emp/\bar{p})))$.

Recall that $\bar{s}=1-s /(r_1 r_2r_3)\in[0,1]$.
Since both $\mathcal{G}$ and $\widehat{\mathcal{G}}$
belong to the set $\mathcal{U}$,
it follows from Lemma \ref{G_sparsity} that
$\|\mathbf{A}\mathbf{g}\|_0\leq s_1m^2p$ and
$\|\widehat{\mathbf{A}}\widehat{\mathbf{g}}\|_0\leq s_2m^2p,$
where $s_1, s_2\in[0,1]$ are the sparsity thresholds related to $\bar{s}$.
Then, we have
\begin{equation}\label{S1S2St}
\begin{aligned}
\|\widehat{\Delta}\|_{0}
=\|\widehat{\mathbf{A}}\widehat{\mathbf{g}}-\mathbf{A}\mathbf{g}\|_{0}
&\leq\|\widehat{\mathbf{A}}\widehat{\mathbf{g}}\|_{0}
+\|\mathbf{A}\mathbf{g}\|_{0}\\
&\leq s_2m^2p+s_1m^2p:=\tilde{s}.
\end{aligned}
\end{equation}
As a consequence, we get
\begin{equation}\label{model_ineq5}
\|\widehat{\Delta}\|_{1}\leq \sqrt{\tilde{s}}\|\widehat{\Delta}\|_{2}.
\end{equation}
In addition,
\begin{equation}\label{model_ineq6}
\begin{aligned}
\left\|(\Delta_{\mathbf{g}})_{\mathbb{S}_{g}}\right\|_1
\leq \sqrt{s}\|\Delta_{\mathbf{g}}\|_{2}
\leq b\sqrt{s}\|\widehat{\Delta}\|_{2},
\end{aligned}
\end{equation}
where
the last inequality follows from
Lemma \ref{HOSVD perturbation bound}.
By H\"older's inequality, we have
\begin{equation}\label{model_ineq7}
\begin{aligned}
\|\operatorname{vec}
(\widehat{\mathcal{G}}_{(1)}^{H})-\operatorname{vec}
(\mathcal{G}_{(1)}^{H})\|_2^2
& \leq\|\operatorname{vec}
(\widehat{\mathcal{G}}_{(1)}^{H})-\operatorname{vec}
(\mathcal{G}_{(1)}^{H})\|_{\infty}\|\operatorname{vec}
(\widehat{\mathcal{G}}_{(1)}^{H})-\operatorname{vec}
(\mathcal{G}_{(1)}^{H})\|_1 \\
& \leq 2 c\|\operatorname{vec}
(\widehat{\mathcal{G}}_{(1)}^{H})-\operatorname{vec}
(\mathcal{G}_{(1)}^{H})\|_1 \\
& =2 c\|\Delta_{\mathbf{g}}\|_1 \\
& \leq 2 c \sqrt{2s}\left\|\Delta_{\mathbf{g}}\right\|_2 \\
& \leq 2 bc \sqrt{2s}\|\widehat{\Delta}\|_2,
\end{aligned}
\end{equation}
where the third inequality follows from the fact that
both $\mathcal{G}$ and $\widehat{\mathcal{G}}$
belong to the set $\mathcal{U}$ and
the last inequality follows from
Lemma \ref{HOSVD perturbation bound}.
Substituting (\ref{model_ineq5}), (\ref{model_ineq6}) and
(\ref{model_ineq7}) into (\ref{model_ineq4}) then gives
\begin{equation}\label{model_ineq8}
\begin{aligned}
\ &\left(\frac{\varpi}{4}+\frac{1}{2}\sqrt{\frac{s\log (m^2p)}{(r_1 r_2 r_3)^2T}}\right)\|\widehat{\Delta}\|_2^2\\
\leq \  &\beta\sqrt{\tilde{s}}\|\widehat{\Delta}\|_{2}
+\sum_{i=1}^{3}2\alpha_i\|\mathbf{L}_{i}\mathbf{A}_{i}\|_F\|\Delta_{\mathbf{A}_{i}}\|_F
+b\beta\sqrt{s}\|\widehat{\Delta}\|_{2}
+2\sqrt{2}bc\sqrt{\frac{s^2\log (m^2p)}{(r_1 r_2 r_3)^2T}}\|\widehat{\Delta}\|_2\\
\leq\ &\beta\sqrt{\tilde{s}}\|\widehat{\Delta}\|_{2}
+\sum_{i=1}^{3}2\alpha_i\|\mathbf{L}_{i}\mathbf{A}_{i}\|_F\|\Delta_{\mathbf{A}_{i}}\|_F
+b\beta\sqrt{s}\|\widehat{\Delta}\|_{2}
+2\sqrt{2}bc\sqrt{\frac{\log (m^2p)}{T}}\|\widehat{\Delta}\|_2
\end{aligned}
\end{equation}
with probability at least $1-2 \exp (-\bar{p}\min(\log mp, \log(21emp/\bar{p})))$,
where the last inequality follows from $s\leq r_1 r_2 r_3$.
Note that
$$
2\alpha_i\|\mathbf{L}_{i}\mathbf{A}_{i}\|_F\|\Delta_{\mathbf{A}_{i}}\|_F
\leq 2\alpha_i\|\mathbf{L}_{i}\|_F\|\mathbf{A}_{i}\|\|\Delta_{\mathbf{A}_{i}}\|_F
\leq 2\alpha_il_i\|\Delta_{\mathbf{A}_{i}}\|_F,
$$
where $l_i:=\|\mathbf{L}_{i}\|_F$.
Using Lemma \ref{HOSVD perturbation bound} further leads to
$$
2\alpha_il_i\|\Delta_{\mathbf{A}_{i}}\|_F\leq 2\alpha_il_i
\frac{2^{\frac{3}{2}}\eta_i}{\rho}\|\widehat{\Delta}\|_2, i=1,2,3.
$$
Taken the above inequality together with (\ref{model_ineq8}) demonstrates that
$$
\begin{aligned}
 \ &\left(\frac{\varpi}{4}+\frac{1}{2}\sqrt{\frac{s\log (m^2p)}{(r_1 r_2 r_3)^2T}}\right)\|\widehat{\Delta}\|_2^2\\
\leq \ &\beta\sqrt{\tilde{s}}\|\widehat{\Delta}\|_{2}
+\sum_{i=1}^{3}\alpha_il_i \frac{2^{\frac{5}{2}}\eta_i}{\rho}\|\widehat{\Delta}\|_2
+b\beta\sqrt{s}\|\widehat{\Delta}\|_{2}
+2\sqrt{2}bc\sqrt{\frac{\log (m^2p)}{T}}\|\widehat{\Delta}\|_2
\end{aligned}
$$
with probability at least $1-2 \exp (-\bar{p}\min(\log mp, \log(21emp/\bar{p})))$.
Consequently, we obtain
\begin{equation}\label{result1}
\begin{aligned}
\|\widehat{\Delta}\|_2
&\leq \frac{4\beta\sqrt{\tilde{s}}
+\sum_{i=1}^{3}2^{\frac{9}{2}}\alpha_il_i \frac{\eta_i}{\rho}
+4b\beta\sqrt{s}
+8\sqrt{2}bc\sqrt{(\log(m^2p))/T}}
{\varpi+2\sqrt{(s\log(m^2p))/((r_1 r_2 r_3)^2T)}}\\
&=\frac{4\beta(\sqrt{\tilde{s}}+b\sqrt{s})
+\hat{a}\kappa
+\bar{c}\sum_{i=1}^{3}\alpha_il_i \eta_i}{\varpi+2\kappa\sqrt{s/(r_1 r_2 r_3)^2}}
\end{aligned}
\end{equation}
with probability at least $1-2 \exp (-\bar{p}\min(\log mp, \log(21emp/\bar{p})))$,
where $\hat{a}=8\sqrt{2}bc$, $\kappa=\sqrt{(\log(m^2p))/T}$ and
$\bar{c}=\frac{2^{\frac{9}{2}}}{\rho}$.

Furthermore, with probability at least $1-2 \exp (-\bar{p}\min(\log mp, \log(21emp/\bar{p})))$,
\[
\begin{aligned}
T^{-1}\sum_{t=1}^{T}\|\widehat{\mathcal{W}}_{(1)}\mathbf{x}_{t}-\mathcal{W}_{(1)}\mathbf{x}_{t}\|_2^2
&=T^{-1}\|\mathbf{X}\widehat{\mathcal{W}}_{(1)}^{H}-\mathbf{X}\mathcal{W}_{(1)}^{H}\|_F^2\\
&=T^{-1}\|\operatorname{vec}(\mathbf{X}\widehat{\mathcal{W}}_{(1)}^{H})
-\operatorname{vec}(\mathbf{X}\mathcal{W}_{(1)}^{H})\|_2^2\\
&=T^{-1}\|\left(\mathbf{I}_m \otimes \mathbf{X}\right)
\widehat{\mathbf{A}}\widehat{\mathbf{g}}-\left(\mathbf{I}_m \otimes \mathbf{X}\right)
\mathbf{A}\mathbf{g}\|_2^2\\
&=T^{-1}\|\left(\mathbf{I}_m \otimes \mathbf{X}\right)
\widehat{\Delta}\|_2^2\\
&=\widehat{\Delta}^H\left(\mathbf{I}_m \otimes (\mathbf{X}^H\mathbf{X}/T)\right)
\widehat{\Delta}\\
&\leq \lambda_{\max}(\mathbf{X}^H\mathbf{X}/T)
\|\widehat{\Delta}\|_2^2\\
&\leq\lambda_{\max}
\left(\boldsymbol{\Sigma}_{\boldsymbol{\varepsilon}}\right) / \mu_{\min }(\mathbf{W})
\|\widehat{\Delta}\|_2^2\\
&\leq
\frac{\lambda_{\max}
\left(\boldsymbol{\Sigma}_{\boldsymbol{\varepsilon}}\right) (4\beta(\sqrt{\tilde{s}}+b\sqrt{s})
+\hat{a}\kappa
+\bar{c}\sum_{i=1}^{3}\alpha_il_i \eta_i)^2}{\mu_{\min }(\mathbf{W})(\varpi+2\kappa\sqrt{s/(r_1 r_2 r_3)^2})^2},
\end{aligned}
\]
where the first equality follows from (\ref{VAR_matrix}),
the third equality follows from (\ref{VAR_vec}),
the fourth equality follows from
$\widehat{\Delta}=\widehat{\mathbf{A}}\widehat{\mathbf{g}}-\mathbf{A}\mathbf{g}$, the first inequality follows from Rayleigh\text{-}Ritz theorem \cite[Theorem 4.2.2]{Horn2013} and $\lambda_{\max}\left(\mathbf{I}_m \otimes (\mathbf{X}^H\mathbf{X}/T)\right)=\lambda_{\max}( (\mathbf{X}^H\mathbf{X}/T))$,
the second inequality follows from the equation (6) of \cite{Masini2022}
and the third inequality follows from (\ref{result1}).
This completes the proof.
\qed

\section*{Appendix C. Gradient Lipschitz Continuity of $\Psi$ in (\ref{DefPsi})}

\begin{lemma}\label{gradient_ele}
Let $f(\mathcal{X}):\mathbb{R}^{n_1 \times n_2 \times n_3}\rightarrow \mathbb{R} $ be a differentiable
function,
and  $\bar{f}(\mathcal{X}_{(1)}):\mathbb{R}^{n_1 \times n_2n_3}\rightarrow \mathbb{R} $ be a differentiable
function with respect to the mode-1 unfolding of $\mathcal{X}$,
where $f(\mathcal{X})=\bar{f}(\mathcal{X}_{(1)})$ for any $\mathcal{X}$.
Then, one has
$$
(\nabla f(\mathcal{X}))_{(1)} = \nabla \bar{f}(\mathcal{X}_{(1)}).
$$
\end{lemma}
\textbf{Proof}.\enspace
For any third-order tensor $\mathcal{X} \in \mathbb{R}^{n_1 \times n_2 \times n_3}$, we get that the $(i,n)$-th element of $\mathcal{X}_{(1)}$ is given by
$$
(\mathcal{X}_{(1)})_{i,n} = \mathcal{X}_{ijl}, \quad \text{where } n = (l-1)n_2 + j, \quad i \in [n_1], j \in [n_2], l \in [n_3].
$$
Note that
$
(\nabla f(\mathcal{X}))_{ijl} = \frac{\partial f}{\partial \mathcal{X}_{ijl}},
$
and
$
((\nabla f(\mathcal{X}))_{(1)})_{i, n} = (\nabla f(\mathcal{X}))_{ijl}, \ \text{where } n = (l-1) n_2 + j.
$
For the gradient of $\bar{f}$ with respect to $ \mathcal{X}_{(1)}$, denoted as $\nabla \bar{f}(\mathcal{X}_{(1)})\in \mathbb{R}^{n_1 \times n_2n_3}$, one has
$$
(\nabla \bar{f}(\mathcal{X}_{(1)}))_{i, n} = \frac{\partial \bar{f}}{\partial (\mathcal{X}_{(1)})_{i, n}}.
$$
Since $(\mathcal{X}_{(1)})_{i, n} = \mathcal{X}_{ijl}$ with $n = (l-1)n_2 + j$ and
$f(\mathcal{X})=\bar{f}(\mathcal{X}_{(1)})$, we get
$$
\frac{\partial \bar{f}}{\partial (\mathcal{X}_{(1)})_{i, n}} = \frac{\partial f}{\partial \mathcal{X}_{ijl}} = (\nabla f(\mathcal{X}))_{ijl}.
$$
Therefore, we obtain that
$
(\nabla f(\mathcal{X}))_{(1)} = \nabla \bar{f}(\mathcal{X}_{(1)}).
$
\qed
\begin{remark}
Let $\bar{Q}(\mathcal{W}_{(1)}):=\frac{1}{2T} \sum_{t=1}^{T}\left\|\mathbf{y}_t-\mathcal{W}_{(1)}\mathbf{x}_{t}\right\|_2^2$ be a
function with respect to $\mathcal{W}_{(1)}$.
Denote the gradient of $\bar{Q}$ with respect to
$\mathcal{W}_{(1)}$ by $\nabla\bar{Q}(\mathcal{W}_{(1)})$.
By Lemma \ref{gradient_ele} and (\ref{DefQ}), we have
\begin{equation}\label{gradientQ_W1}
\nabla\bar{Q}(\mathcal{W}_{(1)})=\frac{1}{T} \sum_{t=1}^{T}\left(\mathcal{W}_{(1)}\mathbf{x}_t-\mathbf{y}_t\right) \mathbf{x}_t^{H}
=(\nabla Q(\llbracket \mathcal{G}; \mathbf{A}_1, \mathbf{A}_2, \mathbf{A}_3\rrbracket))_{(1)},
\end{equation}
where $\nabla Q(\llbracket \mathcal{G}; \mathbf{A}_1, \mathbf{A}_2, \mathbf{A}_3\rrbracket)$ is defined in (\ref{gradient_w}).
\end{remark}

\begin{lemma}\label{lip1}
Suppose that
$\mathbf{A}_i \in \mathbb{R}^{n_i \times r_i}$ with $\mathbf{A}_i\in\mathfrak{B}_{i}$, where $n_1=n_2=m, n_3=p$.
Consequently, $\nabla_\mathcal{G} \Psi(\mathcal{G},\mathbf{A}_1,\mathbf{A}_2,\mathbf{A}_3,\mathbf{U}_1,
\mathbf{U}_2,\mathbf{U}_3)$ is Lipschitz continuous
with Lipschitz constant $L_1:=\frac{1}{T} \sum_{t=1}^T\| \mathbf{x}_t\|_2^2$,
that is,  there exists  a constant $L_1>0$ such that for any $\mathcal{G}^1, \mathcal{G}^2\in \mathbb{R}^{r_1 \times r_2 \times r_3}$, the  following inequality holds:
$$
\left\|\nabla_{\mathcal{G}} \Psi(\mathcal{G}^1,\mathbf{A}_1,\mathbf{A}_2,\mathbf{A}_3,\mathbf{U}_1,
\mathbf{U}_2,\mathbf{U}_3)-\nabla_{\mathcal{G}} \Psi(\mathcal{G}^2,\mathbf{A}_1,\mathbf{A}_2,\mathbf{A}_3,\mathbf{U}_1,
\mathbf{U}_2,\mathbf{U}_3)\right\|_F \leq L_1\left\|\mathcal{G}^1-\mathcal{G}^2\right\|_F.
$$
\end{lemma}
\textbf{Proof}.\enspace
Suppose that $\mathcal{G}^1, \mathcal{G}^2\in \mathbb{R}^{r_1 \times r_2 \times r_3}$, $\mathbf{A}_i\in\mathfrak{B}_{i}, i=1,2,3$, denote
$\mathcal{W}^i=\mathcal{G}^i\times_1 \mathbf{A}_1 \times_2 \mathbf{A}_2 \times_3 \mathbf{A}_3, i=1,2$.
Combining with
(\ref{gradient_w}) and (\ref{gradient_Q}) together, we obtain
\begin{equation}\label{ineq3}
\begin{aligned}
&\left\|\nabla_{\mathcal{G}} Q(\llbracket \mathcal{G}^1; \mathbf{A}_1, \mathbf{A}_2, \mathbf{A}_3\rrbracket)-\nabla_{\mathcal{G}} Q(\llbracket \mathcal{G}^2; \mathbf{A}_1, \mathbf{A}_2, \mathbf{A}_3\rrbracket)\right\|_F\\
=& \left\| \bigg(\frac{1}{T} \sum_{t=1}^{T}(\mathcal{W}_{(1)}^1 \mathbf{x}_t-\mathbf{y}_t)
 \circ \mathbf{X}_t\bigg) \times \mathbf{A}_1^{H} \times \mathbf{A}_2^{H} \times \mathbf{A}_3^{H}
\right.\\
&-\left.
\bigg(\frac{1}{T} \sum_{t=1}^{T}(\mathcal{W}_{(1)}^2 \mathbf{x}_t-\mathbf{y}_t)
\circ \mathbf{X}_t\bigg) \times \mathbf{A}_1^{H} \times \mathbf{A}_2^{H} \times \mathbf{A}_3^{H} \right\|_F \\
=&\left\| \bigg(\frac{1}{T} \sum_{t=1}^{T}\left((\mathcal{W}_{(1)}^1 \mathbf{x}_t-\mathbf{y}_t)
\circ \mathbf{X}_t-( \mathcal{W}_{(1)}^2 \mathbf{x}_t-\mathbf{y}_t)
\circ \mathbf{X}_t\right) \bigg) \times \mathbf{A}_1^{H} \times \mathbf{A}_2^{H}
 \times \mathbf{A}_3^{H}\right\|_F \\
=&\left\| \frac{1}{T} \sum_{t=1}^{T}(\mathcal{W}_{(1)}^1
\mathbf{x}_t-\mathbf{y}_t) \circ \mathbf{X}_t-
\frac{1}{T} \sum_{t=1}^{T}(\mathcal{W}_{(1)}^2
 \mathbf{x}_t-\mathbf{y}_t) \circ \mathbf{X}_t\right\|_F\\
=&\left\| \bigg(\frac{1}{T} \sum_{t=1}^{T}(\mathcal{W}_{(1)}^1
\mathbf{x}_t-\mathbf{y}_t) \circ \mathbf{X}_t\bigg)_{(1)}-
\bigg(\frac{1}{T} \sum_{t=1}^{T}(\mathcal{W}_{(1)}^2
 \mathbf{x}_t-\mathbf{y}_t) \circ \mathbf{X}_t\bigg)_{(1)}\right\|_F,
\end{aligned}
\end{equation}
where the third equality holds since the Frobenius norm is unitarily invariant.
By (\ref{gradientQ_W1}) and (\ref{gradient_w}), we have
\begin{equation}\label{ineq1}
\begin{aligned}
&\left\| \bigg(\frac{1}{T} \sum_{t=1}^{T}(\mathcal{W}_{(1)}^1
\mathbf{x}_t-\mathbf{y}_t) \circ \mathbf{X}_t\bigg)_{(1)}-
\bigg(\frac{1}{T} \sum_{t=1}^{T}(\mathcal{W}_{(1)}^2
 \mathbf{x}_t-\mathbf{y}_t) \circ \mathbf{X}_t\bigg)_{(1)}\right\|_F\\
=\ &\left\| \frac{1}{T} \sum_{t=1}^{T}(\mathcal{W}_{(1)}^1 \mathbf{x}_t-\mathbf{y}_t)
\mathbf{x}_t^{H}-\frac{1}{T}\sum_{t=1}^T (\mathcal{W}_{(1)}^2
\mathbf{x}_t-\mathbf{y}_t)\mathbf{x}_t^{H}\right\|_F \\
= \ &\left\|\frac{1}{T} \sum_{t=1}^T (\mathcal{W}_{(1)}^1 -\mathcal{W}_{(1)}^2
) \mathbf{x}_t\mathbf{x}_t^{H}\right\|_F\\
\leq \ & \left\|\mathcal{W}_{(1)}^1 -\mathcal{W}_{(1)}^2\right\|_F
\left\|\frac{1}{T} \sum_{t=1}^T \mathbf{x}_t \mathbf{x}_t^{H}\right\|_F\\
\leq& \left\|\mathcal{W}_{(1)}^1 -\mathcal{W}_{(1)}^2\right\|_F
\frac{1}{T} \sum_{t=1}^T\left\| \mathbf{x}_t \mathbf{x}_t^{H}\right\|_F\\
= \ & L_1\left\|\mathcal{W}_{(1)}^1 -\mathcal{W}_{(1)}^2\right\|_F,
\end{aligned}
\end{equation}
where the first inequality follows from the Cauchy-Schwarz inequality
and the last equality follows from
$L_1:=\frac{1}{T} \sum_{t=1}^T\| \mathbf{x}_t\|_2^2=\frac{1}{T} \sum_{t=1}^T\| \mathbf{x}_t \mathbf{x}_t^{H}\|_F$.
Observe that
$$
\begin{aligned}
\|\mathcal{W}_{(1)}^1 -\mathcal{W}_{(1)}^2\|_F & =\|\mathcal{W}^1 -\mathcal{W}^2\|_F \\
& =\|\mathcal{G}^1\times_1 \mathbf{A}_1 \times_2 \mathbf{A}_2 \times_3\mathbf{A}_3-
\mathcal{G}^2\times_1 \mathbf{A}_1 \times_2 \mathbf{A}_2 \times_3\mathbf{A}_3\|_F \\
& =\|(\mathcal{G}^1-\mathcal{G}^2) \times_1 \mathbf{A}_1 \times_2 \mathbf{A}_2 \times_3\mathbf{A}_3\|_F \\
& =\|\mathcal{G}^1-\mathcal{G}^2\|_F.
\end{aligned}
$$
Combining the above results and (\ref{ineq1}), (\ref{ineq3}), we have
$$
\|\nabla_{\mathcal{G}} Q(\llbracket \mathcal{G}^1; \mathbf{A}_1, \mathbf{A}_2, \mathbf{A}_3\rrbracket)-\nabla_{\mathcal{G}} Q(\llbracket \mathcal{G}^2; \mathbf{A}_1, \mathbf{A}_2, \mathbf{A}_3\rrbracket)\|_F
\leq L_1\|\mathcal{G}^1-\mathcal{G}^2\|_F.
$$
From (\ref{gradient_H}), we further obtain
$$
\|\nabla_{\mathcal{G}} \Psi(\mathcal{G}^1,\mathbf{A}_1,\mathbf{A}_2,\mathbf{A}_3,\mathbf{U}_1,
\mathbf{U}_2,\mathbf{U}_3)-\nabla_{\mathcal{G}} \Psi(\mathcal{G}^2,\mathbf{A}_1,\mathbf{A}_2,\mathbf{A}_3,\mathbf{U}_1,
\mathbf{U}_2,\mathbf{U}_3)\|_F
\leq L_1\|\mathcal{G}^1-\mathcal{G}^2\|_F.
$$
This completes the proof.
\qed

\begin{lemma}\label{lip2}
Suppose that
$\mathcal{G} \in \mathbb{R}^{r_1 \times r_2 \times r_3}$ with $\mathcal{G}\in\mathcal{D}$, and
$\mathbf{A}_i \in \mathbb{R}^{n_i \times r_i}, i=2,3$ with $\mathbf{A}_i\in\mathfrak{B}_{i}$, where $n_2=m, n_3=p$.
Then $\nabla_{\mathbf{A}_1}\Psi(\mathcal{G},\mathbf{A}_1,\mathbf{A}_2,\mathbf{A}_3,\mathbf{U}_1,
\mathbf{U}_2,\mathbf{U}_3)$ is Lipschitz continuous
with Lipschitz constant $L_2:=\nu^2c_1+\gamma_1$, where $\nu:=\sqrt{r_1r_2r_3}c$ and $c_1:=\frac{1}{T} \sum_{t=1}^T\| \mathbf{x}_t\|_2^2$,
that is, there exists a constant $L_2>0$ such that for
any $\mathbf{A}_1^1, \mathbf{A}_1^2\in \mathbb{R}^{m \times r_1}$, the following inequality holds:
$$
\|\nabla_{\mathbf{A}_1} \Psi(\mathcal{G},\mathbf{A}_1^{1},\mathbf{A}_2,\mathbf{A}_3,\mathbf{U}_1,
\mathbf{U}_2,\mathbf{U}_3)-\nabla_{\mathbf{A}_1} \Psi(\mathcal{G},\mathbf{A}_1^{2},\mathbf{A}_2,\mathbf{A}_3,\mathbf{U}_1,
\mathbf{U}_2,\mathbf{U}_3)\|_F \leq L_2\|\mathbf{A}_1^{1}-\mathbf{A}_1^{2}\|_F.
$$
\end{lemma}
\textbf{Proof}.\enspace
Suppose that $\mathcal{G}\in\mathcal{D}$, $\mathbf{A}_i\in\mathfrak{B}_{i}, i=2,3$.
By (\ref{gradient_w}) and (\ref{gradient_Q}),
for any $\mathbf{A}_1^{1},\mathbf{A}_1^{2}\in\mathbb{R}^{m\times r_1}$,
we have
\begin{equation}\label{lip1_1}
\begin{aligned}
&\left\|\nabla_{\mathbf{A}_1} Q(\llbracket \mathcal{G}; \mathbf{A}_1^1, \mathbf{A}_2, \mathbf{A}_3\rrbracket)-\nabla_{\mathbf{A}_1} Q(\llbracket \mathcal{G}; \mathbf{A}_1^2, \mathbf{A}_2, \mathbf{A}_3\rrbracket)\right\|_F \\
=&\left\|\bigg(\frac{1}{T}\sum_{t=1}^{T}\left(\mathbf{A}_1^{1}\mathcal{G}_{(1)}
\left(\mathbf{A}_3 \otimes \mathbf{A}_2\right)^{H} \mathbf{x}_t-\mathbf{y}_t\right)
\circ \mathbf{X}_t\bigg)_{(1)}\left(\mathbf{A}_3 \otimes \mathbf{A}_2\right)\mathcal{G}_{(1)}^{H}
\right.\\
&-\left.
\bigg(\frac{1}{T}\sum_{t=1}^{T}\left(\mathbf{A}_1^{2}\mathcal{G}_{(1)}
\left(\mathbf{A}_3 \otimes \mathbf{A}_2\right)^{H} \mathbf{x}_t-\mathbf{y}_t\right)
\circ \mathbf{X}_t\bigg)_{(1)}\left(\mathbf{A}_3 \otimes \mathbf{A}_2\right)\mathcal{G}_{(1)}^{H}\right\|_F\\
\leq&\left\|\bigg(\frac{1}{T}\sum_{t=1}^{T}\left(\mathbf{A}_1^{1}\mathcal{G}_{(1)}
\left(\mathbf{A}_3 \otimes \mathbf{A}_2\right)^{H} \mathbf{x}_t-\mathbf{y}_t\right)
\circ \mathbf{X}_t\bigg)_{(1)}
\right.\\
&-\left.
\bigg(\frac{1}{T}\sum_{t=1}^{T}\left(\mathbf{A}_1^{2}\mathcal{G}_{(1)}
\left(\mathbf{A}_3 \otimes \mathbf{A}_2\right)^{H} \mathbf{x}_t-\mathbf{y}_t\right)
\circ \mathbf{X}_t\bigg)_{(1)}\right\|_F\left\|\left(\mathbf{A}_3 \otimes \mathbf{A}_2\right)\mathcal{G}_{(1)}^{H}
\right\|_F\\
=&\left\|\frac{1}{T}\sum_{t=1}^{T}\left(\mathbf{A}_1^{1}\mathcal{G}_{(1)}
\left(\mathbf{A}_3 \otimes \mathbf{A}_2\right)^{H} \mathbf{x}_t-\mathbf{y}_t\right)
\circ \mathbf{X}_t
\right.\\
&-\left.
\frac{1}{T}\sum_{t=1}^{T}\left(\mathbf{A}_1^{2}\mathcal{G}_{(1)}
\left(\mathbf{A}_3 \otimes \mathbf{A}_2\right)^{H} \mathbf{x}_t-\mathbf{y}_t\right)
\circ \mathbf{X}_t\right\|_F
\left\|\left(\mathbf{A}_3 \otimes \mathbf{A}_2\right)\mathcal{G}_{(1)}^{H}
\right\|_F,
\end{aligned}
\end{equation}
where the first inequality follows from the Cauchy-Schwarz inequality.
Proposition \ref{kron_untari} together with \cite[Lemma 2.1]{Chen2021} gives
\begin{equation}\label{ineq2}
\begin{aligned}
\left\|\left(\mathbf{A}_3 \otimes \mathbf{A}_2\right)\mathcal{G}_{(1)}^{H}
\right\|_F
\leq \left\|\mathbf{A}_3 \otimes \mathbf{A}_2\right\|
\left\|\mathcal{G}_{(1)}^{H}
\right\|_F
&\leq
\left\|\mathcal{G}_{(1)}^{H}
\right\|_F=\|\mathcal{G}_{(1)}\|_F=\|\mathcal{G}
\|_F\leq \nu,
\end{aligned}
\end{equation}
where
the last inequality follows from
$\| \mathcal{G}\|_F\leq \sqrt{r_1r_2r_3}\| \mathcal{G}\|_{\infty}\leq \sqrt{r_1r_2r_3}c:=\nu$
for any $\mathcal{G}\in\mathcal{D}$.
In addition,
\begin{equation}\label{lip1_2}
\begin{aligned}
&\left\|\frac{1}{T}\sum_{t=1}^{T}\left(\mathbf{A}_1^{1}\mathcal{G}_{(1)}
\left(\mathbf{A}_3 \otimes \mathbf{A}_2\right)^{H} \mathbf{x}_t-\mathbf{y}_t\right) \circ \mathbf{X}_t
-
\frac{1}{T} \sum_{t=1}^{T}\left(\mathbf{A}_1^{2}\mathcal{G}_{(1)}\left(\mathbf{A}_3
\otimes \mathbf{A}_2\right)^{H} \mathbf{x}_t-\mathbf{y}_t\right) \circ \mathbf{X}_t\right\|_F\\
=\ &\left\|\frac{1}{T} \sum_{t=1}^{T}\left(\mathbf{A}_1^{1}\mathcal{G}_{(1)}\left(\mathbf{A}_3
\otimes \mathbf{A}_2\right)^{H} \mathbf{x}_t-\mathbf{A}_1^{2}\mathcal{G}_{(1)}
\left(\mathbf{A}_3 \otimes \mathbf{A}_2\right)^{H} \mathbf{x}_t\right) \circ \mathbf{X}_t\right\|_F\\
=\ &\left\|\frac{1}{T} \sum_{t=1}^{T}\left(\left(\mathbf{A}_1^{1}-\mathbf{A}_1^{2}\right)
\mathcal{G}_{(1)}\left(\mathbf{A}_3 \otimes \mathbf{A}_2\right)^{H} \mathbf{x}_t\right) \circ
\mathbf{X}_t\right\|_F\\
\leq\ & \frac{1}{T} \sum_{t=1}^{T} \left\|\left(\left(\mathbf{A}_1^{1}-\mathbf{A}_1^{2}\right)
\mathcal{G}_{(1)}\left(\mathbf{A}_3 \otimes \mathbf{A}_2\right)^{H} \mathbf{x}_t\right) \circ \mathbf{X}_t\right\|_F\\
= \ & \frac{1}{T} \sum_{t=1}^{T} \left\|\left(\mathbf{A}_1^{1}-\mathbf{A}_1^{2}\right)\mathcal{G}_{(1)}
\left(\mathbf{A}_3 \otimes \mathbf{A}_2\right)^{H} \mathbf{x}_t \right\|_2 \left\| \mathbf{X}_t\right\|_F\\
\leq \ & \frac{1}{T} \sum_{t=1}^{T} \|\mathbf{A}_1^{1}-\mathbf{A}_1^{2}\|_F
\left\|\mathcal{G}_{(1)}\left(\mathbf{A}_3 \otimes \mathbf{A}_2\right)^{H} \right\|_F
\|\mathbf{x}_t \|_2 \| \mathbf{X}_t\|_F\\
\leq \ & \frac{1}{T} \sum_{t=1}^{T} \nu\|\mathbf{A}_1^{1}-\mathbf{A}_1^{2}\|_F
\|\mathbf{x}_t \|_2 \| \mathbf{X}_t\|_F,
\end{aligned}
\end{equation}
where the last equality holds by
Proposition \ref{out_kronFnorm},
the second inequality holds by the Cauchy-Schwarz inequality
and the last inequality holds by (\ref{ineq2}) and the fact that  $\|\mathcal{G}_{(1)}
\left(\mathbf{A}_3 \otimes \mathbf{A}_2\right)^{H}\|_F=\|
\left(\mathbf{A}_3 \otimes \mathbf{A}_2\right)\mathcal{G}_{(1)}^{H}\|_F$.
From the definitions of $\mathbf{x}_t$ and $\mathbf{X}_t$,
it is easily seen that
\begin{equation}\label{def=def}
\|\mathbf{x}_t\|_2=\left\| \mathbf{X}_t\right\|_F.
\end{equation}
Combining (\ref{lip1_1}), (\ref{ineq2}) and (\ref{lip1_2}), we can easily obtain that
$$
\begin{aligned}
\|\nabla_{\mathbf{A}_1} Q(\llbracket \mathcal{G}; \mathbf{A}_1^1, \mathbf{A}_2, \mathbf{A}_3\rrbracket)-\nabla_{\mathbf{A}_1} Q(\llbracket \mathcal{G}; \mathbf{A}_1^2, \mathbf{A}_2, \mathbf{A}_3\rrbracket)\|_F
&\leq  \nu^2 \frac{1}{T} \sum_{t=1}^{T}
\left\|\mathbf{x}_t \right\|_2^2
\|\mathbf{A}_1^{1}-\mathbf{A}_1^{2}\|_F\\
&=\nu^2c_1
\|\mathbf{A}_1^{1}-\mathbf{A}_1^{2}\|_F,
\end{aligned}
$$
where $c_1:=\frac{1}{T} \sum_{t=1}^{T}
\|\mathbf{x}_t \|_2^2$.
This taken collectively with (\ref{gradient_H}) yields
\begin{equation}\label{lip1_4}
\begin{aligned}
&\|\nabla_{\mathbf{A}_1} \Psi(\mathcal{G},\mathbf{A}_1^{1},\mathbf{A}_2,\mathbf{A}_3,\mathbf{U}_1,
\mathbf{U}_2,\mathbf{U}_3)-\nabla_{\mathbf{A}_1} \Psi(\mathcal{G},\mathbf{A}_1^{2},\mathbf{A}_2,\mathbf{A}_3,\mathbf{U}_1,
\mathbf{U}_2,\mathbf{U}_3)\|_F\\
=\ &\|\nabla_{\mathbf{A}_1} Q(\llbracket \mathcal{G}; \mathbf{A}_1^1, \mathbf{A}_2, \mathbf{A}_3\rrbracket)-\nabla_{\mathbf{A}_1} Q (\llbracket \mathcal{G}; \mathbf{A}_1^2, \mathbf{A}_2, \mathbf{A}_3\rrbracket) +\gamma_1(\mathbf{A}_1^{1}-\mathbf{A}_1^{2})\|_F \\
\leq \ &\|\nabla_{\mathbf{A}_1} Q(\llbracket \mathcal{G}; \mathbf{A}_1^1, \mathbf{A}_2, \mathbf{A}_3\rrbracket)-\nabla_{\mathbf{A}_1} Q(\llbracket \mathcal{G}; \mathbf{A}_1^2, \mathbf{A}_2, \mathbf{A}_3\rrbracket) \|_F+\gamma_1\|\mathbf{A}_1^{1}-\mathbf{A}_1^{2}\|_F  \\
\leq \ &\nu^2c_1\|\mathbf{A}_1^{1}-\mathbf{A}_1^{2}\|_F +\gamma_1\|\mathbf{A}_1^{1}-\mathbf{A}_1^{2}\|_F \\
= \ & (\nu^2c_1+\gamma_1)\|\mathbf{A}_1^{1}-\mathbf{A}_1^{2}\|_F=L_2\|\mathbf{A}_1^{1}-\mathbf{A}_1^{2}\|_F,
\end{aligned}
\end{equation}
where $L_2:=\nu^2c_1+\gamma_1$. This completes the proof.
\qed

\begin{lemma}\label{lip3}
Suppose that
$\mathcal{G} \in \mathbb{R}^{r_1 \times r_2 \times r_3}$ with $\mathcal{G}\in\mathcal{D}$,
$\mathbf{A}_i \in \mathbb{R}^{n_i \times r_i}$ with $\mathbf{A}_i\in\mathfrak{B}_{i}, i=1,3$, where $n_1=m, n_3=p$.
Then $\nabla_{\mathbf{A}_2}\Psi(\mathcal{G},\mathbf{A}_1,\mathbf{A}_2,\mathbf{A}_3,\mathbf{U}_1,
\mathbf{U}_2,\mathbf{U}_3)$ is Lipschitz continuous
with Lipschitz constant $L_3:=\nu^2c_1+\gamma_2$, where $\nu:=\sqrt{r_1r_2r_3}c$ and $c_1:=\frac{1}{T} \sum_{t=1}^T\| \mathbf{x}_t\|_2^2$,
that is, there exists a constant $L_3>0$ such that
for any $\mathbf{A}_2^1, \mathbf{A}_2^2\in \mathbb{R}^{m \times r_2}$, the following inequality holds:
$$
\left\|\nabla_{\mathbf{A}_2} \Psi(\mathcal{G},\mathbf{A}_1,\mathbf{A}_2^{1},\mathbf{A}_3,\mathbf{U}_1,
\mathbf{U}_2,\mathbf{U}_3)-\nabla_{\mathbf{A}_2} \Psi(\mathcal{G},\mathbf{A}_1,\mathbf{A}_2^{2},\mathbf{A}_3,\mathbf{U}_1,
\mathbf{U}_2,\mathbf{U}_3)\right\|_F \leq L_3\left\|\mathbf{A}_2^{1}-\mathbf{A}_2^{2}\right\|_F.
$$
\end{lemma}
\textbf{Proof}.\enspace
Similar to the analysis of $\nabla_{\mathbf{A}_1}Q$,
we have
\begin{equation}\label{lip2_1}
\begin{aligned}
&\left\|\nabla_{\mathbf{A}_2} Q(\llbracket \mathcal{G}; \mathbf{A}_1, \mathbf{A}_2^1, \mathbf{A}_3\rrbracket)-\nabla_{\mathbf{A}_2} Q(\llbracket \mathcal{G}; \mathbf{A}_1, \mathbf{A}_2^2, \mathbf{A}_3\rrbracket)\right\|_F \\
=\ &\left\|\bigg(\frac{1}{T} \sum_{t=1}^{T}\left(\mathbf{A}_1\mathcal{G}_{(1)}
\left(\mathbf{A}_3 \otimes \mathbf{A}_2^{1}\right)^{H} \mathbf{x}_t-\mathbf{y}_t\right) \circ \mathbf{X}_t\bigg)_{(2)}
\left(\mathbf{A}_1 \otimes \mathbf{A}_3\right)\mathcal{G}_{(2)}^{H}
\right.\\
&-\left.
\bigg(\frac{1}{T} \sum_{t=1}^{T}\left(\mathbf{A}_1\mathcal{G}_{(1)}\left(\mathbf{A}_3
\otimes \mathbf{A}_2^{2}\right)^{H} \mathbf{x}_t-\mathbf{y}_t\right) \circ
\mathbf{X}_t\bigg)_{(2)}\left(\mathbf{A}_1 \otimes \mathbf{A}_3\right)\mathcal{G}_{(2)}^{H}\right\|_F\\
\leq \ &\left\|\bigg(\frac{1}{T} \sum_{t=1}^{T}\left(\mathbf{A}_1\mathcal{G}_{(1)}
\left(\mathbf{A}_3 \otimes \mathbf{A}_2^{1}\right)^{H} \mathbf{x}_t-\mathbf{y}_t\right)
\circ \mathbf{X}_t\bigg)_{(2)}
\right.\\
&-\left.
\bigg(\frac{1}{T} \sum_{t=1}^{T}\left(\mathbf{A}_1\mathcal{G}_{(1)}\left(\mathbf{A}_3
\otimes \mathbf{A}_2^{2}\right)^{H} \mathbf{x}_t-\mathbf{y}_t\right) \circ \mathbf{X}_t\bigg)_{(2)}\right\|_F \nu\\
= \ &\left\|\frac{1}{T} \sum_{t=1}^{T}\mathbf{X}_t\otimes \left(\mathbf{A}_1\mathcal{G}_{(1)}
\left(\mathbf{A}_3 \otimes \mathbf{A}_2^{1}\right)^{H} \mathbf{x}_t-\mathbf{y}_t\right)^{H}
\right.\\
&-\left.
\frac{1}{T} \sum_{t=1}^{T}\mathbf{X}_t\otimes \left(\mathbf{A}_1\mathcal{G}_{(1)}
\left(\mathbf{A}_3 \otimes \mathbf{A}_2^{2}\right)^{H} \mathbf{x}_t-\mathbf{y}_t\right)^{H}\right\|_F \nu,
\end{aligned}
\end{equation}
where $\nu:=\sqrt{r_1r_2r_3}c $ and the last equality  follows from Proposition \ref{out_kron1}.
Note that
\begin{equation}\label{lip2_2}
\begin{aligned}
\ &\left\|\frac{1}{T} \sum_{t=1}^{T}\mathbf{X}_t\otimes \left(\mathbf{A}_1\mathcal{G}_{(1)}
\left(\mathbf{A}_3 \otimes \mathbf{A}_2^{1}\right)^{H} \mathbf{x}_t-\mathbf{y}_t\right)^{H}
\right.\\
&-\left.
\frac{1}{T} \sum_{t=1}^{T}\mathbf{X}_t\otimes \left(\mathbf{A}_1\mathcal{G}_{(1)}
\left(\mathbf{A}_3 \otimes \mathbf{A}_2^{2}\right)^{H} \mathbf{x}_t-\mathbf{y}_t\right)^{H}\right\|_F\\
= \ &\left\|\frac{1}{T} \sum_{t=1}^{T}\mathbf{X}_t\otimes \left(\mathbf{A}_1\mathcal{G}_{(1)}
\left(\mathbf{A}_3 \otimes \mathbf{A}_2^{1}\right)^{H} \mathbf{x}_t-\mathbf{A}_1\mathcal{G}_{(1)}
\left(\mathbf{A}_3 \otimes \mathbf{A}_2^{2}\right)^{H} \mathbf{x}_t\right)^{H} \right\|_F\\
= \ &\left\|\frac{1}{T} \sum_{t=1}^{T}\mathbf{X}_t\otimes \left(\mathbf{A}_1\mathcal{G}_{(1)}
\left(\mathbf{A}_3 \otimes (\mathbf{A}_2^{1}- \mathbf{A}_2^{2})\right)^{H} \mathbf{x}_t\right)^{H} \right\|_F\\
\leq\ &\frac{1}{T} \sum_{t=1}^{T}\left\|\mathbf{X}_t\otimes \left(\mathbf{A}_1\mathcal{G}_{(1)}
\left(\mathbf{A}_3 \otimes (\mathbf{A}_2^{1}- \mathbf{A}_2^{2})\right)^{H} \mathbf{x}_t\right)^{H} \right\|_F\\
=\ &\frac{1}{T} \sum_{t=1}^{T}\left\|\mathbf{X}_t \right\|_F \left\|\mathbf{A}_1\mathcal{G}_{(1)}\left(\mathbf{A}_3 \otimes (\mathbf{A}_2^{1}- \mathbf{A}_2^{2})\right)^{H} \mathbf{x}_t \right\|_2,
\end{aligned}
\end{equation}
where the third equality holds by Proposition \ref{kron_pro}.
It follows  from \cite[Example 5.6.6]{Horn2013} that
\begin{equation}\label{lip2_4}
\begin{aligned}
&\left\|\mathbf{A}_1\mathcal{G}_{(1)}\left(\mathbf{A}_3
\otimes (\mathbf{A}_2^{1}- \mathbf{A}_2^{2})\right)^{H} \mathbf{x}_t \right\|_2\\
\leq \ &\left\|\mathbf{A}_1\right\|\left\|\mathcal{G}_{(1)}\left(\mathbf{A}_3
\otimes (\mathbf{A}_2^{1}- \mathbf{A}_2^{2})\right)^{H}\mathbf{x}_t \right\|_2\\
\leq \ &\left\|\mathcal{G}_{(1)}\left(\mathbf{A}_3
\otimes (\mathbf{A}_2^{1}- \mathbf{A}_2^{2})\right)^{H}\mathbf{x}_t \right\|_2\\
\leq \ &\left\|\mathcal{G}_{(1)}\left(\mathbf{A}_3
\otimes (\mathbf{A}_2^{1}- \mathbf{A}_2^{2})\right)^{H}\right\|_F\left\|\mathbf{x}_t \right\|_2\\
\leq \ &\left\|\mathcal{G}_{(1)}\right\|_F\left\|\left(\mathbf{A}_3
\otimes (\mathbf{A}_2^{1}- \mathbf{A}_2^{2})\right)^{H}\right\|\left\|\mathbf{x}_t \right\|_2\\
= \ & \left\|\mathcal{G}\right\|_F \left\|\mathbf{A}_3 \right\|
\left\| (\mathbf{A}_2^{1}- \mathbf{A}_2^{2})^{H}\right\| \left\| \mathbf{x}_t \right\|_2\\
\leq  \ &  \nu \left\|\mathbf{A}_3 \right\|
\left\| (\mathbf{A}_2^{1}- \mathbf{A}_2^{2})^{H}\right\| \left\| \mathbf{x}_t \right\|_2\\
\leq \ & \nu \left\| (\mathbf{A}_2^{1}- \mathbf{A}_2^{2})^{H}
\right\|_F \left\| \mathbf{x}_t \right\|_2,
\end{aligned}
\end{equation}
where the second inequality holds by $\mathbf{A}_1\in\mathfrak{B}_{1}$ with $\left\|\mathbf{A}_1\right\|\leq 1$,
the third inequality holds by the Cauchy-Schwarz inequality,
the fourth inequality holds by \cite[Lemma 2.1]{Chen2021},
the fifth inequality holds by $\mathcal{G}\in\mathcal{D}$ with
$\|\mathcal{G}\|_F\leq \nu$
and the last inequality holds by $\|\mathbf{A}_3\|\leq1$
and any matrix $\mathbf{X}$ satisfies $\|\mathbf{X}\|\leq\|\mathbf{X}\|_F$.
Combining (\ref{lip2_4}) and (\ref{lip2_2}), we have
\begin{equation}\label{lip2_3}
\begin{aligned}
&\left\|\frac{1}{T} \sum_{t=1}^{T}\mathbf{X}_t\otimes \left(\mathbf{A}_1\mathcal{G}_{(1)}
\left(\mathbf{A}_3 \otimes \mathbf{A}_2^{1}\right)^{H} \mathbf{x}_t-\mathbf{y}_t\right)^{H}
\right.\\
&-\left.
\frac{1}{T} \sum_{t=1}^{T}\mathbf{X}_t\otimes \left(\mathbf{A}_1\mathcal{G}_{(1)}
\left(\mathbf{A}_3 \otimes \mathbf{A}_2^{2}\right)^{H} \mathbf{x}_t-\mathbf{y}_t\right)^{H}\right\|_F\\
\leq \ &\frac{1}{T} \sum_{t=1}^{T}\left\|\mathbf{X}_t \right\|_F \nu \left\| (\mathbf{A}_2^{1}- \mathbf{A}_2^{2})^{H}
\right\|_F \left\| \mathbf{x}_t \right\|_2\\
= \ &\frac{1}{T} \sum_{t=1}^{T}\left\| \mathbf{x}_t \right\|_2^2 \nu\left\| \mathbf{A}_2^{1}- \mathbf{A}_2^{2}\right\|_F
=\nu c_1\left\|\mathbf{A}_2^{1}- \mathbf{A}_2^{2}\right\|_F,
\end{aligned}
\end{equation}
where the first equality follows from (\ref{def=def}) and $c_1:=\frac{1}{T} \sum_{t=1}^{T}\| \mathbf{x}_t \|_2^2$.
Then together with (\ref{lip2_1}),
one can obtain that
$$
\|\nabla_{\mathbf{A}_2} Q(\llbracket \mathcal{G}; \mathbf{A}_1, \mathbf{A}_2^1, \mathbf{A}_3\rrbracket)-\nabla_{\mathbf{A}_2} Q(\llbracket \mathcal{G}; \mathbf{A}_1, \mathbf{A}_2^2, \mathbf{A}_3\rrbracket)\|_F \leq
\nu^2c_1\|\mathbf{A}_2^{1}- \mathbf{A}_2^{2}\|_F.
$$
Similar to (\ref{lip1_4}), we have
$$
\|\nabla_{\mathbf{A}_2} \Psi(\mathcal{G},\mathbf{A}_1,\mathbf{A}_2^{1},\mathbf{A}_3,\mathbf{U}_1,
\mathbf{U}_2,\mathbf{U}_3)-\nabla_{\mathbf{A}_2} \Psi(\mathcal{G},\mathbf{A}_1,\mathbf{A}_2^{2},\mathbf{A}_3,\mathbf{U}_1,
\mathbf{U}_2,\mathbf{U}_3)\|_F \leq
L_3\|\mathbf{A}_2^{1}-\mathbf{A}_2^{2}\|_F,
$$
where $L_3:=\nu^2c_1+\gamma_2$.
The proof is completed.
\qed
\begin{lemma}\label{lip4}
Suppose that
$\mathcal{G} \in \mathbb{R}^{r_1 \times r_2 \times r_3}$ with $\mathcal{G}\in\mathcal{D}$,
$\mathbf{A}_i \in \mathbb{R}^{n_i \times r_i}$ with $\mathbf{A}_i\in\mathfrak{B}_{i},i=1,2$,
where $n_1=n_2=m$.
Then $\nabla_{\mathbf{A}_3}\Psi(\mathcal{G},\mathbf{A}_1,\mathbf{A}_2,\mathbf{A}_3,\mathbf{U}_1,
\mathbf{U}_2,\mathbf{U}_3)$ is Lipschitz continuous
with Lipschitz constant $L_4:=\nu^2c_1 +\gamma_3$, where $\nu:=\sqrt{r_1r_2r_3}c$ and $c_1:=\frac{1}{T} \sum_{t=1}^T\| \mathbf{x}_t\|_2^2$,
that is, there exists a constant $L_4>0$ such that
for any $\mathbf{A}_3^1, \mathbf{A}_3^2\in \mathbb{R}^{p \times r_3}$, the following inequality holds:
$$
\|\nabla_{\mathbf{A}_3} \Psi(\mathcal{G},\mathbf{A}_1,\mathbf{A}_2,\mathbf{A}_3^{1},\mathbf{U}_1,
\mathbf{U}_2,\mathbf{U}_3)-\nabla_{\mathbf{A}_3} \Psi(\mathcal{G},\mathbf{A}_1,\mathbf{A}_2,\mathbf{A}_3^{2},\mathbf{U}_1,
\mathbf{U}_2,\mathbf{U}_3)\|_F \leq L_4\|\mathbf{A}_3^{1}-\mathbf{A}_3^{2}\|_F.
$$
\end{lemma}
\textbf{Proof}.\enspace
Similar to the analysis of $\nabla_{\mathbf{A}_1}Q$, by   (\ref{gradient_Q}), we have
\begin{equation}\label{lip3_1}
\begin{aligned}
&\left\|\nabla_{\mathbf{A}_3} Q(\llbracket \mathcal{G}; \mathbf{A}_1, \mathbf{A}_2, \mathbf{A}_3^1\rrbracket)-\nabla_{\mathbf{A}_3} Q(\llbracket \mathcal{G}; \mathbf{A}_1, \mathbf{A}_2, \mathbf{A}_3^2\rrbracket)\right\|_F \\
=&\left\|\bigg(\frac{1}{T} \sum_{t=1}^{T}\left(\mathbf{A}_1\mathcal{G}_{(1)}
\left(\mathbf{A}_3^{1} \otimes \mathbf{A}_2\right)^{H} \mathbf{x}_t-\mathbf{y}_t\right) \circ \mathbf{X}_t\bigg)_{(3)}
\left(\mathbf{A}_2 \otimes \mathbf{A}_1\right)\mathcal{G}_{(3)}^{H}
\right.\\
&-\left.
\bigg(\frac{1}{T} \sum_{t=1}^{T}\left(\mathbf{A}_1\mathcal{G}_{(1)}
\left(\mathbf{A}_3^{2} \otimes \mathbf{A}_2\right)^{H} \mathbf{x}_t-\mathbf{y}_t\right) \circ \mathbf{X}_t\bigg)_{(3)}
\left(\mathbf{A}_2 \otimes \mathbf{A}_1\right)\mathcal{G}_{(3)}^{H}\right\|_F\\
\leq&\left\|\bigg(\frac{1}{T} \sum_{t=1}^{T}\left(\mathbf{A}_1\mathcal{G}_{(1)}
\left(\mathbf{A}_3^{1} \otimes \mathbf{A}_2\right)^{H} \mathbf{x}_t-\mathbf{y}_t\right) \circ \mathbf{X}_t\bigg)_{(3)}
\right.\\
&-\left.
\bigg(\frac{1}{T} \sum_{t=1}^{T}\left(\mathbf{A}_1\mathcal{G}_{(1)}
\left(\mathbf{A}_3^{2} \otimes \mathbf{A}_2\right)^{H} \mathbf{x}_t-\mathbf{y}_t\bigg)
\circ \mathbf{X}_t\right)_{(3)}\right\|_F \nu\\
=&\left\|\frac{1}{T} \sum_{t=1}^{T}\mathbf{X}_t^{H}\otimes \left(\mathbf{A}_1\mathcal{G}_{(1)}\left(\mathbf{A}_3^{1} \otimes \mathbf{A}_2\right)^{H} \mathbf{x}_t-\mathbf{y}_t\right)^{H}
\right.\\
&-\left.
\frac{1}{T} \sum_{t=1}^{T}\mathbf{X}_t^{H}\otimes \left(\mathbf{A}_1\mathcal{G}_{(1)}\left(\mathbf{A}_3^{2} \otimes \mathbf{A}_2\right)^{H} \mathbf{x}_t-\mathbf{y}_t\right)^{H}\right\|_F \nu
\end{aligned}
\end{equation}
where $\nu:=\sqrt{r_1r_2r_3}c $ and  the last equality follows from Proposition \ref{out_kron1}.
Similar
to (\ref{lip2_3}), we have
\begin{equation}\label{TOTB}
\begin{aligned}
&\left\|\frac{1}{T} \sum_{t=1}^{T}\mathbf{X}_t^{H}\otimes \left(\mathbf{A}_1\mathcal{G}_{(1)}\left(\mathbf{A}_3^{1} \otimes \mathbf{A}_2\right)^{H} \mathbf{x}_t-\mathbf{y}_t\right)^{H}
\right.\\
&-\left.
\frac{1}{T} \sum_{t=1}^{T}\mathbf{X}_t^{H}\otimes \left(\mathbf{A}_1\mathcal{G}_{(1)}\left(\mathbf{A}_3^{2} \otimes \mathbf{A}_2\right)^{H} \mathbf{x}_t-\mathbf{y}_t\right)^{H}\right\|_F\\
\leq &\frac{1}{T} \sum_{t=1}^{T}\left\|\mathbf{X}_t \right\|_F \left\|\left(\mathbf{A}_1\mathcal{G}_{(1)}
\left((\mathbf{A}_3^{1}- \mathbf{A}_3^{2}) \otimes \mathbf{A}_2\right)^{H} \mathbf{x}_t\right)^{H} \right\|_F
\leq \nu c_1 \|\mathbf{A}_3^{1}-\mathbf{A}_3^{2}\|_F,
\end{aligned}
\end{equation}
where $c_1:=\frac{1}{T} \sum_{t=1}^{T}\| \mathbf{x}_t \|_2^2$.
Combining  (\ref{lip3_1}) and (\ref{TOTB}), we get that
$$
\|\nabla_{\mathbf{A}_3} Q(\llbracket \mathcal{G}; \mathbf{A}_1, \mathbf{A}_2, \mathbf{A}_3^1\rrbracket)-\nabla_{\mathbf{A}_3} Q(\llbracket \mathcal{G}; \mathbf{A}_1, \mathbf{A}_2, \mathbf{A}_3^2\rrbracket)\|_F \leq  \nu^2c_1\|\mathbf{A}_3^{1}-\mathbf{A}_3^{2}\|_F.
$$
Furthermore, similar to (\ref{lip1_4}), we have
$$
\|\nabla_{\mathbf{A}_3} \Psi(\mathcal{G},\mathbf{A}_1,\mathbf{A}_2,\mathbf{A}_3^{1},\mathbf{U}_1,
\mathbf{U}_2,\mathbf{U}_3)-\nabla_{\mathbf{A}_3} \Psi(\mathcal{G},\mathbf{A}_1,\mathbf{A}_2,\mathbf{A}_3^{2},\mathbf{U}_1,
\mathbf{U}_2,\mathbf{U}_3)\|_F \leq L_4\left\|\mathbf{A}_3^{1}-\mathbf{A}_3^{2}\right\|_F,
$$
where $L_4:=\nu^2c_1 +\gamma_3$.
This completes the proof.
\qed

It can be easily verified from (\ref{gradient_H}) that
\begin{equation}\label{lip_U}
\begin{aligned}
\|\nabla_{\mathbf{U}_1} \Psi(\mathcal{G},\mathbf{A}_1,\mathbf{A}_2,\mathbf{A}_3,\mathbf{U}_1^{1},
\mathbf{U}_2,\mathbf{U}_3)-\nabla_{\mathbf{U}_1} \Psi(\mathcal{G},\mathbf{A}_1,\mathbf{A}_2,\mathbf{A}_3,\mathbf{U}_1^{2},
\mathbf{U}_2,\mathbf{U}_3)\|_F
\leq \gamma_1\|\mathbf{U}_1^{1}-\mathbf{U}_1^{2}\|_F,\\
\|\nabla_{\mathbf{U}_2} \Psi(\mathcal{G},\mathbf{A}_1,\mathbf{A}_2,\mathbf{A}_3,\mathbf{U}_1,
\mathbf{U}_2^{1},\mathbf{U}_3)-\nabla_{\mathbf{U}_2} \Psi(\mathcal{G},\mathbf{A}_1,\mathbf{A}_2,\mathbf{A}_3,\mathbf{U}_1,
\mathbf{U}_2^{2},\mathbf{U}_3)\|_F
\leq \gamma_2\|\mathbf{U}_2^{1}-\mathbf{U}_2^{2}\|_F,\\
\|\nabla_{\mathbf{U}_3} \Psi(\mathcal{G},\mathbf{A}_1,\mathbf{A}_2,\mathbf{A}_3,\mathbf{U}_1,
\mathbf{U}_2,\mathbf{U}_3^{1})-\nabla_{\mathbf{U}_3} \Psi(\mathcal{G},\mathbf{A}_1,\mathbf{A}_2,\mathbf{A}_3,\mathbf{U}_1,
\mathbf{U}_2,\mathbf{U}_3^{2})\|_F
\leq \gamma_3\|\mathbf{U}_3^{1}-\mathbf{U}_3^{2}\|_F.
\end{aligned}
\end{equation}

\section*{Appendix D. Proof of Lemma \ref{suffcond}}
It can be seen from (\ref{palm}) that
\begin{equation}\label{G1}
\begin{aligned}
&\beta\|\mathcal{G}^{k+1}\|_1 +\delta_{\mathcal{D}}(\mathcal{G}^{k+1})+\left\langle
\nabla_{\mathcal{G}} \Psi(\mathcal{G}^{k}, \mathbf{A}_1^{k}, \mathbf{A}_2^{k}, \mathbf{A}_3^{k},\mathbf{U}_1^{k},
\mathbf{U}_2^{k},\mathbf{U}_3^{k}),\mathcal{G}^{k+1}- \mathcal{G}^k\right\rangle+ \frac{\rho_1}{2}\|\mathcal{G}^{k+1}-\mathcal{G}^k\|_F^2\\
&\leq\beta\|\mathcal{G}^{k}\|_1 +\delta_{\mathcal{D}}(\mathcal{G}^{k}).
\end{aligned}
\end{equation}
In addition, applying Lemma \ref{lip1}
and \cite[Lemma 1]{Bolte2014}, we have that
\begin{equation}\label{G2}
\begin{aligned}
&\Psi\left(\mathcal{G}^{k+1}, \mathbf{A}_1^k, \mathbf{A}_2^k, \mathbf{A}_3^k,\mathbf{U}_1^k,
\mathbf{U}_2^k,\mathbf{U}_3^k\right)
\leq \Psi\left(\mathcal{G}^{k}, \mathbf{A}_1^k, \mathbf{A}_2^k, \mathbf{A}_3^k,\mathbf{U}_1^k,
\mathbf{U}_2^k,\mathbf{U}_3^k\right)\\
&+\left\langle \nabla_{\mathcal{G}} \Psi(\mathcal{G}^{k}, \mathbf{A}_1^{k}, \mathbf{A}_2^{k}, \mathbf{A}_3^{k},\mathbf{U}_1^{k},
\mathbf{U}_2^{k},\mathbf{U}_3^{k}),\mathcal{G}^{k+1}-
\mathcal{G}^k\right\rangle+\frac{L_1}{2}\|\mathcal{G}^{k+1}- \mathcal{G}^k\|_F^2.
\end{aligned}
\end{equation}
Combining (\ref{G1}) and (\ref{G2}), we obtain
$$
\begin{aligned}
&\Psi\left(\mathcal{G}^{k+1}, \mathbf{A}_1^k, \mathbf{A}_2^k, \mathbf{A}_3^k,\mathbf{U}_1^k,
\mathbf{U}_2^k,\mathbf{U}_3^k\right)
+\beta\|\mathcal{G}^{k+1}\|_1 +\delta_{\mathcal{D}}(\mathcal{G}^{k+1})\\
&+\left\langle
\nabla_{\mathcal{G}} \Psi(\mathcal{G}^{k}, \mathbf{A}_1^{k}, \mathbf{A}_2^{k}, \mathbf{A}_3^{k},\mathbf{U}_1^{k},
\mathbf{U}_2^{k},\mathbf{U}_3^{k}),\mathcal{G}^{k+1}- \mathcal{G}^k\right\rangle
+\frac{\rho_1}{2}\|\mathcal{G}^{k+1}- \mathcal{G}^k\|_F^2\\
& \leq \Psi\left(\mathcal{G}^{k}, \mathbf{A}_1^k, \mathbf{A}_2^k, \mathbf{A}_3^k,\mathbf{U}_1^k,
\mathbf{U}_2^k,\mathbf{U}_3^k\right)+\left\langle \nabla_{\mathcal{G}} \Psi(\mathcal{G}^{k}, \mathbf{A}_1^{k}, \mathbf{A}_2^{k}, \mathbf{A}_3^{k},\mathbf{U}_1^{k},
\mathbf{U}_2^{k},\mathbf{U}_3^{k}),\mathcal{G}^{k+1}-
\mathcal{G}^k\right\rangle\\
&+\frac{L_1}{2}\|\mathcal{G}^{k+1}- \mathcal{G}^k\|_F^2
+\beta\|\mathcal{G}^{k}\|_1 +\delta_{\mathcal{D}}(\mathcal{G}^{k}),
\end{aligned}
$$
which yields
\begin{equation}\label{suff_condition1}
\begin{aligned}
&\Psi\left(\mathcal{G}^{k+1}, \mathbf{A}_1^k, \mathbf{A}_2^k, \mathbf{A}_3^k,\mathbf{U}_1^k,
\mathbf{U}_2^k,\mathbf{U}_3^k\right)
+\beta\|\mathcal{G}^{k+1}\|_1 +\delta_{\mathcal{D}}(\mathcal{G}^{k+1})+ \frac{\rho_1-L_1}{2}\|\mathcal{G}^{k+1}- \mathcal{G}^k\|_F^2\\
& \leq \Psi\left(\mathcal{G}^{k}, \mathbf{A}_1^k, \mathbf{A}_2^k, \mathbf{A}_3^k,\mathbf{U}_1^k,
\mathbf{U}_2^k,\mathbf{U}_3^k\right)
+\beta\|\mathcal{G}^{k}\|_1 +\delta_{\mathcal{D}}(\mathcal{G}^{k}).
\end{aligned}
\end{equation}
By repeating the similar ideas and techniques, and combining Lemma \ref{lip2}-\ref{lip4} with (\ref{lip_U}) and (\ref{palm}), we can obtain
\begin{equation}\label{suff_condition2}
\begin{aligned}
&\delta_{\mathfrak{B}_1}(\mathbf{A}_1^{k+1})+\Psi\left(\mathcal{G}^{k+1}, \mathbf{A}_1^{k+1}, \mathbf{A}_2^k,
\mathbf{A}_3^k,\mathbf{U}_1^k,
\mathbf{U}_2^k,\mathbf{U}_3^k\right)+ \frac{\rho_2-L_2}{2}\|\mathbf{A}_1^{k+1}-\mathbf{A}_1^{k}\|_F^2\\
&\leq\delta_{\mathfrak{B}_1}(\mathbf{A}_1^{k})+\Psi\left(\mathcal{G}^{k+1}, \mathbf{A}_1^{k}, \mathbf{A}_2^k,
\mathbf{A}_3^k,\mathbf{U}_1^k,
\mathbf{U}_2^k,\mathbf{U}_3^k\right),\\
&\delta_{\mathfrak{B}_2}(\mathbf{A}_2^{k+1})+\Psi\left(\mathcal{G}^{k+1}, \mathbf{A}_1^{k+1}, \mathbf{A}_2^{k+1},
\mathbf{A}_3^k,\mathbf{U}_1^k,
\mathbf{U}_2^k,\mathbf{U}_3^k\right)+ \frac{\rho_3-L_3}{2}\|\mathbf{A}_2^{k+1}-\mathbf{A}_2^{k}\|_F^2\\
&\leq\delta_{\mathfrak{B}_2}(\mathbf{A}_2^{k})+\Psi\left(\mathcal{G}^{k+1}, \mathbf{A}_1^{k+1}, \mathbf{A}_2^k,
\mathbf{A}_3^k,\mathbf{U}_1^k,
\mathbf{U}_2^k,\mathbf{U}_3^k\right),\\
&\delta_{\mathfrak{B}_3}(\mathbf{A}_3^{k+1})+\Psi\left(\mathcal{G}^{k+1}, \mathbf{A}_1^{k+1}, \mathbf{A}_2^{k+1},
\mathbf{A}_3^{k+1},\mathbf{U}_1^k,
\mathbf{U}_2^k,\mathbf{U}_3^k\right)+ \frac{\rho_4-L_4}{2}\|\mathbf{A}_3^{k+1}-\mathbf{A}_3^{k}\|_F^2\\
&\leq \delta_{\mathfrak{B}_3}(\mathbf{A}_3^{k})+\Psi\left(\mathcal{G}^{k+1}, \mathbf{A}_1^{k+1}, \mathbf{A}_2^{k+1},
\mathbf{A}_3^k,\mathbf{U}_1^k,
\mathbf{U}_2^k,\mathbf{U}_3^k\right),\\
&\alpha_1\operatorname{tr}\left((\mathbf{U}_1^{k+1})^{H}\mathbf{L}_1 \mathbf{U}_1^{k+1}\right)+\Psi\left(\mathcal{G}^{k+1},
\mathbf{A}_1^{k+1}, \mathbf{A}_2^{k+1}, \mathbf{A}_3^{k+1},\mathbf{U}_1^{k+1},
\mathbf{U}_2^k,\mathbf{U}_3^k\right)+\frac{\rho_5-\gamma_1}{2}\|\mathbf{U}_1^{k+1}- \mathbf{U}_1^{k}\|_F^2\\
& \leq \alpha_1\operatorname{tr}\left((\mathbf{U}_1^{k})^{H}\mathbf{L}_1 \mathbf{U}_1^{k}\right)+\Psi\left(\mathcal{G}^{k+1},
\mathbf{A}_1^{k+1}, \mathbf{A}_2^{k+1}, \mathbf{A}_3^{k+1},\mathbf{U}_1^k,
\mathbf{U}_2^k,\mathbf{U}_3^k\right),\\
&\alpha_2\operatorname{tr}\left((\mathbf{U}_2^{k+1})^{H}\mathbf{L}_2 \mathbf{U}_2^{k+1}\right)+\Psi\left(\mathcal{G}^{k+1},
\mathbf{A}_1^{k+1}, \mathbf{A}_2^{k+1}, \mathbf{A}_3^{k+1},\mathbf{U}_1^{k+1},
\mathbf{U}_2^{k+1},\mathbf{U}_3^k\right)\\
&+\frac{\rho_6-\gamma_2}{2}\|\mathbf{U}_2^{k+1}- \mathbf{U}_2^{k}\|_F^2\\
& \leq \alpha_2\operatorname{tr}\left((\mathbf{U}_2^{k})^{H}\mathbf{L}_2 \mathbf{U}_2^{k}\right)+\Psi\left(\mathcal{G}^{k+1},
\mathbf{A}_1^{k+1}, \mathbf{A}_2^{k+1}, \mathbf{A}_3^{k+1},\mathbf{U}_1^{k+1},
\mathbf{U}_2^k,\mathbf{U}_3^k\right),\\
&\alpha_3\operatorname{tr}\left((\mathbf{U}_3^{k+1})^{H}\mathbf{L}_3 \mathbf{U}_3^{k+1}\right)+\Psi\left(\mathcal{G}^{k+1},
\mathbf{A}_1^{k+1}, \mathbf{A}_2^{k+1}, \mathbf{A}_3^{k+1},\mathbf{U}_1^{k+1},
\mathbf{U}_2^{k+1},\mathbf{U}_3^{k+1}\right)\\
&+\frac{\rho_7-\gamma_3}{2}
\|\mathbf{U}_3^{k+1}- \mathbf{U}_3^{k}\|_F^2\\
& \leq \alpha_3\operatorname{tr}\left((\mathbf{U}_3^{k})^{H}\mathbf{L}_3 \mathbf{U}_3^{k}\right)+\Psi\left(\mathcal{G}^{k+1},
\mathbf{A}_1^{k+1}, \mathbf{A}_2^{k+1}, \mathbf{A}_3^{k+1},\mathbf{U}_1^{k+1},
\mathbf{U}_2^{k+1},\mathbf{U}_3^k\right).
\end{aligned}
\end{equation}
This together with (\ref{suff_condition1}) leads to
$$
\begin{aligned}
&\beta\|\mathcal{G}^{k+1}\|_1 +\delta_{\mathcal{D}}(\mathcal{G}^{k+1})
+\Psi\left(\mathcal{G}^{k+1}, \mathbf{A}_1^{k+1}, \mathbf{A}_2^{k+1}, \mathbf{A}_3^{k+1},\mathbf{U}_1^{k+1},
\mathbf{U}_2^{k+1},\mathbf{U}_3^{k+1}\right)
+\sum_{i=1}^3\delta_{\mathfrak{B}_i}(\mathbf{A}_i^{k+1})\\
+&\frac{\rho_1-L_1}{2}\|\mathcal{G}^{k+1}- \mathcal{G}^k\|_F^2+
\frac{\rho_2-L_2}{2}\|\mathbf{A}_1^{k+1}-\mathbf{A}_1^{k}\|_F^2
+\frac{\rho_3-L_3}{2}\|\mathbf{A}_2^{k+1}-\mathbf{A}_2^{k}\|_F^2\\
+&\frac{\rho_4-L_4}{2}\|\mathbf{A}_3^{k+1}-\mathbf{A}_3^{k}\|_F^2
+\frac{\rho_5-\gamma_1}{2}\|\mathbf{U}_1^{k+1}- \mathbf{U}_1^{k}\|_F^2
+\frac{\rho_6-\gamma_2}{2}\|\mathbf{U}_2^{k+1}- \mathbf{U}_2^{k}\|_F^2\\
+&\frac{\rho_7-\gamma_3}{2}
\|\mathbf{U}_3^{k+1}- \mathbf{U}_3^{k}\|_F^2
+\sum_{i=1}^3\alpha_i\operatorname{tr}\left((\mathbf{U}_i^{k+1})^{H}\mathbf{L}_i \mathbf{U}_i^{k+1}\right)\\
\leq \ &\beta\|\mathcal{G}^{k}\|_1 +\delta_{\mathcal{D}}(\mathcal{G}^{k})
+\Psi\left(\mathcal{G}^{k}, \mathbf{A}_1^{k}, \mathbf{A}_2^{k}, \mathbf{A}_3^{k},\mathbf{U}_1^{k},
\mathbf{U}_2^{k},\mathbf{U}_3^{k}\right)
+\sum_{i=1}^3\delta_{\mathfrak{B}_i}(\mathbf{A}_i^{k})
+\sum_{i=1}^3\alpha_i\operatorname{tr}\left((\mathbf{U}_i^{k})^{H}\mathbf{L}_i \mathbf{U}_i^{k}\right).
\end{aligned}
$$
Denote $\bar{\rho}:=\min\{\rho_1-L_1,\rho_2-L_2,\rho_3-L_3,
\rho_4-L_4,\rho_5-\gamma_1,\rho_6-\gamma_2,\rho_7-\gamma_3\}$,
then one has
$$
\begin{aligned}
&\beta\|\mathcal{G}^{k+1}\|_1 +\delta_{\mathcal{D}}(\mathcal{G}^{k+1})
+\Psi\left(\mathcal{G}^{k+1}, \mathbf{A}_1^{k+1}, \mathbf{A}_2^{k+1}, \mathbf{A}_3^{k+1},\mathbf{U}_1^{k+1},
\mathbf{U}_2^{k+1},\mathbf{U}_3^{k+1}\right)
+\sum_{i=1}^3\delta_{\mathfrak{B}_i}(\mathbf{A}_i^{k+1})\\
+&\frac{\bar{\rho}}{2}\left(\|\mathcal{G}^{k+1}- \mathcal{G}^k\|_F^2+
\|\mathbf{A}_1^{k+1}-\mathbf{A}_1^{k}\|_F^2
+\|\mathbf{A}_2^{k+1}-\mathbf{A}_2^{k}\|_F^2
+\|\mathbf{A}_3^{k+1}-\mathbf{A}_3^{k}\|_F^2\right.\\
+&\left.\|\mathbf{U}_1^{k+1}- \mathbf{U}_1^{k}\|_F^2
+\|\mathbf{U}_2^{k+1}- \mathbf{U}_2^{k}\|_F^2
+\|\mathbf{U}_3^{k+1}- \mathbf{U}_3^{k}\|_F^2\right)
+\sum_{i=1}^3\alpha_i\operatorname{tr}\left((\mathbf{U}_i^{k+1})^{H}\mathbf{L}_i \mathbf{U}_i^{k+1}\right)\\
\leq \ &\beta\|\mathcal{G}^{k}\|_1 +\delta_{\mathcal{D}}(\mathcal{G}^{k})
+\Psi\left(\mathcal{G}^{k}, \mathbf{A}_1^{k}, \mathbf{A}_2^{k}, \mathbf{A}_3^{k},\mathbf{U}_1^{k},
\mathbf{U}_2^{k},\mathbf{U}_3^{k}\right)
+\sum_{i=1}^3\delta_{\mathfrak{B}_i}(\mathbf{A}_i^{k})
+\sum_{i=1}^3\alpha_i\operatorname{tr}\left((\mathbf{U}_i^{k})^{H}\mathbf{L}_i \mathbf{U}_i^{k}\right).
\end{aligned}
$$
Recalling the definition of $F$ in (\ref{objecF}) and
the definition of $\Psi$ in (\ref{DefPsi}),
the above inequality is equivalent to
\begin{equation}\label{suff}
\begin{aligned}
&F\left(\mathcal{G}^{k+1}, \mathbf{A}_1^{k+1}, \mathbf{A}_2^{k+1}, \mathbf{A}_3^{k+1},\mathbf{U}_1^{k+1},
\mathbf{U}_2^{k+1},\mathbf{U}_3^{k+1}\right)
+\frac{\bar{\rho}}{2}\left(\|\mathcal{G}^{k+1}- \mathcal{G}^k\|_F^2+
\|\mathbf{A}_1^{k+1}-\mathbf{A}_1^{k}\|_F^2
\right.\\
+&\left.
\|\mathbf{A}_2^{k+1}-\mathbf{A}_2^{k}\|_F^2
+\|\mathbf{A}_3^{k+1}-\mathbf{A}_3^{k}\|_F^2
+\|\mathbf{U}_1^{k+1}- \mathbf{U}_1^{k}\|_F^2
+\|\mathbf{U}_2^{k+1}- \mathbf{U}_2^{k}\|_F^2
+\|\mathbf{U}_3^{k+1}- \mathbf{U}_3^{k}\|_F^2\right)
\\
\leq \ &
F\left(\mathcal{G}^{k}, \mathbf{A}_1^{k}, \mathbf{A}_2^{k}, \mathbf{A}_3^{k},\mathbf{U}_1^{k},
\mathbf{U}_2^{k},\mathbf{U}_3^{k}\right).
\end{aligned}
\end{equation}
This completes the proof.

\section*{Appendix E. Proof of Lemma \ref{Relacond}}
Firstly, we provide the first-order optimality conditions for the PALM algorithm in (\ref{palm}):
\begin{equation}\label{opti}
\left\{
\begin{aligned}
0& \in  \partial(\beta\|\mathcal{G}^{k+1}\|_1 +\delta_{\mathcal{D}}(\mathcal{G}^{k+1}))+
\nabla_{\mathcal{G}}\Psi(\mathcal{G}^{k}, \mathbf{A}_1^{k}, \mathbf{A}_2^{k}, \mathbf{A}_3^{k},\mathbf{U}_1^{k},
\mathbf{U}_2^{k},\mathbf{U}_3^{k})+\rho_1(\mathcal{G}^{k+1}- \mathcal{G}^k),\\
0& \in  \partial \delta_{\mathfrak{B}_1}(\mathbf{A}_1^{k+1}) +  \nabla_{\mathbf{A}_1} \Psi(\mathcal{G}^{k+1}, \mathbf{A}_1^{k}, \mathbf{A}_2^{k}, \mathbf{A}_3^{k},\mathbf{U}_1^{k},
\mathbf{U}_2^{k},\mathbf{U}_3^{k})+
\rho_2(\mathbf{A}_1^{k+1}-\mathbf{A}_1^{k}), \\
0& \in  \partial \delta_{\mathfrak{B}_2}(\mathbf{A}_2^{k+1}) +  \nabla_{\mathbf{A}_2} \Psi(\mathcal{G}^{k+1}, \mathbf{A}_1^{k+1}, \mathbf{A}_2^{k}, \mathbf{A}_3^{k},\mathbf{U}_1^{k},
\mathbf{U}_2^{k},\mathbf{U}_3^{k})+
\rho_3(\mathbf{A}_2^{k+1}-\mathbf{A}_2^{k}), \\
0& \in  \partial \delta_{\mathfrak{B}_3}(\mathbf{A}_3^{k+1}) + \nabla_{\mathbf{A}_3} \Psi(\mathcal{G}^{k+1}, \mathbf{A}_1^{k+1}, \mathbf{A}_2^{k+1}, \mathbf{A}_3^{k},\mathbf{U}_1^{k},
\mathbf{U}_2^{k},\mathbf{U}_3^{k})+
\rho_4(\mathbf{A}_3^{k+1}-\mathbf{A}_3^{k}), \\
0& =  2 \alpha_1\mathbf{L}_1 \mathbf{U}_1^{k+1} + \nabla_{\mathbf{U}_1} \Psi(\mathcal{G}^{k+1}, \mathbf{A}_1^{k+1}, \mathbf{A}_2^{k+1}, \mathbf{A}_3^{k+1},\mathbf{U}_1^{k},
\mathbf{U}_2^{k},\mathbf{U}_3^{k})
+\rho_5(\mathbf{U}_1^{k+1}-\mathbf{U}_1^{k}), \\
0& =  2 \alpha_2\mathbf{L}_2 \mathbf{U}_2^{k+1} + \nabla_{\mathbf{U}_2} \Psi(\mathcal{G}^{k+1}, \mathbf{A}_1^{k+1}, \mathbf{A}_2^{k+1}, \mathbf{A}_3^{k+1},\mathbf{U}_1^{k+1},
\mathbf{U}_2^{k},\mathbf{U}_3^{k})+
\rho_6(\mathbf{U}_2^{k+1}-\mathbf{U}_2^{k}), \\
0& =  2 \alpha_3\mathbf{L}_3 \mathbf{U}_3^{k+1} + \nabla_{\mathbf{U}_3}\Psi(\mathcal{G}^{k+1}, \mathbf{A}_1^{k+1}, \mathbf{A}_2^{k+1}, \mathbf{A}_3^{k+1},\mathbf{U}_1^{k+1},
\mathbf{U}_2^{k+1},\mathbf{U}_3^{k})+
\rho_7(\mathbf{U}_3^{k+1}-\mathbf{U}_3^{k}).
\end{aligned}
\right.
\end{equation}
Denote
\begin{equation}\label{gradient}
\begin{aligned}
N_1= \ & \nabla_\mathcal{G} \Psi(\mathcal{G}^{k+1}, \mathbf{A}_1^{k+1},\mathbf{A}_2^{k+1},
\mathbf{A}_3^{k+1},\mathbf{U}_1^{k+1},
\mathbf{U}_2^{k+1},\mathbf{U}_3^{k+1})\\
&-\nabla_\mathcal{G} \Psi(\mathcal{G}^{k}, \mathbf{A}_1^{k},\mathbf{A}_2^{k},
\mathbf{A}_3^{k},\mathbf{U}_1^{k},
\mathbf{U}_2^{k},\mathbf{U}_3^{k})
-\rho_1(\mathcal{G}^{k+1}-\mathcal{G}^k), \\
 N_2= \ &\nabla_{\mathbf{A}_1}\Psi(\mathcal{G}^{k+1}, \mathbf{A}_1^{k+1},\mathbf{A}_2^{k+1},
\mathbf{A}_3^{k+1},\mathbf{U}_1^{k+1},
\mathbf{U}_2^{k+1},\mathbf{U}_3^{k+1})\\
&-\nabla_{\mathbf{A}_1} \Psi(\mathcal{G}^{k+1}, \mathbf{A}_1^{k},\mathbf{A}_2^{k},
\mathbf{A}_3^{k},\mathbf{U}_1^{k},
\mathbf{U}_2^{k},\mathbf{U}_3^{k})
-\rho_2(\mathbf{A}_1^{k+1}-\mathbf{A}_1^k), \\
 N_3= \ &\nabla_{\mathbf{A}_2}\Psi(\mathcal{G}^{k+1}, \mathbf{A}_1^{k+1},\mathbf{A}_2^{k+1},
\mathbf{A}_3^{k+1},\mathbf{U}_1^{k+1},
\mathbf{U}_2^{k+1},\mathbf{U}_3^{k+1})\\
&-\nabla_{\mathbf{A}_2} \Psi(\mathcal{G}^{k+1}, \mathbf{A}_1^{k+1},\mathbf{A}_2^{k},
\mathbf{A}_3^{k},\mathbf{U}_1^{k},
\mathbf{U}_2^{k},\mathbf{U}_3^{k})
-\rho_3(\mathbf{A}_2^{k+1}-\mathbf{A}_2^k), \\
 N_4= \ &\nabla_{\mathbf{A}_3}\Psi(\mathcal{G}^{k+1}, \mathbf{A}_1^{k+1},\mathbf{A}_2^{k+1},
\mathbf{A}_3^{k+1},\mathbf{U}_1^{k+1},
\mathbf{U}_2^{k+1},\mathbf{U}_3^{k+1})\\
&-\nabla_{\mathbf{A}_3} \Psi(\mathcal{G}^{k+1}, \mathbf{A}_1^{k+1},\mathbf{A}_2^{k+1},
\mathbf{A}_3^{k},\mathbf{U}_1^{k},
\mathbf{U}_2^{k},\mathbf{U}_3^{k})
-\rho_4(\mathbf{A}_3^{k+1}-\mathbf{A}_3^k), \\
 N_5= \ &\nabla_{\mathbf{U}_1}\Psi(\mathcal{G}^{k+1}, \mathbf{A}_1^{k+1},\mathbf{A}_2^{k+1},
\mathbf{A}_3^{k+1},\mathbf{U}_1^{k+1},
\mathbf{U}_2^{k+1},\mathbf{U}_3^{k+1})\\
&-\nabla_{\mathbf{U}_1} \Psi(\mathcal{G}^{k+1}, \mathbf{A}_1^{k+1},\mathbf{A}_2^{k+1},
\mathbf{A}_3^{k+1},\mathbf{U}_1^{k},
\mathbf{U}_2^{k},\mathbf{U}_3^{k})
-\rho_5(\mathbf{U}_1^{k+1}-\mathbf{U}_1^k), \\
 N_6= \ &\nabla_{\mathbf{U}_2}\Psi(\mathcal{G}^{k+1}, \mathbf{A}_1^{k+1},\mathbf{A}_2^{k+1},
\mathbf{A}_3^{k+1},\mathbf{U}_1^{k+1},
\mathbf{U}_2^{k+1},\mathbf{U}_3^{k+1})\\
&-\nabla_{\mathbf{U}_2} \Psi(\mathcal{G}^{k+1}, \mathbf{A}_1^{k+1},\mathbf{A}_2^{k+1},
\mathbf{A}_3^{k+1},\mathbf{U}_1^{k+1},
\mathbf{U}_2^{k},\mathbf{U}_3^{k})
-\rho_6(\mathbf{U}_2^{k+1}-\mathbf{U}_2^k), \\
 N_7= \ &\nabla_{\mathbf{U}_3}\Psi(\mathcal{G}^{k+1}, \mathbf{A}_1^{k+1},\mathbf{A}_2^{k+1},
\mathbf{A}_3^{k+1},\mathbf{U}_1^{k+1},
\mathbf{U}_2^{k+1},\mathbf{U}_3^{k+1})\\
&-\nabla_{\mathbf{U}_3} \Psi(\mathcal{G}^{k+1}, \mathbf{A}_1^{k+1},\mathbf{A}_2^{k+1},
\mathbf{A}_3^{k+1},\mathbf{U}_1^{k+1},
\mathbf{U}_2^{k+1},\mathbf{U}_3^{k})
-\rho_7(\mathbf{U}_3^{k+1}-\mathbf{U}_3^k).
\end{aligned}
\end{equation}
Then we get
$$
\begin{aligned}
& N_1 \in \partial_{\mathcal{G}}F(\mathcal{G}^{k+1},\mathbf{A}_1^{k+1},\mathbf{A}_2^{k+1},\mathbf{A}_3^{k+1},
\mathbf{U}_1^{k+1},
\mathbf{U}_2^{k+1},\mathbf{U}_3^{k+1}), \\
& N_2 \in \partial_{\mathbf{A}_1}F(\mathcal{G}^{k+1},\mathbf{A}_1^{k+1},\mathbf{A}_2^{k+1},\mathbf{A}_3^{k+1},
\mathbf{U}_1^{k+1},
\mathbf{U}_2^{k+1},\mathbf{U}_3^{k+1}), \\
& N_3 \in \partial_{\mathbf{A}_2}F(\mathcal{G}^{k+1},\mathbf{A}_1^{k+1},\mathbf{A}_2^{k+1},\mathbf{A}_3^{k+1},
\mathbf{U}_1^{k+1},
\mathbf{U}_2^{k+1},\mathbf{U}_3^{k+1}), \\
& N_4 \in \partial_{\mathbf{A}_3}F(\mathcal{G}^{k+1},\mathbf{A}_1^{k+1},\mathbf{A}_2^{k+1},\mathbf{A}_3^{k+1},
\mathbf{U}_1^{k+1},
\mathbf{U}_2^{k+1},\mathbf{U}_3^{k+1}), \\
& N_5 \in \partial_{\mathbf{U}_1}F(\mathcal{G}^{k+1},\mathbf{A}_1^{k+1},\mathbf{A}_2^{k+1},\mathbf{A}_3^{k+1},
\mathbf{U}_1^{k+1},
\mathbf{U}_2^{k+1},\mathbf{U}_3^{k+1}), \\
& N_6 \in \partial_{\mathbf{U}_2}F(\mathcal{G}^{k+1},\mathbf{A}_1^{k+1},\mathbf{A}_2^{k+1},\mathbf{A}_3^{k+1},
\mathbf{U}_1^{k+1},
\mathbf{U}_2^{k+1},\mathbf{U}_3^{k+1}), \\
& N_7 \in \partial_{\mathbf{U}_3}F(\mathcal{G}^{k+1},\mathbf{A}_1^{k+1},\mathbf{A}_2^{k+1},\mathbf{A}_3^{k+1},
\mathbf{U}_1^{k+1},
\mathbf{U}_2^{k+1},\mathbf{U}_3^{k+1}).
\end{aligned}
$$
By \cite[Proposition 2.1]{Attouch2010}, we obtain that
$$
N^{k+1}:=\left(N_1, N_2, N_3, N_4, N_5, N_6, N_7\right) \in \partial F (\mathcal{G}^{k+1},\mathbf{A}_1^{k+1},\mathbf{A}_2^{k+1},\mathbf{A}_3^{k+1},
\mathbf{U}_1^{k+1},
\mathbf{U}_2^{k+1},\mathbf{U}_3^{k+1}).
$$

Note that
\begin{equation}\label{N1_12}
\begin{aligned}
&\|N_1\|_F\\
=&\left\|\nabla_\mathcal{G} \Psi(\mathcal{G}^{k+1}, \mathbf{A}_1^{k+1},\mathbf{A}_2^{k+1},
\mathbf{A}_3^{k+1},\mathbf{U}_1^{k+1},
\mathbf{U}_2^{k+1},\mathbf{U}_3^{k+1})\right.\\
&\left.-\nabla_\mathcal{G} \Psi(\mathcal{G}^{k}, \mathbf{A}_1^{k},\mathbf{A}_2^{k},
\mathbf{A}_3^{k},\mathbf{U}_1^{k},
\mathbf{U}_2^{k},\mathbf{U}_3^{k})
-\rho_1(\mathcal{G}^{k+1}-\mathcal{G}^k)\right\|_F\\
\leq&\left\|\nabla_\mathcal{G} \Psi(\mathcal{G}^{k+1}, \mathbf{A}_1^{k+1},\mathbf{A}_2^{k+1},
\mathbf{A}_3^{k+1},\mathbf{U}_1^{k+1},
\mathbf{U}_2^{k+1},\mathbf{U}_3^{k+1})\right.\\
&\left.-\nabla_\mathcal{G} \Psi(\mathcal{G}^{k}, \mathbf{A}_1^{k},\mathbf{A}_2^{k},
\mathbf{A}_3^{k},\mathbf{U}_1^{k},
\mathbf{U}_2^{k},\mathbf{U}_3^{k})\right\|_F
+\rho_1\|\mathcal{G}^{k+1}-\mathcal{G}^k\|_F\\
=&\left\|\nabla_\mathcal{G} Q(\llbracket \mathcal{G}^{k+1}; \mathbf{A}_1^{k+1}, \mathbf{A}_2^{k+1}, \mathbf{A}_3^{k+1}\rrbracket)
-\nabla_\mathcal{G} Q(\llbracket \mathcal{G}^{k}; \mathbf{A}_1^{k}, \mathbf{A}_2^{k}, \mathbf{A}_3^{k}\rrbracket)\right\|_F+\rho_1\|\mathcal{G}^{k+1}-\mathcal{G}^k\|_F\\
=&\left\|\left(\nabla_\mathcal{G} Q(\llbracket \mathcal{G}^{k+1}; \mathbf{A}_1^{k+1}, \mathbf{A}_2^{k+1}, \mathbf{A}_3^{k+1}\rrbracket)\right)_{(1)}
-\left(\nabla_\mathcal{G} Q(\llbracket \mathcal{G}^{k}; \mathbf{A}_1^{k}, \mathbf{A}_2^{k}, \mathbf{A}_3^{k}\rrbracket)\right)_{(1)}\right\|_F+\rho_1\|\mathcal{G}^{k+1}-\mathcal{G}^k\|_F,
\end{aligned}
\end{equation}
where the second equality follows from (\ref{gradient_H}).
Notice that (\ref{DefQ}) can be rewritten as
\begin{equation}\label{QGtf}
\begin{aligned}
Q(\llbracket \mathcal{G}; \mathbf{A}_1, \mathbf{A}_2, \mathbf{A}_3\rrbracket)
&=\frac{1}{2T} \sum_{t=1}^{T}\left\|\mathbf{y}_t-\mathbf{A}_1\mathcal{G}_{(1)}
\left(\mathbf{A}_3 \otimes \mathbf{A}_2\right)^{H}\mathbf{x}_{t}\right\|_2^2.
\end{aligned}
\end{equation}
Combining (\ref{QGtf}) with Lemma \ref{gradient_ele}, simple calculation gives
\begin{equation}\label{N1_7}
\begin{aligned}
&\left\|\left(\nabla_\mathcal{G} Q(\llbracket \mathcal{G}^{k+1}; \mathbf{A}_1^{k+1}, \mathbf{A}_2^{k+1}, \mathbf{A}_3^{k+1}\rrbracket)\right)_{(1)}
-\left(\nabla_\mathcal{G} Q(\llbracket \mathcal{G}^{k}; \mathbf{A}_1^{k}, \mathbf{A}_2^{k}, \mathbf{A}_3^{k}\rrbracket)\right)_{(1)}\right\|_F\\
=& \left\| \frac{1}{T} \sum_{t=1}^T (\mathbf{A}_1^{k+1})^{H}
\left(\mathbf{A}_1^{k+1}\mathcal{G}_{(1)}^{k+1}
(\mathbf{A}_3^{k+1} \otimes \mathbf{A}_2^{k+1})^{H}
\mathbf{x}_t-\mathbf{y}_t\right)\left((\mathbf{A}_3^{k+1} \otimes \mathbf{A}_2^{k+1})^{H}
\mathbf{x}_t\right)^{H} \right.\\
&~~~-\left. \frac{1}{T} \sum_{t=1}^T (\mathbf{A}_1^{k})^{H}
\left(\mathbf{A}_1^{k}\mathcal{G}_{(1)}^{k}
(\mathbf{A}_3^{k} \otimes \mathbf{A}_2^{k})^{H} \mathbf{x}_t-\mathbf{y}_t\right)
\left((\mathbf{A}_3^{k} \otimes \mathbf{A}_2^{k})^{H} \mathbf{x}_t\right)^{H}
\right\|_F\\
=& \left\| \frac{1}{T} \sum_{t=1}^T \left(\mathcal{G}_{(1)}^{k+1}
(\mathbf{A}_3^{k+1} \otimes \mathbf{A}_2^{k+1})^{H}
\mathbf{x}_t-(\mathbf{A}_1^{k+1})^{H}\mathbf{y}_t\right)\left((\mathbf{A}_3^{k+1} \otimes
\mathbf{A}_2^{k+1})^{H} \mathbf{x}_t\right)^{H} \right.\\
-&\left. \frac{1}{T} \sum_{t=1}^T \left(\mathcal{G}_{(1)}^{k}
(\mathbf{A}_3^{k} \otimes \mathbf{A}_2^{k})^{H}
\mathbf{x}_t-(\mathbf{A}_1^{k})^{H}\mathbf{y}_t\right)\left((\mathbf{A}_3^{k} \otimes \mathbf{A}_2^{k})^{H} \mathbf{x}_t\right)^{H}\right\|_F\\
=& \left\| \frac{1}{T} \sum_{t=1}^T \mathcal{G}_{(1)}^{k+1}
(\mathbf{A}_3^{k+1} \otimes \mathbf{A}_2^{k+1})^{H}
\mathbf{x}_t\mathbf{x}_t^{H}(\mathbf{A}_3^{k+1} \otimes \mathbf{A}_2^{k+1})
\right.\\
&~~~-\left.\frac{1}{T} \sum_{t=1}^T(\mathbf{A}_1^{k+1})^{H}\mathbf{y}_t
\mathbf{x}_t^{H}(\mathbf{A}_3^{k+1} \otimes \mathbf{A}_2^{k+1})
+\frac{1}{T} \sum_{t=1}^T(\mathbf{A}_1^{k})^{H}\mathbf{y}_t\mathbf{x}_t^{H}
(\mathbf{A}_3^{k} \otimes \mathbf{A}_2^{k})
\right.\\
&~~~-\left.\frac{1}{T} \sum_{t=1}^T \mathcal{G}_{(1)}^{k}
(\mathbf{A}_3^{k} \otimes \mathbf{A}_2^{k})^{H} \mathbf{x}_t\mathbf{x}_t^{H}
(\mathbf{A}_3^{k} \otimes \mathbf{A}_2^{k})
\right\|_F\\
\leq & \left\| \frac{1}{T} \sum_{t=1}^T \mathcal{G}_{(1)}^{k+1}
(\mathbf{A}_3^{k+1} \otimes \mathbf{A}_2^{k+1})^{H} \mathbf{x}_t\mathbf{x}_t^{H}
(\mathbf{A}_3^{k+1} \otimes \mathbf{A}_2^{k+1})
\right.\\
-&\left.
\frac{1}{T} \sum_{t=1}^T \mathcal{G}_{(1)}^{k}
(\mathbf{A}_3^{k} \otimes \mathbf{A}_2^{k})^{H} \mathbf{x}_t\mathbf{x}_t^{H}
(\mathbf{A}_3^{k} \otimes \mathbf{A}_2^{k})\right\|_F
\\
+&\left\|\frac{1}{T} \sum_{t=1}^T (\mathbf{A}_1^{k})^{H}\mathbf{y}_t\mathbf{x}_t^{H}
(\mathbf{A}_3^{k} \otimes \mathbf{A}_2^{k})-\frac{1}{T}
\sum_{t=1}^T (\mathbf{A}_1^{k+1})^{H}\mathbf{y}_t\mathbf{x}_t^{H}(\mathbf{A}_3^{k+1}
\otimes \mathbf{A}_2^{k+1})\right\|_F\\
\leq \  & \frac{1}{T} \sum_{t=1}^T \left\|  \mathcal{G}_{(1)}^{k+1}
(\mathbf{A}_3^{k+1} \otimes \mathbf{A}_2^{k+1})^{H} \mathbf{x}_t\mathbf{x}_t^{H}
(\mathbf{A}_3^{k+1} \otimes \mathbf{A}_2^{k+1})
- \mathcal{G}_{(1)}^{k}
(\mathbf{A}_3^{k} \otimes \mathbf{A}_2^{k})^{H} \mathbf{x}_t\mathbf{x}_t^{H}
(\mathbf{A}_3^{k} \otimes \mathbf{A}_2^{k})\right\|_F
\\
&~~~+\frac{1}{T} \sum_{t=1}^T \left\|(\mathbf{A}_1^{k})^{H}\mathbf{y}_t\mathbf{x}_t^{H}
(\mathbf{A}_3^{k} \otimes \mathbf{A}_2^{k})-(\mathbf{A}_1^{k+1})^{H}\mathbf{y}_t\mathbf{x}_t^{H}(\mathbf{A}_3^{k+1}
\otimes \mathbf{A}_2^{k+1})\right\|_F,
\end{aligned}
\end{equation}
where the second equality results from
$(\mathbf{A}_1^{k+1})^{H}\mathbf{A}_1^{k+1}=\mathbf{I}_{r_i}$.
Next, we start by analyzing the second term on the right side of (\ref{N1_7}).
Note that
\begin{equation}\label{N1_9}
\begin{aligned}
&\left\|(\mathbf{A}_1^{k})^{H}\mathbf{y}_t\mathbf{x}_t^{H}
(\mathbf{A}_3^{k} \otimes \mathbf{A}_2^{k})-
(\mathbf{A}_1^{k+1})^{H}\mathbf{y}_t\mathbf{x}_t^{H}(\mathbf{A}_3^{k+1}
\otimes \mathbf{A}_2^{k+1})\right\|_F\\
\leq \ & \left\|
(\mathbf{A}_1^{k})^{H}\mathbf{y}_t\mathbf{x}_t^{H}(\mathbf{A}_3^{k}
\otimes \mathbf{A}_2^{k})
-(\mathbf{A}_1^{k+1})^{H}\mathbf{y}_t\mathbf{x}_t^{H}(\mathbf{A}_3^{k}
\otimes \mathbf{A}_2^{k})
\right.\\
&~~+\left.
(\mathbf{A}_1^{k+1})^{H}\mathbf{y}_t\mathbf{x}_t^{H}(\mathbf{A}_3^{k}
\otimes \mathbf{A}_2^{k})-
(\mathbf{A}_1^{k+1})^{H}\mathbf{y}_t\mathbf{x}_t^{H}(\mathbf{A}_3^{k}
\otimes \mathbf{A}_2^{k+1})
\right.\\
&~~+\left.(\mathbf{A}_1^{k+1})^{H}\mathbf{y}_t\mathbf{x}_t^{H}(\mathbf{A}_3^{k}
\otimes \mathbf{A}_2^{k+1})-
(\mathbf{A}_1^{k+1})^{H}\mathbf{y}_t\mathbf{x}_t^{H}(\mathbf{A}_3^{k+1}
\otimes \mathbf{A}_2^{k+1})\right\|_F\\
\leq \  & \left\|
(\mathbf{A}_1^{k})^{H}\mathbf{y}_t\mathbf{x}_t^{H}(\mathbf{A}_3^{k}
\otimes \mathbf{A}_2^{k})
-(\mathbf{A}_1^{k+1})^{H}\mathbf{y}_t\mathbf{x}_t^{H}(\mathbf{A}_3^{k}
\otimes \mathbf{A}_2^{k})\right\|_F\\
&~~+\left\|
(\mathbf{A}_1^{k+1})^{H}\mathbf{y}_t\mathbf{x}_t^{H}(\mathbf{A}_3^{k}
\otimes \mathbf{A}_2^{k})-
(\mathbf{A}_1^{k+1})^{H}\mathbf{y}_t\mathbf{x}_t^{H}(\mathbf{A}_3^{k}
\otimes \mathbf{A}_2^{k+1})\right\|_F\\
&~~+\left\|(\mathbf{A}_1^{k+1})^{H}\mathbf{y}_t\mathbf{x}_t^{H}
(\mathbf{A}_3^{k} \otimes \mathbf{A}_2^{k+1})-
(\mathbf{A}_1^{k+1})^{H}\mathbf{y}_t\mathbf{x}_t^{H}(\mathbf{A}_3^{k+1}
 \otimes \mathbf{A}_2^{k+1})\right\|_F.
\end{aligned}
\end{equation}
By the Cauchy-Schwarz inequality, one can get that
\begin{equation}\label{N1_13}
\begin{aligned}
&\left\|
(\mathbf{A}_1^{k})^{H}\mathbf{y}_t\mathbf{x}_t^{H}(\mathbf{A}_3^{k} \otimes \mathbf{A}_2^{k})
-(\mathbf{A}_1^{k+1})^{H}\mathbf{y}_t\mathbf{x}_t^{H}(\mathbf{A}_3^{k} \otimes \mathbf{A}_2^{k})\right\|_F\\
= \ & \left\|
(\mathbf{A}_1^{k}-\mathbf{A}_1^{k+1})^{H}\mathbf{y}_t\mathbf{x}_t^{H}(\mathbf{A}_3^{k} \otimes \mathbf{A}_2^{k})\right\|_F\\
\leq \ & \left\|
(\mathbf{A}_1^{k}-\mathbf{A}_1^{k+1})^{H}\right\|_F
\left\|\mathbf{y}_t\mathbf{x}_t^{H}(\mathbf{A}_3^{k} \otimes \mathbf{A}_2^{k})\right\|_F\\
\leq \ & \|
\mathbf{A}_1^{k}-\mathbf{A}_1^{k+1}\|_F
\|\mathbf{y}_t\mathbf{x}_t^{H}\|_F
\|\mathbf{A}_3^{k} \otimes \mathbf{A}_2^{k}\|\\
\leq \ & \|
\mathbf{A}_1^{k}-\mathbf{A}_1^{k+1}\|_F
\|\mathbf{y}_t\mathbf{x}_t^{H}\|_F,
\end{aligned}
\end{equation}
where the second inequality holds by \cite[Lemma 2.1]{Chen2021}
and the last inequality holds by Proposition \ref{kron_untari}.
One can similarly obtain that
$$
\begin{aligned}
&\left\|
(\mathbf{A}_1^{k+1})^{H}\mathbf{y}_t\mathbf{x}_t^{H}(\mathbf{A}_3^{k} \otimes
\mathbf{A}_2^{k})-
(\mathbf{A}_1^{k+1})^{H}\mathbf{y}_t\mathbf{x}_t^{H}(\mathbf{A}_3^{k} \otimes
\mathbf{A}_2^{k+1})\right\|_F
\leq \|\mathbf{y}_t\mathbf{x}_t^{H}\|_F\|
\mathbf{A}_2^{k}-\mathbf{A}_2^{k+1}\|_F,\\
&\left\|(\mathbf{A}_1^{k+1})^{H}\mathbf{y}_t\mathbf{x}_t^{H}
(\mathbf{A}_3^{k} \otimes \mathbf{A}_2^{k+1})-
(\mathbf{A}_1^{k+1})^{H}\mathbf{y}_t\mathbf{x}_t^{H}(\mathbf{A}_3^{k+1}
\otimes \mathbf{A}_2^{k+1})\right\|_F
\leq \|\mathbf{y}_t\mathbf{x}_t^{H}\|_F\|
\mathbf{A}_3^{k}-\mathbf{A}_3^{k+1}\|_F.
\end{aligned}
$$
Taking this together with (\ref{N1_9}) and (\ref{N1_13}), we can deduce
\begin{equation}\label{N1_14}
\begin{aligned}
&\left\|(\mathbf{A}_1^{k})^{H}\mathbf{y}_t\mathbf{x}_t^{H}
(\mathbf{A}_3^{k} \otimes \mathbf{A}_2^{k})-
(\mathbf{A}_1^{k+1})^{H}\mathbf{y}_t\mathbf{x}_t^{H}(\mathbf{A}_3^{k+1} \otimes \mathbf{A}_2^{k+1})\right\|_F\\
\leq&
\|\mathbf{y}_t\mathbf{x}_t^{H}\|_F
\left(\|\mathbf{A}_1^{k}-\mathbf{A}_1^{k+1}\|_F+
\|
\mathbf{A}_2^{k}-\mathbf{A}_2^{k+1}\|_F
+\|
\mathbf{A}_3^{k}-\mathbf{A}_3^{k+1}\|_F\right).
\end{aligned}
\end{equation}
Summing up the above inequality over $t=1,\ldots,T$,
we get that
\begin{equation}\label{N1_8}
\begin{aligned}
 \ &\frac{1}{T} \sum_{t=1}^T \left\|(\mathbf{A}_1^{k})^{H}\mathbf{y}_t\mathbf{x}_t^{H}
(\mathbf{A}_3^{k} \otimes \mathbf{A}_2^{k})-
(\mathbf{A}_1^{k+1})^{H}\mathbf{y}_t\mathbf{x}_t^{H}(\mathbf{A}_3^{k+1} \otimes \mathbf{A}_2^{k+1})\right\|_F\\
\leq \ &
\frac{1}{T} \sum_{t=1}^T \|\mathbf{y}_t\mathbf{x}_t^{H}\|_F
\left(\|\mathbf{A}_1^{k}-\mathbf{A}_1^{k+1}\|_F+
\|
\mathbf{A}_2^{k}-\mathbf{A}_2^{k+1}\|_F
+\|
\mathbf{A}_3^{k}-\mathbf{A}_3^{k+1}\|_F\right)\\
=\ &
c_2
\left(\|\mathbf{A}_1^{k}-\mathbf{A}_1^{k+1}\|_F+
\|
\mathbf{A}_2^{k}-\mathbf{A}_2^{k+1}\|_F
+\|
\mathbf{A}_3^{k}-\mathbf{A}_3^{k+1}\|_F\right),
\end{aligned}
\end{equation}
where
\begin{equation}\label{DegiC2}	
c_2:=\frac{1}{T} \sum_{t=1}^T \|\mathbf{y}_t\mathbf{x}_t^{H}\|_F.
\end{equation}

Now we analyze the first term on the right side of (\ref{N1_7}).
It is easily seen that
\begin{equation}\label{N1_15}
\begin{aligned}
\ & \left\| \mathcal{G}_{(1)}^{k+1}
(\mathbf{A}_3^{k+1} \otimes \mathbf{A}_2^{k+1})^{H}
\mathbf{x}_t\mathbf{x}_t^{H}(\mathbf{A}_3^{k+1} \otimes \mathbf{A}_2^{k+1})
-
\mathcal{G}_{(1)}^{k}
(\mathbf{A}_3^{k} \otimes \mathbf{A}_2^{k})^{H}
\mathbf{x}_t\mathbf{x}_t^{H}(\mathbf{A}_3^{k} \otimes \mathbf{A}_2^{k})\right\|_F\\
= \ &\left\|\mathcal{G}_{(1)}^{k+1}
(\mathbf{A}_3^{k+1} \otimes \mathbf{A}_2^{k+1})^{H}
\mathbf{x}_t\mathbf{x}_t^{H}(\mathbf{A}_3^{k+1} \otimes \mathbf{A}_2^{k+1})
-\mathcal{G}_{(1)}^{k}
(\mathbf{A}_3^{k+1} \otimes \mathbf{A}_2^{k+1})^{H}
\mathbf{x}_t\mathbf{x}_t^{H}(\mathbf{A}_3^{k+1} \otimes \mathbf{A}_2^{k+1})
\right.\\
 \ & +\left.
\mathcal{G}_{(1)}^{k}
(\mathbf{A}_3^{k+1} \otimes \mathbf{A}_2^{k+1})^{H}
\mathbf{x}_t\mathbf{x}_t^{H}(\mathbf{A}_3^{k+1} \otimes \mathbf{A}_2^{k+1})
-\mathcal{G}_{(1)}^{k}
(\mathbf{A}_3^{k} \otimes \mathbf{A}_2^{k})^{H}
\mathbf{x}_t\mathbf{x}_t^{H}(\mathbf{A}_3^{k} \otimes \mathbf{A}_2^{k})\right\|_F\\
\leq \  & \left\| \mathcal{G}_{(1)}^{k+1}
(\mathbf{A}_3^{k+1} \otimes \mathbf{A}_2^{k+1})^{H}
\mathbf{x}_t\mathbf{x}_t^{H}(\mathbf{A}_3^{k+1} \otimes \mathbf{A}_2^{k+1})
\right. -\left.
\mathcal{G}_{(1)}^{k}
(\mathbf{A}_3^{k+1} \otimes \mathbf{A}_2^{k+1})^{H}
\mathbf{x}_t\mathbf{x}_t^{H}(\mathbf{A}_3^{k+1} \otimes \mathbf{A}_2^{k+1})\right\|_F\\
 \ & +\left\| \mathcal{G}_{(1)}^{k}
(\mathbf{A}_3^{k+1} \otimes \mathbf{A}_2^{k+1})^{H}
\mathbf{x}_t\mathbf{x}_t^{H}(\mathbf{A}_3^{k+1} \otimes \mathbf{A}_2^{k+1})
-\mathcal{G}_{(1)}^{k}
(\mathbf{A}_3^{k} \otimes \mathbf{A}_2^{k})^{H} \mathbf{x}_t\mathbf{x}_t^{H}(\mathbf{A}_3^{k} \otimes \mathbf{A}_2^{k})\right\|_F.
\end{aligned}
\end{equation}
Applying the Cauchy-Schwarz inequality to (\ref{N1_15}) yields
\begin{equation}\label{N1_10}
\begin{aligned}
 \ &  \left\| \mathcal{G}_{(1)}^{k+1}
(\mathbf{A}_3^{k+1} \otimes \mathbf{A}_2^{k+1})^{H} \mathbf{x}_t\mathbf{x}_t^{H}(\mathbf{A}_3^{k+1} \otimes \mathbf{A}_2^{k+1})
-
\mathcal{G}_{(1)}^{k}
(\mathbf{A}_3^{k+1} \otimes \mathbf{A}_2^{k+1})^{H} \mathbf{x}_t\mathbf{x}_t^{H}(\mathbf{A}_3^{k+1} \otimes \mathbf{A}_2^{k+1})\right\|_F\\
= \ &\left\|(\mathcal{G}_{(1)}^{k+1}-\mathcal{G}_{(1)}^{k})
(\mathbf{A}_3^{k+1} \otimes \mathbf{A}_2^{k+1})^{H} \mathbf{x}_t\mathbf{x}_t^{H}(\mathbf{A}_3^{k+1} \otimes \mathbf{A}_2^{k+1})\right\|_F\\
\leq \ & \|\mathcal{G}_{(1)}^{k+1}-\mathcal{G}_{(1)}^{k}\|_F
\left\|(\mathbf{A}_3^{k+1} \otimes \mathbf{A}_2^{k+1})^{H} \mathbf{x}_t\mathbf{x}_t^{H}(\mathbf{A}_3^{k+1} \otimes \mathbf{A}_2^{k+1})\right\|_F\\
= \ &\|\mathcal{G}^{k+1}-\mathcal{G}^{k}\|_F
\left\|(\mathbf{A}_3^{k+1} \otimes \mathbf{A}_2^{k+1})^{H} \mathbf{x}_t\right\|_2^2\\
\leq \ &\|\mathcal{G}^{k+1}-\mathcal{G}^{k}\|_F
\left(\|\mathbf{A}_3^{k+1} \otimes \mathbf{A}_2^{k+1}\|
\left\|\mathbf{x}_t\right\|_2\right)^2\\
\leq \ &\|\mathcal{G}^{k+1}-\mathcal{G}^{k}\|_F \left\|\mathbf{x}_t\right\|_2^2,
\end{aligned}
\end{equation}
where the second inequality follows from \cite[Example 5.6.6]{Horn2013}.

On the other hand, for the second term of the right hand inequality in (\ref{N1_15}), we  obtain that
\begin{equation}\label{N1}
\begin{aligned}
&\left\| \mathcal{G}_{(1)}^{k}
(\mathbf{A}_3^{k+1} \otimes \mathbf{A}_2^{k+1})^{H}
\mathbf{x}_t\mathbf{x}_t^{H}(\mathbf{A}_3^{k+1} \otimes \mathbf{A}_2^{k+1})
-\mathcal{G}_{(1)}^{k}
(\mathbf{A}_3^{k} \otimes \mathbf{A}_2^{k})^{H}
\mathbf{x}_t\mathbf{x}_t^{H}(\mathbf{A}_3^{k} \otimes \mathbf{A}_2^{k})\right\|_F\\
\leq \ &\| \mathcal{G}_{(1)}^{k}\|_F
\left\|(\mathbf{A}_3^{k+1} \otimes \mathbf{A}_2^{k+1})^{H}
\mathbf{x}_t\mathbf{x}_t^{H}(\mathbf{A}_3^{k+1} \otimes \mathbf{A}_2^{k+1})
-
(\mathbf{A}_3^{k} \otimes \mathbf{A}_2^{k})^{H}
\mathbf{x}_t\mathbf{x}_t^{H}(\mathbf{A}_3^{k} \otimes \mathbf{A}_2^{k})\right\|_F\\
= \ &\| \mathcal{G}_{(1)}^{k}\|_F
\left\|(\mathbf{A}_3^{k+1} \otimes \mathbf{A}_2^{k+1})^{H}
\mathbf{x}_t\mathbf{x}_t^{H}(\mathbf{A}_3^{k+1} \otimes \mathbf{A}_2^{k+1})
-(\mathbf{A}_3^{k+1} \otimes \mathbf{A}_2^{k+1})^{H} \mathbf{x}_t\mathbf{x}_t^{H}
(\mathbf{A}_3^{k+1} \otimes \mathbf{A}_2^{k})
\right.\\
\ &+\left.
(\mathbf{A}_3^{k+1} \otimes \mathbf{A}_2^{k+1})^{H} \mathbf{x}_t\mathbf{x}_t^{H}
(\mathbf{A}_3^{k+1} \otimes \mathbf{A}_2^{k})
-(\mathbf{A}_3^{k+1} \otimes \mathbf{A}_2^{k+1})^{H} \mathbf{x}_t\mathbf{x}_t^{H}
(\mathbf{A}_3^{k} \otimes \mathbf{A}_2^{k})
\right.\\
 \ &+\left.
(\mathbf{A}_3^{k+1} \otimes \mathbf{A}_2^{k+1})^{H} \mathbf{x}_t\mathbf{x}_t^{H}
(\mathbf{A}_3^{k} \otimes \mathbf{A}_2^{k})
-(\mathbf{A}_3^{k+1} \otimes \mathbf{A}_2^{k})^{H} \mathbf{x}_t\mathbf{x}_t^{H}
(\mathbf{A}_3^{k} \otimes \mathbf{A}_2^{k})
\right.\\
 \ &+\left.
(\mathbf{A}_3^{k+1} \otimes \mathbf{A}_2^{k})^{H} \mathbf{x}_t\mathbf{x}_t^{H}
(\mathbf{A}_3^{k} \otimes \mathbf{A}_2^{k})
-(\mathbf{A}_3^{k} \otimes \mathbf{A}_2^{k})^{H} \mathbf{x}_t\mathbf{x}_t^{H}
(\mathbf{A}_3^{k} \otimes \mathbf{A}_2^{k})\right\|_F\\
\leq \ &\nu
\left\|(\mathbf{A}_3^{k+1} \otimes \mathbf{A}_2^{k+1})^{H}
\mathbf{x}_t\mathbf{x}_t^{H}(\mathbf{A}_3^{k+1} \otimes \mathbf{A}_2^{k+1})
-(\mathbf{A}_3^{k+1} \otimes \mathbf{A}_2^{k+1})^{H}
 \mathbf{x}_t\mathbf{x}_t^{H}
(\mathbf{A}_3^{k+1} \otimes \mathbf{A}_2^{k})\right\|_F\\
\ &+
\nu\left\|(\mathbf{A}_3^{k+1} \otimes \mathbf{A}_2^{k+1})^{H}
\mathbf{x}_t\mathbf{x}_t^{H}
(\mathbf{A}_3^{k+1} \otimes \mathbf{A}_2^{k})
-(\mathbf{A}_3^{k+1} \otimes \mathbf{A}_2^{k+1})^{H}
\mathbf{x}_t\mathbf{x}_t^{H}
(\mathbf{A}_3^{k} \otimes \mathbf{A}_2^{k})\right\|_F\\
\ &+
\nu\left\|
(\mathbf{A}_3^{k+1} \otimes \mathbf{A}_2^{k+1})^{H}
\mathbf{x}_t\mathbf{x}_t^{H}
(\mathbf{A}_3^{k} \otimes \mathbf{A}_2^{k})
-(\mathbf{A}_3^{k+1} \otimes \mathbf{A}_2^{k})^{H}
\mathbf{x}_t\mathbf{x}_t^{H}
(\mathbf{A}_3^{k} \otimes \mathbf{A}_2^{k})\right\|_F\\
\ &+
\nu\left\|(\mathbf{A}_3^{k+1} \otimes \mathbf{A}_2^{k})^{H}
 \mathbf{x}_t\mathbf{x}_t^{H}
(\mathbf{A}_3^{k} \otimes \mathbf{A}_2^{k})
-(\mathbf{A}_3^{k} \otimes \mathbf{A}_2^{k})^{H}
 \mathbf{x}_t\mathbf{x}_t^{H}
(\mathbf{A}_3^{k} \otimes \mathbf{A}_2^{k})\right\|_F,
\end{aligned}
\end{equation}
where the second inequality holds the triangle inequality and
$\| \mathcal{G}_{(1)}^{k}\|_F\leq \nu:=\sqrt{r_1r_2r_3}c$.
Similar to (\ref{lip2_4}), we have
\begin{equation}\label{N1_1}
\begin{aligned}
& \left\|(\mathbf{A}_3^{k+1} \otimes \mathbf{A}_2^{k+1})^{H}
 \mathbf{x}_t\mathbf{x}_t^{H}(\mathbf{A}_3^{k+1} \otimes \mathbf{A}_2^{k+1})
-(\mathbf{A}_3^{k+1} \otimes \mathbf{A}_2^{k+1})^{H}
 \mathbf{x}_t\mathbf{x}_t^{H}(\mathbf{A}_3^{k+1} \otimes
  \mathbf{A}_2^{k})\right\|_F\\
=&\left\|(\mathbf{A}_3^{k+1} \otimes \mathbf{A}_2^{k+1})^{H}
 \mathbf{x}_t \mathbf{x}_t^{H}\left(\mathbf{A}_3^{k+1} \otimes
(\mathbf{A}_2^{k+1}-\mathbf{A}_2^k)\right)\right\|_F\\
\leq&\left\|\mathbf{A}_3^{k+1} \otimes \mathbf{A}_2^{k+1}\right\|
\left\|\mathbf{x}_t \mathbf{x}_t^{H}\left(\mathbf{A}_3^{k+1}
\otimes(\mathbf{A}_2^{k+1}-\mathbf{A}_2^k)\right)\right\|_F \\
\leq &\left\| \mathbf{x}_t \mathbf{x}_t^{H}\left(\mathbf{A}_3^{k+1}
 \otimes(\mathbf{A}_2^{k+1}-\mathbf{A}_2^k) \right)\right\|_F\\
\leq& \left\| \mathbf{x}_t \mathbf{x}_t^{H}\right\|_F
\|\mathbf{A}_3^{k+1} \otimes(\mathbf{A}_2^{k+1}-\mathbf{A}_2^k)\|\\
=&\left\| \mathbf{x}_t \mathbf{x}_t^{H}\right\|_F
\|\mathbf{A}_3^{k+1}\| \|\mathbf{A}_2^{k+1}-\mathbf{A}_2^k\|\\
\leq& \left\| \mathbf{x}_t \mathbf{x}_t^{H}\right\|_F \|\mathbf{A}_2^{k+1}-\mathbf{A}_2^k\|_F.
\end{aligned}
\end{equation}
Similarly, we can get that
\begin{equation}\label{N1_2}
\begin{aligned}
&\left\|(\mathbf{A}_3^{k+1} \otimes \mathbf{A}_2^{k+1})^{H}
\mathbf{x}_t\mathbf{x}_t^{H}(\mathbf{A}_3^{k+1} \otimes \mathbf{A}_2^{k})
-(\mathbf{A}_3^{k+1} \otimes \mathbf{A}_2^{k+1})^{H}
\mathbf{x}_t\mathbf{x}_t^{H}(\mathbf{A}_3^{k} \otimes \mathbf{A}_2^{k})\right\|_F\\
\leq & \left\| \mathbf{x}_t \mathbf{x}_t^{H}\right\|_F
\|\mathbf{A}_3^{k+1}-\mathbf{A}_3^k\|_F,\\
&\left\|(\mathbf{A}_3^{k+1} \otimes \mathbf{A}_2^{k+1})^{H}
\mathbf{x}_t\mathbf{x}_t^{H}(\mathbf{A}_3^{k} \otimes \mathbf{A}_2^{k})
-(\mathbf{A}_3^{k+1} \otimes \mathbf{A}_2^{k})^{H}
\mathbf{x}_t\mathbf{x}_t^{H}(\mathbf{A}_3^{k} \otimes \mathbf{A}_2^{k})\right\|_F\\
\leq &\left\| \mathbf{x}_t \mathbf{x}_t^{H}\right\|_F
\|\mathbf{A}_2^{k+1}-\mathbf{A}_2^k\|_F,\\
&\left\|
(\mathbf{A}_3^{k+1} \otimes \mathbf{A}_2^{k})^{H}
\mathbf{x}_t\mathbf{x}_t^{H}(\mathbf{A}_3^{k} \otimes
\mathbf{A}_2^{k})
-(\mathbf{A}_3^{k} \otimes \mathbf{A}_2^{k})^{H}
\mathbf{x}_t\mathbf{x}_t^{H}(\mathbf{A}_3^{k} \otimes \mathbf{A}_2^{k})\right\|_F\\
\leq& \left\| \mathbf{x}_t \mathbf{x}_t^{H}\right\|_F
\|\mathbf{A}_3^{k+1}-\mathbf{A}_3^k\|_F.
\end{aligned}
\end{equation}
%
Substituting (\ref{N1_1}) and (\ref{N1_2}) into (\ref{N1}), we immediately establish
\begin{equation}\label{NN1_11}
\begin{aligned}
\ &\left\| \mathcal{G}_{(1)}^{k}
(\mathbf{A}_3^{k+1} \otimes \mathbf{A}_2^{k+1})^{H}
\mathbf{x}_t\mathbf{x}_t^{H}(\mathbf{A}_3^{k+1} \otimes \mathbf{A}_2^{k+1})
-\mathcal{G}_{(1)}^{k}
(\mathbf{A}_3^{k} \otimes \mathbf{A}_2^{k})^{H}
\mathbf{x}_t\mathbf{x}_t^{H}(\mathbf{A}_3^{k} \otimes \mathbf{A}_2^{k})\right\|_F\\
\leq \  & 2\nu \left\| \mathbf{x}_t \mathbf{x}_t^{H}\right\|_F
\left(\|\mathbf{A}_2^{k+1}-\mathbf{A}_2^k\|_F
+\|\mathbf{A}_3^{k+1}-\mathbf{A}_3^k\|_F\right).
\end{aligned}
\end{equation}
Plugging (\ref{NN1_11}) and (\ref{N1_10}) into (\ref{N1_15}), we have
\begin{equation}\label{N1_11}
\begin{aligned}
&\left\| \mathcal{G}_{(1)}^{k+1}
(\mathbf{A}_3^{k+1} \otimes \mathbf{A}_2^{k+1})^{H}
\mathbf{x}_t\mathbf{x}_t^{H}(\mathbf{A}_3^{k+1} \otimes \mathbf{A}_2^{k+1})
-
\mathcal{G}_{(1)}^{k}
(\mathbf{A}_3^{k} \otimes \mathbf{A}_2^{k})^{H}
\mathbf{x}_t\mathbf{x}_t^{H}(\mathbf{A}_3^{k} \otimes \mathbf{A}_2^{k})\right\|_F\\
\leq \ &\|\mathcal{G}^{k+1}-\mathcal{G}^{k}\|_F \left\|\mathbf{x}_t\right\|_2^2+2\nu \left\| \mathbf{x}_t \mathbf{x}_t^{H}\right\|_F
\left(\|\mathbf{A}_2^{k+1}-\mathbf{A}_2^k\|_F
+\|\mathbf{A}_3^{k+1}-\mathbf{A}_3^k\|_F\right).
\end{aligned}
\end{equation}
Summing up the above inequality over $t=1,\ldots,T$,
one can obtain that
\begin{equation}\label{N1_3}
\begin{aligned}
\ &\frac{1}{T} \sum_{t=1}^T \left\| \mathcal{G}_{(1)}^{k+1}
(\mathbf{A}_3^{k+1} \otimes \mathbf{A}_2^{k+1})^{H}
\mathbf{x}_t\mathbf{x}_t^{H}(\mathbf{A}_3^{k+1} \otimes \mathbf{A}_2^{k+1})
\right.\\
&~~~~~~~~~~~~-\left.
\mathcal{G}_{(1)}^{k}
(\mathbf{A}_3^{k} \otimes \mathbf{A}_2^{k})^{H}
\mathbf{x}_t\mathbf{x}_t^{H}(\mathbf{A}_3^{k} \otimes \mathbf{A}_2^{k})\right\|_F\\
\leq \  & \frac{1}{T} \sum_{t=1}^T \|\mathcal{G}^{k+1}-\mathcal{G}^{k}\|_F \|\mathbf{x}_t\|_2^2+2\nu \frac{1}{T} \sum_{t=1}^T  \| \mathbf{x}_t\|_2^2
\|\mathbf{A}_2^{k+1}-\mathbf{A}_2^k\|_F\\
&+2\nu \frac{1}{T} \sum_{t=1}^T  \| \mathbf{x}_t\|_2^2
\|\mathbf{A}_3^{k+1}-\mathbf{A}_3^k\|_F\\
= \ & c_1\|\mathcal{G}^{k+1}-\mathcal{G}^{k}\|_F+2\nu c_1\|\mathbf{A}_2^{k+1}-\mathbf{A}_2^k\|_F
+2\nu c_1\|\mathbf{A}_3^{k+1}-\mathbf{A}_3^k\|_F,
\end{aligned}
\end{equation}
where
\begin{equation}\label{DefC1}
c_1:= \frac{1}{T} \sum_{t=1}^T \| \mathbf{x}_t\|_2^2.
\end{equation}
Substituting (\ref{N1_3}) and (\ref{N1_8}) into (\ref{N1_7}) establishes
$$
\begin{aligned}
&\left\|\left(\nabla_\mathcal{G} Q(\llbracket \mathcal{G}^{k+1}; \mathbf{A}_1^{k+1}, \mathbf{A}_2^{k+1}, \mathbf{A}_3^{k+1}\rrbracket)\right)_{(1)}
-\left(\nabla_\mathcal{G} Q(\llbracket \mathcal{G}^{k}; \mathbf{A}_1^{k}, \mathbf{A}_2^{k}, \mathbf{A}_3^{k}\rrbracket)\right)_{(1)}\right\|_F\\
\leq \ &c_2
\left(\|\mathbf{A}_1^{k}-\mathbf{A}_1^{k+1}\|_F+
\|\mathbf{A}_2^{k}-\mathbf{A}_2^{k+1}\|_F
+\|\mathbf{A}_3^{k}-\mathbf{A}_3^{k+1}\|_F\right)+
c_1\|\mathcal{G}^{k+1}-\mathcal{G}^{k}\|_F\\
&+2\nu c_1\|\mathbf{A}_2^{k+1}-\mathbf{A}_2^k\|_F
+2\nu c_1\|\mathbf{A}_3^{k+1}-\mathbf{A}_3^k\|_F.
\end{aligned}
$$
This taken together with (\ref{N1_12}) leads to
\begin{equation}\label{N1_4}
\begin{aligned}
\|N_1\|_F&\leq c_2
\left(\|\mathbf{A}_1^{k}-\mathbf{A}_1^{k+1}\|_F+
\|\mathbf{A}_2^{k}-\mathbf{A}_2^{k+1}\|_F
+\|\mathbf{A}_3^{k}-\mathbf{A}_3^{k+1}\|_F\right)+
c_1\|\mathcal{G}^{k+1}-\mathcal{G}^{k}\|_F\\
&~~~~+2\nu c_1\|\mathbf{A}_2^{k+1}-\mathbf{A}_2^k\|_F
+2\nu c_1\|\mathbf{A}_3^{k+1}-\mathbf{A}_3^k\|_F
+\rho_1\|\mathcal{G}^{k+1}-\mathcal{G}^k\|_F\\
&=(\rho_1+c_1)\|\mathcal{G}^{k+1}-\mathcal{G}^k\|_F+
c_2\|\mathbf{A}_1^{k+1}-\mathbf{A}_1^{k}\|_F+
(c_2+2\nu c_1)\|\mathbf{A}_2^{k+1}-\mathbf{A}_2^{k}\|_F\\
&~~~~+(c_2+2\nu c_1)\|\mathbf{A}_3^{k+1}-\mathbf{A}_3^{k}\|_F.
\end{aligned}
\end{equation}

Similar to the above analysis about the upper bound of $N_1$, we can easily obtain
that
\begin{equation}\label{N2_7}
\begin{aligned}
\|N_2\|_F=&\left\|\nabla_{\mathbf{A}_1} \Psi\left(\mathcal{G}^{k+1}, \mathbf{A}_1^{k+1},\mathbf{A}_2^{k+1},
\mathbf{A}_3^{k+1},\mathbf{U}_1^{k+1},
\mathbf{U}_2^{k+1},\mathbf{U}_3^{k+1}\right)\right.\\
&\left.-\nabla_{\mathbf{A}_1} \Psi\left(\mathcal{G}^{k+1}, \mathbf{A}_1^{k},\mathbf{A}_2^{k},
\mathbf{A}_3^{k},\mathbf{U}_1^{k},
\mathbf{U}_2^{k},\mathbf{U}_3^{k}\right)
-\rho_2(\mathbf{A}_1^{k+1}-\mathbf{A}_1^{k})\right\|_F\\
\leq&\left\|\nabla_{\mathbf{A}_1} \Psi\left(\mathcal{G}^{k+1}, \mathbf{A}_1^{k+1},\mathbf{A}_2^{k+1},
\mathbf{A}_3^{k+1},\mathbf{U}_1^{k+1},
\mathbf{U}_2^{k+1},\mathbf{U}_3^{k+1}\right)\right.\\
&\left.-\nabla_{\mathbf{A}_1} \Psi\left(\mathcal{G}^{k+1}, \mathbf{A}_1^{k},\mathbf{A}_2^{k},
\mathbf{A}_3^{k},\mathbf{U}_1^{k},
\mathbf{U}_2^{k},\mathbf{U}_3^{k}\right)\right\|_F
+\rho_2\|\mathbf{A}_1^{k+1}-\mathbf{A}_1^{k}\|_F\\
=&\left\|\nabla_{\mathbf{A}_1} Q(\llbracket \mathcal{G}^{k+1}; \mathbf{A}_1^{k+1}, \mathbf{A}_2^{k+1}, \mathbf{A}_3^{k+1}\rrbracket)-\gamma_1(\mathbf{U}_1^{k+1}-\mathbf{A}_1^{k+1})
\right.\\
&\left.-
\nabla_{\mathbf{A}_1} Q(\llbracket \mathcal{G}^{k+1}; \mathbf{A}_1^{k}, \mathbf{A}_2^{k}, \mathbf{A}_3^{k}\rrbracket)+\gamma_1(\mathbf{U}_1^{k}-\mathbf{A}_1^{k})\right\|_F
+\rho_2\|\mathbf{A}_1^{k+1}-\mathbf{A}_1^{k}\|_F\\
\leq
&\left\|\nabla_{\mathbf{A}_1} Q(\llbracket \mathcal{G}^{k+1}; \mathbf{A}_1^{k+1}, \mathbf{A}_2^{k+1}, \mathbf{A}_3^{k+1}\rrbracket)
-\nabla_{\mathbf{A}_1} Q(\llbracket \mathcal{G}^{k+1}; \mathbf{A}_1^{k}, \mathbf{A}_2^{k}, \mathbf{A}_3^{k}\rrbracket)\right\|_F \\
& +\gamma_1\|\mathbf{U}_1^{k+1}-\mathbf{U}_1^k\|_F
+(\gamma_1+\rho_2)\|\mathbf{A}_1^{k+1}-\mathbf{A}_1^k\|_F.
\end{aligned}
\end{equation}
where the second equality follows from (\ref{gradient_H}) and the last inequality follows from the triangle inequality.
Taking (\ref{gradient_Q}) and (\ref{gradientQ_W1}) together yields
\begin{equation}\label{N2_1}
\begin{aligned}
&\left\|\nabla_{\mathbf{A}_1} Q(\llbracket \mathcal{G}^{k+1}; \mathbf{A}_1^{k+1}, \mathbf{A}_2^{k+1}, \mathbf{A}_3^{k+1}\rrbracket)
-\nabla_{\mathbf{A}_1} Q(\llbracket \mathcal{G}^{k+1}; \mathbf{A}_1^{k}, \mathbf{A}_2^{k}, \mathbf{A}_3^{k}\rrbracket)\right\|_F\\
=&\left\|(\nabla Q(\llbracket \mathcal{G}^{k+1}; \mathbf{A}_1^{k+1}, \mathbf{A}_2^{k+1}, \mathbf{A}_3^{k+1}\rrbracket))_{(1)}(\mathbf{A}_3^{k+1} \otimes \mathbf{A}_2^{k+1})(\mathcal{G}_{(1)}^{k+1})^{H}
\right.\\
&~~-\left.
(\nabla Q(\llbracket \mathcal{G}^{k+1}; \mathbf{A}_1^{k}, \mathbf{A}_2^{k}, \mathbf{A}_3^{k}\rrbracket))_{(1)}(\mathbf{A}_3^{k} \otimes \mathbf{A}_2^{k})(\mathcal{G}_{(1)}^{k+1})^{H}\right\|_F\\
=& \left\| \frac{1}{T} \sum_{t=1}^T \left(\mathbf{A}_1^{k+1}\mathcal{G}_{(1)}^{k+1}
(\mathbf{A}_3^{k+1} \otimes \mathbf{A}_2^{k+1})^{H} \mathbf{x}_t-\mathbf{y}_t\right)
\mathbf{x}_t^{H}
(\mathbf{A}_3^{k+1} \otimes \mathbf{A}_2^{k+1})(\mathcal{G}_{(1)}^{k+1})^{H}\right.\\
&~~-\left.
\frac{1}{T} \sum_{t=1}^T \left(\mathbf{A}_1^{k}\mathcal{G}_{(1)}^{k+1}
(\mathbf{A}_3^{k} \otimes \mathbf{A}_2^{k})^{H} \mathbf{x}_t-\mathbf{y}_t\right)
\mathbf{x}_t^{H}
(\mathbf{A}_3^{k} \otimes \mathbf{A}_2^{k})(\mathcal{G}_{(1)}^{k+1})^{H}\right\|_F\\
\leq \ & \frac{1}{T}
 \sum_{t=1}^T \left\|\mathbf{A}_1^{k+1}\mathcal{G}_{(1)}^{k+1}
(\mathbf{A}_3^{k+1} \otimes \mathbf{A}_2^{k+1})^{H} \mathbf{x}_t
\mathbf{x}_t^{H}
(\mathbf{A}_3^{k+1} \otimes \mathbf{A}_2^{k+1})(\mathcal{G}_{(1)}^{k+1})^{H}
\right.\\
&~~-\left.
\mathbf{y}_t\mathbf{x}_t^{H}(\mathbf{A}_3^{k+1} \otimes \mathbf{A}_2^{k+1})
(\mathcal{G}_{(1)}^{k+1})^{H}-\mathbf{A}_1^{k}\mathcal{G}_{(1)}^{k+1}
(\mathbf{A}_3^{k} \otimes \mathbf{A}_2^{k})^{H} \mathbf{x}_t
\mathbf{x}_t^{H}
(\mathbf{A}_3^{k} \otimes \mathbf{A}_2^{k})(\mathcal{G}_{(1)}^{k+1})^{H}
\right.\\
&~~+\left.
\mathbf{y}_t\mathbf{x}_t^{H}
(\mathbf{A}_3^{k} \otimes \mathbf{A}_2^{k})(\mathcal{G}_{(1)}^{k+1})^{H}\right\|_F\\
\leq \ & \frac{1}{T} \sum_{t=1}^T \left\| \mathbf{A}_1^{k+1}\mathcal{G}_{(1)}^{k+1}
(\mathbf{A}_3^{k+1} \otimes \mathbf{A}_2^{k+1})^{H} \mathbf{x}_t
\mathbf{x}_t^{H}
(\mathbf{A}_3^{k+1} \otimes \mathbf{A}_2^{k+1})(\mathcal{G}_{(1)}^{k+1})^{H}
\right.\\
&~~-\left.
\mathbf{A}_1^{k}\mathcal{G}_{(1)}^{k+1}
(\mathbf{A}_3^{k} \otimes \mathbf{A}_2^{k})^{H} \mathbf{x}_t
\mathbf{x}_t^{H}
(\mathbf{A}_3^{k} \otimes \mathbf{A}_2^{k})(\mathcal{G}_{(1)}^{k+1})^{H}\right\|_F\\
&~~+\frac{1}{T} \sum_{t=1}^T \left\|
\mathbf{y}_t\mathbf{x}_t^{H}
\left((\mathbf{A}_3^{k+1} \otimes \mathbf{A}_2^{k+1})(\mathcal{G}_{(1)}^{k+1})^{H}
-(\mathbf{A}_3^{k} \otimes \mathbf{A}_2^{k})(\mathcal{G}_{(1)}^{k+1})^{H}\right)\right\|_F,
\end{aligned}
\end{equation}
where the second equality follows from Lemma \ref{gradient_ele}.
First, let us consider the second term on the right-hand side of the above inequality.
Note that
\begin{equation}\label{N2_9}
\begin{aligned}
&\left\|
\mathbf{y}_t\mathbf{x}_t^{H}
\left((\mathbf{A}_3^{k+1} \otimes \mathbf{A}_2^{k+1})(\mathcal{G}_{(1)}^{k+1})^{H}
-(\mathbf{A}_3^{k} \otimes \mathbf{A}_2^{k})(\mathcal{G}_{(1)}^{k+1})^{H}\right)\right\|_F\\
=&\left\|
\mathbf{y}_t\mathbf{x}_t^{H}
\left((\mathbf{A}_3^{k+1} \otimes \mathbf{A}_2^{k+1})
-(\mathbf{A}_3^{k+1} \otimes \mathbf{A}_2^{k})
+(\mathbf{A}_3^{k+1} \otimes \mathbf{A}_2^{k})
-(\mathbf{A}_3^{k} \otimes\mathbf{A}_2^{k})\right)
\mathcal{G}_{(1)}^{k+1})^{H}\right\|_F\\
\leq&
\left\|\mathbf{y}_t\mathbf{x}_t^{H}
\left((\mathbf{A}_3^{k+1} \otimes \mathbf{A}_2^{k+1})
-(\mathbf{A}_3^{k+1} \otimes \mathbf{A}_2^{k}))\right)
\mathcal{G}_{(1)}^{k+1})^{H}\right\|_F\\
&+\left\|
\mathbf{y}_t\mathbf{x}_t^{H}
\left((\mathbf{A}_3^{k+1} \otimes \mathbf{A}_2^{k})
-(\mathbf{A}_3^{k} \otimes \mathbf{A}_2^{k}))\right)
\mathcal{G}_{(1)}^{k+1})^{H}\right\|_F\\
=&\left\|\mathbf{y}_t\mathbf{x}_t^{H}
\left(\mathbf{A}_3^{k+1} \otimes (\mathbf{A}_2^{k+1}-\mathbf{A}_2^{k})\right)
(\mathcal{G}_{(1)}^{k+1})^{H}\right\|_F+\left\|\mathbf{y}_t\mathbf{x}_t^{H}
\left((\mathbf{A}_3^{k+1}-\mathbf{A}_3^{k})\otimes\mathbf{A}_2^{k}\right)
(\mathcal{G}_{(1)}^{k+1})^{H}\right\|_F.
\end{aligned}
\end{equation}
By using the same discussion in (\ref{lip2_4}), we can obtain that
\begin{equation}\label{N2_8}
\begin{aligned}
&\left\|\mathbf{y}_t\mathbf{x}_t^{H}
\left(\mathbf{A}_3^{k+1} \otimes (\mathbf{A}_2^{k+1}-\mathbf{A}_2^{k})\right)
(\mathcal{G}_{(1)}^{k+1})^{H}\right\|_F\leq \nu\left\|\mathbf{y}_t\mathbf{x}_t^{H}\right\|_F
\| \mathbf{A}_2^{k+1}-\mathbf{A}_2^{k}\|_F,\\
&\left\|\mathbf{y}_t\mathbf{x}_t^{H}
\left((\mathbf{A}_3^{k+1}-\mathbf{A}_3^{k})\otimes\mathbf{A}_2^{k}\right)
(\mathcal{G}_{(1)}^{k+1})^{H}\right\|_F\leq \nu\left\|\mathbf{y}_t\mathbf{x}_t^{H}\right\|_F
\| \mathbf{A}_3^{k+1}-\mathbf{A}_3^{k}\|_F.
\end{aligned}
\end{equation}
Substituting (\ref{N2_8}) into (\ref{N2_9}) establishes
$$
\begin{aligned}
&\left\|
\mathbf{y}_t\mathbf{x}_t^{H}
\left((\mathbf{A}_3^{k+1} \otimes \mathbf{A}_2^{k+1})(\mathcal{G}_{(1)}^{k+1})^{H}
-(\mathbf{A}_3^{k} \otimes \mathbf{A}_2^{k})(\mathcal{G}_{(1)}^{k+1})^{H}\right)\right\|_F\\
\leq\ &\nu\left\|\mathbf{y}_t\mathbf{x}_t^{H}\right\|_F
\| \mathbf{A}_2^{k+1}-\mathbf{A}_2^{k}\|_F
+\nu\left\|\mathbf{y}_t\mathbf{x}_t^{H}\right\|_F
\| \mathbf{A}_3^{k+1}-\mathbf{A}_3^{k}\|_F,
\end{aligned}
$$
which indicates that
\begin{equation}\label{N2_2}
\begin{aligned}
&\frac{1}{T} \sum_{t=1}^T \left\|
\mathbf{y}_t\mathbf{x}_t^{H}
\left((\mathbf{A}_3^{k+1} \otimes \mathbf{A}_2^{k+1})(\mathcal{G}_{(1)}^{k+1})^{H}
-(\mathbf{A}_3^{k} \otimes \mathbf{A}_2^{k})(\mathcal{G}_{(1)}^{k+1})^{H}\right)\right\|_F\\
\leq \ &\frac{1}{T} \sum_{t=1}^T
\nu\left\|\mathbf{y}_t\mathbf{x}_t^{H}\right\|_F
\| \mathbf{A}_2^{k+1}-\mathbf{A}_2^{k}\|_F
+\frac{1}{T} \sum_{t=1}^T \nu\left\|\mathbf{y}_t\mathbf{x}_t^{H}\right\|_F
\| \mathbf{A}_3^{k+1}-\mathbf{A}_3^{k}\|_F\\
= \ &c_2\nu\|\mathbf{A}_2^{k+1}-\mathbf{A}_2^{k}\|_F
+c_2\nu\|\mathbf{A}_3^{k+1}-\mathbf{A}_3^{k}\|_F,
\end{aligned}
\end{equation}
where $c_2$ is defined in (\ref{DegiC2}).

For the first term on the right-hand side of (\ref{N2_1}), a similar technique as that of
(\ref{N1_15}) yields
$$
\begin{aligned}
\ & \left\|\mathbf{A}_1^{k+1}\mathcal{G}_{(1)}^{k+1}
(\mathbf{A}_3^{k+1} \otimes \mathbf{A}_2^{k+1})^{H} \mathbf{x}_t
\mathbf{x}_t^{H}
(\mathbf{A}_3^{k+1} \otimes \mathbf{A}_2^{k+1})(\mathcal{G}_{(1)}^{k+1})^{H}
\right.\\
&~~-\left.
\mathbf{A}_1^{k}\mathcal{G}_{(1)}^{k+1}
(\mathbf{A}_3^{k} \otimes \mathbf{A}_2^{k})^{H} \mathbf{x}_t
\mathbf{x}_t^{H}
(\mathbf{A}_3^{k} \otimes \mathbf{A}_2^{k})(\mathcal{G}_{(1)}^{k+1})^{H}\right\|_F\\
\leq\ &\|\mathbf{A}_1^{k+1} -\mathbf{A}_1^k\|_F
\nu^2 \left\| \mathbf{x}_t \mathbf{x}_t^{H}\right\|_F+\nu^2 \left\| \mathbf{x}_t \mathbf{x}_t^{H}\right\|_F
\left(2\|\mathbf{A}_2^{k+1}-\mathbf{A}_2^k\|_F
+2\|\mathbf{A}_3^{k+1}-\mathbf{A}_3^k\|_F\right).
\end{aligned}
$$
Summing up the above inequality over $t=1,\ldots,T$,
one can easily obtain that
\begin{equation}\label{N2_6}
\begin{aligned}
& \frac{1}{T} \sum_{t=1}^T \left\|\mathbf{A}_1^{k+1}\mathcal{G}_{(1)}^{k+1}
(\mathbf{A}_3^{k+1} \otimes \mathbf{A}_2^{k+1})^{H} \mathbf{x}_t
\mathbf{x}_t^{H}
(\mathbf{A}_3^{k+1} \otimes \mathbf{A}_2^{k+1})(\mathcal{G}_{(1)}^{k+1})^{H}
\right.\\
&~~-\left.
\mathbf{A}_1^{k}\mathcal{G}_{(1)}^{k+1}
(\mathbf{A}_3^{k} \otimes \mathbf{A}_2^{k})^{H} \mathbf{x}_t
\mathbf{x}_t^{H}
(\mathbf{A}_3^{k} \otimes \mathbf{A}_2^{k})(\mathcal{G}_{(1)}^{k+1})^{H}\right\|_F\\
\leq \ &\frac{1}{T} \sum_{t=1}^T \|\mathbf{A}_1^{k+1} -\mathbf{A}_1^k\|_F
\nu^2\left\| \mathbf{x}_t \mathbf{x}_t^{H}\right\|_F
+\frac{1}{T} \sum_{t=1}^T 2\nu^2 \left\| \mathbf{x}_t \mathbf{x}_t^{H}\right\|_F
\|\mathbf{A}_2^{k+1}-\mathbf{A}_2^k\|_F\\
&+\frac{1}{T} \sum_{t=1}^T 2 \nu^2 \left\| \mathbf{x}_t \mathbf{x}_t^{H}\right\|_F\|\mathbf{A}_3^{k+1}-\mathbf{A}_3^k\|_F\\
=\ &\nu^2c_1\|\mathbf{A}_1^{k+1}-\mathbf{A}_1^k\|_F+
2\nu^2c_1\|\mathbf{A}_2^{k+1}-\mathbf{A}_2^k\|_F
+2\nu^2c_1\|\mathbf{A}_3^{k+1}-\mathbf{A}_3^k\|_F,
\end{aligned}
\end{equation}
where $c_1$ is defined in (\ref{DefC1}).
Plugging (\ref{N2_6}) and (\ref{N2_2}) into (\ref{N2_1}),
we conclude that
$$
\begin{aligned}
\ &\left\|\nabla_{\mathbf{A}_1} Q(\llbracket \mathcal{G}^{k+1}; \mathbf{A}_1^{k+1}, \mathbf{A}_2^{k+1}, \mathbf{A}_3^{k+1}\rrbracket)
-\nabla_{\mathbf{A}_1} Q(\llbracket \mathcal{G}^{k+1}; \mathbf{A}_1^{k}, \mathbf{A}_2^{k}, \mathbf{A}_3^{k}\rrbracket)\right\|_F\\
\leq \ &\nu c_2\|\mathbf{A}_2^{k+1}-\mathbf{A}_2^{k}\|_F
+\nu c_2\|\mathbf{A}_3^{k+1}-\mathbf{A}_3^{k}\|_F+
\nu^2c_1\|\mathbf{A}_1^{k+1} -\mathbf{A}_1^k\|_F\\
 \ &+2\nu^2c_1\|\mathbf{A}_2^{k+1}-\mathbf{A}_2^k\|_F
+2\nu^2c_1\|\mathbf{A}_3^{k+1}-\mathbf{A}_3^k\|_F.
\end{aligned}
$$
Combining the above inequality and (\ref{N2_7}), we
arrive at
\begin{equation}\label{N2}
\begin{aligned}
\|N_2\|_F=\ &\left\|\nabla_{\mathbf{A}_1} \Psi(\mathcal{G}^{k+1}, \mathbf{A}_1^{k+1},\mathbf{A}_2^{k+1},
\mathbf{A}_3^{k+1},\mathbf{U}_1^{k+1},
\mathbf{U}_2^{k+1},\mathbf{U}_3^{k+1})\right.\\
&\left.-\nabla_{\mathbf{A}_1} \Psi(\mathcal{G}^{k+1}, \mathbf{A}_1^{k},\mathbf{A}_2^{k},
\mathbf{A}_3^{k},\mathbf{U}_1^{k},
\mathbf{U}_2^{k},\mathbf{U}_3^{k})
-\rho_2(\mathbf{A}_1^{k+1}-\mathbf{A}_1^{k})\right\|_F\\
\leq \ &\nu c_2\|\mathbf{A}_2^{k+1}-\mathbf{A}_2^{k}\|_F+
\nu c_2\|\mathbf{A}_3^{k+1}-\mathbf{A}_3^{k}\|_F+
\nu ^2c_1\|\mathbf{A}_1^{k+1} -\mathbf{A}_1^k\|_F\\
&+2\nu ^2c_1\|\mathbf{A}_2^{k+1}-\mathbf{A}_2^k\|_F
+2\nu ^2c_1\|\mathbf{A}_3^{k+1}-\mathbf{A}_3^k\|_F
+\gamma_1\|\mathbf{U}_1^{k+1}-\mathbf{U}_1^k\|_F\\
&+(\rho_2+\gamma_1)\|\mathbf{A}_1^{k+1}-\mathbf{A}_1^k\|_F\\
= \ &(\nu^2c_1+\rho_2+\gamma_1)\|\mathbf{A}_1^{k+1}-\mathbf{A}_1^k\|_F
+(\nu c_2+2\nu^2c_1)\|\mathbf{A}_2^{k+1}-\mathbf{A}_2^k\|_F\\
&+(\nu c_2+2\nu^2c_1)\|\mathbf{A}_3^{k+1}-\mathbf{A}_3^{k}\|_F
+\gamma_1\|\mathbf{U}_1^{k+1}-\mathbf{U}_1^k\|_F.
\end{aligned}
\end{equation}

Similar to the analysis of $N_2$, we have
\begin{equation}\label{N3_1}
\begin{aligned}
\|N_3\|_F=\ &\left\|\nabla_{\mathbf{A}_2} \Psi(\mathcal{G}^{k+1}, \mathbf{A}_1^{k+1},\mathbf{A}_2^{k+1},
\mathbf{A}_3^{k+1},\mathbf{U}_1^{k+1},
\mathbf{U}_2^{k+1},\mathbf{U}_3^{k+1})\right.\\
&\left.-\nabla_{\mathbf{A}_2} \Psi(\mathcal{G}^{k+1}, \mathbf{A}_1^{k+1},\mathbf{A}_2^{k},
\mathbf{A}_3^{k},\mathbf{U}_1^{k},
\mathbf{U}_2^{k},\mathbf{U}_3^{k})
-\rho_3(\mathbf{A}_2^{k+1}-\mathbf{A}_2^{k})\right\|_F\\
\leq \ & \left\|\nabla_{\mathbf{A}_2} Q(\llbracket \mathcal{G}^{k+1}; \mathbf{A}_1^{k+1}, \mathbf{A}_2^{k+1}, \mathbf{A}_3^{k+1}\rrbracket)
-\nabla_{\mathbf{A}_2} Q(\llbracket \mathcal{G}^{k+1}; \mathbf{A}_1^{k+1}, \mathbf{A}_2^{k}, \mathbf{A}_3^{k}\rrbracket)\right\|_F \\
& +\gamma_2\|\mathbf{U}_2^{k+1}-\mathbf{U}_2^k\|_F
+(\gamma_2+\rho_3)\|\mathbf{A}_2^{k+1}-\mathbf{A}_2^k\|_F.
\end{aligned}
\end{equation}
Note that (\ref{gradient_Q}) together with (\ref{gradient_w}) yields
\begin{equation}\label{N3_2}
\begin{aligned}
\ &\left\|\nabla_{\mathbf{A}_2} Q(\llbracket \mathcal{G}^{k+1}; \mathbf{A}_1^{k+1}, \mathbf{A}_2^{k+1}, \mathbf{A}_3^{k+1}\rrbracket)
-\nabla_{\mathbf{A}_2} Q(\llbracket \mathcal{G}^{k+1}; \mathbf{A}_1^{k+1}, \mathbf{A}_2^{k}, \mathbf{A}_3^{k}\rrbracket)\right\|_F\\
=\ &\left\|\bigg(\frac{1}{T} \sum_{t=1}^T (\mathbf{A}_1^{k+1}\mathcal{G}_{(1)}^{k+1}
(\mathbf{A}_3^{k+1} \otimes \mathbf{A}_2^{k+1})^{H} \mathbf{x}_t-\mathbf{y}_t)\circ \mathbf{X}_t\bigg)_{(2)}(\mathbf{A}_1^{k+1} \otimes \mathbf{A}_3^{k+1})(\mathcal{G}_{(2)}^{k+1})^{H}
\right.\\
\ &~~-\left.
\bigg(\frac{1}{T} \sum_{t=1}^T (\mathbf{A}_1^{k+1}\mathcal{G}_{(1)}^{k+1}
(\mathbf{A}_3^{k} \otimes \mathbf{A}_2^{k})^{H} \mathbf{x}_t-\mathbf{y}_t)\circ \mathbf{X}_t\bigg)_{(2)}(\mathbf{A}_1^{k+1} \otimes \mathbf{A}_3^{k})(\mathcal{G}_{(2)}^{k+1})^{H}\right\|_F\\
=\ & \left\| \frac{1}{T} \sum_{t=1}^T \mathbf{X}_t \otimes \left(\mathbf{A}_1^{k+1}
\mathcal{G}_{(1)}^{k+1}
(\mathbf{A}_3^{k+1} \otimes \mathbf{A}_2^{k+1})^{H} \mathbf{x}_t-
\mathbf{y}_t\right)^{H}
(\mathbf{A}_1^{k+1} \otimes \mathbf{A}_3^{k+1})(\mathcal{G}_{(2)}^{k+1})^{H}\right.\\
\ &~~-\left.
\frac{1}{T} \sum_{t=1}^T \mathbf{X}_t \otimes \left(\mathbf{A}_1^{k+1}\mathcal{G}_{(1)}^{k+1}
(\mathbf{A}_3^{k} \otimes \mathbf{A}_2^{k})^{H} \mathbf{x}_t-\mathbf{y}_t\right)^{H}
(\mathbf{A}_1^{k+1} \otimes \mathbf{A}_3^{k})(\mathcal{G}_{(2)}^{k+1})^{H}\right\|_F\\
\leq \ & \frac{1}{T} \sum_{t=1}^T \left\|  \mathbf{X}_t \otimes
\left(\mathbf{A}_1^{k+1}\mathcal{G}_{(1)}^{k+1}
(\mathbf{A}_3^{k+1} \otimes \mathbf{A}_2^{k+1})^{H}
\mathbf{x}_t-\mathbf{y}_t\right)^{H}
(\mathbf{A}_1^{k+1} \otimes \mathbf{A}_3^{k+1})(\mathcal{G}_{(2)}^{k+1})^{H}\right.\\
\ &~~ -\left.
\mathbf{X}_t \otimes \left(\mathbf{A}_1^{k+1}\mathcal{G}_{(1)}^{k+1}
(\mathbf{A}_3^{k} \otimes \mathbf{A}_2^{k})^{H} \mathbf{x}_t-\mathbf{y}_t\right)^{H}
(\mathbf{A}_1^{k+1}\otimes \mathbf{A}_3^{k})(\mathcal{G}_{(2)}^{k+1})^{H}\right\|_F\\
=\ &\frac{1}{T} \sum_{t=1}^T \left\|
\mathbf{X}_t \otimes\left(\mathbf{A}_1^{k+1}\mathcal{G}_{(1)}^{k+1}
(\mathbf{A}_3^{k+1} \otimes \mathbf{A}_2^{k+1})^{H}
\mathbf{x}_t-\mathbf{y}_t\right)^{H}\mathbf{A}_1^{k+1} \otimes \mathbf{A}_3^{k+1})
(\mathcal{G}_{(2)}^{k+1})^{H}
\right.\\
\ &~~ -\left.
\mathbf{X}_t \otimes\left(\mathbf{A}_1^{k+1}\mathcal{G}_{(1)}^{k+1}
(\mathbf{A}_3^{k+1} \otimes \mathbf{A}_2^{k})^{H}
\mathbf{x}_t-\mathbf{y}_t\right)^{H}(\mathbf{A}_1^{k+1} \otimes \mathbf{A}_3^{k+1})
(\mathcal{G}_{(2)}^{k+1})^{H}
\right.\\
\ &~~ +\left.
\mathbf{X}_t \otimes\left(\mathbf{A}_1^{k+1}\mathcal{G}_{(1)}^{k+1}
(\mathbf{A}_3^{k+1} \otimes \mathbf{A}_2^{k})^{H}
\mathbf{x}_t-\mathbf{y}_t\right)^{H}(\mathbf{A}_1^{k+1} \otimes
\mathbf{A}_3^{k+1})(\mathcal{G}_{(2)}^{k+1})^{H}
\right.\\
 \ &~~-\left.
\mathbf{X}_t \otimes\left(\mathbf{A}_1^{k+1}\mathcal{G}_{(1)}^{k+1}
(\mathbf{A}_3^{k} \otimes \mathbf{A}_2^{k})^{H}
\mathbf{x}_t-\mathbf{y}_t\right)^{H}(\mathbf{A}_1^{k+1} \otimes
\mathbf{A}_3^{k+1})(\mathcal{G}_{(2)}^{k+1})^{H}
\right.\\
\ &~~ +\left.
\mathbf{X}_t \otimes\left(\mathbf{A}_1^{k+1}\mathcal{G}_{(1)}^{k+1}
(\mathbf{A}_3^{k} \otimes \mathbf{A}_2^{k})^{H}
\mathbf{x}_t-\mathbf{y}_t\right)^{H}(\mathbf{A}_1^{k+1} \otimes
\mathbf{A}_3^{k+1})(\mathcal{G}_{(2)}^{k+1})^{H}
\right.\\
\ &~~-\left.
\mathbf{X}_t \otimes\left(\mathbf{A}_1^{k+1}\mathcal{G}_{(1)}^{k+1}
(\mathbf{A}_3^{k} \otimes \mathbf{A}_2^{k})^{H}
\mathbf{x}_t-\mathbf{y}_t\right)^{H}(\mathbf{A}_1^{k+1} \otimes
 \mathbf{A}_3^{k})(\mathcal{G}_{(2)}^{k+1})^{H}\right\|_F,
\end{aligned}
\end{equation}
where the second equality holds by Proposition \ref{out_kron1}.
Therefore, we get
\begin{equation}\label{QA2Gt}
	\begin{aligned}
\ &\left\|\nabla_{\mathbf{A}_2} Q(\llbracket \mathcal{G}^{k+1}; \mathbf{A}_1^{k+1}, \mathbf{A}_2^{k+1}, \mathbf{A}_3^{k+1}\rrbracket)
-\nabla_{\mathbf{A}_2} Q(\llbracket \mathcal{G}^{k+1}; \mathbf{A}_1^{k+1}, \mathbf{A}_2^{k}, \mathbf{A}_3^{k}\rrbracket)\right\|_F\\
\leq\ &\frac{1}{T} \sum_{t=1}^T \left\|
\mathbf{X}_t \otimes\left(\mathbf{A}_1^{k+1}\mathcal{G}_{(1)}^{k+1}
(\mathbf{A}_3^{k+1} \otimes \mathbf{A}_2^{k+1})^{H}
\mathbf{x}_t-\mathbf{y}_t\right)^{H}\mathbf{A}_1^{k+1} \otimes \mathbf{A}_3^{k+1})
(\mathcal{G}_{(2)}^{k+1})^{H}
\right.\\
\ &~~ -\left.
\mathbf{X}_t \otimes\left(\mathbf{A}_1^{k+1}\mathcal{G}_{(1)}^{k+1}
(\mathbf{A}_3^{k+1} \otimes \mathbf{A}_2^{k})^{H}
\mathbf{x}_t-\mathbf{y}_t\right)^{H}(\mathbf{A}_1^{k+1} \otimes \mathbf{A}_3^{k+1})
(\mathcal{G}_{(2)}^{k+1})^{H}\right\|_F\\
&~~+\frac{1}{T} \sum_{t=1}^T \left\|
\mathbf{X}_t \otimes\left(\mathbf{A}_1^{k+1}\mathcal{G}_{(1)}^{k+1}
(\mathbf{A}_3^{k+1} \otimes \mathbf{A}_2^{k})^{H}
\mathbf{x}_t-\mathbf{y}_t\right)^{H}(\mathbf{A}_1^{k+1} \otimes
\mathbf{A}_3^{k+1})(\mathcal{G}_{(2)}^{k+1})^{H}
\right.\\
\ &~~ -\left.
\mathbf{X}_t \otimes\left(\mathbf{A}_1^{k+1}\mathcal{G}_{(1)}^{k+1}
(\mathbf{A}_3^{k} \otimes \mathbf{A}_2^{k})^{H}
\mathbf{x}_t-\mathbf{y}_t\right)^{H}(\mathbf{A}_1^{k+1} \otimes
\mathbf{A}_3^{k+1})(\mathcal{G}_{(2)}^{k+1})^{H}\right\|_F\\
\ &~~ +\frac{1}{T} \sum_{t=1}^T \left\|
\mathbf{X}_t \otimes\left(\mathbf{A}_1^{k+1}\mathcal{G}_{(1)}^{k+1}
(\mathbf{A}_3^{k} \otimes \mathbf{A}_2^{k})^{H}
\mathbf{x}_t-\mathbf{y}_t\right)^{H}(\mathbf{A}_1^{k+1} \otimes
\mathbf{A}_3^{k+1})(\mathcal{G}_{(2)}^{k+1})^{H}
\right.\\
\ &~~ -\left.
\mathbf{X}_t \otimes\left(\mathbf{A}_1^{k+1}\mathcal{G}_{(1)}^{k+1}
(\mathbf{A}_3^{k} \otimes \mathbf{A}_2^{k})^{H}
\mathbf{x}_t-\mathbf{y}_t\right)^{H}(\mathbf{A}_1^{k+1} \otimes
\mathbf{A}_3^{k})(\mathcal{G}_{(2)}^{k+1})^{H}\right\|_F.
\end{aligned}
\end{equation}

Now we focus on the first term on the right-hand side of the above inequality.
$$
\begin{aligned}
\ &\left\|
\mathbf{X}_t \otimes\left(\mathbf{A}_1^{k+1}\mathcal{G}_{(1)}^{k+1}
(\mathbf{A}_3^{k+1} \otimes \mathbf{A}_2^{k+1})^{H}
\mathbf{x}_t-\mathbf{y}_t\right)^{H}(\mathbf{A}_1^{k+1} \otimes
\mathbf{A}_3^{k+1})(\mathcal{G}_{(2)}^{k+1})^{H}
\right.\\
\ &~~-\left.
\mathbf{X}_t \otimes\left(\mathbf{A}_1^{k+1}\mathcal{G}_{(1)}^{k+1}
(\mathbf{A}_3^{k+1} \otimes \mathbf{A}_2^{k})^{H}
\mathbf{x}_t-\mathbf{y}_t\right)^{H}(\mathbf{A}_1^{k+1} \otimes
\mathbf{A}_3^{k+1})(\mathcal{G}_{(2)}^{k+1})^{H}\right\|_F\\
\leq \ &
\left\|\mathbf{X}_t \otimes\left(\mathbf{A}_1^{k+1}\mathcal{G}_{(1)}^{k+1}
(\mathbf{A}_3^{k+1} \otimes \mathbf{A}_2^{k+1})^{H}
 \mathbf{x}_t-\mathbf{y}_t\right)^{H}
\right.\\
\ &~~-\left.
\mathbf{X}_t \otimes\left(\mathbf{A}_1^{k+1}\mathcal{G}_{(1)}^{k+1}
(\mathbf{A}_3^{k+1} \otimes \mathbf{A}_2^{k})^{H}
\mathbf{x}_t-\mathbf{y}_t\right)^{H}\right\|_F\left\|(\mathbf{A}_1^{k+1} \otimes
\mathbf{A}_3^{k+1})(\mathcal{G}_{(2)}^{k+1})^{H}\right\|_F\\
\leq \ &
\left\|\mathbf{X}_t \otimes\left(\mathbf{A}_1^{k+1}\mathcal{G}_{(1)}^{k+1}
(\mathbf{A}_3^{k+1} \otimes \mathbf{A}_2^{k+1})^{H} \mathbf{x}_t-\mathbf{y}_t\right)^{H}
\right.\\
\ &~~ -\left.
\mathbf{X}_t \otimes\left(\mathbf{A}_1^{k+1}\mathcal{G}_{(1)}^{k+1}
(\mathbf{A}_3^{k+1} \otimes \mathbf{A}_2^{k})^{H}
\mathbf{x}_t-\mathbf{y}_t\right)^{H}\right\|_F \nu\\
=\ &
\nu\left\|\mathbf{X}_t \otimes\left(\mathbf{x}_t^{H}(\mathbf{A}_3^{k+1}
 \otimes(\mathbf{A}_2^{k+1}-\mathbf{A}_2^k))
(\mathcal{G}_{(1)}^{k+1})^{H}
(\mathbf{A}_1^{k+1})^{H}\right)
\right\|_F,
\end{aligned}
$$
where the second inequality follows from (\ref{ineq2}).
Using (\ref{lip2_4}) further leads to
$$
\begin{aligned}
\ &
\nu\left\|\mathbf{X}_t \otimes\left(\mathbf{x}_t^{H}\left(\mathbf{A}_3^{k+1}
\otimes(\mathbf{A}_2^{k+1}-\mathbf{A}_2^k)\right)
(\mathcal{G}_{(1)}^{k+1})^{H}
(\mathbf{A}_1^{k+1})^{H}\right)
\right\|_F\\
= \ &\nu\left\|\mathbf{X}_t\right\|_F\left\|\left(\mathbf{x}_t^{H}(\mathbf{A}_3^{k+1}
 \otimes(\mathbf{A}_2^{k+1}-\mathbf{A}_2^k))
(\mathcal{G}_{(1)}^{k+1})^{H}
(\mathbf{A}_1^{k+1})^{H}\right)
\right\|_F\\
\leq\ &\nu\left\|\mathbf{X}_t\right\|_F\left\|\mathbf{x}_t\right\|_2
\nu\|\mathbf{A}_2^{k+1}-\mathbf{A}_2^k\|_F\\
= \ &\nu^2\left\|\mathbf{x}_t\right\|_2^2
\|\mathbf{A}_2^{k+1}-\mathbf{A}_2^k\|_F,
\end{aligned}
$$
where the first equality holds by Proposition \ref{kron_pro} and the second equality follows from the fact that
$\|\mathbf{x}_t\|_2=\left\|\mathbf{X}_t\right\|_F$.
Summing the aforementioned inequality from $t+1$ to $T$ yields
\begin{equation}\label{N3_3}
\begin{aligned}
\ &
\frac{1}{T} \sum_{t=1}^T \nu\left\|\mathbf{X}_t \otimes\left(\mathbf{x}_t^{H}\left(\mathbf{A}_3^{k+1}
\otimes(\mathbf{A}_2^{k+1}-\mathbf{A}_2^k)\right)
(\mathcal{G}_{(1)}^{k+1})^{H}
(\mathbf{A}_1^{k+1})^{H}\right)
\right\|_F\\
\leq\ &\frac{1}{T} \sum_{t=1}^T \nu^2\left\|\mathbf{x}_t\right\|_2^2
\|\mathbf{A}_2^{k+1}-\mathbf{A}_2^k\|_F
= \nu^2c_1\|\mathbf{A}_2^{k+1}-\mathbf{A}_2^k\|_F,
\end{aligned}
\end{equation}
where $c_1$ is defined in (\ref{DefC1}).
Similarly, we can get that
\begin{equation}\label{N3_4}
\begin{aligned}
&\frac{1}{T} \sum_{t=1}^T \left\|
\mathbf{X}_t \otimes\left(\mathbf{A}_1^{k+1}\mathcal{G}_{(1)}^{k+1}
(\mathbf{A}_3^{k+1} \otimes \mathbf{A}_2^{k})^{H} \mathbf{x}_t-\mathbf{y}_t\right)^{H}(\mathbf{A}_1^{k+1} \otimes \mathbf{A}_3^{k+1})(\mathcal{G}_{(2)}^{k+1})^{H}
\right.\\
-&\left.
\mathbf{X}_t \otimes\left(\mathbf{A}_1^{k+1}\mathcal{G}_{(1)}^{k+1}
(\mathbf{A}_3^{k} \otimes \mathbf{A}_2^{k})^{H}
\mathbf{x}_t-\mathbf{y}_t\right)^{H}(\mathbf{A}_1^{k+1} \otimes
\mathbf{A}_3^{k+1})(\mathcal{G}_{(2)}^{k+1})^{H}\right\|_F
\leq \nu^2c_1
\|\mathbf{A}_3^{k+1}-\mathbf{A}_3^k\|_F,\\
&\frac{1}{T} \sum_{t=1}^T \left\|\mathbf{X}_t \otimes\left(\mathbf{A}_1^{k+1}\mathcal{G}_{(1)}^{k+1}
(\mathbf{A}_3^{k} \otimes \mathbf{A}_2^{k})^{H} \mathbf{x}_t-\mathbf{y}_t\right)^{H}(\mathbf{A}_1^{k+1} \otimes \mathbf{A}_3^{k+1})(\mathcal{G}_{(2)}^{k+1})^{H}
\right.\\
-&\left.
\mathbf{X}_t \otimes\left(\mathbf{A}_1^{k+1}\mathcal{G}_{(1)}^{k+1}
(\mathbf{A}_3^{k} \otimes \mathbf{A}_2^{k})^{H}
\mathbf{x}_t-\mathbf{y}_t\right)^{H}(\mathbf{A}_1^{k+1} \otimes
\mathbf{A}_3^{k})(\mathcal{G}_{(2)}^{k+1})^{H}\right\|_F
\leq
\nu^2c_1\|\mathbf{A}_3^{k+1}-\mathbf{A}_3^k\|_F.
\end{aligned}
\end{equation}
Plugging (\ref{N3_4}) and (\ref{N3_3}) into (\ref{QA2Gt}) yields
\begin{equation}\label{N3_5}
\begin{aligned}
\ &\left\|\nabla_{\mathbf{A}_2} Q(\llbracket \mathcal{G}^{k+1}; \mathbf{A}_1^{k+1}, \mathbf{A}_2^{k+1}, \mathbf{A}_3^{k+1}\rrbracket)
-\nabla_{\mathbf{A}_2} Q(\llbracket \mathcal{G}^{k+1}; \mathbf{A}_1^{k+1}, \mathbf{A}_2^{k}, \mathbf{A}_3^{k}\rrbracket)\right\|_F\\
\leq \
&\nu^2c_1\|\mathbf{A}_2^{k+1}-\mathbf{A}_2^k\|_F+
2\nu^2c_1\|\mathbf{A}_3^{k+1}-\mathbf{A}_3^k\|_F.
\end{aligned}
\end{equation}
Combining (\ref{N3_5}) and (\ref{N3_1}), we obtain that
\begin{equation}\label{N3}
\begin{aligned}
\|N_3\|_F
\leq \ &\nu^2c_1\|\mathbf{A}_2^{k+1}-\mathbf{A}_2^k\|_F+
2\nu^2c_1\|\mathbf{A}_3^{k+1}-\mathbf{A}_3^k\|_F
+\gamma_2\|\mathbf{U}_2^{k+1}-\mathbf{U}_2^k\|_F\\
\ &+(\gamma_2+\rho_3)\|\mathbf{A}_2^{k+1}-\mathbf{A}_2^k\|_F\\
= \  &(\nu^2c_1+\gamma_2+\rho_3)\|\mathbf{A}_2^{k+1}-\mathbf{A}_2^k\|_F+
2\nu^2c_1\|\mathbf{A}_3^{k+1}-\mathbf{A}_3^k\|_F
+\gamma_2\|\mathbf{U}_2^{k+1}-\mathbf{U}_2^k\|_F.
\end{aligned}
\end{equation}

Regarding $N_4$ defined in (\ref{gradient}),
we have
\begin{equation}\label{N4_1}
\begin{aligned}
\|N_4\|_F=\ &\left\|\nabla_{\mathbf{A}_3} \Psi(\mathcal{G}^{k+1}, \mathbf{A}_1^{k+1},\mathbf{A}_2^{k+1},
\mathbf{A}_3^{k+1},\mathbf{U}_1^{k+1},
\mathbf{U}_2^{k+1},\mathbf{U}_3^{k+1})\right.\\
\ &~~\left.-\nabla_{\mathbf{A}_3} \Psi(\mathcal{G}^{k+1}, \mathbf{A}_1^{k+1},\mathbf{A}_2^{k+1},
\mathbf{A}_3^{k},\mathbf{U}_1^{k},
\mathbf{U}_2^{k},\mathbf{U}_3^{k})
-\rho_4(\mathbf{A}_3^{k+1}-\mathbf{A}_3^{k})\right\|_F\\
\leq\  & \left\|\nabla_{\mathbf{A}_3} Q(\llbracket \mathcal{G}^{k+1}; \mathbf{A}_1^{k+1}, \mathbf{A}_2^{k+1}, \mathbf{A}_3^{k+1}\rrbracket)
-\nabla_{\mathbf{A}_3} Q(\llbracket \mathcal{G}^{k+1}; \mathbf{A}_1^{k+1}, \mathbf{A}_2^{k+1}, \mathbf{A}_3^{k}\rrbracket)\right\|_F \\
\ &~~ +\gamma_3\|\mathbf{U}_3^{k+1}-\mathbf{U}_3^k\|_F
+(\gamma_3+\rho_4)\|\mathbf{A}_3^{k+1}-\mathbf{A}_3^k\|_F.
\end{aligned}
\end{equation}
By (\ref{gradient_Q}), (\ref{gradient_w}) and Proposition \ref{out_kron1},
we can get
\begin{equation}\label{N4_2}
\begin{aligned}
\ &\left\|\nabla_{\mathbf{A}_3} Q(\llbracket \mathcal{G}^{k+1}; \mathbf{A}_1^{k+1}, \mathbf{A}_2^{k+1}, \mathbf{A}_3^{k+1}\rrbracket)
-\nabla_{\mathbf{A}_3} Q(\llbracket \mathcal{G}^{k+1}; \mathbf{A}_1^{k+1}, \mathbf{A}_2^{k+1}, \mathbf{A}_3^{k}\rrbracket)\right\|_F\\
=\ &\left\|\bigg(\frac{1}{T} \sum_{t=1}^T (\mathbf{A}_1^{k+1}\mathcal{G}_{(1)}^{k+1}
(\mathbf{A}_3^{k+1} \otimes \mathbf{A}_2^{k+1})^{H} \mathbf{x}_t-\mathbf{y}_t)\circ \mathbf{X}_t\bigg)_{(3)}(\mathbf{A}_2^{k+1} \otimes \mathbf{A}_1^{k+1})(\mathcal{G}_{(3)}^{k+1})^{H}
\right.\\
\ &~~ -\left.
\bigg(\frac{1}{T} \sum_{t=1}^T (\mathbf{A}_1^{k+1}\mathcal{G}_{(1)}^{k+1}
(\mathbf{A}_3^{k} \otimes \mathbf{A}_2^{k+1})^{H} \mathbf{x}_t-\mathbf{y}_t)\circ \mathbf{X}_t\bigg)_{(3)}(\mathbf{A}_2^{k+1} \otimes \mathbf{A}_1^{k+1})(\mathcal{G}_{(3)}^{k+1})^{H}\right\|_F\\
=\ & \left\| \frac{1}{T} \sum_{t=1}^T \mathbf{X}_t^{H} \otimes \left(\mathbf{A}_1^{k+1}\mathcal{G}_{(1)}^{k+1}
(\mathbf{A}_3^{k+1} \otimes \mathbf{A}_2^{k+1})^{H} \mathbf{x}_t-\mathbf{y}_t\right)^{H}
(\mathbf{A}_2^{k+1} \otimes \mathbf{A}_1^{k+1})(\mathcal{G}_{(3)}^{k+1})^{H}\right.\\
\ &~~-\left.
\frac{1}{T} \sum_{t=1}^T \mathbf{X}_t^{H} \otimes \left(\mathbf{A}_1^{k+1}\mathcal{G}_{(1)}^{k+1}
(\mathbf{A}_3^{k} \otimes \mathbf{A}_2^{k+1})^{H} \mathbf{x}_t-\mathbf{y}_t\right)^{H}
(\mathbf{A}_2^{k+1} \otimes \mathbf{A}_1^{k+1})(\mathcal{G}_{(3)}^{k+1})^{H}\right\|_F\\
\leq \ & \frac{1}{T} \sum_{t=1}^T \left\|  \mathbf{X}_t^{H} \otimes \left(\mathbf{A}_1^{k+1}\mathcal{G}_{(1)}^{k+1}
(\mathbf{A}_3^{k+1} \otimes \mathbf{A}_2^{k+1})^{H} \mathbf{x}_t-\mathbf{y}_t\right)^{H}
(\mathbf{A}_2^{k+1} \otimes \mathbf{A}_1^{k+1})(\mathcal{G}_{(3)}^{k+1})^{H}\right.\\
\ &~~ -\left.
\mathbf{X}_t^{H} \otimes \left(\mathbf{A}_1^{k+1}\mathcal{G}_{(1)}^{k+1}
(\mathbf{A}_3^{k} \otimes \mathbf{A}_2^{k+1})^{H} \mathbf{x}_t-\mathbf{y}_t\right)^{H}
(\mathbf{A}_2^{k+1} \otimes \mathbf{A}_1^{k+1})(\mathcal{G}_{(3)}^{k+1})^{H}\right\|_F\\
\leq \ & \frac{1}{T} \sum_{t=1}^T \left\|  \mathbf{X}_t^{H} \otimes \left(\mathbf{A}_1^{k+1}\mathcal{G}_{(1)}^{k+1}
(\mathbf{A}_3^{k+1} \otimes \mathbf{A}_2^{k+1})^{H} \mathbf{x}_t-\mathbf{y}_t\right)^{H}
\right.\\
\ &~~ -\left.
\mathbf{X}_t^{H} \otimes \left(\mathbf{A}_1^{k}\mathcal{G}_{(1)}^{k+1}
(\mathbf{A}_3^{k} \otimes \mathbf{A}_2^{k+1})^{H} \mathbf{x}_t-\mathbf{y}_t\right)^{H}\right\|_F
\left\|(\mathbf{A}_2^{k+1} \otimes \mathbf{A}_1^{k+1})(\mathcal{G}_{(3)}^{k+1})^{H}\right\|_F\\
\leq \ & \frac{1}{T} \sum_{t=1}^T \left\|  \mathbf{X}_t^{H} \otimes \left(\mathbf{A}_1^{k+1}\mathcal{G}_{(1)}^{k+1}
(\mathbf{A}_3^{k+1} \otimes \mathbf{A}_2^{k+1})^{H} \mathbf{x}_t-\mathbf{y}_t\right)^{H}
\right.\\
\ &~~-\left.
\mathbf{X}_t^{H} \otimes \left(\mathbf{A}_1^{k+1}\mathcal{G}_{(1)}^{k+1}
(\mathbf{A}_3^{k} \otimes \mathbf{A}_2^{k+1})^{H} \mathbf{x}_t-\mathbf{y}_t\right)^{H}\right\|_F \nu,
\end{aligned}
\end{equation}
where the last equality holds by (\ref{ineq2}).
Using similar techniques as those of (\ref{lip2_4}) to obtain
$$
\begin{aligned}
\ &\left\|  \mathbf{X}_t^{H} \otimes (\mathbf{A}_1^{k+1}\mathcal{G}_{(1)}^{k+1}
(\mathbf{A}_3^{k+1} \otimes \mathbf{A}_2^{k+1})^{H} \mathbf{x}_t-\mathbf{y}_t)^{H}
-
\mathbf{X}_t^{H} \otimes (\mathbf{A}_1^{k+1}\mathcal{G}_{(1)}^{k+1}
(\mathbf{A}_3^{k} \otimes \mathbf{A}_2^{k+1})^{H} \mathbf{x}_t-\mathbf{y}_t)^{H}\right\|_F\\
= &\left\|
\mathbf{X}_t^{H}  \otimes(\mathbf{A}_1^{k+1}\mathcal{G}_{(1)}^{k+1}
((\mathbf{A}_3^{k+1}-\mathbf{A}_3^{k}) \otimes \mathbf{A}_2^{k+1})^{H} \mathbf{x}_t)^{H}\right\|_F\\
\leq  &\left\| \mathbf{x}_t \right\|_2^2 \nu\|\mathbf{A}_3^{k+1}- \mathbf{A}_3^{k}\|_F.
\end{aligned}
$$
Combining the above inequality with (\ref{N4_2}) to get
\begin{equation}\label{N4_3}
\begin{aligned}
\ &\left\|\nabla_{\mathbf{A}_3} Q(\llbracket \mathcal{G}^{k+1}; \mathbf{A}_1^{k+1}, \mathbf{A}_2^{k+1}, \mathbf{A}_3^{k+1}\rrbracket)
-\nabla_{\mathbf{A}_3} Q(\llbracket \mathcal{G}^{k+1}; \mathbf{A}_1^{k+1}, \mathbf{A}_2^{k+1}, \mathbf{A}_3^{k}\rrbracket)\right\|_F\\
\leq \ &\frac{1}{T} \sum_{t=1}^T \left\| \mathbf{x}_t \right\|_2^2
\nu^2\|\mathbf{A}_3^{k+1}- \mathbf{A}_3^{k}\|_F
=\nu^2c_1\|\mathbf{A}_3^{k+1}- \mathbf{A}_3^{k}\|_F.
\end{aligned}
\end{equation}
Substituting (\ref{N4_3}) into (\ref{N4_1}) yields
\begin{equation}\label{N4}
\begin{aligned}
\|N_4\|_F
\leq\ &\nu^2c_1\|\mathbf{A}_3^{k+1}- \mathbf{A}_3^{k}\|_F
+\gamma_3\|\mathbf{U}_3^{k+1}-\mathbf{U}_3^k\|_F
+(\gamma_3+\rho_4)\|\mathbf{A}_3^{k+1}-\mathbf{A}_3^k\|_F\\
=\ &(\nu^2c_1+\gamma_3+\rho_4)\|\mathbf{A}_3^{k+1}- \mathbf{A}_3^{k}\|_F
+\gamma_3\|\mathbf{U}_3^{k+1}-\mathbf{U}_3^k\|_F.
\end{aligned}
\end{equation}

For $N_5$ in (\ref{gradient}), it follows from (\ref{gradient_H}) that
\begin{equation}\label{N5}
\begin{aligned}
\|N_5\|_F= \ &\left\|\nabla_{\mathbf{U}_1} \Psi(\mathcal{G}^{k+1}, \mathbf{A}_1^{k+1},\mathbf{A}_2^{k+1},
\mathbf{A}_3^{k+1},\mathbf{U}_1^{k+1},
\mathbf{U}_2^{k+1},\mathbf{U}_3^{k+1})\right.\\
\ &~~ \left.-\nabla_{\mathbf{U}_1}\Psi(\mathcal{G}^{k+1}, \mathbf{A}_1^{k+1},\mathbf{A}_2^{k+1},
\mathbf{A}_3^{k+1},\mathbf{U}_1^{k},
\mathbf{U}_2^{k},\mathbf{U}_3^{k})
-\rho_5(\mathbf{U}_1^{k+1}-\mathbf{U}_1^{k})\right\|_F\\
=\  & \left\|\gamma_1(\mathbf{U}_1^{k+1}-\mathbf{A}_1^{k+1})
-\gamma_1(\mathbf{U}_1^{k}-\mathbf{A}_1^{k+1})
-\rho_5(\mathbf{U}_1^{k+1}-\mathbf{U}_1^{k})\right\|_F\\
=\  & \left\|\gamma_1 (\mathbf{U}_1^{k+1}-\mathbf{U}_1^{k})
-\rho_5(\mathbf{U}_1^{k+1}-\mathbf{U}_1^{k})\right\|_F\\
\leq\  & \|\gamma_1(\mathbf{U}_1^{k+1}-\mathbf{U}_1^{k})\|_F
+\|\rho_5(\mathbf{U}_1^{k+1}-\mathbf{U}_1^{k})\|_F\\
=\ &(\gamma_1+\rho_5)\|\mathbf{U}_1^{k+1}-\mathbf{U}_1^{k}\|_F.
\end{aligned}
\end{equation}
Similarly, we can easily obtain that
\begin{equation}\label{N6}
\begin{aligned}
\|N_6\|_F=&\left\|\nabla_{\mathbf{U}_2} \Psi(\mathcal{G}^{k+1}, \mathbf{A}_1^{k+1},\mathbf{A}_2^{k+1},
\mathbf{A}_3^{k+1},\mathbf{U}_1^{k+1},
\mathbf{U}_2^{k+1},\mathbf{U}_3^{k+1})\right.\\
&\left.-\nabla_{\mathbf{U}_2} \Psi(\mathcal{G}^{k+1}, \mathbf{A}_1^{k+1},\mathbf{A}_2^{k+1},
\mathbf{A}_3^{k+1},\mathbf{U}_1^{k+1},
\mathbf{U}_2^{k},\mathbf{U}_3^{k})
-\rho_6(\mathbf{U}_2^{k+1}-\mathbf{U}_2^{k})\right\|_F\\
\leq\ &(\gamma_2+\rho_6)\|\mathbf{U}_2^{k+1}-\mathbf{U}_2^{k}\|_F,\\
\|N_7\|_F=&\left\|\nabla_{\mathbf{U}_3} \Psi(\mathcal{G}^{k+1}, \mathbf{A}_1^{k+1},\mathbf{A}_2^{k+1},
\mathbf{A}_3^{k+1},\mathbf{U}_1^{k+1},
\mathbf{U}_2^{k+1},\mathbf{U}_3^{k+1})\right.\\
&\left.-\nabla_{\mathbf{U}_3} \Psi(\mathcal{G}^{k+1}, \mathbf{A}_1^{k+1},\mathbf{A}_2^{k+1},
\mathbf{A}_3^{k+1},\mathbf{U}_1^{k+1},
\mathbf{U}_2^{k+1},\mathbf{U}_3^{k})
-\rho_7(\mathbf{U}_3^{k+1}-\mathbf{U}_3^{k})\right\|_F\\
\leq\ &(\gamma_3+\rho_7)\|\mathbf{U}_3^{k+1}-\mathbf{U}_3^{k}\|_F.
\end{aligned}
\end{equation}

By taking (\ref{N1_4}) collectively with (\ref{N2}), (\ref{N3}), (\ref{N4}), (\ref{N5}) and
(\ref{N6}), we get
$$
\begin{aligned}
\ &\|N^{k+1}\|_F=\|\left(N_1, N_2, N_3, N_4, N_5, N_6, N_7\right)\|_F\\
\leq \ &\|N_1\|_F+\| N_2\|_F+\| N_3\|_F+\| N_4\|_F+\|  N_5\|_F+\|  N_6\|_F+\|  N_7\|_F\\
\leq \ &(\rho_1+c_1)\|\mathcal{G}^{k+1}-\mathcal{G}^k\|_F+
c_2\|\mathbf{A}_1^{k}-\mathbf{A}_1^{k+1}\|_F
+(c_2+2\nu c_1)\|\mathbf{A}_2^{k}-\mathbf{A}_2^{k+1}\|_F\\
&+(c_2+2\nu c_1)\|\mathbf{A}_3^{k}-\mathbf{A}_3^{k+1}\|_F+
(\nu^2c_1+\rho_2+\gamma_1)\|\mathbf{A}_1^{k+1}-\mathbf{A}_1^k\|_F\\
&+(\nu c_2+2\nu^2c_1)\|\mathbf{A}_2^{k+1}-\mathbf{A}_2^k\|_F
+(\nu c_2+2\nu^2c_1)\|\mathbf{A}_3^{k+1}-\mathbf{A}_3^{k}\|_F
+\gamma_1\|\mathbf{U}_1^{k+1}-\mathbf{U}_1^k\|_F\\
&+(\nu^2c_1+\gamma_2+\rho_3)\|\mathbf{A}_2^{k+1}-\mathbf{A}_2^k\|_F
+2\nu^2c_1\|\mathbf{A}_3^{k+1}-\mathbf{A}_3^k\|_F
+\gamma_2\|\mathbf{U}_2^{k+1}-\mathbf{U}_2^k\|_F\\
&+(\nu^2c_1+\gamma_3+\rho_4)\|\mathbf{A}_3^{k+1}- \mathbf{A}_3^{k}\|_F
+\gamma_3\|\mathbf{U}_3^{k+1}-\mathbf{U}_3^k\|_F\\
&+(\gamma_1+\rho_5)\|\mathbf{U}_1^{k+1}-\mathbf{U}_1^{k}\|_F+
(\gamma_2+\rho_6)\|\mathbf{U}_2^{k+1}-\mathbf{U}_2^{k}\|_F+
(\gamma_3+\rho_7)\|\mathbf{U}_3^{k+1}-\mathbf{U}_3^{k}\|_F\\
= \ &(\rho_1+c_1)\|\mathcal{G}^{k+1}-\mathcal{G}^k\|_F+
(c_2+\nu^2c_1+\rho_2+\gamma_1)\|\mathbf{A}_1^{k}-\mathbf{A}_1^{k+1}\|_F\\
&+(c_1+2\nu c_1+\nu c_2+3\nu^2c_1+\gamma_2+\rho_3)\|
\mathbf{A}_2^{k}-\mathbf{A}_2^{k+1}\|_F\\
&+(c_2+2\nu c_1+\nu c_2+5\nu^2c_1+\gamma_3+\rho_4)\|
\mathbf{A}_3^{k}-\mathbf{A}_3^{k+1}\|_F\\
&+(2\gamma_1+\rho_5)\|\mathbf{U}_1^{k+1}-\mathbf{U}_1^{k}\|_F+
(2\gamma_2+\rho_6)\|\mathbf{U}_2^{k+1}-\mathbf{U}_2^{k}\|_F+
(2\gamma_3+\rho_7)\|\mathbf{U}_3^{k+1}-\mathbf{U}_3^{k}\|_F\\
\leq \ &\bar{\delta}\left(\|\mathcal{G}^{k+1}-\mathcal{G}^k\|_F
+\|\mathbf{A}_1^{k+1}-\mathbf{A}_1^k\|_F
+\|\mathbf{A}_2^{k+1}-\mathbf{A}_2^k\|_F
\right.\\
&+\left.
\|\mathbf{A}_3^{k+1}-\mathbf{A}_3^k\|_F
+\|\mathbf{U}_1^{k+1}-\mathbf{U}_1^{k}\|_F
+\|\mathbf{U}_2^{k+1}-\mathbf{U}_2^{k}\|_F
+\|\mathbf{U}_3^{k+1}-\mathbf{U}_3^{k}\|_F\right),
\end{aligned}
$$
where $\bar{\delta}:=\max\{\rho_1+c_1,c_2+\nu^2c_1+\rho_2+\gamma_1,c_1+2\nu c_1+\nu c_2+3\nu^2c_1+\gamma_2+\rho_3,
c_2+\nu c_2+2\nu c_1+5\nu^2c_1+\gamma_3+\rho_4,2\gamma_1+\rho_5,2\gamma_2+\rho_6,2\gamma_3+\rho_7\}$.
This completes the proof.

\bibliographystyle{abbrv}

\bibliography{reference}

\begin{thebibliography}{10}

\bibitem{Attouch2010}
H.~Attouch, J.~Bolte, P.~Redont, and A.~Soubeyran.
\newblock Proximal alternating minimization and projection methods for
  nonconvex problems: An approach based on the kurdyka-Łojasiewicz inequality.
\newblock {\em Math. Oper. Res.}, 35(2):438--457, 2010.

\bibitem{Attouch2013}
H.~Attouch, J.~Bolte, and B.~F. Svaiter.
\newblock Convergence of descent methods for semi-algebraic and tame problems:
  proximal algorithms, forward{\textendash}backward splitting, and regularized
  {Gauss-Seidel} methods.
\newblock {\em Math. Program.}, 137(1-2):91--129, 2013.

\bibitem{bahadori2014fast}
M.~T. Bahadori, Q.~R. Yu, and Y.~Liu.
\newblock Fast multivariate spatio-temporal analysis via low rank tensor
  learning.
\newblock {\em Advances in Neural Information Processing Systems}, 27, 2014.

\bibitem{basu2019low}
S.~Basu, X.~Li, and G.~Michailidis.
\newblock Low rank and structured modeling of high-dimensional vector
  autoregressions.
\newblock {\em IEEE Trans. Signal Process.}, 67(5):1207--1222, 2019.

\bibitem{Basu2015}
S.~Basu and G.~Michailidis.
\newblock Regularized estimation in sparse high-dimensional time series models.
\newblock {\em Ann. Statist.}, 43(4):1535--1567, 2015.

\bibitem{Beck2017}
A.~Beck.
\newblock {\em First-order methods in optimization}.
\newblock Society for Industrial and Applied Mathematics, Philadelphia, 2017.

\bibitem{Belkin2002}
M.~Belkin and P.~Niyogi.
\newblock Laplacian eigenmaps and spectral techniques for embedding and
  clustering.
\newblock In {\em Proceedings of the 15th International Conference on Neural
  Information Processing Systems}, page 585–591, 2001.

\bibitem{Bolte2014}
J.~Bolte, S.~Sabach, and M.~Teboulle.
\newblock Proximal alternating linearized minimization for nonconvex and
  nonsmooth problems.
\newblock {\em Math. Program.}, 146(1-2):459--494, 2014.

\bibitem{Cai2011}
D.~Cai, X.~He, J.~Han, and T.~S. Huang.
\newblock Graph regularized nonnegative matrix factorization for data
  representation.
\newblock {\em IEEE Trans. Pattern Anal. Mach. Intell.}, 33(8):1548--1560,
  2011.

\bibitem{chen2022factor}
R.~Chen, D.~Yang, and C.-H. Zhang.
\newblock Factor models for high-dimensional tensor time series.
\newblock {\em J. Amer. Statist. Assoc.}, 117(537):94--116, 2022.

\bibitem{chen2025forecasting}
X.~Chen, C.~Zhang, X.-L. Zhao, N.~Saunier, and L.~Sun.
\newblock Forecasting sparse movement speed of urban road networks with
  nonstationary temporal matrix factorization.
\newblock {\em Transportation Sci.}, 59(3):670--687, 2025.

\bibitem{Chen2021}
Y.~Chen, Y.~Chi, J.~Fan, and C.~Ma.
\newblock Spectral methods for data science: A statistical perspective.
\newblock {\em Foundations and Trends® in Machine Learning}, 14(5):566--806,
  2021.

\bibitem{chung1997spectral}
F.~R. Chung.
\newblock {\em Spectral Graph Theory}, volume~92.
\newblock Providence: American Mathematical Society, 1997.

\bibitem{DeLathauwer2000}
L.~De~Lathauwer, B.~De~Moor, and J.~Vandewalle.
\newblock A multilinear singular value decomposition.
\newblock {\em SIAM J. Matrix Anal. Appl.}, 21(4):1253--1278, 2000.

\bibitem{de2008forecasting}
C.~De~Mol, D.~Giannone, and L.~Reichlin.
\newblock Forecasting using a large number of predictors: Is {B}ayesian
  shrinkage a valid alternative to principal components?
\newblock {\em J. Econometrics}, 146(2):318--328, 2008.

\bibitem{duchon2012time}
C.~Duchon and R.~Hale.
\newblock {\em Time Series Analysis in Meteorology and Climatology: An
  Introduction}.
\newblock Wiley-Blackwell, Oxford, 2012.

\bibitem{golub2013matrix}
G.~H. Golub and C.~F. Van~Loan.
\newblock {\em Matrix Computations, 4th Edition}.
\newblock Johns Hopkins University Press, Philadelphia, PA, 2013.

\bibitem{gorrostieta2012investigating}
C.~Gorrostieta, H.~Ombao, P.~B{\'e}dard, and J.~N. Sanes.
\newblock Investigating brain connectivity using mixed effects vector
  autoregressive models.
\newblock {\em NeuroImage}, 59(4):3347--3355, 2012.

\bibitem{Gueridi2022}
D.~Gueridi and F.~Kittaneh.
\newblock Inequalities for the kronecker product of matrices.
\newblock {\em Ann. Funct. Anal.}, 13(3):50, 2022.

\bibitem{hampton2013quantifying}
S.~E. Hampton, E.~E. Holmes, L.~P. Scheef, M.~D. Scheuerell, S.~L. Katz, D.~E.
  Pendleton, and E.~J. Ward.
\newblock Quantifying effects of abiotic and biotic drivers on community
  dynamics with multivariate autoregressive (mar) models.
\newblock {\em Ecology}, 94(12):2663--2669, 2013.

\bibitem{han2015direct}
F.~Han, H.~Lu, and H.~Liu.
\newblock A direct estimation of high dimensional stationary vector
  autoregressions.
\newblock {\em J. Mach. Learn. Res.}, 16(1):3115--3150, 2015.

\bibitem{Han2022}
R.~Han, R.~Willett, and A.~R. Zhang.
\newblock An optimal statistical and computational framework for generalized
  tensor estimation.
\newblock {\em Ann. Statist.}, 50(1):1--29, 2022.

\bibitem{harris2021time}
K.~D. Harris, A.~Aravkin, R.~Rao, and B.~W. Brunton.
\newblock Time-varying autoregression with low-rank tensors.
\newblock {\em SIAM J. Appl. Dyn. Syst.}, 20(4):2335--2358, 2021.

\bibitem{hillar2013most}
C.~J. Hillar and L.-H. Lim.
\newblock Most tensor problems are {NP}-hard.
\newblock {\em J. ACM}, 60(6):45, 2013.

\bibitem{Matthew2021}
M.~Hirn.
\newblock Lecture 01: Introduction to spectral graph theory.
\newblock Technical report, 2021.

\bibitem{hitchcock1927expression}
F.~L. Hitchcock.
\newblock The expression of a tensor or a polyadic as a sum of products.
\newblock {\em J. Math. Phys.}, 6(1-4):164--189, 1927.

\bibitem{Horn2013}
R.~A. Horn and C.~R. Johnson.
\newblock {\em Matrix Analysis}.
\newblock Cambridge University Press, New York, 2nd edition, 2013.

\bibitem{kock2025data}
A.~B. Kock, R.~S. Pedersen, and J.~R.-V. S{\o}rensen.
\newblock Data-driven tuning parameter selection for high-dimensional vector
  autoregressions.
\newblock {\em J. Amer. Statist. Assoc.}, pages 1--19, 2025.

\bibitem{kolda2009tensor}
T.~G. Kolda and B.~W. Bader.
\newblock Tensor decompositions and applications.
\newblock {\em SIAM Rev.}, 51(3):455--500, 2009.

\bibitem{Koop2013}
G.~M. Koop.
\newblock Forecasting with medium and large bayesian vars.
\newblock {\em J. Appl. Econometrics}, 28(2):177--203, 2013.

\bibitem{Masini2022}
R.~P. Masini, M.~C. Medeiros, and E.~F. Mendes.
\newblock Regularized estimation of high\text{-}dimensional vector
  autoregressions with weakly dependent innovations.
\newblock {\em J. Time Series Anal.}, 43(4):532--557, 2022.

\bibitem{miao2023high}
K.~Miao, P.~C. Phillips, and L.~Su.
\newblock High-dimensional {VAR}s with common factors.
\newblock {\em J. Econometrics}, 233(1):155--183, 2023.

\bibitem{Negahban2011}
S.~Negahban and M.~J. Wainwright.
\newblock Estimation of (near) low-rank matrices with noise and
  high-dimensional scaling.
\newblock {\em Ann. Statist.}, 39(2):1069--1097, 2011.

\bibitem{Qiu2024}
D.~Qiu, B.~Yang, and X.~Zhang.
\newblock Robust tensor completion via dictionary learning and generalized
  nonconvex regularization for visual data recovery.
\newblock {\em IEEE Trans. Circuits Syst. Video Technol.}, 34(11):11026--11039,
  2024.

\bibitem{Qiu2022}
Y.~Qiu, G.~Zhou, Y.~Wang, Y.~Zhang, and S.~Xie.
\newblock A generalized graph regularized non-negative tucker decomposition
  framework for tensor data representation.
\newblock {\em IEEE Trans. Cybernet.}, 52(1):594--607, 2022.

\bibitem{raskutti2019convex}
G.~Raskutti, M.~Yuan, and H.~Chen.
\newblock Convex regularization for high-dimensional multiresponse tensor
  regression.
\newblock {\em Ann. Statist.}, 47(3):1554--1584, 2019.

\bibitem{reinsel1983some}
G.~Reinsel.
\newblock Some results on multivariate autoregressive index models.
\newblock {\em Biometrika}, 70(1):145--156, 1983.

\bibitem{Rockafellar1998}
R.~T. Rockafellar and R.~J.~B. Wets.
\newblock {\em Variational Analysis}.
\newblock Springer, New York, 1998.

\bibitem{Roughan2012}
M.~Roughan, Y.~Zhang, W.~Willinger, and L.~Qiu.
\newblock Spatio-temporal compressive sensing and internet traffic matrices
  (extended version).
\newblock {\em IEEE/ACM Trans. Netw.}, 20(3):662--676, 2012.

\bibitem{samadi2024reduced}
S.~Y. Samadi and H.~W.~B. Herath.
\newblock Reduced-rank envelope vector autoregressive model.
\newblock {\em J. Bus. Econom. Statist.}, 42(3):918--932, 2024.

\bibitem{sims1980macroeconomics}
C.~A. Sims.
\newblock Macroeconomics and reality.
\newblock {\em Econometrica}, 48(1):1--48, 1980.

\bibitem{Stewart1990}
G.~W. Stewart and J.-g. Sun.
\newblock {\em Matrix Perturbation Theory}.
\newblock Academic Press, New York, 1990.

\bibitem{tsay2005analysis}
R.~S. Tsay.
\newblock {\em Analysis of Financial Time Series}.
\newblock John Wiley \& Sons, New York, 2005.

\bibitem{wang2023rate}
D.~Wang and R.~S. Tsay.
\newblock Rate-optimal robust estimation of high-dimensional vector
  autoregressive models.
\newblock {\em Ann. Statist.}, 51(2):846--877, 2023.

\bibitem{Wang2022a}
D.~Wang, X.~Zhang, G.~Li, and R.~Tsay.
\newblock High-dimensional vector autoregression with common response and
  predictor factors.
\newblock {\em arXiv:2203.15170}, 2022.

\bibitem{wang2024high}
D.~Wang, Y.~Zheng, and G.~Li.
\newblock High-dimensional low-rank tensor autoregressive time series modeling.
\newblock {\em J. Econometrics}, 238(1):105544, 2024.

\bibitem{Wang2022}
D.~Wang, Y.~Zheng, H.~Lian, and G.~Li.
\newblock High-dimensional vector autoregressive time series modeling via
  tensor decomposition.
\newblock {\em J. Amer. Statist. Assoc.}, 117(539):1338--1356, 2022.

\bibitem{Wang2024}
H.~Wang, T.~Liu, R.~Li, M.~X. Cheng, T.~Zhao, and J.~Gao.
\newblock Roselora: Row and column-wise sparse low-rank adaptation of
  pre-trained language model for knowledge editing and fine-tuning.
\newblock In {\em Proceedings of the 2024 Conference on Empirical Methods in
  Natural Language Processing}, pages 996--1008. Association for Computational
  Linguistics, 2024.

\bibitem{wang2014multi}
H.~Wang, C.~Weng, and J.~Yuan.
\newblock Multi-feature spectral clustering with minimax optimization.
\newblock In {\em Proceedings of the IEEE Conference on Computer Vision and
  Pattern Recognition}, pages 4106--4113, 2014.

\bibitem{wang2023regularized}
Z.~Wang, A.~Safikhani, Z.~Zhu, and D.~S. Matteson.
\newblock Regularized estimation in high-dimensional vector auto-regressive
  models using spatio-temporal information.
\newblock {\em Statist. Sinica}, 33:1271--1294, 2023.

\bibitem{wong2020lasso}
K.~C. Wong, Z.~Li, and A.~Tewari.
\newblock Lasso guarantees for $\beta$-mixing heavy-tailed time series.
\newblock {\em Ann. Statist.}, 48(2):1124--1142, 2020.

\bibitem{Yang2019}
K.~Yang, X.~Dong, and L.~Toni.
\newblock Error analysis on graph laplacian regularized estimator.
\newblock {\em arXiv:1902.03720}, 2019.

\bibitem{Yu2015}
Y.~Yu, T.~Wang, and R.~J. Samworth.
\newblock A useful variant of the davis–kahan theorem for statisticians.
\newblock {\em Biometrika}, 102(2):315--323, 2015.

\bibitem{zhu2020nonconcave}
X.~Zhu.
\newblock Nonconcave penalized estimation in sparse vector autoregression
  model.
\newblock {\em Electron. J. Stat.}, 14(1):1413--1448, 2020.

\bibitem{Zou2006}
H.~Zou, T.~Hastie, and R.~Tibshirani.
\newblock Sparse principal component analysis.
\newblock {\em J. Comput. Graph. Statist.}, 15(2):265--286, 2006.

\end{thebibliography}

\end{document}